\documentclass[11pt,letterpaper]{amsart}
              
\usepackage{amsmath,amssymb,latexsym, amsbsy, amsfonts, amscd, amsthm, mathrsfs, xy,comment, xcolor}
\usepackage{multirow, tabularx}
\usepackage{hyperref}
\usepackage{enumitem, empheq, url}
 \usepackage[bibtex-style]{amsrefs}
\xyoption{all}
\usepackage{tikz-cd}
\usepackage{rotating, scalerel}

\theoremstyle{plain}

\newtheorem{theorem}{Theorem}[section]

\newtheorem{prop}[theorem]{Proposition}
\newtheorem{corollary}[theorem]{Corollary}

\newtheorem{lemma}[theorem]{Lemma}

\newcommand{\longtwoheadrightarrow}{\relbar\joinrel\twoheadrightarrow}

\theoremstyle{definition}
\newtheorem{definition}[theorem]{Definition}
\newtheorem{remark}[theorem]{Remark}
\newtheorem{example}[theorem]{Example}

\long\def\symbolfootnote[#1]#2{\begingroup
\def\thefootnote{\fnsymbol{footnote}}\footnote[#1]{#2}\endgroup}

\def\lra{\longrightarrow}

\DeclareMathOperator{\GL}{GL}
\def\ZG{\mathbf{Z}[G]}

\def\D{\mathbf{D}}

\def\A{\mathbf{A}}
\def\cG{\mathcal{G}}
\def\cI{\mathcal{I}}
\def\cA{\mathcal{A}}

\def\sgn{\mathrm{sgn}}
\def\N{\mathrm{N}}
\def\1{\mf{1}}

\DeclareMathOperator{\lrar}{Tor}
\DeclareMathOperator{\id}{id}

\DeclareMathOperator{\Fitt}{Fitt}
\DeclareMathOperator{\Ann}{Ann}
\DeclareMathOperator{\Sym}{Sym}
\DeclareMathOperator{\Ind}{Ind}
\DeclareMathOperator{\cts}{cts}

\DeclareMathOperator{\Frac}{Frac}

\DeclareMathOperator{\res}{res}

\DeclareMathOperator{\cyc}{cyc}

\DeclareMathOperator{\cond}{cond}

\DeclareMathOperator{\rec}{rec} 
\DeclareMathOperator{\Hom}{Hom} \DeclareMathOperator{\End}{End}

\DeclareMathOperator{\coker}{coker} 
\DeclareMathOperator{\ord}{ord}

\DeclareMathOperator{\cusps}{cusps}
 \DeclareMathOperator{\real}{Re}
 
 \DeclareMathOperator{\Cl}{Cl}
 \DeclareMathOperator{\red}{red}
 
 \DeclareMathOperator{\Frob}{Frob}
\DeclareMathOperator{\Spec}{Spec}
\DeclareMathOperator{\Ext}{Ext}
\DeclareMathOperator{\Sel}{Sel}
\DeclareMathOperator{\tr}{tr}

\DeclareMathOperator{\im}{im}

\newcommand{\bpsi}{\pmb{\psi}}

\newcommand{\mat}[4]{\begin{pmatrix} {#1} & {#2} \\ {#3} & {#4} \end{pmatrix}}

\newcommand{\mf}{\mathfrak }

\def\fa{\mathfrak{a}}
\def\fn{\mathfrak{n}}
\def\fp{\mathfrak{p}}
\def\fq{\mathfrak{q}}

\def\fm{\mathfrak{m}}
\def\ft{\mathfrak{t}}

\def\fl{\mathfrak{l}}
\def\fP{\mathfrak{P}}

\def\fm{\mathfrak{m}}

\def\fd{\mathfrak{d}}

\def\T{\mathbf{T}}
\def\Z{\mathbf{Z}}
\def\F{\mathbf{F}}
\def\Q{\mathbf{Q}}
\def\C{\mathbf{C}}
\def\G{\mathbf{G}}
\def\B{\mathbf{B}}
\def\R{\mathbf{R}}

\def\bdf{\begin{defn}}
\def\edf{\end{defn}}
\def\cM{\mathcal{M}}

\def\cD{\mathcal{D}}
\def\cO{\mathcal{O}}

\def\cF{\mathcal{F}}
\def\cB{\mathcal{B}}
\def\cP{\mathcal{P}}

\def\cA{\mathcal{A}}
\def\cC{\mathcal{C}}
\def\cV{\mathcal{V}}

\def\cG{\mathcal{G}}

\def\fb{\mathfrak{b}}

\def\Gal{{\rm Gal}}

\def\ab{{\rm ab}}

\def\loc{\text{loc}}
\def\ram{\text{ram}}
\def\ab{\text{ab}}

\begin{document}

\title{The Brumer--Stark Conjecture over $\Z$}
\author{Samit Dasgupta, Mahesh Kakde, Jesse Silliman, \\ and Jiuya Wang}

\maketitle

{\small 
\paragraph{Abstract.}
 In this paper we give a complete proof of the Brumer--Stark conjecture over $\Z$.  In previous work, the first two named authors proved the result over $\Z[1/2]$. There are two conceptual difficulties that occur when attempting to generalize to the localization at $p=2$.  
  The first issue concerns the Ritter--Weiss module  $\nabla_{\Sigma}^{\Sigma'}(H)$ studied in our previous work.  We define a new generalized Ritter--Weiss module  $\nabla_{\Sigma}^{\Sigma'}(H)$ in which the depletion set 
$\Sigma$ need not contain the archimedean primes. The connection between our Ritter--Weiss modules and Galois cohomology is handled in a more conceptual way than in our previous work, and requires new results on class formations.  The second issue is the need for a general version of Ribet's Lemma that applies in the residually indistinguishable case.  This requires the application of results from matrix invariant theory and the use of rational cohomology.
}

\tableofcontents

 \section{Introduction}

Let $F$ be a totally real field and let $H$ be a finite abelian CM extension of $F$.  Let $\fp$ denote a prime ideal of $F$ that splits completely in $H$.  The Brumer--Stark conjecture states the existence of a $\fp$-unit in $H$ whose valuations at the primes above $\fp$ are related to the values $L(\chi, 0)$ as $\chi$ ranges over the characters of the extension $H/F$.  See Theorem~\ref{t:bs}  for a precise statement.  In previous work, the first two named authors proved the Brumer--Stark conjecture over $\Z[1/2]$, i.e.~we proved the existence of an element in the $\fp$-unit group tensored with $\Z[1/2]$ with the desired properties.  

In this paper, we give a complete proof of the Brumer--Stark conjecture over $\Z$.  In fact we obtain a strong refinement of the Brumer--Stark conjecture that yields the Fitting ideal of a certain Ritter--Weiss module.
In a companion paper  to this one, we show how to use this result to deduce
 the minus part of the Equivariant Tamagawa Number Conjecture (ETNC) for the Tate motive associated to $H/F$ \cite{etnc}. Again, we obtain this integrally over $\Z$ and not just over $\Z[1/2]$.  The proof of ETNC, which is obtained by applying an idea of Bullach--Burns--Daoud--Seo \cite{bbds} to the results on Brumer--Stark proven here, yields many important new corollaries. These include Rubin's higher rank Brumer--Stark conjecture, the integral Gross--Stark conjecture, and the higher rank version due to Popescu.
 One should also obtain
 the classical Main Conjecture of Iwasawa Theory over totally real fields at the prime $p=2$, though there does not appear to be a precise statement of the Main Conjecture at $p=2$ in the literature (but see \cite{bks}*{Conjecture 3.1}).  Our proof of ETNC is rather formal, with the main arithmetic input being Theorem~\ref{t:main} from this paper.

There are two conceptual difficulties that occur when attempting to generalize our previous work on the Brumer--Stark conjecture to the prime $p=2$.   
The first of these concerns the Ritter--Weiss module $\nabla_{\Sigma}^{\Sigma'}(H)$ appearing in our previous work.
It turns out that this module is not what appears in the Galois representations that we study locally at 2.  We must alter the module
 $\nabla_{\Sigma}^{\Sigma'}(H)$ by modifying the local conditions at the archimedean places.  On the analytic side this modification is reflected in the power of 2 in the denominator of the constant terms of Eisenstein series for $F$. 
 
These modifications have an important effect on our constructions with group-ring valued Hilbert modular forms.  In the case that the extension $H/F$ is not ramified at any primes above 2, the set $\Sigma$ that one would naturally like to take is empty; however it is an important assumption in our study of Ritter--Weiss modules that $\Sigma$ is non-empty.  One is therefore led to introduce an auxiliary prime $\fl$ to the set $\Sigma$ that does not divide the conductor of the extension $H/F$.  The inclusion of this prime leads to the presence of oldforms in our constructions with Hilbert modular forms, which in turn yields a Hecke algebra $\T$ that is not reduced.  Now working with a non-reduced Hecke algebra, we show that the Galois representations we construct have the desired properties.  The connection between our generalized Ritter--Weiss modules and Galois cohomology is developed in a more conceptual way than in \cite{dk}, and requires new results on class formations.

The next issue at $p=2$ is in the application of Ribet's method, which provides a strategy for constructing nontrivial extensions of a $p$-adic Galois representation $\rho_1$ by another such representation $\rho_2$.  Our previous proof of Brumer--Stark over $\Z_p$ for $p$ odd is an application of a  general version of Ribet's method.  However, an important assumption that occurs throughout the literature is that the representations $\rho_i$ are residually distinguishable, i.e.\ that 
\[ \overline{\rho}_1 \not\cong \overline{\rho}_2 \pmod{p}. \] 
In the setting of Brumer--Stark, $\rho_1$ is  the trivial character and $\rho_2$ is a totally odd character $\chi$.  When $p$ is odd, $\chi(c) = -1 \not\equiv 1\pmod{p}$, where $c$ denotes complex conjugation.  Ribet's method therefore carries through.  However when $p=2$, we may have $\chi \equiv 1 \pmod{2}$, and a new construction is necessary.  
In this paper, we prove a strong version of Ribet's Lemma that applies in the residually indistinguishable case (see Theorem~\ref{t:rl}). 

In the rest of the introduction, we state our results precisely.  First we recall the statement of the Brumer--Stark conjecture,
and next we describe our results on Ritter--Weiss modules and Ribet's Lemma.

\subsection{The Brumer--Stark Conjecture}

Recall that we have fixed a totally real field $F$ of degree $n$ and a finite abelian extension $H/F$ that is a CM field.  Write $G = \Gal(H/F)$.   Let $S$ and $T$ denote finite  disjoint sets of places of $F$.
Associated to any character $\chi \colon G \longrightarrow \C^*$  one has the Artin $L$-function
\begin{equation} \label{e:ls} 
L_S(\chi, s) = \prod_{\substack{\fp \not \in S \\ \fp \text{ finite}}} \frac{1}{1 - \chi(\fp) \N\fp^{-s}}, \qquad \real(s) > 1, 
\end{equation}
and its ``$T$-smoothed" version
\begin{equation} \label{e:lst}
L_{S,T}(\chi, s) = L_S(\chi, s) \prod_{\substack{\fp \in T \\ \fp \text{ finite}}} (1 - \chi(\fp)\N\fp^{1-s}). 
\end{equation}

Assume now that $S$ contains the set $S_{\infty}$ of real places and the set $S_{\ram}$ of finite primes ramifying in $H$. 
Assume that $T$ is torsion-free: the only root of unity in $H$
congruent to $1$ at every prime of $H$ above $T$ is $1$.
Since $S\cap T=\emptyset$ and
$S\supset S_\infty\cup S_{\ram}$, the primes in $T$ are
finite and unramified in $H/F$.

The familiar Deligne--Ribet condition---that
$T$ contain two primes of different residue characteristics,
or one prime of residue characteristic larger than $n+1$---is
sufficient for torsion-freeness
(see \cite{dk}*{(1) and the following paragraph}).
We do not impose this stronger numerical condition.

\begin{theorem} [Brumer--Stark Conjecture] \label{t:bs}  
 Let $\fp \not\in S \cup T$ be a prime of $F$ that splits completely in $H$.  Fix a prime $\fP$ of $H$ above $\fp$.  There exists an element  $u \in H^*$ 
 satisfying the following.
 \begin{itemize}
\item  We have $|u|_w = 1$ for all places $w$ of $H$ not lying above $\fp$, including the complex places. 
\item We have
\begin{equation} \label{e:bs}
\sum_{\sigma \in G} \chi(\sigma) \ord_{\sigma^{-1}(\fP)}(u) = L_{S,T}(\chi, 0)
 \end{equation}
for all characters $\chi\colon G \longrightarrow \C^*$. 
\item We have $u \equiv 1 \pmod{\fq\cO_H}$ for all $\fq \in T$.
\end{itemize}
\end{theorem}
In previous work, the first two named authors proved that there exists $u \in \cO_H[1/\fp]^* \otimes \Z[1/2]$ satisfying the properties of Theorem~\ref{t:bs}.  In this paper, we complete the proof of the theorem.

\subsection{Ritter--Weiss Modules} \label{s:rwm}

The {\em Stickelberger element} $\Theta_{S,T} \in \C[G]$ is the unique element such that 
\begin{equation}  \chi(\Theta_{S,T}) = L_{S,T}(\chi^{-1}, 0)
 \end{equation}
for all characters $\chi$ of $G$.  A theorem of Deligne--Ribet \cite{dr} and Cassou--Nogu\`es states that $\Theta_{S,T}$ lies in $\Z[G]$.
As is well-known, up to issues at $p=2$, the Brumer--Stark conjecture is equivalent to the statement that the Stickelberger element $\Theta_{S,T}$ annihilates the $T$-smoothed class group $\Cl^T(H)$.
Ritter and Weiss defined a certain module $\nabla_S^T(H)$ that is related to $\Cl^T(H)$ but (on the algebraic side) has a more precise relationship with Galois cohomology, and (on the analytic side) has a more precise relationship with the Stickelberger element  
$\Theta_{S,T}$ (for details see Theorem \ref{t:sbs}, in particular, we note that $\Theta_{S,T}$ annihilates the full group $\Cl^T(H)$ and not just its quotient $\Cl^T_S(H)$, defined below).  In order to complete the proof of Brumer--Stark, we must generalize the construction of Ritter and Weiss.

% For the moment, let $F$ be an arbitrary number field, and $H$ a finite abelian extension of $F$ with $G = \Gal(H/F)$.  
 Let $\Sigma, \Sigma'$ denote disjoint finite sets of places of $F$. {For a set of places $S$ of $F$, denote by $S_H$ the set of primes of $H$ lying above those in $S$.} 
{Define
\[ \begin{split}  \cO_{H,S,\Sigma'}^* =  \{ u \in H^* \colon & |u|_\fq =1 \text{ for all finite } \fq \not\in S_H, \text{ and }
 \\  & u \equiv 1 \!\!\! \pmod{\fq} \text{ for all finite } \fq \in \Sigma'_H\}.
\end{split} \]
}
 {In particular, we have
 \[ \cO_{H, S_{\infty}, \Sigma'}^* = \{u \in \cO_H^* : u \equiv 1 \!\!\! \pmod{\fq} \text{ for all finite } \fq \in \Sigma'_H\}. \]   
 
 \begin{itemize}
\item The pair $(\Sigma, \Sigma')$ is said to be an {\bf RW pair} if \begin{equation} \label{e:rwdef}
 \Sigma \cup \Sigma' \supset S_\infty.\end{equation}
 In \S4 we define a generalized
Ritter--Weiss module $\nabla_\Sigma^{\Sigma'}(H)$ for any RW pair.
\item We say that $\Sigma'$ is {\bf torsion-free} if  \[ (\cO_{H, S_{\infty}, \Sigma'}^*)_{\rm tor} = \{1\}. \]
\item For a fixed prime $p$, we say that an RW pair $(\Sigma,\Sigma')$ is {\bf $p$-admissible} if
$\Sigma'$ is torsion-free, $\Sigma\neq\emptyset$, and
\[
\Sigma\supset \{v\in S_{\mathrm{ram}}:v\mid p\},\qquad
\Sigma'\supset \{v\in S_{\mathrm{ram}}:v\nmid p\}.
\]
\end{itemize}
See Remark~\ref{r:conditions} below for a discussion of these conditions.}

Our definition of the $\Z[G]$-module $\nabla_{\Sigma}^{\Sigma'}(H)$ {in (\ref{e:rwm}) below} generalizes the construction of Ritter--Weiss \cite{rw}, which did not consider the smoothing set $\Sigma'$, and our previous work \cite{dk}, which imposed an additional condition $\Sigma \supset S_\infty$.  {When $(\Sigma,\Sigma')$ is $p$-admissible,} the module  $\nabla_{\Sigma}^{\Sigma'}(H)$ sits in a short exact sequence of $\Z[G]$-modules
\[ 0 \longrightarrow \Cl_\Sigma^{\Sigma'}(H) \longrightarrow \nabla_{\Sigma}^{\Sigma'}(H) \longrightarrow X_{H, \Sigma} \longrightarrow 0. \]
Here $X_{H, \Sigma}$ denotes the group of degree zero divisors {on $\Sigma_H$,} and $\Cl_\Sigma^{\Sigma'}(H)$ is the $\Sigma$-depleted, $\Sigma'$-smoothed class group of $H$ (for a definition see \cite{dk}*{page 292}).

%We now specialize the situation and impose further conditions on $\Sigma, \Sigma'$.  Assume that $F$ is totally real and that $H$ is a CM field.  
Denote by $c \in G$ the complex conjugation, and define
\[ \Z[G]_- = \Z[G]/(1 + c). \] 
For any $\Z[G]$-module $M$, we define 
\[ M_- = M \otimes_{\Z[G]} \Z[G]_- = M/(1+c)M. \]

We write \[ \nabla_{\Sigma}^{\Sigma'}(H)_{p, -} = \nabla_{\Sigma}^{\Sigma'}(H) \otimes_{\Z[G]} \Z_p[G]_-. \]
{When $(\Sigma, \Sigma')$ is $p$-admissible}, we prove that $ \nabla_{\Sigma}^{\Sigma'}(H)_{p, -} $ is quadratically presented over $\Z_p[G]_-$, and hence that its 0th Fitting ideal is principal.  The following is our main result {and is proven in \S\ref{ss:ie}.
The statement is the natural adaptation to the RW pair $(\Sigma, \Sigma')$ of a conjecture of Burns, Kurihara, and Sano (see Remark~\ref{r:bks}).}
{
\begin{theorem}  \label{t:main}
If $(\Sigma,\Sigma')$ is $p$-admissible, we have
\[ \Fitt_{\Z_p[G]_-}(\nabla_{\Sigma}^{\Sigma'}(H)_{p,-}) = (\Theta_{\Sigma, \Sigma'}/2^t),
\]
where $\Fitt$ denotes the 0th Fitting ideal and $t = |S_\infty \cap \Sigma'|$.
\end{theorem}
}

Here $\Theta_{\Sigma,\Sigma'}/2^t$  denotes the unique element of $\Z_p[G]_-$ such that
\begin{equation} \label{e:td}
 \chi(\Theta_{\Sigma, \Sigma'}/2^t) = L_{\Sigma,\Sigma'}(\chi^{-1}, 0)/2^t 
 \end{equation}
for all odd characters $\chi \in \hat{G}$.  The congruence of Deligne--Ribet \cite{dr}*{Corollary 8.9} together with the argument of Kurihara recalled in \cite{dk}*{Lemma 3.4} implies that there exists an element satisfying (\ref{e:td}) in $\frac{c-1}{2}\Z_p[G]$.  The existence of  $\Theta_{\Sigma,\Sigma'}/2^t \in \Z_p[G]_-$ follows. We thank Kurihara for pointing this out to us.

Theorem~\ref{t:main} was proven for $p$ odd in \cite{dk}, so it suffices to prove the result for $p=2$.  However, since most of the arguments presented here work for all $p$, we do not in general set $p=2$ in what follows (though we do provide certain arguments only for $p=2$ when the analogous results for $p$ odd were already proven in \cite{dk}).

The fact that Theorem~\ref{t:main} implies the Brumer--Stark conjecture (Theorem~\ref{t:bs}) is proven in \S\ref{s:bs}.

\begin{remark}  \label{r:conditions}
 We indicate where the conditions in $p$-admissibility enter.
The RW condition is needed to define the generalized Ritter--Weiss module.
The torsion-free condition is used to obtain projective presentations.

The condition $\Sigma\neq\emptyset$ is technical but essential: one element
of $\Sigma$ is used to rigidify the construction and obtain an explicit
presentation in \S4.6. Note for example that $X_{H,\Sigma}$ collapses when $\Sigma$ is empty (more importantly the degree map $Y_{H,\Sigma} \rightarrow \Z$ used to define $X_{H,\Sigma}$ is not surjective).

Finally, the distribution of finite ramified primes
between $\Sigma$ and $\Sigma'$ is chosen to match the local conditions
arising from the $p$-adic Galois representation constructed from
group-ring valued Hilbert modular forms.
\end{remark}

\begin{remark} \label{r:bks}  When $S \supset S_\infty \cup S_{\ram}$, the statement
\begin{equation} \label{e:bks}
 \Fitt_{\Z[G]_-}(\nabla_{S}^{T}(H)_{-}) = (\Theta_{S, T})
\end{equation}
was conjectured by Burns, Kurihara, and Sano in \cite{bks}.  Our initial goal was to prove this result directly. However, in our modular forms constructions we found it necessary to modify the $S$ and $T$ sets based on the prime $p$. This led to the statement of Theorem~\ref{t:main}---essentially the same statement but with the $S, T$ sets replaced by $\Sigma, \Sigma'$ (and a power of 2 included to account for relaxing the assumption $S_\infty \subset \Sigma$).  Note that while (\ref{e:bks}) can be stated over $\Z$, the sets $\Sigma, \Sigma'$ depend on $p$ and Theorem~\ref{t:main} only holds over $\Z_p$.  Comment added in proof: In work subsequent to this one \cite{etnc}, we proved that Theorem~\ref{t:main} for all $p$ actually implies (\ref{e:bks}), however this was not known to us when this paper was initially written.
\end{remark}

\begin{remark} \label{r:pdef}
As we discuss in \S\ref{s:sr}, the $\Z_p[G]_-$-module $\nabla_{\Sigma}^{\Sigma'}(H)_{p,-}$ has a canonical transpose denoted  $\Sel_{\Sigma}^{\Sigma'}(H)_{p,-}$, and Theorem~\ref{t:main} is equivalent to the statement
\[ \Fitt_{\Z_p[G]_-}(\Sel_{\Sigma}^{\Sigma'}(H)_{p,-}) = (\Theta_{\Sigma, \Sigma'}^\#/2^t), \]
where $\#$ is the involution on $\Z_p[G]_-$ induced by $g \mapsto g^{-1}$.
\end{remark}

\subsection{Ribet's Method} \label{ss:ribetmethod}

{
We first consider the case $t = n$, i.e.\ $S_\infty \subset \Sigma'$.}
To prove Theorem~\ref{t:main}, we construct a $\Z_p[G]_-$-module $M$ and a surjective $\Z_p[G]_-$-module map
\begin{equation}
\label{e:surjmap}
 \nabla_{\Sigma}^{\Sigma'}(H)_{p,-} \longtwoheadrightarrow M 
 \end{equation}
 such that
 \begin{equation} \label{e:fittinc}
  \Fitt_{\Z_p[G]_-}(M) \subset (\Theta_{\Sigma, \Sigma'}/2^n). \end{equation}
  General properties of Fitting ideals imply that (\ref{e:surjmap}) and (\ref{e:fittinc}) together yield
\begin{equation} \label{e:fittinc2}
      \Fitt_{\Z_p[G]_-}( \nabla_{\Sigma}^{\Sigma'}(H)_{p,-}) \subset (\Theta_{\Sigma, \Sigma'}/2^n).  
  \end{equation}  
The important and surprisingly delicate result that the inclusion (\ref{e:fittinc2}) is sufficient to imply an equality (i.e. Theorem~\ref{t:main}) is proven in \S\ref{ss:ie}.

In order to construct (\ref{e:surjmap}), we provide in \S\ref{s:gc} a cohomological interpretation of the Ritter--Weiss module $\nabla_{\Sigma}^{\Sigma'}(H)$.  We show that the existence of (\ref{e:surjmap}) is equivalent to 
the existence of a surjective cohomology class $\kappa \in H^1(F, M)$ satisfying certain local conditions.  
Such a cohomology class can be viewed as an extension of the canonical character $G_F \lra G$ by the trivial character.
We are therefore led to a fundamental problem in algebraic number theory: constructing non-trivial extensions of $p$-adic Galois representations with prescribed local properties (i.e. elements of Selmer groups).

Ribet's method, first established in the groundbreaking paper \cite{ribet}, 
provides a strategy for solving this problem.  Ribet's method has played a central role in Iwasawa theory. 
 After Ribet's original proof of the converse to Herbrand's theorem, his method was used by Mazur--Wiles to prove the Main Conjecture of Iwasawa theory over $\Q$ \cite{mw} and by Wiles to prove the Main Conjecture over arbitrary totally real fields (with some restrictions at $p=2$) \cite{wiles}.  An important innovation introduced in Mazur--Wiles is the use of Fitting ideals. More recently, generalizations of Ribet's method have been used by Skinner--Urban to prove the main conjecture for $p$-ordinary elliptic curves \cite{su} and by the first two authors of this paper to prove the Brumer--Stark conjecture away from $p=2$ \cite{dk}.

Ribet's idea is that the generator of the ideal of interest, in our case {
$\Theta_{\Sigma, \Sigma'}^\#/2^n$,} can be realized as the constant term of a certain Eisenstein series.  In our setting, these Eisenstein series are examples of {\em group-ring valued Hilbert modular forms}, as defined in \cite{dk}.  In \S\ref{s:modforms}, we 
construct  a group-ring family $F_k$ of cusp forms that is congruent to a family of Eisenstein series modulo $\Theta_{\Sigma, \Sigma'}^\#/2^n$. 
%  Here $\#$ denotes the $\Z$-linear map from $\Z[G]$ to itself given by $g \mapsto g^{-1}$ for all $g\in G$.
If we define  $\T$ as the Hecke algebra of the space of group-ring valued Hilbert modular cusp forms as  in \cite{dk}, then the cusp form $F_k$ gives rise to a $\Z_p[G]_-$-algebra homomorphism
{
\[ \varphi \colon \T \longrightarrow \Z_p[G]_-/(\Theta_{\Sigma, \Sigma'}^\#/2^n). \]
}
Let $I = \ker \varphi$.

 The Galois representations associated to Hilbert modular cusp forms yield a representation
 \[ \rho \colon G_F \longrightarrow \GL_2(K), 
 \]
{where $K = \Frac(\T)$ is the total ring of fractions of $\T$.  The representation $\rho$} has the property that the characteristic polynomial $P_{\rho(g)}(x)$ lies in $\T[x]$ for all $g \in G_F$ and 
 \begin{equation} \label{e:redii}
 P_{\rho(g)}(x) \equiv (x - \chi(g))(x - \psi(g)) \pmod{I} 
 \end{equation}
 for two characters $\chi, \psi \colon G_F \lra \T^*$.  The congruence (\ref{e:redii}) follows from the congruence between $F_k$ and the Eisenstein series.
Our construction of the Galois representation $\rho$ requires special care in the case that ${\T}$ is not reduced, a new feature which arises here from the insertion of an auxiliary prime $\fl$ into the set 
$\Sigma$.

What remains is to construct the module $M$ and cohomology class $\kappa$ from the homomorphism $\rho$.  This is the subject of Ribet's Lemma, which is discussed in the next subsection.  {
We mention in passing that moving from the case $t=n$ to arbitrary $t$ is a purely algebraic procedure conducted at the end---as described in \S\ref{s:minimal}, moving real places from $\Sigma'$ to $\Sigma$ multiplies the Fitting ideal of $\nabla_{\Sigma}^{\Sigma'}(H)_{p,-}$ by 2.
}

\subsection{Ribet's Lemma}

Suppose we are working over a local ring $(\T, \fm)$, and we are attempting to construct an extension of a $p$-adic Galois representation $\rho_1$ by another representation $\rho_2$.  An important assumption that occurs throughout the literature is that the representations $\rho_i$ are residually distinguishable, i.e.\ that 
\begin{equation} \label{e:rd}
 \overline{\rho}_1 \not\cong \overline{\rho}_2 \pmod{\fm}. \end{equation}
The seminal book on Ribet's method by Bellaiche--Chenevier \cite{bc}, which works in the language of pseudorepresentations, 
assumes throughout that the pseudorepresentations considered are ``residually multiplicity free'', i.e.\ that the pseudorepresentation analog of (\ref{e:rd}) holds.
Here we prove a very general version of Ribet's Lemma for $\GL_2$ where we do not impose the assumption that the associated characters are residually distinguished, since this is what we need to prove Brumer--Stark at $p=2$. The following is a simplified version; {see Theorem \ref{t:rl} for the full version}.

\begin{theorem} \label{t:rli} Let $\T$ be a complete Noetherian reduced  local ring,  let $I \subset \T$ be an ideal, and let $\cG$ be a group.  Suppose we have a representation
\[ \rho \colon \cG \longrightarrow \GL_2(K), \quad K=\Frac(\T), \]
such that the characteristic polynomial $P_{\rho(g)}(x)$ lies in $\T[x]$ for all $g \in \cG$ and 
\begin{equation} \label{e:redi}
 P_{\rho(g)}(x) \equiv (x - \chi(g))(x - \psi(g)) \pmod{I} 
 \end{equation}
for two characters $\chi, \psi \colon \cG \longrightarrow \T^*$.
  Suppose  that for every projection of $K$ onto one of its  field factors $K \rightarrow k$, the projection of $\rho$ to 
  $\GL_2(k)$ is an irreducible representation of $\cG$ over $k$.  Finally, assume that the extension of $\rho$ to a $\T$-algebra homomorphism $\T[\cG] \longrightarrow M_2(K)$ has image that is finitely generated as a $\T$-module.
  
Then there exists a finitely generated $\T$-module $M$ and a surjective cohomology class \[ \kappa \in H^1(\cG, M(\chi\psi^{-1})) \] such that
\begin{equation}  \label{e:fitmi}
 \Fitt_{\T}(M) \subset I.
 \end{equation}
\end{theorem}
Here $M(\chi\psi^{-1})$ denotes the $\T$-module $M$ endowed with a $\cG$-action via $\chi\psi^{-1}$.
A cohomology class is called surjective if for every representative cocycle $\kappa$, the elements $\kappa(g)$ for $g \in \cG$ generate $M$ as a $\T$-module.   

The novelty of Theorem~\ref{t:rli} is that we do not assume that $\chi \not\equiv \psi \pmod{\fm}$.  
The residually indistinguishable case $\chi \equiv \psi \pmod{\fm}$ was studied recently by Ophir and Weiss, who obtained a version of Ribet's Lemma when working with representations over a DVR \cite{ow}.  As shown by Hajjar Mu\~noz in \cite{rafah}, one can deduce Theorem~\ref{t:rli} when $\T$ is a DVR from the results of Ophir--Weiss.  However, in typical arithmetic applications, the local rings $\T$ that occur are Hecke algebras that are rarely DVRs.  In this paper we establish a new technique to handle the  general case. 

Before describing the proof of Theorem~\ref{t:rli}, it is illuminating to 
first consider the residually distinguishable case $\chi \not\equiv \psi \pmod{\fm}$.  In this setting, Theorem~\ref{t:rli} can be proven following Mazur--Wiles (see \cite{mw}*{Chapter 5, \S 5, Proposition 1}).  Fix $\tau \in \cG$ such that $\chi(\tau) \not\equiv \psi(\tau) \pmod{\fm}$.  By Hensel's lemma, $P_{\rho(\tau)}(x)$ has two distinct roots in $\T$, 
congruent to $\chi(\tau)$ and $\psi(\tau) \pmod{\fm}$ respectively.
Choose a basis for $\rho$ consisting of the associated eigenvectors for $\rho(\tau)$.  Write \[ \rho(\sigma) = \mat{a(\sigma)}{b(\sigma)}{c(\sigma)}{d(\sigma)} \] and let $B$ denote the $\T$-module generated by the elements $b(\sigma) \in K$ for $\sigma \in \cG$.  The function 
\[\kappa(\sigma) = \psi^{-1}(\sigma) \overline{b}(\sigma) \in B/IB \]
defines a cocycle representing a surjective class in $H^1(\cG, B/IB(\chi \psi^{-1}))$.  Furthermore, the irreducibility assumption on $\rho$ implies that $B$ is a faithful $\T$-module which in turn implies that $\Fitt_\T(B/IB) \subset I$.  This concludes our sketch of  the proof in the residually distinguished case.

In the residually indistinguishable case, we define another canonical $\T$-module $M$ and a surjective class $\kappa \in H^1(\cG, M(\chi\psi^{-1}))$.  The main effort is to establish the inclusion $\Fitt_{\T}(M) \subset I$.  To prove this inclusion, we apply {\em matrix invariant theory} and {\em rational cohomology}.  Let us describe these ingredients in greater detail.

By definition, $\Fitt_{\T}(M)$ is the ideal generated by the determinants of all square matrices $D$ of relations occurring in a finite presentation of $M$ over $\T$.
In \S\ref{s:ti} we prove that  assumption (\ref{e:redi}) implies that certain expressions involving the traces and determinants of the matrices $\rho(g)$ lie in the ideal $I$.  It therefore suffices to prove that for all $D$, the element $\det(D) \in \T$ can be expressed in terms of these traces and determinants.  While this is possible to show ``by hand" using some combinatorics in small situations (see for instance \S\ref{s:example}), a general proof requires a more conceptual approach.  

We first pass to a ring of formal variables $R$ that is naturally endowed with a homomorphism $\pi \colon R \lra K$ and prove that it suffices to show that certain expressions in $R$ can be expressed modulo $\ker(\pi)$ in terms of traces and determinants of certain matrices taking values in $R$.  The advantage of this is that the ring $R$ is naturally endowed with an action of the algebraic group  $\G = \GL_2$ over $\Z$.  We then apply a theorem of De Concini  and Procesi {stated in Theorem \ref{t:invariants} below and known as the fundamental theorem of matrix invariant theory \cite{cp}*{Theorem 1.10}.  This result identifies} the subring $A \subset R$ generated by traces and determinants as the subring of $R$ invariant under the action of $\G$.

Now, the element  $x \in R$ naturally mapping to $\det(D) \in \T \subset K$ is not invariant under the action of $\G$, but we prove that its image in $R/J$ is invariant for an appropriate $\G$-invariant ideal $J \subset \ker(\pi)$.  In order to apply the theorem of De Concini and Procesi, we must then show that $x$ lies in the image of the natural map
\[ A = R^{\G} \lra (R/J)^{\G}. \]
The cohomology theory of algebraic group actions goes by the name {\em rational cohomology} (or {\em Hochschild cohomology}).  Here ``rational" refers to actions via rational maps rather than the rational numbers---indeed, it is essential that  we work integrally over $\Z$.  Our goal  is to show that the class $\beta \in H^1(\G, J)$ associated to $x$  under the connecting homomorphism $(R/J)^{\G} \lra H^1(\G, J)$ vanishes.

To do this, we first recall an important result in rational cohomology: the restriction to the lower triangular Borel $\B \subset \G$ induces an isomorphism 
\[ H^i(\G, J) \cong H^i(\B, J) \quad \text{ for all } i \ge 0. \]
Hence it suffices to prove the vanishing of the image of $\beta$ in  $H^1(\B, J)$.  For this, we define a $\B$-invariant subideal $J' \subset J$ such that $\beta$ is in the image of the canonical push-forward  
\begin{equation}
\label{e:pf}
 H^1(\B, J') \lra H^1(\B, J). \end{equation}
We prove that  (\ref{e:pf}) is actually the zero map, and hence that $\beta$  vanishes.  We stress that the ideal $J'$ is only endowed with an action of $\B$ (not the full group $\G$), so restriction to the Borel is a crucial part of our argument.  If the algebraic groups $\B$ and $\G$ are replaced by their $\Z$-valued points, then
the restriction map $H^1(\G(\Z), J) \lra H^1(\B(\Z), J)$ on ordinary group cohomology is in general {\em not} an isomorphism (or even injective).  Our use of rational cohomology in place of group cohomology is therefore essential.

The proof of the vanishing of (\ref{e:pf}) is intricate. We construct
a bounded exact complex $\cC^\bullet$ of $\B$-modules ending
in $\cC^0=J'$, together with a morphism to a complex
$\cD^\bullet$ ending in $\cD^0=J$, whose degree-zero map is
the inclusion $J'\hookrightarrow J$. The positive-degree terms
of $\cD^\bullet$ are acyclic $\B$-modules. The construction
uses Koszul complexes and their generalization by Buchsbaum--Rim
to determinantal ideals. Theorems~\ref{t:push} and~\ref{t:comm}
show that this setup forces the push-forward (\ref{e:pf}) to vanish.

In the main text we  prove a stronger version of Theorem~\ref{t:rli} that  relaxes the assumption that $\T$ is reduced and that establishes local conditions for the cohomology class we construct.  See Theorem~\ref{t:rl} below.
It is the local conditions that require the use of the Buchsbaum--Rim complex; for the global picture in Theorem~\ref{t:rli}, the classical Koszul complex suffices.  The reader who wishes to understand the main aspects of our argument in a simplified setting (e.g.\ where we ignore local conditions) is encouraged to consult the announcement \cite{das}.

\bigskip

We would like to thank Rafah  Hajjar Mu\~noz, who visited the first named author at Duke University from UPC (Barcelona) to write an undergraduate senior thesis in the 2021--2022 academic year.  Rafah did computer examples generalizing \S\ref{s:example}, and it was by analyzing the formulas he produced that we realized the role of the Koszul complex in the proof of  the vanishing of $\beta \in H^1(\G, J)$.  We thank 
Jaclyn Lang and Aniruddha Sudarshan for pointing out the discussion of Remark~\ref{remark:top}. 
We would like to thank Frank Calegari, David Loeffler, and Sug Woo Shin for helpful discussions regarding the constructions of \S\ref{s:modforms}, particularly the results in \S\ref{ss:gr}.
We thank Robert Boltje, Brian Conrad, Corrado De Concini, Masato Kurihara, Andreas Nickel, Cristian Popescu, Claudio Procesi, and Geordie Williamson for helpful discussions.   
We thank the anonymous referee for a very careful reading of the manuscript and for suggesting several changes that greatly improved the exposition.

The first named author is supported by NSF grant DMS--2200787. The second named author is supported by DST-SERB grant SB/SJF/2020-21/11, SERB MATRICS grant MTR/2020/000215, SERB SUPRA grant SPR/2019/ 000422, DST FIST program - 2021 [TPN - 700661], and Tata Education and Development Trust. The fourth named author is supported by NSF grant DMS--2201346. 

\bigskip {\em AI Acknowledgement.} After submission of the paper and multiple rounds of human refereeing, we ran the paper through ChatGPT 5.6-Pro to flag any errors.  We made the necessary corrections and passed the updated version through the AI again.  We repeated this process until the AI certified that there were no remaining issues.  The last few passages used ChatGPT 6-Pro, which was released while this was taking place.  A majority of the issues flagged by the AI were typographical or minor mathematical issues with simple local repairs.  However, the AI did catch several important errors as well, such as missing or misapplied hypotheses in lemmas.  In cases where ChatGPT did not itself supply a simple fix simultaneously with identifying the issue, the authors were able to supply the necessary repairs. These corrections were subsequently checked by ChatGPT. Overall, for this paper, ChatGPT played the role of a very thorough referee.

\section{Ritter--Weiss Modules and Brumer--Stark} \label{s:sr}

In this section we  first state precisely the properties satisfied by the generalized Ritter--Weiss modules $\nabla_{\Sigma}^{\Sigma'}(H)$ that we construct in \S\ref{s:constr-rw}.
We then explain how Theorem~\ref{t:main} implies the Brumer--Stark conjecture.  
Next, we will establish some preliminaries and outline the proof of Theorem~\ref{t:main}.

\subsection{Properties of Ritter--Weiss Modules}\label{s:properties-rw}

Throughout \S\ref{s:properties-rw}, $H/F$  is an abelian extension of number fields  with $G = \Gal(H/F)$, and
 $(S,T)$ is an RW pair of finite disjoint sets of places of $F$ (see (\ref{e:rwdef})). In \S\ref{s:constr-rw}, we will construct a $\Z[G]$-module $\nabla_S^T(H)$, the \emph{generalized Ritter--Weiss module}, as the cokernel of a map of finitely generated $\Z[G]$-modules.  In this section we state the main properties of the modules  
 $\nabla_S^T(H)$ which will be proved in \S\ref{s:constr-rw}.
 The module $\nabla_S^T(H)$ is defined for every
RW pair, without any assumptions on the signatures of $F$ or $H$.

For general signatures, the congruences defining
$\cO_{H,S,T}^*$ are imposed at the finite places in $T_H$;
we also require positivity at each real place of $H$ above $T$.
Complex places in $T_H$ impose no additional condition on
$\cO_{H,S,T}^*$. Correspondingly, $\Cl_S^T(H)$ is the ray
class group for the modulus consisting of the finite primes
in $T_H$ together with these real places, depleted at the
finite places in $S_H$.

We will impose additional hypotheses that we now define only where they are needed.
\begin{itemize}
\item We say that an RW pair $(S,T)$ is {\bf complete} if
\[
S\cup T \supset S_{\mathrm{ram}}.
\]
\item
We say that $T$ is {\bf tame} if for every finite place $v\in T$ and $w$ a place of $H$ above $v$, the extension $H_w/F_v$ is tamely ramified (the wild
inertia subgroup of $\Gal(H_w/F_v)$ is trivial).
\item For a prime $p$, we say that $T$ is {\bf $p$-tame} if for every finite
place $v\in T$ and $w$ a place of $H$ above $v$, the extension $H_w/F_v$ is $p$-tamely ramified
(the wild inertia subgroup of $\Gal(H_w/F_v)$ has size coprime to $p$).
\end{itemize}
We remind the reader that $T$ is called torsion-free if $(\cO_{H, S_\infty, T}^*)_{\rm tor} = \{1\}$.
Note that if $(\Sigma, \Sigma')$ is a $p$-admissible RW pair, then
$(\Sigma, \Sigma')$ is complete and $\Sigma'$ is $p$-tame and torsion-free.

\begin{remark}
For orientation: RW pairs are used to define $\nabla_S^T(H)$; complete
RW pairs give the clean class-group exact sequence and the Galois-cohomology
interpretation; torsion-free tame RW pairs give projective presentations;
and $p$-admissible pairs are the complete $p$-local pairs used in the
statement of the main theorem.
\end{remark}

The theorems below often distinguish between results for general extensions $H/F$ where $S$ contains all the infinite places, versus results for CM extensions $H/F$ where $S \cup T$ contains all the infinite places. The former results have largely been proven in \cite{rw} and \cite{dk}, and we reprove them for the sake of completeness.

\begin{theorem} \label{t:nablasequence} {
For any RW pair $(S,T)$, we have:}
\begin{enumerate}
\item The module $\nabla_S^T(H)$ contains $\Cl_S^T(H)$ as a $\Z[G]$-submodule, giving an exact sequence
\begin{equation} \label{e:cnx} 0 \lra \Cl_S^T(H) \lra \nabla_S^T(H) \lra X \lra 0. 
\end{equation}
Under the additional assumption that $(S,T)$ is complete, $X=X_{H,S}$ {as defined in \S\ref{s:rwm}}.\footnote{{If $(S,T)$ is not complete, then $X$ is the kernel of a canonical map
  \[ Y_{H,S} \oplus \bigoplus_{v \in S'_{\ram}} W_v^* \lra \Z,\] where $S'_{\ram}$ is the set of ramified primes not in $S \cup T$ 
  and $W_v^* := \Ind_{G_w}^G W_w^*$ is the dual of the local Ritter--Weiss module defined in \ref{ss:rwc}; see equation (143), the definition of $Z$ preceding it, and (147) in \cite{dk} for details (note we write $Z^\theta$ instead of $X$ in {\em loc.\ cit.}).}
 }
 \item If $F$ is totally real, $H/F$ is an abelian CM extension,
$(S,T)$ is complete, and $S$ contains a place $v$ such that
$c\in G_v$, then the sequence \eqref{e:cnx} remains exact after
tensoring with $\Z[G]_-$:
\[
0\longrightarrow
\Cl_S^T(H)_-
\longrightarrow
\nabla_S^T(H)_-
\longrightarrow
(X_{H,S})_-
\longrightarrow0.
\]
  \end{enumerate}
\end{theorem}

These results are proven in Lemma \ref{l:nablaseq} and Proposition \ref{p:exactness}.

\begin{theorem} \label{thm:quad} {
 Let $(S,T)$ be an RW pair such that $T$ is torsion-free and tame.}
\begin{enumerate}
\item If $S \supset S_\infty$, {there is a map of finitely generated projective $\Z[G]$ modules $V^{\theta} \lra B^{\theta}$ and} an exact sequence
\[ 0 \lra \cO_{H,S,T}^* \lra V^{\theta}\lra B^{\theta} \lra \nabla_S^T(H) \lra 0. \] If in addition $(S,T)$ is complete, {the last three terms in the exact sequence give a locally quadratic presentation of $\nabla_S^T(H)$, i.e.\ the modules $V^{\theta}$ and $B^{\theta}$ are locally free of the same constant rank}.
\item {If $H/F$ is a CM extension and $S \neq \emptyset$,} there is a map of finitely-generated projective $\Z[G]_{-}$-modules $V' \lra (B^{\theta})_{-}$  and an exact sequence
\begin{equation}\label{e:nabla-presentation} 0 \lra (\cO_{H,S,T}^*)^{-}   \lra V'  \lra (B^{\theta})_{-}  \lra \nabla_S^T(H)_{-} \lra 0. \end{equation}
{If in addition $(S,T)$ is complete, the last three terms give a  locally quadratic presentation of $\nabla_S^T(H)_{-}$.}
\end{enumerate}
{If we assume $T$ is only $p$-tame instead of tame, analogous results hold after tensoring over $\Z$ with  $\Z_{(p)}$.}
\end{theorem}

{These results are proven in  Corollary \ref{c:uvbn} (which leads to (\ref{e:fourterm}) and (\ref{e:vp})), Proposition \ref{prop:quad-res}, and Proposition~\ref{prop:quad-resCM}. The notation $V^\theta$ and $B^\theta$ will become clear in \S\ref{s:qp} when these modules are defined as the kernels of maps denoted $\theta$ whose codomains are certain modules $V$ and $B$, respectively.}

{\begin{remark} \label{r:neeqp}
If $(\Sigma, \Sigma')$ is a $p$-admissible RW pair, then all the conditions of Theorem~\ref{thm:quad} (2) hold, and hence
$\nabla_{\Sigma}^{\Sigma'}(H)_{p,-}$ is a quadratically presented $\Z_p[G]_-$-module.  Indeed, since $\Z_p[G]_-$ is a semi-local ring,  the fact that $V'$ and 
$(B^{\theta})_{-}$ are locally free of the same constant rank implies that the $p$-adic completions $(V')_p$ and $(B^{\theta})_{p,-}$ are free of the same rank over 
$\Z_p[G]_-$.
\end{remark}
}

{
\begin{remark} \label{r:quad-res} In the appendix to \cite{etnc}, it is proven using the results of this paper that in Case 2 of Theorem \ref{thm:quad}, $\nabla^T_S(H)_-$ is in fact quadratically presented over $\Z[G]_-$ (not just locally).  This result is not used in this paper.
\end{remark}
}

We will prove that $\nabla_S^T(H)$ has the following relationship to Galois cohomology {in Theorem \ref{t:cocycle} below.}

\begin{theorem}\label{t:cocycle-nabla}
{Let $(S,T)$ be a complete RW pair with $S\neq\emptyset$,
and fix $v_0\in S$. Let $N$ be a prodiscrete $\Z[G]$-module, regarded as a
$G_F$-module by inflation.
Assume either that $T$ is tame, or that $T$ is
$p$-tame and $N$ is a $p$-group.}
Suppose we are given a {continuous} 1-cocycle $\kappa \colon G_F \lra N$ and $x_v \in N$ for $v \in S \setminus \{v_0\}$ such that 

\begin{enumerate}
\item $\kappa(G_{F_{v_0}}) = 0$
\item $\kappa|_{G_{F_v}}$ is unramified for $v \notin (S \cup T)$
\item $\kappa|_{G_{F_v}}$ is tamely ramified for finite $v \in T$
\item for $v \in S \setminus \{v_0\}$, $\kappa(\sigma_v) = (\sigma_v - 1)x_v$ for all $\sigma_v \in G_{F_v}$.
\end{enumerate}
Then there exists a $G$-equivariant map $\nabla_S^T(H) \lra N$ whose image is the $\Z[G]$-submodule generated by $\kappa( G_{F})$ and the $x_v$.
\end{theorem}

\subsection{Brumer--Stark} \label{s:bs}

In \S\ref{s:bs}--\S\ref{s:summary} we return to the
Brumer--Stark setting: $F$ is totally real and $H/F$ is a finite
abelian CM extension.  For simplicity, we set $S=S_\infty\cup S_{\ram}$, as it is well-known that the result for general $S \supset S_\infty \cup S_{\ram}$ follows from this setting.
The set $T$ is a finite torsion-free set of finite primes disjoint from $S$.
In particular, every prime in $T$ is unramified in $H/F$.

As mentioned in the introduction, we will prove in a forthcoming paper \cite{etnc} that Theorem~\ref{t:main} implies the Equivariant Tamagawa Number Conjecture for the minus part of the Tate motive associated to $H/F$.  This is a powerful statement with many corollaries, including the Brumer--Stark conjecture.  Because of its singular importance, however, we would like to prove here that the Brumer--Stark conjecture follows directly from Theorem~\ref{t:main}.

We begin by introducing a slight variant of the ray class group that is adapted to the
anti-unit formulation of the Brumer--Stark conjecture.  Let $I^T(H)$ denote
the group of fractional ideals of $H$ relatively prime to $T_H$, and as usual write
\[
 I^T(H)^-=\ker(1+c:I^T(H)\longrightarrow I^T(H)).
\]
Let
\[
 H_T^{\times,-}
 =
 \{u\in H^\times:u\equiv 1\pmod{\fq}\text{ for all }\fq\in T_H,
 \quad u^{1+c}=1\}.
\]
Since the ideal generated by every element of $H_T^{\times,-}$ belongs to
$I^T(H)^-$, we may define
\begin{equation}\label{e:clau}
 \Cl_{\mathrm{au}}^T(H)
 =
 I^T(H)^-/\{(u):u\in H_T^{\times,-}\},
\end{equation}
where the subscript au stands for anti-unit.  This is naturally a $\Z[G]_-$-module.

\begin{lemma}\label{l:clau}
The natural map
\[
\rho\colon \Cl_{\mathrm{au}}^T(H)\longrightarrow \Cl^T(H)^-
\]
has kernel and cokernel that are finite and annihilated by $2$. In particular,
$\Cl_{\mathrm{au}}^T(H)$ is finite, and $\rho$ becomes an
isomorphism after inverting $2$.
\end{lemma}

\begin{proof}
An ideal $\fa\in I^T(H)^-$ maps to an anti-invariant ray class, so
the natural map $I^T(H)\to\Cl^T(H)$ induces $\rho$.

Suppose that a class in $\ker(\rho)$ is represented by
$\fa=(x)$, where $x\equiv1\pmod T$. Since
$\fa^{1+c}=(1)$, we have $xx^c\in\cO_{H,S_\infty,T}^*.$
Torsion-freeness of $T$ implies that $c$ acts trivially on
$\cO_{H,S_\infty,T}^*$. Hence
\[
[\fa]\longmapsto
xx^c \pmod{(\cO_{H,S_\infty,T}^*)^2}
\]
defines an injection
\[
\ker(\rho)\hookrightarrow
\cO_{H,S_\infty,T}^*/(\cO_{H,S_\infty,T}^*)^2.
\]
Indeed, if $xx^c=\epsilon^2$, then $u=x/\epsilon$ satisfies
$u\equiv1\pmod T$, $u^{1+c}=1$, and $(u)=\fa$. Thus
$\ker(\rho)$ is annihilated by $2$.

Finally, if $x\in\Cl^T(H)^-$ is represented by $\fa\in I^T(H)$,
then
\[
\fa^{1-c}\in I^T(H)^-
\]
and its image under $\rho$ is
\[
(1-c)x=2x.
\]
Therefore $\coker(\rho)$ is also annihilated by $2$. The remaining
assertions follow.
\end{proof}

Our goal in this subsection is to prove the following result assuming Theorem~\ref{t:main}, and then to deduce the Brumer--Stark conjecture (with an improved power of 2) from this result.
\begin{theorem}\label{t:sbs}
Assume Theorem~\ref{t:main}.  With $(S,T)$ as above, we have
\begin{equation}\label{e:sbs-au}
 \frac{\Theta_{S,T}}{2^{n-1}}
 \in
 \Ann_{\Z[G]_-}\bigl(\Cl_{\mathrm{au}}^T(H)\bigr).
\end{equation}
\end{theorem}

Theorem~\ref{t:sbs} will be proved by relating $\Cl_{\mathrm{au}}^T(H)$ to the transpose of a certain generalized
Ritter--Weiss module. Fix a real place $w$ of $F$ and put
\[
 T^\dagger=T\cup(S_\infty\setminus\{w\}).
\]
Let $\Sel_{\{w\}}^{T^\dagger}(H)_{-}$ denote the transpose of $\nabla_{\{w\}}^{T^\dagger}(H)_{-}$ over $\Z[G]_{-}$ associated to the presentation (\ref{e:nabla-presentation}) defining this latter $\Z[G]_{-}$-module, as described in Lemma~\ref{l:transprop}.

\begin{lemma}\label{l:causel}
Fix a real place $w$ of $F$, and put
\[
T^\dagger=T\cup(S_\infty\setminus\{w\}).
\]
There is a natural surjective homomorphism of $\Z[G]_-$-modules
\[
\Sel_{\{w\}}^{T^\dagger}(H)_-
\longtwoheadrightarrow
\Cl_{\mathrm{au}}^T(H)^\vee.
\]
\end{lemma}

We defer the proof of Lemma~\ref{l:causel} until \S\ref{ss:changing-sets}.

\begin{lemma} \label{l:tas}
Assume Theorem~\ref{t:main}.  Let $(S,T)$ be as above. 
We have \[
 \left(\frac{\Theta_{S,T}}{2^{n-1}}\right)^\#
 \in
 \Ann_{\Z_2[G]_-}
 \left(\Sel_{\{w\}}^{T^\dagger}(H)_{2,-}\right).
\]
\end{lemma}
\begin{proof}
Define
\[ \Sigma = \{w\} \cup  \{ v \in S_{\ram} \colon v \mid 2  \}, \qquad \Sigma' = T \cup \{ v \in S_{\ram}, v \nmid 2\}  \cup S_\infty \setminus \{w\}.
\]
Then $(\Sigma,\Sigma')$ is $2$-admissible: it is an RW pair
by construction, $\Sigma\neq\emptyset$ since $w\in\Sigma$,
the finite ramified primes are distributed as required, and
$\Sigma'$ is torsion-free because it contains the torsion-free
set $T$. Therefore Theorem~\ref{t:main} applies:
\[ \Fitt_{\Z_2[G]_-}(\nabla_{\Sigma}^{\Sigma'}(H)_{2,-}) = (\Theta_{\Sigma, \Sigma'}/2^{n-1}). \]
Let $\Sel_{\Sigma}^{\Sigma'}(H)_{2,-}$ denote the transpose of $\nabla_{\Sigma}^{\Sigma'}(H)_{2,-}$ over $\Z_2[G]_{-}$ associated to the presentation (\ref{e:nabla-presentation}).

Let $s = |\{ v \in S_{\ram}, v \nmid 2\}|$.  The argument of \cite{dk}*{Theorem 3.7} has a projective-presentation analogue in our setting, which yields
\begin{align*}
 \Fitt_{\Z_2[G]_-}(\Sel_{\Sigma}^{T^\dagger}(H)_{2,-}) &= \ \Fitt_{\Z_2[G]_-}^s(\nabla_{\Sigma}^{T^\dagger}(H)_{2,-})^\# \\
&= \ (\Theta_{\Sigma, T}/2^{n-1})^\# \prod_{v \in S_{\ram}, \ v \nmid 2} \!\!\!\!\! (\N I_v, 1 - \sigma_v e_v).
 \end{align*}
Here $\Fitt^s_{\Z_2[G]_-}$ denotes the $s$th Fitting ideal over the ring $\Z_2[G]_-$,
\[ \N I_v = \sum_{\tau \in I_v} \tau, \qquad e_v = \frac{1}{|I_v|} \N I_v \in \Q[G],\]
and $\sigma_v$ denotes a Frobenius lift modulo $I_v$.

  Let $s' = |S_{\ram}|$.
 The arguments of \cite{dk}*{Appendix B}  yield
 \begin{align*}
 \Fitt_{\Z_2[G]_-}(\Sel_{\{w\}}^{T^\dagger}(H)_{2,-}) &= \ \Fitt_{\Z_2[G]_-}^{s'}(\nabla_{\{w\}}^{T^\dagger}(H)_{2,-})^\# \\
&= \ (\Theta_{\emptyset, T}/2^{n-1})^\# \prod_{v \in S_{\ram}} \!\! (\N I_v, 1 - \sigma_v e_v) \\
&\supset (\Theta_{S,T}/2^{n-1})^\#. 
 \end{align*}
 It follows that $\Theta_{S,T}^\#/2^{n-1}$ annihilates $\Sel_{\{w\}}^{T^\dagger}(H)_{2,-}$.
\end{proof}

\begin{proof}[Proof of Theorem~\ref{t:sbs}]
Lemma~\ref{l:clau} yields  an isomorphism
\[
 \Cl_{\mathrm{au}}^T(H)\otimes\Z[1/2]
 \cong
 \Cl^T(H)^-\otimes\Z[1/2].
\]
Thus the prime-to-$2$ part of (\ref{e:sbs-au}) follows from the Strong
Brumer--Stark theorem proved in \cite{dk}:
 \[  \Theta^\#_{S,T} \in \Fitt_{\Z[1/2][G]_-}(\Cl^T(H)^\vee \otimes_{\Z[G]} \Z[1/2][G]_-).
 \]

It remains to work over $\Z_2$. 
Lemmas~\ref{l:causel} and~\ref{l:tas} imply that
$\Theta_{S,T}^{\#}/2^{n-1}$ annihilates
$(\Cl_{\mathrm{au}}^T(H)_2)^\vee$.
By Pontryagin duality, $\Theta_{S,T}/2^{n-1}$ annihilates
$\Cl_{\mathrm{au}}^T(H)_2$, as required.
\end{proof}

The Brumer--Stark conjecture now follows immediately; in fact we again
obtain the stronger statement with $\Theta_{S,T}$ replaced by
$\Theta_{S,T}/2^{n-2}$.

\begin{proof}[Proof of Theorem~\ref{t:bs}]
Assume first that $n\geq2$.  Let $\fP$ be as in the statement and put
\[
 \alpha=\frac{\Theta_{S,T}}{2^{n-1}}.
\]
The anti-invariant ideal $\fP^{1-c}$
defines an element of $\Cl_{\mathrm{au}}^T(H)$.  By
Theorem~\ref{t:sbs},
\[
 \fP^{\alpha(1-c)}=(u)
\]
for some $u\in H^\times$ satisfying
\[
 u\equiv1\pmod T,\qquad u^{1+c}=1.
\]
Since $\alpha c=-\alpha$, we have
\[
 \alpha(1-c)=2\alpha
 =\frac{\Theta_{S,T}}{2^{n-2}}, 
\qquad \text{hence}
 (u) = \fP^{\frac{\Theta_{S,T}}{2^{n-2}}}.
\]
Thus $u$ is a $\fp$-unit, it is congruent to $1$ modulo $T$, and
$u^{1+c}=1$ implies that $u$ has absolute value $1$ at every complex
embedding.  The equivariant valuation identity above is precisely
(\ref{e:bs}) with $\Theta_{S,T}$ replaced by
$\Theta_{S,T}/2^{n-2}$.

The result for $n=1$ (i.e. $F=\Q$) is Stickelberger's classical theorem.   In \cite{etnc} we show how our results here can be used to prove that the Brumer--Stark conjecture holds with $\Theta_{S,T}$ replaced by $\Theta_{S,T}/2^{n-1}$ when $H/F$ is ramified at any finite prime.  In particular, we recover the (known) result for $F=\Q$.
\end{proof}

\subsection{Minimal Sets} \label{s:minimal}

{
Fix a prime $p$. In this subsection we record how the Ritter--Weiss
modules $\nabla$ and Stickelberger elements $\Theta$, and hence the
statement of Theorem~\ref{t:main}, change as a $p$-admissible pair
$(\Sigma,\Sigma')$ varies. This will allow us, in the modular form
constructions below, to reduce to certain pairs with $\Sigma$ minimal. The
needed comparison result, proved in \S\ref{ss:changing-sets}, is the following.}

\begin{lemma} \label{l:addprimes} {Let $(\Sigma,\Sigma')$ be $p$-admissible.} We  have the following:
\begin{itemize} \item Let $w \in \Sigma'$ be a real place.  Then \[ \Fitt_{\Z_p[G]_-}(\nabla_{\Sigma \cup \{w\}}^{\Sigma' \setminus{w}}(H)_{p,-}) = 
2 \Fitt_{\Z_p[G]_-}(\nabla_{\Sigma}^{\Sigma'}(H)_{p,-}). \]
\item Let $\fl \not \in \Sigma \cup \Sigma'$, so by the assumptions, $\fl$ is a finite prime unramified in $H$. Put $\sigma_{\fl}$ for the associated Frobenius element. Then
\[  \Fitt_{\Z_p[G]_-}(\nabla_{\Sigma \cup \{\fl \}}^{\Sigma'}(H)_{p,-}) = 
(1 - \sigma_{\fl}^{-1}) \Fitt_{\Z_p[G]_-}(\nabla_{\Sigma}^{\Sigma'}(H)_{p,-})
\]
and 
\[ \Theta_{\Sigma \cup \{\fl\}, \Sigma'} = (1 - \sigma_{\fl}^{-1})\Theta_{\Sigma, \Sigma'}. \]
Similarly we have \[  \Fitt_{\Z_p[G]_-}(\nabla_{\Sigma}^{\Sigma' \cup \{\fl \}}(H)_{p,-}) = 
(1 - \sigma_{\fl}^{-1} \N\fl) \Fitt_{\Z_p[G]_-}(\nabla_{\Sigma}^{\Sigma'}(H)_{p,-})
\]
and 
\[ \Theta_{\Sigma, \Sigma' \cup \{\fl \}} = (1 - \sigma_{\fl}^{-1}\N\fl)\Theta_{\Sigma, \Sigma'}. \]
\end{itemize}
\end{lemma}

{Thus Theorem~\ref{t:main} for a $p$-admissible pair follows from Theorem~\ref{t:main} for
any other $p$-admissible pair obtained by deleting unnecessary finite
unramified primes and moving real places between $\Sigma$ and
$\Sigma'$, provided $\Sigma$ remains nonempty. We now choose minimal
pairs for this reduction.

As in \cite{dk}*{(8.2.1) and (8.2.2)}, it is convenient to establish two cases that will appear throughout the paper.
\begin{itemize}\item {\bf Case 1.} There are no primes above $p$  ramified in $H/F$.
\item {\bf Case 2.} There exists a prime above $p$  ramified in $H/F$.
\end{itemize}

\begin{lemma}\label{l:minimal} Suppose that the statement of Theorem~\ref{t:main} holds for the following specific $p$-admissible pairs 
$(\Sigma, \Sigma')$:
\begin{itemize}\item {\bf Case 1.}  \[  \Sigma = \{\fl\},  \quad \Sigma' = T \cup S_\infty \cup S_{\ram}, \] where $\fl\notin\Sigma'$ is one fixed prime of $F$ whose
Frobenius $\sigma_\fl\in G$ is the complex conjugation $c$.
\item {\bf Case 2.}  \[ \Sigma = \{v \in S_{\ram} \colon v \mid p \}, 
 \Sigma' = T \cup S_\infty \cup \{v \in S_{\ram} \colon v \nmid p \}. \]
 \end{itemize}
 Then Theorem~\ref{t:main} holds for all $p$-admissible pairs  $(\Sigma, \Sigma')$ with $\Sigma' \supset T$.
\end{lemma}

\begin{proof}
In Case 1, any $p$-admissible pair  with minimal $\Sigma$ set and $\Sigma'$ set containing $T$ will either have the form
  \begin{align*}
 \Sigma_1 &= \{\fl'\} , \\
 \Sigma'_1 &= T \cup S_\infty \cup S_{\ram}
 \end{align*}
for a finite prime $\fl'$ or have the form
  \begin{align*}
 \Sigma_2 &= \{v\} , \\
 \Sigma'_2 &= T \cup (S_\infty \setminus \{v\}) \cup S_{\ram}
 \end{align*}
for a real place $v$. If $\fl'=\fl$, the pair
$(\Sigma_1,\Sigma'_1)$ is the assumed pair and no comparison
is needed. Otherwise, with $(\Sigma,\Sigma')$ the fixed pair
in Case 1, Lemma~\ref{l:addprimes} yields
\begin{align*} \Fitt_{\Z_p[G]_-}(\nabla_{ \{\fl,\fl'\}}^{\Sigma'_1}(H)_{p,-}) &=
 (1 - \sigma_\fl^{-1})  \Fitt_{\Z_p[G]_-}(\nabla_{\Sigma_1}^{\Sigma'_1}(H)_{p,-}) \\
  &= 2  \Fitt_{\Z_p[G]_-}(\nabla_{\Sigma_1}^{\Sigma'_1}(H)_{p,-}).
\end{align*}
The second equality above follows since $\sigma_\fl$ equals complex conjugation, which acts as $-1$ on the minus side.
Therefore
\begin{align} \Fitt_{\Z_p[G]_-}(\nabla_{\Sigma_1}^{\Sigma'_1}(H)_{p,-}) &=  \frac{1}{2}\Fitt_{\Z_p[G]_-}(\nabla_{ \{\fl,\fl'\}}^{\Sigma'_1}(H)_{p,-}) \nonumber \\
  &= \frac{1 - \sigma_{\fl'}^{-1}}{2}  \Fitt_{\Z_p[G]_-}(\nabla_{\Sigma}^{\Sigma'}(H)_{p,-})
\label{e:nssp}
\end{align}
 by another application of Lemma~\ref{l:addprimes}.
A similar calculation shows
 \begin{equation} \label{e:tssp}
 \Theta_{\Sigma_1, \Sigma'_1} = \frac{1}{2} \Theta_{ \{\fl',\fl\}, \Sigma_1'} = \frac{1 - \sigma_{\fl'}^{-1}}{2} \Theta_{\Sigma, \Sigma'}. \end{equation}
Comparing (\ref{e:nssp}) and (\ref{e:tssp}), it follows that the statement of Theorem~\ref{t:main} for $(\Sigma, \Sigma')$ implies the result for $(\Sigma_1, \Sigma_1')$.

A similar calculation holds for $(\Sigma_2, \Sigma_2')$: by Lemma~\ref{l:addprimes}, moving an infinite place from $\Sigma'$ to $\Sigma$ multiplies the Fitting ideal of $\nabla$ by 2.  We leave details to the reader.  Therefore, in Case 1, proving Theorem~\ref{t:main} for the particular pair $(\Sigma, \Sigma')$ implies the result for all $p$-admissible pairs with $\Sigma' \supset T$.

The argument in Case 2 is simpler since there is a unique minimal $\Sigma$ as $(\Sigma,\Sigma')$ ranges over $p$-admissible pairs.
The result follows directly from Lemma~\ref{l:addprimes}, since $\Fitt_{\Z_p[G]_-}(\nabla_\Sigma^{\Sigma'}(H)_{p,-})$ and
$(\Theta_{\Sigma,\Sigma'}/2^t)$, where $t=|S_\infty\cap\Sigma'|$, have the same transformation laws upon adding finite unramified primes or moving real places.
\end{proof}
}

For the proof of Theorem~\ref{t:main}, it suffices to consider
$p$-admissible pairs $(\Sigma,\Sigma')$ for which $\Sigma'$
contains a finite torsion-free set $T$ of primes unramified
in $H/F$ and not lying above $p$.

Indeed, given any $p$-admissible pair $(\Sigma,\Sigma')$,
choose a finite torsion-free set $T_0$ of primes unramified
in $H/F$, disjoint from $\Sigma\cup\Sigma'$ and from the
primes above $p$. Such a set can be chosen outside any
prescribed finite set of primes. The enlarged pair
$(\Sigma,\Sigma'\cup T_0)$ is again $p$-admissible.
Put
\[
E_{T_0}=\prod_{v\in T_0}(1-\N v\,\sigma_v^{-1})
\in\Z_p[G]_-,
\qquad
t=|S_\infty\cap\Sigma'|.
\]
By
Lemma~\ref{l:addprimes} and the corresponding Euler-factor
identity, we have
\begin{align*}
\Fitt_{\Z_p[G]_-}
\bigl(\nabla_\Sigma^{\Sigma'\cup T_0}(H)_{p,-}\bigr)
&=
E_{T_0}\,
\Fitt_{\Z_p[G]_-}
\bigl(\nabla_\Sigma^{\Sigma'}(H)_{p,-}\bigr),\\
\frac{\Theta_{\Sigma,\Sigma'\cup T_0}}{2^t}
&=
E_{T_0}\frac{\Theta_{\Sigma,\Sigma'}}{2^t}.
\end{align*}
Each factor of $E_{T_0}$ is a nonzerodivisor in
$\Z_p[G]_-$: its value at every odd character is nonzero,
since $\N v>1$. Cancelling $E_{T_0}$, the equality in
Theorem~\ref{t:main} for the enlarged pair implies the
equality for the original pair.

\subsection{Character Group Rings} \label{s:cgr}

Let $\cO$ denote the ring of integers of a finite extension of $\Q_p$ containing the values of all characters $\chi \in \hat{G}$.  There is a $\Z_p$-algebra injection
\[ \Z_p[G] \longrightarrow \prod_{\chi \in \hat{G}} \cO, \qquad x \mapsto (\chi(x)).
\]
In our constructions with modular forms, we would like to work with local rings (as opposed to the semilocal ring $\Z_p[G]_-$).  Furthermore, it will be convenient if the Stickelberger element $\Theta_{\Sigma, \Sigma'}$ is a nonzerodivisor in our ring.  To this end, for every $\Gal(\overline{\Q}_p/\Q_p)$-stable set of characters $\Psi \subset \hat{G}$ we define the associated {\em character group ring}
\[ R_\Psi = \text{Image}( \Z_p[G] \longrightarrow \prod_{\chi \in \Psi} \cO), \qquad x \mapsto (\chi(x)). \]
{We define 
\[ 
\Psi = \{ \chi \in \hat{G} \colon \chi \text{ is odd and } \chi(G_v) \neq 1 \text{ for all } v \in \Sigma \}. 
\]
The image of the Stickelberger element  $\Theta_{\Sigma, \Sigma'}$ in $R_\Psi$ is a nonzerodivisor, since if $\chi$ is an odd character, $L_{\Sigma, \Sigma'}(\chi, 0) = 0$ if and only if there exists $v \in \Sigma$ such that $\chi(G_v) =1$.}   If $M$ is any $\Z_p[G]$-module, we write $M_\Psi = M \otimes_{\Z_p[G]} R_\Psi$.
\bigskip

Write $G = G_p \times G'$ where $G_p$ is a $p$-group and $G'$ has order relatively prime to $p$.  For each $\Gal(\overline{\Q}_p/\Q_p)$-conjugacy class of characters $\Phi$ of $G'$, we define

\begin{equation} \label{e:phi0def}
{\Phi_0 = \{ \chi \in \Psi \colon \chi|_{G'} \in \Phi\}.}
 \end{equation}
Then, for $\Phi_0\neq\emptyset$, the ring $R_{\Phi_0}$ is
local and $\Theta_{\Sigma,\Sigma'}$ is a nonzerodivisor in it.
We have a decomposition
 \begin{equation} \label{e:proddecomp}
 {R_\Psi \cong \prod_{\Phi} R_{\Phi_0},} 
  \end{equation}
where $\Phi$ runs through the
$\Gal(\overline{\Q}_p/\Q_p)$-conjugacy classes of characters
of $G'$ for which $\Phi_0\neq\emptyset$.
All subsequent constructions with the local rings
$R_{\Phi_0}$ are restricted to these nonempty components.
 
 The following result was proven in \cite{dk}*{Lemmas 2.4 and 2.5}.
 
 \begin{lemma} \label{l:size}  Let $M$ be a quadratically presented module over a character group ring $R_\Psi$ such that $\Fitt_{R_\Psi}(M) = (x)$.  Then $M$ is finite if and only if $x$ is a nonzerodivisor, and in this case we have
 \[ |M| = | R_\Psi/(x)| =  \Big|\Z_p \ / \prod_{\psi \in \Psi} \psi(x)\Big|. \]
 \end{lemma}
 {Note that the product $\prod_{\psi \in \Psi} \psi(x)$ is $\Z_p$-valued as it is Galois invariant.}
 
\subsection{Inclusion Implies Equality} \label{ss:ie}
{
We remain in the Brumer--Stark setting fixed in \S\ref{s:bs}. Fix a prime $p$, and
put
\begin{equation} \label{e:bestsigma}
\Sigma=S_\infty\cup \{v\in S_{\mathrm{ram}}:v\mid p\},\qquad
\Sigma'=T\cup \{v\in S_{\mathrm{ram}}:v\nmid p\}.
\end{equation}
This pair is $p$-admissible by the standing assumptions of \S\ref{s:bs}.}
Our goal in this section is to prove:
\begin{theorem} \label{t:incl}
{
Suppose that for every
$\Gal(\overline{\mathbf Q}_p/\mathbf Q_p)$-conjugacy class $\Phi$ of characters of $G'$, we have}
\begin{equation} \label{e:incl}
 \Fitt_{R_{\Phi_0}}(\nabla_{\Sigma}^{\Sigma'}(H)_{\Phi_0}) \subset (\Theta_{\Sigma, \Sigma'}), \end{equation}
 where $\Phi_0$ is as in (\ref{e:phi0def}).
 Moreover, suppose that the same inclusion holds for each CM subfield  $H' \subset H$ containing $F$, with $\Sigma$ and $\Sigma'$ defined by (\ref{e:bestsigma}) for the extension $H'/F$.
 
Then every such inclusion is actually an equality, and 
\begin{equation} \label{e:equality}
 \Fitt_{\Z_p[G]_-}(\nabla_{\Sigma}^{\Sigma'}(H)_{p,-}) = (\Theta_{\Sigma, \Sigma'}). \end{equation}
\end{theorem}

 {
 In \S\ref{s:minimal} we defined certain minimal sets $(\Sigma, \Sigma')$.  
By the proof of Lemma~\ref{l:minimal}---simply replace the relevant equalities by inclusions---the inclusion
\begin{equation} \label{e:minincl}
 \Fitt_{R_{\Phi_0}}(\nabla_{\Sigma}^{\Sigma'}(H)_{\Phi_0}) \subset (\Theta_{\Sigma, \Sigma'}/2^t) \end{equation}
for the minimal pairs defined in \S\ref{s:minimal} implies the inclusion (\ref{e:incl}) for the pairs defined in (\ref{e:bestsigma}).
By Theorem~\ref{t:incl}, this implies the equality (\ref{e:equality}).
Now the process can be reversed, and (\ref{e:equality}) implies that
\[ \Fitt_{\Z_p[G]_-}(\nabla_{\Sigma}^{\Sigma'}(H)_{p,-}) = (\Theta_{\Sigma, \Sigma'}/2^t) \]
for the minimal pairs $(\Sigma, \Sigma')$.  By Lemma~\ref{l:minimal}, this implies 
Theorem~\ref{t:main} for all $p$-admissible pairs. In summary, the following is a corollary of  Theorem~\ref{t:incl}.

\begin{corollary} \label{c:incl}
In order to prove Theorem~\ref{t:main}, it suffices to prove the inclusion (\ref{e:minincl}) for the minimal pairs $(\Sigma, \Sigma')$ chosen in \S\ref{s:minimal}, for $H$ and all of its CM subfields containing $F$.
\end{corollary}
}

{Recall} \[ \Psi = \{ \chi \in \hat{G} \colon \chi \text{ is odd and } \chi(G_v) \neq 1 \text{ for all } v \in \Sigma\} \] and let $R_\Psi$ be the associated character group ring.  In view of the product decomposition (\ref{e:proddecomp}), the following lemma implies that the first assertion of Theorem~\ref{t:incl} implies the second.

\begin{lemma}  \label{l:zeroproj}
With notation as above, if
\[ \Fitt_{R_\Psi}(\nabla_{\Sigma}^{\Sigma'}(H)_{\Psi}) = (\Theta_{\Sigma, \Sigma'}), \]
then 
\[  \Fitt_{\Z_p[G]_-}(\nabla_{\Sigma}^{\Sigma'}(H)_{p,-}) = (\Theta_{\Sigma, \Sigma'}). \]
\end{lemma}

 \begin{proof} Lemma~\ref{l:zeroproj} is proven in \cite{dk}*{Lemma 7.1}.  
 The fact that $p$ is always assumed odd in that paper is irrelevant to the proof, but for completeness we reproduce the proof here.
 The key point is that if $\overline{\Psi}$ denotes the set of odd characters of $G$ not lying in $\Psi$, then
\[  \Fitt_{R_{\overline{\Psi}}}(\nabla_{\Sigma}^{\Sigma'}(H)_{\overline{\Psi}}) =  0 = \Theta_{\Sigma, \Sigma'} R_{\overline{\Psi}}.\]
Consider the canonical injection
\[
\iota : \Z_p[G]_- \longrightarrow R_{\Psi} \times R_{\overline{\Psi}}, \qquad \qquad \text{denoted } \iota(x) = (\iota_1(x), \iota_2(x)).
\]
By {Remark~\ref{r:neeqp},} the $\Z_p[G]_-$-module $\nabla_{\Sigma}^{\Sigma'}(H)_{p,-}$ is quadratically presented, so we can write \[ \Fitt_{\Z_p[G]_-}(\nabla_{\Sigma}^{\Sigma'}(H)_{p,-}) = (x) \] for some $x \in \Z_p[G]_-$. Then we have
\[
\iota_2(x) = 0 = \iota_2(\Theta_{\Sigma, \Sigma'}).
\]
The given equality implies that there exists $y \in R_{\Psi}^*$ such that \[ \iota_1(x) = y \cdot \iota_1(\Theta_{\Sigma, \Sigma'}),\] {as $\iota_1(x)$ and $\iota_1(\Theta_{\Sigma, \Sigma'})$ are nonzerodivisors by the definition of $\Psi$.} 
{As \[ 
\Z_p[G]_- \lra R_{\Psi}  \] is a surjective map of semilocal rings,} $y$ has a lift $\widetilde{y}$ in $\Z_p[G]_-^*$. We then have that $x - \widetilde{y} \Theta_{\Sigma, \Sigma'}$ vanishes under both $\iota_1$ and $\iota_2$. It follows that $x = \widetilde{y} \Theta_{\Sigma, \Sigma'}$, giving the desired equality. 
\end{proof}

Theorem~\ref{t:incl} was proven for $p$ odd in \cite{dk}, {so for the remainder of \S\ref{ss:ie} we will assume that $p = 2$
and prove the first assertion of Theorem~\ref{t:incl}. 

By  Remark~\ref{r:neeqp},} the module $\nabla_{\Sigma}^{\Sigma'}(H)_{\Psi}$ is  quadratically presented and therefore its Fitting ideal is principal. Write \[  \Fitt_{R_\Psi}(\nabla_{\Sigma}^{\Sigma'}(H)_{\Psi}) = (x). \]  Since $R_\Psi = \prod_{\Phi} R_{\Phi_0}$,  (\ref{e:incl}) implies that $(x) \subset (\Theta_{\Sigma, \Sigma'})$ in $R_\Psi$ and hence it suffices to prove 
that 
\[ |R_\Psi/(x)| = |R_\Psi/(\Theta_{\Sigma, \Sigma'})|. \]
 Now it follows from Lemma~\ref{l:size} that
\[ |\nabla_{\Sigma}^{\Sigma'}(H)_{\Psi}| =  |R_\Psi/(x)| = \Big|\Z_2 \ \Big/ \prod_{\alpha \in \Psi} \alpha(x)\Big|\] and 
\[  |R_\Psi/(\Theta_{\Sigma, \Sigma'})| = \Big|\Z_2 \ \Big/ \prod_{\alpha \in \Psi} L_{\Sigma, \Sigma'}(\alpha^{-1}, 0) \Big| = 2\text{-part of } \prod_{\alpha \in \Psi} L_{\Sigma, \Sigma'}(\alpha^{-1}, 0). \]
{
Note here that the product on the right side is an algebraic integer and invariant under $\Gal(\overline{\Q}_2/\Q_2)$, and hence lies in $\Z_2$.}
We must therefore show
\begin{equation} \label{e:nsize}
 |\nabla_{\Sigma}^{\Sigma'}(H)_{\Psi}| \doteq \prod_{\alpha \in \Psi} L_{\Sigma, \Sigma'}(\alpha^{-1}, 0),
\end{equation}
where $\doteq$ indicates equality up to a $p$-adic unit.
This will be achieved through the analytic class number formula. We prove (\ref{e:nsize}) by partitioning $\Psi$ into $\Gal(\overline{\Q}_2/\Q_2)$ conjugacy classes.  For such a conjugacy class {$\Omega$}, we have
\[ | \nabla_{\Sigma}^{\Sigma'}(H)_{\Omega}| \doteq \prod_{\alpha \in \Omega} \alpha(x). \]
{
In particular, we have
\[ |\nabla_{\Sigma}^{\Sigma'}(H)_{\Psi}| = \prod_{\Omega} | \nabla_{\Sigma}^{\Sigma'}(H)_{\Omega}|, \]
so (\ref{e:nsize}) will follow from the following result whose proof will take up the remainder of \S\ref{ss:ie}.}
\begin{theorem} \label{t:nsize2}
With notation as above, under the assumption of Theorem~\ref{t:incl}, we have
\begin{equation}
  | \nabla_{\Sigma}^{\Sigma'}(H)_{\Omega}| \doteq \prod_{\alpha \in \Omega} L_{\Sigma, \Sigma'}(\alpha^{-1}, 0). \label{e:phiprod}
  \end{equation}
\end{theorem}
Taking the product of (\ref{e:phiprod}) over all possible  $\Gal(\overline{\Q}_2/\Q_2)$ conjugacy classes $\Omega \subset \Psi$ yields (\ref{e:nsize}) and completes the proof of Theorem~\ref{t:incl}.  In the remainder of the section, we prove Theorem~\ref{t:nsize2}.
Let $H_{\Omega}$ denote the fixed field in $H$ of the kernel of any character in $\Omega$. 
 
 \begin{lemma} \label{l:descent} We have 
\begin{equation} 
  \nabla_{\Sigma}^{\Sigma'}(H)_{\Omega} \cong  \nabla_{\Sigma}^{\Sigma'}(H_{\Omega})_{\Omega}. \end{equation}
  \end{lemma}
  
  \begin{proof}
Since $\nabla_{\Sigma}^{\Sigma'}(H)_{\Omega} \cong (\nabla_{\Sigma}^{\Sigma'}(H)_{\Gal(H/H_{\Omega})})_{\Omega}$, this follows from Lemma \ref{lem:coinv}.  Note that the RW pair $(\Sigma, \Sigma')$ chosen at the beginning of this subsection is complete, and that $\Sigma'$ is 2-tame.
 \end{proof}

Write $G_{\Omega} = \Gal(H_{\Omega}/F)$.  {The conjugacy class $\Omega$ can be written 
\[ \Omega = \Omega_0 \Omega_1 \] where $\Omega_0$ is a conjugacy class of 2-power order characters and $\Omega_1$ is a conjugacy class of odd order characters.  This corresponds to the decomposition \[ G_{\Omega} = G_{\Omega,2} \times G'_{\Omega} \] of $G_{\Omega}$ as a product of a 2-group and an odd order group.  Since $G_\Omega$ is cyclic,  $G_{\Omega,2}$ is as well.  Hence 
the conjugacy class $\Omega_0$ is the set of {\em all} odd characters of the cyclic 2-group $G_{\Omega,2}$.  Indeed, all the characters $\psi$ of $\Z/2^m\Z$ such that $\psi(2^{m-1}) = -1$ are conjugate.    

As a result, any conjugacy class of odd characters of $G_{\Omega}$ can be written $\Upsilon = \Omega_0 \Upsilon_1$ where $\Upsilon_1$ is a conjugacy class of characters of $G'_{\Omega}$.
We therefore have a decomposition
\begin{equation} \label{e:gpprod}
 \Z_2[G_{\Omega}]_- \cong \prod_{\Upsilon} R_{\Upsilon}, \end{equation}
as $\Upsilon$ ranges over the conjugacy classes of odd characters of $G_{\Omega}$.  This is the key point where $p=2$ is being used in this proof (this decomposition does not hold in general for $p$ odd, and explains why a different argument was used in \cite{dk}*{\S5}).
 
In view of (\ref{e:gpprod}), we have
\begin{equation} \label{e:mupsilon}
 M_-  \cong \prod_{\Upsilon} M_{\Upsilon}
\end{equation}
for any $\Z_2[G_{\Omega}]$-module $M$.}

{Let $H_{\Upsilon}$ denote the fixed field of the kernel of any element in $\Upsilon$.}
 To avoid trivial zeroes, we need to consider
\[ \Sigma_c = \{v \in \Sigma \colon c \in G_v \} \supset S_\infty. \]
Let 
\[ \Sigma_{\Upsilon}  = \Sigma_c \cup \{ v \in S_{\ram}(H_{\Upsilon}/F) \colon v \mid 2 \}. \]
The pair $(\Sigma_\Upsilon,\Sigma')$ is complete for
$H_\Upsilon/F$: every finite place ramified in $H_\Upsilon/F$
above $2$ lies in $\Sigma_\Upsilon$, while every such place away
from $2$ lies in $\Sigma'$. Moreover,
\[
S_\infty\subset\Sigma_c\subset\Sigma_\Upsilon,
\]
so $\Sigma_\Upsilon$ contains a place whose decomposition group
contains $c$. Therefore Theorem~\ref{t:nablasequence} gives an exact
sequence over $\Z_2[G_\Upsilon]_-$. Projecting to the factor
$R_\Upsilon$ in \eqref{e:gpprod}, we obtain
\[
0\longrightarrow
\Cl_{\Sigma_\Upsilon}^{\Sigma'}(H_\Upsilon)_\Upsilon
\longrightarrow
\nabla_{\Sigma_\Upsilon}^{\Sigma'}(H_\Upsilon)_\Upsilon
\longrightarrow
(X_{H_\Upsilon,\Sigma_\Upsilon})_\Upsilon
\longrightarrow0.
\]
 The inclusion (\ref{e:incl}) for the field $H_\Upsilon$ implies that
\begin{equation} \label{e:div}
\prod_{\alpha \in \Upsilon} L_{\Sigma_\Upsilon, \Sigma'}(\alpha^{-1}, 0) \ \Big\vert \ |\nabla_{\Sigma_\Upsilon}^{\Sigma'}(H_\Upsilon)_\Upsilon |
 =  | \Cl_{\Sigma_{\Upsilon}}^{\Sigma'}(H_{\Upsilon})_\Upsilon | \cdot | (X_{H_\Upsilon, \Sigma_\Upsilon})_\Upsilon |,
 \end{equation}
where the divisibility occurs in $\Z_2$.
\begin{lemma} \label{l:sigmac}
We have
\begin{itemize}
\item $L_{\Sigma_\Upsilon, \Sigma'}(\alpha^{-1}, 0) = L_{\Sigma_c, \Sigma'}(\alpha^{-1}, 0)$.
\item {$ \Cl_{\Sigma_{\Upsilon}}^{\Sigma'}(H_{\Upsilon})_\Upsilon \cong  \Cl_{\Sigma_c}^{\Sigma'}(H_{\Upsilon})_\Upsilon \cong 
\Cl_{\Sigma_c}^{\Sigma'}(H_{\Omega})_\Upsilon.$}
\item $(X_{H_\Upsilon, \Sigma_\Upsilon})_\Upsilon \cong (X_{H_\Upsilon, \Sigma_c})_\Upsilon \cong (X_{H_\Omega, \Sigma_c})_\Upsilon$
\end{itemize}
 \end{lemma}
\begin{proof} The first bullet point is clear, since by definition, if $v \in \Sigma_\Upsilon \setminus \Sigma_c$ then any $\alpha \in \Upsilon$ is ramified at $v$.  

For the second bullet point, the kernel of \[ \Cl_{\Sigma_c}^{\Sigma'}(H_{\Upsilon}) \longtwoheadrightarrow  \Cl_{\Sigma_{\Upsilon}}^{\Sigma'}(H_{\Upsilon}) \] is the subgroup generated by the primes lying above  $v \in \Sigma_\Upsilon \setminus \Sigma_c$.  Such a prime is fixed by the inertia subgroup $I_v$, but since any $\alpha \in \Upsilon$ is ramified at $v$ there exists a $\sigma \in I_v$ such that $\alpha(\sigma) \neq 1 $ for all $\alpha \in \Upsilon$.  Since $c \not \in G_v$, $I_v$ has odd order and 
$\alpha(\sigma)$ is an odd order root of unity.  Thus $1 - \alpha(\sigma)$ is a 2-adic unit, and the image of the class represented by $v$ in $ \Cl_{\Sigma_c}^{\Sigma'}(H_{\Upsilon})_\Upsilon$ is trivial.  This proves the first isomorphism.  {The second isomorphism follows from the proof of \cite{dk}*{Lemma 5.2}.}

The proof of the first isomorphism in the third bullet point is similar to the first isomorphism of the previous bullet point {(the essential point again being that $1 - \alpha(\sigma)$ is a 2-adic unit)}. The second isomorphism is elementary.
\end{proof}

Combining (\ref{e:div}) with Lemma~\ref{l:sigmac}, taking the product over all $\Upsilon$, and applying (\ref{e:mupsilon}), we obtain
\begin{equation} \label{e:div2}
\prod_{\alpha \in \hat{G}_{\Omega}, \ \alpha \text{ odd}} 
L_{\Sigma_c, \Sigma'}(\alpha^{-1}, 0) \ \Big\vert \ 
| \Cl_{\Sigma_c}^{\Sigma'}(H_{\Omega})_- | \cdot 
| (X_{H_{\Omega},\Sigma_c})_- |.  
 \end{equation}

\begin{lemma} \label{l:diveq}
The divisibility in (\ref{e:div2}) is an equality up to sign.
\end{lemma}

\begin{proof}
By the Artin formalism for $L$-functions, the left side of (\ref{e:div2}) may be written
\[ \prod_{\alpha \in \hat{G}_{\Omega}, \ \alpha \text{ odd}} L_{\Sigma_c, \Sigma'}(\alpha^{-1}, 0) = L_{\Sigma_c, \Sigma'}(H_{\Omega}/H_{\Omega}^+, \epsilon, 0)
 \]
where $H_{\Omega}^+$ denotes the maximal totally real subfield of the CM field $H_{\Omega}$, and $\epsilon$ is the non-trivial character of $\Gal(H_{\Omega}/H_{\Omega}^+)$. Dividing the analytic class number formulas for $\zeta_{H_{\Omega}, \Sigma_c, \Sigma'}(s)$ and $\zeta_{H_{\Omega}^+, \Sigma_c, \Sigma'}(s)$ at $s=0$ {as in \cite{dk}*{Equation (25)}} yields
\begin{equation} \label{e:cn}
  L_{\Sigma_c, \Sigma'}(H_{\Omega}/H_{\Omega}^+, \epsilon, 0) = \pm \frac{|\Cl_{\Sigma_c}^{\Sigma'}(H_{\Omega})|}{|\Cl_{\Sigma_c}^{\Sigma'}(H_{\Omega}^+)|} \cdot 
\frac{R_{\Sigma_c}^{\Sigma'}(H_{\Omega})}{R_{\Sigma_c}^{\Sigma'}(H_{\Omega}^+)}
\end{equation}
by \cite{pcmi}*{Equation (16)}.
The groups associated to the indicated regulators are the same, namely 
\begin{equation} \label{e:equalunits}
 \cO_{H_{\Omega}, \Sigma_c, \Sigma'}^* = \cO_{H_{\Omega}^+, \Sigma_c, \Sigma'}^*. \end{equation}
{To see this equality, one inclusion is clear.  For the reverse, let $u \in \cO_{H_{\Omega}, \Sigma_c, \Sigma'}^*$ and consider $x = u/u^c$.  Then $x$ is a $\Sigma_c$-unit, and since $c \in G_v$ for any $v \in \Sigma_c$, we have $|x|_v = 1$ for all $v \in \Sigma_c$ as well.  Therefore $x$ has absolute value 1 at every place and is therefore a root of unity.  Since $\Sigma'$ is torsion-free, $x=1$, so $u = u^c$ and $u \in \cO_{H_{\Omega}^+, \Sigma_c, \Sigma'}^*$ as desired.}

The normalization of the absolute values used in the definition of the regulators $R_{\Sigma_c}^{\Sigma'}(H_{\Omega})$ and $R_{\Sigma_c}^{\Sigma'}(H_{\Omega}^+)$
 is off by a factor of 2 in each entry of the regulator matrices, {since $H_\Omega$ is totally complex and $H_{\Omega}^+$ is totally real. 
 The regulator matrices have rank $\#(\Sigma_c)_{H_{\Omega}} -1$}.  As a result, the ratio of the regulators on the right side of (\ref{e:cn}) is $2^{\#(\Sigma_c)_{H_{\Omega}} - 1}$. {To complete the proof, we will show that}
\begin{equation} \label{e:cnratio}
\frac{|\Cl_{\Sigma_c}^{\Sigma'}(H_{\Omega})|}{|\Cl_{\Sigma_c}^{\Sigma'}(H_{\Omega}^+)|} =  |\Cl_{\Sigma_c}^{\Sigma'}(H_{\Omega})_- |
\end{equation}
and
\begin{equation} \label{e:xsize}
|(X_{H_{\Omega}, \Sigma_c})_-| = 2^{|(\Sigma_c)_{H_{\Omega}}| - 1}.
\end{equation}
The second of these is straightforward and follows from the short exact sequence
\[ 0 \longrightarrow (X_{H_{\Omega},\Sigma_c})_- \longrightarrow (Y_{H_{\Omega},\Sigma_c})_- \longrightarrow \Z_- = \Z/2\Z \longrightarrow 0. 
\]
(The left exactness of this sequence can be established directly, and also follows from Lemma~\ref{l:tor}, which implies $\lrar_1^{\Z[G]}(\Z, \Z[G]_-) = \Z_+[2] = 0$.)
Indeed, for $v \in \Sigma_c$ we have $c \in G_v$ so the module $(Y_{H_{\Omega},v})_- \cong \Z[G_{\Omega}/G_v]/2$ has size $2^{|\{w \mid v\}|}$, where the exponent  is the number of places $w$ of $H_{\Omega}$ dividing $v$.   This proves (\ref{e:xsize}).

To prove (\ref{e:cnratio}), we note there is a short exact sequence
\begin{equation} \label{e:cgs}
\begin{tikzcd}
 1 \ar[r] & \Cl_{\Sigma_c}^{\Sigma'}(H_{\Omega})^-  \ar[r] & \Cl_{\Sigma_c}^{\Sigma'}(H_{\Omega})  \ar[r,"\N"] &
\Cl_{\Sigma_c}^{\Sigma'}(H_{\Omega}^+)  \ar[r] &  1,
\end{tikzcd}
\end{equation}
where the labelled arrow is given by the norm.  Here the {\em superscript} on $ \Cl_{\Sigma_c}^{\Sigma'}(H_{\Omega})^-$ denotes the largest {\em subgroup} (rather than quotient, which is indicated by a subscript) on which $c$ acts as $-1$.  It is clear that $\Cl_{\Sigma_c}^{\Sigma'}(H_{\Omega})^-$ contains the kernel of the norm but a small observation is needed to deduce equality: suppose $\fa \fa^c = (x)\fb$ where $x \in (H_{\Omega})^*_{\Sigma'}$ and $\fb$ is in the subgroup of fractional ideals generated by the primes above $\Sigma_c$.  Dividing this equation by its conjugate, we see that $(x/x^c) = (1)$.  Therefore \[ x/x^c \in \cO_{H_{\Omega}, S_\infty, \Sigma'}^*. \] {We now argue as in the proof of (\ref{e:equalunits}): $x/x^c$ is a unit with absolute value 1 in every complex embedding, so it must be a root of unity. Therefore $x/x^c=1$ since $\Sigma'$ is torsion-free}, so $x = x^c$, i.e.\ $x \in (H_{\Omega}^+)^*_{\Sigma'}$.  Thus $\fa$ lies in the kernel of the norm.

We should also comment on the surjectivity of the norm map in (\ref{e:cgs}).  By class field theory, the cokernel of this map is identified with the Galois group of the largest subextension of $H_{\Omega}/H_{\Omega}^+$ in which the primes in $\Sigma_c$ split completely.  Since this includes the infinite primes, the largest such subextension is trivial (the point here is that in the definition of $ \Cl_{\Sigma_c}^{\Sigma'}(H_{\Omega}^+)$ one takes the quotient by principal ideals generated by all elements of $(H_{\Omega}^+)^*_{\Sigma'}$, not just the totally positive elements). 

To conclude the proof, we just need to note the equality \[  | \Cl_{\Sigma_c}^{\Sigma'}(H_{\Omega})^-| = |  \Cl_{\Sigma_c}^{\Sigma'}(H_{\Omega})_- |, \] 
which follows from the tautological  exact sequence
\[  \begin{tikzcd}
1 \ar[r] &   \Cl_{\Sigma_c}^{\Sigma'}(H_{\Omega})^-  \ar[r] &  \Cl_{\Sigma_c}^{\Sigma'}(H_{\Omega}) \ar[r,"1+c"] 
 &  \Cl_{\Sigma_c}^{\Sigma'}(H_{\Omega})  \ar[dll, out=-10, in=170] & \  \\ 
 &  \Cl_{\Sigma_c}^{\Sigma'}(H_{\Omega})_- \ar[r] \ar[r] & 1. 
 \end{tikzcd} \]
\end{proof}

In view of Lemma~\ref{l:diveq}, it follows that all the divisibilities (\ref{e:div}) used to deduce (\ref{e:div2}) are equalities up to 2-adic units.  For $\Upsilon = \Omega$ this says
\begin{equation}
  | \nabla_{\Sigma_{\Omega}}^{\Sigma'}(H_{\Omega})_{\Omega} | \doteq \prod_{\alpha \in \Omega} L_{\Sigma_{\Omega}, \Sigma'}(\alpha^{-1}, 0). \label{e:phiprod2}
  \end{equation}
  Now, the primes in $\Sigma\setminus\Sigma_{\Omega}$ are unramified for the characters in $H_{\Omega}$.  Adding unramified primes to the depletion set simply multiplies the Fitting ideal of $\nabla$ by the Euler factor $(1 - \sigma_v^{-1})$ at those primes {(see Lemma \ref{l:addprimes})}, so we immediately obtain from (\ref{e:phiprod2})
\[
  | \nabla_{\Sigma}^{\Sigma'}(H_{\Omega})_{\Omega}| \doteq \prod_{\alpha \in \Omega} L_{\Sigma, \Sigma'}(\alpha^{-1}, 0). \]
In view of Lemma~\ref{l:descent}, this concludes the proof of Theorem~\ref{t:nsize2} and hence of Theorem~\ref{t:incl}.

\subsection{Summary} \label{s:summary}

The following result summarizes the results obtained in \S\ref{s:sr} so far. {We assume Theorems~\ref{t:nablasequence}, \ref{thm:quad}, and \ref{t:cocycle-nabla}; these will be proven in \S\ref{s:constr-rw}.}

\begin{theorem} \label{t:incleq}
Fix a prime $p$. In Case 1 (no primes above $p$ are ramified in $H$) let
\[ \Sigma = \{\fl\} \]
for some finite unramified prime $\fl \nmid p $ whose Frobenius in $G$ equals the complex conjugation $c$.  In Case 2, let 
\[  \Sigma = \{v \in S_{\ram} \colon v \mid p \}. \]
In both cases let 
\[  \Sigma' = T \cup \{v \in S_{\ram} \colon v \nmid p \} \cup S_\infty, \]
where $T$ is a finite torsion-free set of finite primes
disjoint from \[S=S_{\ram}\cup S_\infty.\]
For each $\Gal(\overline{\Q}_p / \Q_p)$-conjugacy class of characters $\Phi$ of $G'$, let
 \[ \Phi_0 = \{ \chi \in \hat{G} \colon \chi \text{ is odd, } \chi(G_v) \neq 1 \text{ for all } v \in \Sigma, \text{ and } \chi|_{G'} \in \Phi \} \] 
and let $R = R_{\Phi_0}$ be the associated character group ring.
Suppose that for each $\Phi$ we have a finitely generated $R$-module $M$ and a cocycle $\kappa \in Z^1(G_F, M)$ satisfying the conditions of Theorem~\ref{t:cocycle-nabla} such that 
$M$ is generated over $R$ by $\kappa(G_F)$ and the $x_v$ for $v\in\Sigma\setminus\{v_0\}$, and \begin{equation} \label{e:ninc}
 \Fitt_R(M) \subset (\Theta_{\Sigma, \Sigma'}/2^n) R. 
 \end{equation}
Moreover, suppose the same result holds for every CM subfield $H' \subset H$ containing $F$.
Then Theorem~\ref{t:main} holds.
\end{theorem}

\begin{proof} {
Since $(\Sigma,\Sigma')$ is $p$-admissible, the hypotheses of Theorem~\ref{t:cocycle-nabla} are satisfied: the pair is complete, $\Sigma\neq\emptyset$,
and $\Sigma'$ is $p$-tame. Applying the theorem, we} obtain a surjective map $\nabla_{\Sigma}^{\Sigma'}(H)_R \longrightarrow M$, which in conjunction with the inclusion (\ref{e:ninc}) yields
\[ \Fitt_R(\nabla_{\Sigma}^{\Sigma'}(H)_R) \subset   (\Theta_{\Sigma, \Sigma'}/2^n)R.
\]
By assumption, we get the same inclusion for all CM subfields of $H$ containing $F$ as well.
Corollary~\ref{c:incl} then implies that Theorem~\ref{t:main} holds.
\end{proof}

\subsection{Ribet's Lemma} \label{ss:rl}
Theorem~\ref{t:incleq} demonstrates the importance of constructing Galois cohomology classes that generate large modules, where in our context ``large'' is defined precisely by (\ref{e:ninc}).  Our technique for constructing these classes is Ribet's Lemma.  In this paper, we prove the following general version of Ribet's Lemma for $\GL_2$.
  Of central importance is the fact that the condition of residual distinguishability
 $\chi \not\equiv \psi \pmod{\fm}$ is  not assumed, where $\fm$ is the maximal ideal of the ring $\T$ in the theorem below.

\begin{theorem} \label{t:rl}
 Let $\T \subset \widetilde{\T}$ be an inclusion of Noetherian commutative  rings, with $(\T, \fm)$ local.  
 % Suppose that $\T$ and $\widetilde{\T}$ are complete with respect to the maximal ideal of $\T$. 
  Let $\widetilde{I} \subset \widetilde{\T}$ be a nontrivial ideal and let $I = \widetilde{I} \cap \T$.  Let $K = \Frac(\widetilde{\T})$ be the total ring of fractions of $\widetilde{\T}$.
 Assume that $K$ is a product of local rings and that the maximal ideals of $K$ are principal.
  Let $\cG$ be a  group and suppose we are given a representation
\[ \rho \colon \cG \longrightarrow \GL_2(K) 
\]
satisfying the following conditions.
\begin{itemize}
\item  
 For each $\sigma \in \cG,$  the characteristic polynomial $P_{\rho(\sigma)}(x)$ lies in $\T[x]$.  Furthermore we have
\begin{equation} \label{e:redpoly}
P_{\rho(\sigma)}(x) \equiv (x - \chi(\sigma))(x - \psi(\sigma)) \pmod{I} 
\end{equation}
for two characters $\chi, \psi \colon \cG \longrightarrow \T^*$ such that $\chi \equiv \psi \pmod{\fm}$.
\item   Let $K_0 = 
\red(K)$ denote the maximal reduced quotient of $K$. Write $K_0 = \prod_{i=1}^{m} k_i$ as a product of fields.  For every projection \[ K \rightarrow K_0 \rightarrow k_i,\] the projection of $\rho$ to $\GL_2(k_i)$ is an irreducible representation of $\cG$ over $k_i$.
\item The extension of $\rho$ to a $\T$-algebra homomorphism $\T[\cG] \longrightarrow M_2(K)$ has image that is finitely generated as a $\T$-module.
\item We are given a  set of subgroups $\cG_v \subset \cG$, indexed by a disjoint union of finite sets $\Sigma \sqcup \cP$ (we allow $\Sigma, \cP$ to be empty). For each $v \in \Sigma \sqcup \cP$ there exists a basis in which the restriction of $\rho$ has the form
\begin{equation} \label{e:rholocal}
 \rho|_{\cG_v} \cong \mat{\eta_v}{0}{*}{\xi_v}
\end{equation}
for two characters $\xi_v, \eta_v \colon \cG_v \longrightarrow \widetilde{\T}^*.$  
\item  For each $v \in \Sigma$ we have the congruence 
  \begin{equation} \xi_v \equiv \psi|_{\cG_v} \pmod{\widetilde{I}}. \label{e:Sigmaassumption}
   \end{equation}
\item For each $v \in \cP$ we are given a subgroup $\cI_v \subset \cG_v$ such that 
\begin{equation} \label{e:Passumption}
(\xi_v)|_{\cI_v} \equiv \chi|_{\cI_v} \pmod{\widetilde{I}}. \end{equation}
   \end{itemize}
 If $\Sigma$ is nonempty,  fix $v_0 \in \Sigma$.  Choose an element $\sigma_v \in \cG_v$ for each $v \in \cP$.  Then there exists a finitely generated $\T$-module $M$ and a 
 cocycle \[ \kappa \in Z^1(\cG, M(\chi\psi^{-1})) \] 
satisfying the following conditions.
\begin{itemize}
\item For each $v \in \cP$, we have $\kappa(\sigma) = 0$ for all $\sigma \in \cI_v$.
\item $\kappa(\ker \rho \cap \ker \psi) = 0$.
\item If $\Sigma$ is nonempty, we have $\kappa(\cG_{v_0}) = 0$ and for each $v \in \Sigma \setminus \{v_0\}$, there exists $y_v \in M$ such that
\[ \kappa(\sigma) = (\chi\psi^{-1}(\sigma) - 1) y_v \]
for all $\sigma \in \cG_v.$
\item $M$ is generated over $\T$ by $\kappa(\cG)$ and $\{y_v \colon v \in \Sigma \setminus \{v_0\}\}$.

\item We have \begin{equation} \label{e:nincl}
 \left(\prod_{v \in \cP} (\xi_v(\sigma_v) - \chi(\sigma_v))\right)\Fitt_{\T}(M) \subset \widetilde{I}. \end{equation}
\end{itemize}
\end{theorem}
Theorem~\ref{t:rl} will be proved in \S\ref{s:cocycle}--\ref{s:cc}.
We conclude this section by combining Theorems~\ref{t:incleq} and~\ref{t:rl} to show that
in order to prove our main result (Theorem~\ref{t:main}), it suffices to construct a Galois representation with the desired properties.
In proving Theorem~\ref{t:main}, the auxiliary-set reduction in
\S\ref{s:minimal} allows us to choose the torsion-free unramified
set $T$ away from the primes above $p$. For the minimal pairs
of Theorem~\ref{t:incleq}, the smoothing set $\Sigma'$
therefore also contains no primes above $p$.

\begin{theorem} \label{t:galrep}
Let the notation be as in Theorem~\ref{t:incleq}. Suppose that
for every nonempty character component $R=R_{\Phi_0}$ on which
$\Theta_{\Sigma,\Sigma'}^\#/2^n$ is a nonunit, and likewise
for every CM subfield of $H$ containing $F$, the following
data exist (with notation relative to that field).
There is an inclusion of commutative Noetherian $R$-algebras
$\T\subset\widetilde\T$, with $\T$ local, such that both rings
are complete with respect to the maximal ideal of $\T$ and
$K=\Frac(\widetilde\T)$ has finite dimension over $\Q_p$.
There is a continuous representation
$\rho:G_F\longrightarrow\GL_2(K)$ satisfying the hypotheses
of Theorem~\ref{t:rl}, using also its residually distinguishable
version stated in the first remark of \S\ref{s:cocycle}, where:
\begin{itemize}
\item The set
 $\Sigma$ is as in Theorem~\ref{t:incleq} and $\cP$ is equal to the set of primes above $p$ that are not in $\Sigma$.
\item We have  $\psi = \bpsi$, where $\bpsi$ is the canonical character \[ \bpsi\colon G_F \lra G \lra \T^*.\]  We have $\chi \equiv 1 \pmod{p^m}$ for some positive integer $m$ large enough that $(\Theta^\#_{\Sigma, \Sigma'}/2^n)$ divides $p^m$ in $R$.
\item The representation $\rho$ is unramified outside $\Sigma \cup \Sigma' \cup \cP$.
\item $\T/I$ is cyclic as an $R$-module, i.e.\ the structure map $R \lra \T/I$ is surjective.
 \item  If $y \in R$ and $(\prod_{v \in \cP} (\xi_v(\sigma_v) - \chi(\sigma_v))) y \in \widetilde{I}$, then
$y \in (\Theta^\#_{\Sigma, \Sigma'}/2^n)R$.
\end{itemize}
Then Theorem~\ref{t:main} follows.
\end{theorem}

\begin{proof} On a character component $R$ on which
$\Theta_{\Sigma,\Sigma'}/2^n$ is a unit, take $M=0$,
$\kappa=0$, and $x_v=0$ for
$v\in\Sigma\setminus\{v_0\}$ in Theorem~\ref{t:incleq}.
All the local and generation conditions hold, and
\[
\Fitt_R(0)=R=(\Theta_{\Sigma,\Sigma'}/2^n)R.
\]
The group-ring involution identifies such components with
the inverse character components on which
$\Theta_{\Sigma,\Sigma'}^\#/2^n$ is a unit.
The same argument applies to every CM subfield of $H$
containing $F$. We may therefore restrict to a component
$R$ on which $\Theta_{\Sigma,\Sigma'}^\#/2^n$ is a nonunit
and use the data in the hypothesis.

First we comment on the finite generation of $\rho(\T[G_F])$
as a $\T$-module.
The image $\rho(G_F)$ is compact and hence bounded in
$M_2(K)$. Its $\Z_p$-span is contained in a finite-rank
$\Z_p$-lattice and is therefore finitely generated over $\Z_p$.
Its $\T$-span, which is $\rho(\T[G_F])$, is consequently
finitely generated over $\T$.

Now, Theorem~\ref{t:rl}  provides a $\T$-module $M$ and a cocycle 
\[ \kappa \in Z^1(G_F, M(\chi\bpsi^{-1})).\]  
It will be clear from our construction that the continuity of $\rho$ implies the continuity of $\kappa$ (see Remarks~\ref{remark:top} and \ref{remark:cont}).
We consider $\overline{M} = M/(IM, p^mM)$.  The projection of $\kappa$ to $\overline{M}$ can be viewed as a cocycle 
 \[ \overline{\kappa} \in Z^1(G_F, \overline{M}(\bpsi^{-1})) \] since $\chi \equiv 1 \pmod{p^m}$.

 Since $\T/I$ is a cyclic $R$-module, the $\T$-module generators of $M$ reduce to $R$-module generators of $\overline{M}$.  All of the properties required by Theorem~\ref{t:cocycle-nabla} are directly seen to be satisfied by the corresponding properties in Theorem~\ref{t:rl}, except for possibly the discussion of ramification.   The Galois representation $\rho$ is unramified outside $\Sigma \cup \Sigma' \cup \cP$, as is the character $\psi$,  so $\kappa$ and $\overline{\kappa}$ are unramified outside this set by the second bullet point of the conclusion of Theorem~\ref{t:rl}.  Furthermore, we are given in Theorem~\ref{t:rl} that $\overline{\kappa}$ is unramified at $\cP$, and that the local conditions  at $\Sigma$ hold.  Finally, $\overline{\kappa}$ is tamely ramified at the finite primes in $\Sigma'$ since our module is pro-$p$ and $\Sigma'$ contains no primes above $p$.

Let
\[
\Phi_0^{-1}=\{\alpha^{-1}:\alpha\in\Phi_0\},\qquad
R^\iota=R_{\Phi_0^{-1}}.
\]
The group-ring involution induces an isomorphism
$\#:R^\iota\xrightarrow{\sim}R$.
With its $G_F$-action through $\bpsi^{-1}$, the module
$\overline M$ is naturally an $R^\iota$-module via this
isomorphism. Thus it suffices to prove
\begin{equation}\label{e:fittnbar}
\Fitt_R(\overline M)
 \subset (\Theta_{\Sigma,\Sigma'}^\#/2^n).
\end{equation}
Indeed, this is equivalent to the required Fitting containment
over $R^\iota$, and inversion permutes the character components
which occur in Theorem~\ref{t:incleq}.
Put $A=\T/I$ and $\overline A=A/p^mA$.
The maps $R\to\overline A$ and $\T\to\overline A$ are
surjective. The module $\overline M$ is the base change of
$M$ from $\T$ to $\overline A$; it is also the base change of
itself from $R$ to $\overline A$.
Base change for Fitting ideals therefore identifies the images
of $\Fitt_R(\overline M)$ and $\Fitt_\T(M)$ in $\overline A$.
Consequently, for $y\in\Fitt_R(\overline M)$ there are
$i\in I$ and $r\in R$ such that, after applying the structure
map $R\to\T$,
\[
y+i+p^m r\in\Fitt_\T(M).
\]
Here the coefficient of $p^m$ can be chosen in $R$ because
$R\to A$ is surjective. Equation~\ref{e:nincl} now gives
\[\prod_{v \in \cP} (\xi_v(\sigma_v) - \chi(\sigma_v)) ( y + p^mr) \in \widetilde{I}. \]
The last assumption of the theorem then yields \[ y + p^m r \in (\Theta_{\Sigma, \Sigma'}^\#/2^n) \] and since $\Theta_{\Sigma, \Sigma'}^\#/2^n$ divides $p^m$, we obtain $y \in (\Theta_{\Sigma, \Sigma'}^\#/2^n)$ as desired.
\end{proof}

In \S\ref{s:modforms}, we will  use group-ring valued Hilbert modular forms to construct the Galois representation required by Theorem~\ref{t:galrep}.

\section{Modular Forms} \label{s:modforms}

Let $\Sigma, \Sigma'$, and $R = R_{\Phi_0}$ be as in the statement of Theorem~\ref{t:incleq}. In particular, $(\Sigma, \Sigma')$ are the minimal pair defined in Lemma~\ref{l:minimal}. This notation will be fixed throughout \S\ref{s:modforms}.
Throughout this section, assume that $R\neq 0$ and that
$\Theta_{\Sigma,\Sigma'}^\#/2^n$ is a nonunit in $R$.
By Theorem~\ref{t:galrep} and its proof, it is enough to
construct the required data on these components.
Recall that $\Sigma'$ contains no primes above $p$ (see the comment preceding Theorem~\ref{t:galrep}). 
 
The goal of this section is to describe certain Hecke algebras $\T_\fm \subset \widetilde{\T}_{\fm}$, which will be Noetherian $R$-algebras that are complete with respect to the maximal ideal $\fm \subset \T$.
We will associate to this algebra:
\begin{enumerate}
\item An $R$-algebra homomorphism $\varphi \colon \widetilde{\T}_{{\fm}} \longrightarrow W$ where $W$ is another $R$-algebra, satisfying the following:
\begin{enumerate}
\item The homomorphism $\varphi$ is Eisenstein in the sense that
\[ \varphi(T_\fq) = \epsilon_{\cyc}^{k-1}(\fq)+ \bpsi(\fq) \] for primes $\fq \not\in \Sigma \cup \Sigma', \fq \nmid p$.
\item The restriction of $\varphi$ to $\T_\fm$ induces an isomorphism \[ \varphi(\T_\fm) \cong R/(x \Theta_{\Sigma, \Sigma'}^\#/2^n)\] where $x$ is a certain nonzerodivisor in $R$. 
\item For a specific element $U \in \widetilde{\T}_\fm$ to be defined later, if $y \in R$ and $\varphi(Uy) = 0$, then $y \in (\Theta_{\Sigma, \Sigma'}^\#/2^n)R$. 
\end{enumerate}
 \item A continuous Galois representation
\[ \rho \colon G_F \longrightarrow \GL_2(K), \qquad K = \Frac(\widetilde{\T}_\fm) \]
such that $\tr(\rho(\Frob_\fq)) = T_\fq$ for $\fq \not\in \Sigma \cup \Sigma', \fq \nmid p$.
\end{enumerate}

The homomorphism $\varphi$ and representation $\rho$ will satisfy certain other local properties
as required in Theorem~\ref{t:galrep}.   The construction of $\varphi$ follows  the argument in \cite{dk}*{\S 8.5} with one exception---the introduction of the prime $\fl \in \Sigma$ in Case 1 requires an altered definition of the necessary cusp form (Corollary \ref{c:cuspc1}).  Similarly in the construction of the Galois representation $\rho$, subtleties arising from forms old at $\fl$ arise in Case 1.  In this section we will focus only on these new aspects, while stating without proof the results carried over directly from \cite{dk}.
For simplicity, we will simply write $\Theta^\#$ for $\Theta_{\Sigma, \Sigma'}^\#$ in this section.

\subsection{Hecke Algebra} \label{ss:hecke}

We use the definitions and notations concerning {\em Hilbert modular forms} and {\em group-ring valued Hilbert modular forms} as given in \cite{dk}*{\S7.2--7.3}.  We will follow closely the construction of {\em loc.\ cit.} \S8.   Define
\begin{align*}
\fn &= \cond(H/F) \prod_{\fq \mid \Sigma'} \fq, \\
\fP &= \gcd(p^{\infty}, \fn), \\
\fP' &= \prod_{\fp \mid p, \fp \nmid \fP} \fp. 
\end{align*}

Recall also the cases that were defined in \S\ref{s:minimal} above---Case 1 is when $\fP = 1$ and Case 2 is when $\fP \neq 1$.
In Case 1, choose once and for all a prime $\fl\nmid\fn p$
with the weight-independent nonscalarity property of
Lemma~\ref{l:nonscalarprime}. In particular, its Frobenius in
$G$ is $c$. This choice is made before choosing the weight
$k$. We define
\[
\widetilde{\fn} = \left\{ \begin{array}{cc} \fn \fl  & \text{ in Case 1} \\ \fn & \text{ in Case 2.}\end{array} \right.
\]
We will be considering group-ring valued Hilbert modular forms of level $\widetilde{\fn} \fP'$ and weight $k$, where $k$ is a large positive integer congruent to $1$ modulo $(p-1)p^N$ for $N$ sufficiently large. We denote by 
\[ \bpsi \colon G_F \lra G \lra R^* \] the canonical character and consider the space of $p$-ordinary cusp forms $S_k(\widetilde{\fn} \fP', R, \bpsi)^{p-\ord}$, {
as defined in \cite{dk}*{\S7.2.8--\S7.2.9}.}

{
We define the Hecke algebra ${\T}$ to be the $R$-subalgebra of  \[ \End_R(S_k(\widetilde{\fn} \fP', R, \bpsi)^{p-\ord}) \] generated by the following Hecke operators:}
\begin{itemize}
\item  $T_{\fq}$ for $\fq \nmid \widetilde{\fn} \fP'$, 
\item $U_{\fp}$ for $\fp \mid \fP$, 
\item the diamond operators $S(\fa)$ for $(\fa, \widetilde{\fn} \fP') = 1$.
\end{itemize}
 Let $\widetilde{\T} \supset \T$ denote algebra generated over $\T$ by the operators $U_\fp$ for $\fp \mid \fP'$, and in Case 1, the operator $U_\fl$.

\bigskip

In the next subsection, we will construct a modular form \[ F_k(\bpsi) \in S_k(\widetilde{\fn} \fP', R, \bpsi)^{p-\ord}\] that is congruent to an Eisenstein series modulo the maximal ideal of $R$ (in fact we prove a much stronger congruence).  The existence of this form implies the following.

\begin{lemma} \label{l:einmax}  Let $k$ denote the residue field of the local ring $R$.  There is an $R$-algebra homomorphism $\overline{\varphi} : \T \lra k$ given by 
\begin{itemize}
\item $\overline{\varphi}(T_{\fq}) = 1 + \overline{\bpsi}(\fq)$ for $\fq \nmid \widetilde{\fn} p$; 
\item $\overline{\varphi}(U_{\fp}) = 1$ for $\fp \mid \fP$;
\item $\overline{\varphi}(S(\fa)) = \overline{\bpsi}(\fa)$, {for $(\fa, \widetilde{\fn} \fP') = 1$.}
\end{itemize}
\end{lemma}
Denote by $\mathfrak{m} \subset \T$ the kernel of $\overline{\varphi}$.  Denote by $\T_{\mathfrak{m}}$ and $\widetilde{\T}_{\mathfrak{m}} = \widetilde{\T} \otimes_{\T} \T_{\mathfrak{m}}$ the $\mathfrak{m}$-adic completions of $\T$ and $\widetilde{\T}$, respectively.

We now describe these Hecke algebras more concretely.
Denote by {$\mathcal M$} the set of $p$-ordinary cuspidal newforms of weight $k$, level dividing $ \widetilde{\fn} \fP'$, and nebentypus $\psi$ {for some $\psi \in\Phi_0$}. For each $f \in \cM$,  denote by $f_p$ the ordinary stabilization of $f$ with respect to all primes $\fp \mid p$.  Let $E = E_f$ denote the finite extension of $\Q_p$ generated by the Fourier coefficients of $f$, let $\cO$ be the ring of integers of $E$, and let $\pi$ be a uniformizer for $E$.
Let 
\[
\overline{\cM} = \{ f \in \cM : c(1, (f_p)|_t) \equiv \overline{\varphi}(t) \pmod{\pi} \text{ for all } t \in \T \}.
\]
 We can state this another way. For each $f \in {\cM}$, 
 let $\lambda_f:\T\longrightarrow E_f$ be the eigenvalue
homomorphism, $\lambda_f(t)=c(1,(f_p)|_t)$, and define the prime ideal
$\mathcal P_f=\ker(\lambda_f)$.
 Then $\overline{\cM}$ is the set of $f \in \cM$ such that $\mathcal{P}_f \subset \fm$.

For $f \in \overline{\cM}$, write $K = \Frac(\widetilde{\T}_\fm)$ for the total ring of fractions of $\widetilde{\T}_\fm$. Write $K_f$ for the total ring of fractions of the localization $\widetilde{\T}_{\mathcal{P}_f} = \T_{\mathcal{P}_f} \otimes_{\T} \widetilde{\T}$. 
Then we have
\begin{equation} \label{e:kdec}
K = \prod_{\substack{f \in \overline{\cM} \\ \mathrm{distinct} \ \cP_f}} \!\!\!\!\! K_f.
\end{equation}
Here the product ranges over $\Gal(\overline{\Q}_p/\Q_p)$-orbits of forms $f \in  \overline{\cM}$, i.e. one representative $f$ is taken for each associated prime ideal $\mathcal{P}_f \subset \T$.

We can be explicit about $K_f$.
First note that $\Frac(\T_{\mathcal{P}_f}) = E = $ the finite extension of $\Q_p$ generated by the Fourier coefficients of $f$.
Next we note that $U_\fp \in \Frac(\T)$ for $\fp \mid \fP'$, which results from the discussion of \cite{dk}*{\S8.5}. 
However, the same is not necessarily true in Case 1 for $U_\fl$.  As we outline below, there are two situations that  occur.

\begin{itemize}
\item We are in Case 2, or in Case 1 and the newform $f$ has level divisible by $\fl$.  Then $K_f = E$.
\item We are in Case 1 and the newform $f$ has level not divisible by $\fl$.
 Then \begin{equation} \label{e:kfdef}
 K_f = E[U_\fl] / g(U_\fl) \end{equation}
 where  
\[ g(x) =  x^2 - c(f, \fl)x+ \psi(\fl) \N\fl^{k-1} \]
is the Hecke polynomial at $\fl$ of $f$.  As always there are three possibilities for the quadratic polynomial $g(x)$: it may be irreducible over $E$, in which case $K_f$ is a quadratic field extension of $E$; it may split as a product of distinct linear factors, in which case $K_f \cong E \times E$; it may factor as a square, in which case $K_f \cong E[x]/x^2$ is the ring of dual numbers over $E$ (and is in particular not reduced but has a principal maximal ideal).
\end{itemize}

\subsection{Cusp form and homomorphism} \label{ss:cusp}

We now define the modular forms whose existence yields the desired homomorphism $\varphi$ on $\widetilde{\T}_\fm$. 
We will work in a weight $k > 1$ such that $k \equiv 1 \pmod{(p-1)p^N}$ for $N$ sufficiently large.
For now, it will be convenient to ensure that \begin{equation} \label{e:epscong}
\epsilon_{\cyc}^{k-1} \equiv 1 \pmod{\Theta^\#/2^n}, \end{equation}
where $\epsilon_{\cyc}$ is the $p$-adic cyclotomic character. 
  Therefore we choose a positive integer $m$ large enough such that {$\Theta^\#/2^n$} divides $p^m$ in $R$ and ensure that $N \ge m$, so that $\epsilon^{k-1}_{\cyc} \equiv 1 \pmod{p^m}$.  The integer $N$ will need to be made further sufficiently large in the results below.

We first recall the following theorem of Silliman \cite{silliman}, generalizing an earlier construction of Hida--Wiles. 
\begin{theorem} \label{t:hasse} Let $m$ be a positive integer. For positive integers $k \equiv 0 \pmod{(p-1)p^N}$ with $N$ sufficiently large, there is a modular form $V_k \in M_k(1, \Z_p, 1)$ such that $V_k \equiv 1 \pmod{p^m}$, and such that the normalized constant term $c_{\cA}(0,V_k)$ for each cusp $[\cA] \in \cusps(1)$ is congruent to $1 \pmod{p^m}$. 
\end{theorem}
We next refer the reader to \cite{dk}*{\S8.3, Eqn.~(102)}, where a certain group-ring family of Eisenstein series $W_k(\bpsi, 1) \in M_k(\fn, R, \bpsi)$ was defined.
In the (easier) Case 2, the construction of \cite{dk}*{Theorem 8.18} works without change for our application.
We only remark that the factor of $2^n$ was ignored in {\em loc.\ cit.} since this factor is a unit for $p$ odd, but we must be careful about recording this factor now that we are including $p=2$.

\begin{theorem} \label{t:fdefcase2}
Suppose we are in Case 2. For an integer $k > 1$ with $k \equiv 1 \pmod{(p-1)p^N}$ and $N$ sufficiently large, there exists a group ring form $H_k(\bpsi) \in M_k(\fn, R, \bpsi)$ such that 
\[
\widetilde{F}_k(\bpsi) = e_{\fP}^{\ord}\Big(W_1(\bpsi, 1)V_{k-1} - (\Theta^{\#}/2^n) H_k(\bpsi) \Big)
\]
lies in $S_k(\fn, R, \bpsi)$. 
\end{theorem}

In Case 1, we need to modify the construction of \cite{dk} to account for the new prime $\fl$.  

\begin{lemma} \label{l:vkl} Let the  notation be as in Theorem \ref{t:hasse}, with $N \ge m$. For any prime $\fl \nmid p$, the form 
\[
V_{k} - V_{k}|_{\fl} \in M_k(\fl, \Z_p, 1)
\]
has constant terms divisible by $p^m$ at each cusp. {Here $|_{\fl}$ denotes the level-raising operator as in \cite{dk}*{7.2.7}.}
\end{lemma} 
\begin{proof} We freely borrow notation from \cite{dkaux} here. If we are given a cusp $\mathcal{A} \in C_{\infty}(\fl)$ (see \cite{dkaux}*{\S3.5}), then by \cite{dkaux}*{Lemma 3.13} the constant term of $V_{k}|_{\fl}$ at $\mathcal{A}$ is equal to the constant term of $V_{k}$ at the cusp denoted $\mathcal{A}'$ in {\em loc.\ cit.} The result is then immediate since the normalized constant term of $V_{k}$ at every cusp is congruent to $1 \pmod{p^m}$. 

As we will explain, when $\mathcal{A} \in C_0(\fl)$, a similar calculation shows that the constant term of $V_{k}|_{\fl}$ at $\mathcal{A}$ equals $\N \fl^{-k}$ times the constant term of $V_{k}$ at $\mathcal{A}'$.   This will give the desired result since $k \equiv 0 \pmod{p^N(p-1)}$ and $N \ge m$, so $\N \fl^{-k} \equiv 1 \pmod{p^m}$.

Let $\mathcal{A} \in C_0(\fl)$.  A direct calculation using the definitions \cite{dkaux}*{equations (4) and (11)} reduces the statement above  to proving $\mathfrak{b}_{\mathcal{A}'} = \mathfrak{b}_{\mathcal{A}} \fl^{-1} \alpha_{\mu}$. We have
\begin{align}
\mathfrak{b}_{\mathcal{A}'} &= \alpha \alpha_{\mu} \mathcal{O}_F + \gamma (\ft_{\mu} \fd)^{-1} \nonumber \\
& =  \alpha \alpha_{\mu} \mathcal{O}_F + \gamma \alpha_{\mu} (\ft_{\lambda}\fd)^{-1} \fl^{-1}  \nonumber \\
& = (\alpha \fl   + \gamma (\ft_{\lambda} \fd)^{-1}) \alpha_{\mu} \fl^{-1}.  \label{e:baeqn}
\end{align}
By definition of $\fb_{\mathcal{A}}$, we have $\mathcal{O}_{F} = \alpha \mathfrak{b}_{\mathcal{A}}^{-1} + \gamma (\ft_{\lambda} \fd)^{-1} \mathfrak{b}_{\mathcal{A}}^{-1}$. Since $\mathcal{A} \in C_0(\fl)$, the ideal $\gamma (\ft_{\lambda} \fd)^{-1} \mathfrak{b}_{\mathcal{A}}^{-1}$ is coprime to $\fl$. Therefore 
\[
\mathcal{O}_{F} =  \alpha \fl \mathfrak{b}_{\mathcal{A}}^{-1} + \gamma (\ft_{\lambda} \fd)^{-1} \mathfrak{b}_{\mathcal{A}}^{-1},
\]
or equivalently,  $\fb_{\mathcal{A}} = \alpha \fl + \gamma (\ft_{\lambda} \fd)^{-1}$.  Substituting into (\ref{e:baeqn})
gives the desired result $\fb_{\mathcal{A}'} = \fb_{\mathcal{A}} \fl^{-1} \alpha_{\mu}$.
\end{proof}

In Case 1, put
\[
\Theta_0=\Theta_{\emptyset,\Sigma'},\qquad
D=\Theta_0^\#/2^n\in\Frac(R).
\]
Since $\bpsi(\fl)=-1$, our notation $\Theta=\Theta_{\Sigma,\Sigma'}$
with $\Sigma=\{\fl\}$ gives $\Theta=2\Theta_0$.
We define
\[
x=\Theta_{S_\infty}(1-k)/\Theta_{S_\infty}(0)\in R,
\]
which is a nonzerodivisor for $k$ sufficiently close to $1$, as in
\cite{dk}*{Lemma 8.16}. Here
$\chi(\Theta_{S_\infty}(1-k))=L(\chi^{-1},1-k)$, the usual primitive $L$-function with no depletion or smoothing sets.
For an ideal $\fa$, the subscript in $\psi_\fa$ means that the
character is considered with modulus divisible by $\fa$.
The level-raising identity is
\begin{equation}\label{e:fh}
W_j(\bpsi_\fl,1)
 =W_j(\bpsi,1)-\bpsi(\fl)W_j(\bpsi,1)|_\fl
 =W_j(\bpsi,1)+W_j(\bpsi,1)|_\fl.
\end{equation}

\begin{corollary}\label{c:cuspc1}
Suppose we are in Case 1, and fix the auxiliary prime $\fl$.
For $k>1$ and $k\equiv1\pmod{(p-1)p^N}$ with $N$ sufficiently
large, there is a group ring form
$H_k(\bpsi_\fl)\in M_k(\fn\fl,R,\bpsi)$ such that
\[
\widetilde F_k(\bpsi_\fl)
=xW_1(\bpsi_\fl,1)V_{k-1}-W_k(\bpsi_\fl,1)
 -(x\Theta^\#/2^n)H_k(\bpsi_\fl)
\]
lies in $S_k(\fn\fl,R,\bpsi)$.
\end{corollary}

\begin{proof}
We explain the argument at $p=2$, which is the case requiring the
normalization; the same argument works when $p$ is odd.
Let $G_k\in M_k(\fn,R,\bpsi)$ be the form of
\cite{dk}*{Theorem 8.9}. Its normalized constant term is
\[
\gamma_{\mathcal A}
=\operatorname{sgn}(\N\alpha)\,
  \bpsi^{-1}(\alpha\fb_{\mathcal A}^{-1})
\]
at $[\mathcal A]\in C_\infty(\fn)$, and is zero at the other
cusps. Such forms exist integrally at $2$ by the surjectivity
of the constant-term map in \cite{silliman}*{Theorem 1}.
As in \cite{dk}*{Proposition 8.11}, for any prescribed $M$ the form
\[
f_k=W_1(\bpsi,1)V_{k-1}-x^{-1}W_k(\bpsi,1)-DG_k
\]
has all normalized constant terms divisible by $2^M$ after
specialization at each character, provided $N$ is sufficiently
large. This constant-term assertion is made over $\Frac(R)$;
we do not assert that $D$ is integral.

Consider the constant-term vector
\[
b_k=\tfrac12\operatorname{const}(G_k+G_k|_\fl)
\]
at level $\fn\fl$. This vector is integral over $R$.
Indeed, over a cusp of $C_\infty(\fn)$ its entries are
\[
\begin{cases}
\gamma_{\mathcal A},&[\mathcal A]\in C_\infty(\fl,\fn\fl),\\
\dfrac{1-\N\fl^{-k}}2\gamma_{\mathcal A},
 &[\mathcal A]\in C_0(\fl,\fn\fl),
\end{cases}
\]
and they are zero over the other cusps.
For this calculation, level raising multiplies a normalized
constant term by $1$ or $\N\fl^{-k}$, respectively, while the
corresponding $\gamma$ is multiplied by $1$ or
$\bpsi(\fl)^{-1}=-1$. These are the same cusp transformations
as in Lemma~\ref{l:vkl} and \cite{dkaux}*{Lemma 3.13}.
The integrality follows because $\N\fl$ is odd. The vector is in the
$\bpsi$-isotypic part of the constant-term module: this holds
over $\Frac(R)$, and the integral constant-term module is
$\Z_2$-torsion-free.

Since we are in Case 1 and are working at $p=2$, the level
$\fn$ is prime to $2$. Our auxiliary prime also satisfies
$\fl\nmid2$, so the $2$-primary part of $\fn\fl$ is
the unit ideal. Hence every cusp at this level is
$2$-unramified in the sense of
\cite{silliman}*{Definition~4.15}.
The $\bpsi$-isotypic surjectivity in
\cite{silliman}*{Theorem~1} therefore supplies, for $k$
sufficiently large, a form
\[
G_k^\fl\in M_k(\fn\fl,R,\bpsi)
\]
whose constant-term vector at all cusps is $b_k$.
Set
\[
f_k^\fl=W_1(\bpsi_\fl,1)V_{k-1}
       -x^{-1}W_k(\bpsi_\fl,1)-2DG_k^\fl.
\]
Using \eqref{e:fh}, we obtain
\begin{align*}
f_k^\fl={}&f_k+f_k|_\fl
 +D(G_k+G_k|_\fl-2G_k^\fl)\\
 &+W_1(\bpsi,1)|_\fl\,(V_{k-1}-V_{k-1}|_\fl).
\end{align*}
The term $G_k+G_k|_\fl-2G_k^\fl$ has zero constant terms.
Lemma~\ref{l:vkl}, with its precision increased to absorb the
fixed denominators in the constant terms of $W_1$, shows that
the last term has arbitrarily highly $2$-divisible constant
terms. Thus the constant terms of $f_k^\fl$ are divisible by
$2^M$ at every character, for $N$ sufficiently large.

The same characterwise-to-group-ring passage as in
\cite{dk}*{Theorem 8.15} gives an integral group ring form
$h_k\in M_k(\fn\fl,R,\bpsi)$ such that
\[
f_k^\fl-\frac{2^M}{|G|}h_k
\]
is cuspidal. Choose $M$ large enough that
\[
\frac{2^M}{|G|\,2D}\in R
\]
and put
\[
H_k(\bpsi_\fl)
 =G_k^\fl+\frac{2^M}{|G|\,2D}h_k.
\]
Multiplying the preceding cuspidal form by $x$ proves the
asserted identity, since $2D=\Theta^\#/2^n$.
Its nonconstant Fourier coefficients are integral as in
\cite{dk}*{Theorem 8.17}: the only possible denominators in
$W_1V_{k-1}$ come from the constant terms of $W_1$, and these
are absorbed by the highly divisible nonconstant coefficients
of $V_{k-1}$. Its constant terms vanish.
\end{proof}

\begin{remark}
The symbols $W_j(\bpsi_\fl,1)$ refer to the imprimitive
Eisenstein series in \eqref{e:fh};
$H_k(\bpsi_\fl)$ and $\widetilde F_k(\bpsi_\fl)$ are the forms
constructed in the proof. We divide the constant-term vector
of $G_k+G_k|_\fl$ by $2$, not the whole modular form.
\end{remark}

Recall that our ring $R$ is associated to a $\Gal(\overline{\Q}_p / \Q_p)$-conjugacy class of characters $\Phi$ of $G' \subset G$.  Let $\chi \in \Phi$ and define  \[ \fP'' = \prod_{\fp \mid \fP', \chi(\fp) \neq 1} \!\!\!\!\!\! \fp. \]
This is of course independent of the choice of $\chi \in \Phi$.
 The proof of Corollaries 8.19 and 8.21 of \cite{dk} applies directly to deduce the following from Corollary~\ref{c:cuspc1}.
\begin{corollary} \label{c:cuspfk} Suppose we are in Case 1. For $k > 1$ and $ k \equiv 1 \pmod{(p-1)p^N}$ with $N$ sufficiently large, there exists a cuspidal group ring family $F_k(\bpsi_{\fl}) \in S_k(\fn \fP' \fl, R, \bpsi)^{p-\ord}$ such that 
\begin{equation} \label{e:fcong}
F_k(\bpsi_{\fl}) \equiv  \begin{cases}
xW_1(\bpsi_{\fl}, 1) - W_{k}(\bpsi_{\fl}, 1_p)  \pmod{x \Theta^{\#}/2^n} & \text{if } \fP' \neq \fP'',
 \\
W_1(\bpsi_{\fl p}, 1)  \pmod{\Theta^{\#}/2^n} & \text{if } \fP' = \fP''.
\end{cases}
\end{equation}
In Case 2 there exists a cuspidal group ring family $F_k(\bpsi) \in S_k(\fn \fP', R, \bpsi)^{p-\ord}$ such that 
\[
F_k(\bpsi) \equiv W_1(\bpsi_{\fP\fP''}, 1) \pmod{\Theta^{\#}/2^n}.
\]
\end{corollary}

Lemma~\ref{l:einmax} stated above follows directly from Corollary~\ref{c:cuspfk} as in \cite{dk}*{Lemma 8.22}. {
We are now ready to construct the homomorphism $\varphi$ described at the beginning of \S\ref{s:modforms}.}

\begin{theorem} \label{t:heckehom} In both Cases 1 and 2, there exists a non-zerodivisor $x \in R$, an $R/(x\Theta^{\#}/2^n)$-algebra $W$, and a surjective $R$-algebra homomorphism $\varphi: \widetilde{\T}_{\mathfrak{m}} \lra W$ satisfying the following properties.  We have $x = 1$ unless we are in Case 1 and $\fP' \neq \fP''$.
\begin{itemize}
\item The structure map $R/(x\Theta^{\#}/2^n) \lra W$ is an injection. 
\item The restriction of $\varphi$ to $\T_{\mathfrak{m}}$ takes values in $R/(x\Theta^{\#}/2^n) \subset W$. More precisely, 
\begin{align}
\varphi(S(\fa)) &= \bpsi(\fa) \text{ for } {(\fa, \widetilde{\fn} \fP') = 1,} \nonumber \\
\varphi(T_{\fq}) &= \epsilon_{\cyc}^{k-1}(\fq) + \bpsi(\fq) \text{ for } \fq \nmid \tilde{\fn} p, \label{e:qcong} \\
\varphi(U_{\fp}) &= 1 \text{ for } \fp \mid \fP.\label{e:pcong}
\end{align}
\item In Case 1 we have
\begin{equation} \label{e:lcong} \varphi(U_{\fl}) = \begin{cases}
\epsilon_{\cyc}^{k-1}(\fl) & \text{if } \fP' \neq \fP'' \\
1 & \text{if } \fP' = \fP''.
\end{cases}
\end{equation}
\item Let 
\[
U = \prod_{\fp \mid \fP'} (U_{\fp} - \bpsi(\fp)) \in \widetilde{\T}_{\mathfrak{m}}.
\]
If $y \in R$ and $\varphi(U)y = 0$ in $W$, then $y \in (\Theta^{\#}/2^n)$. 
\end{itemize}
\end{theorem}

\begin{proof}
This result is proved using the forms $F_k(\bpsi_\fl)$ and $F_k(\bpsi)$ exactly as in the proof of 
\cite{dk}*{Theorem 8.23}, {with the parity of $p$ playing no role}.  The only extra feature is that in Case 1, we must verify that 
\[ F_k(\bpsi_\fl)|_{U_\fl} \equiv \begin{cases} 
\N\fl^{k-1} F_k(\bpsi_\fl)  \pmod{x\Theta^\#/2^n} & \text{if } \fP' \neq \fP'' \\
F_k(\bpsi_\fl) \pmod{\Theta^\#/2^n} & \text{if } \fP' = \fP''.
\end{cases} \]
This follows directly from (\ref{e:fcong}).  Note that the first case uses (\ref{e:epscong}).
\end{proof}

Define
\[
\widetilde I
 = \ker\!\left(
     \varphi:\widetilde{\T}_{\fm}\longrightarrow W
   \right),
\qquad
I=\widetilde I\cap\T_{\fm}.
\]
By Theorem~\ref{t:heckehom}, the restriction of $\varphi$
induces an isomorphism of $R$-algebras
\[
\T_{\fm}/I \xrightarrow{\sim} R/(x\Theta^\#/2^n).
\]
In particular, the structure map $R\longrightarrow\T_{\fm}/I$
is surjective.

\subsection{Galois Representation} \label{ss:gr}

It remains now to construct the Galois representation $\rho \colon G_F \longrightarrow \GL_2(K).$ We do this on each factor of $K = \prod_{f \in \overline{\cM}} K_f$ in the decomposition (\ref{e:kdec}).
To each $f$ (say with nebentypus $\psi$),  Hida and Wiles \cite{wilesordrep}*{Theorems 1 and 2} attach a continuous irreducible Galois representation 
\begin{equation} \label{e:rhofdef}
\rho_f\colon G_F \longrightarrow \GL_2(E)
\end{equation}
satisfying the following properties:
\begin{itemize}
\item $\rho_f$ is unramified outside $\widetilde{\fn} p$. 
\item For all primes $\fq \nmid \widetilde{\fn} p$, the characteristic polynomial of $\rho(\Frob_{\fq})$, where $\Frob_\fq$ denotes arithmetic Frobenius, is given by 
\begin{equation} \label{e:charpoly}
\text{char}(\rho_f(\Frob_{\fq}))(x) = x^2 - c(\fq, f) x + \psi(\fq) \epsilon_{\cyc}^{k-1}(\fq),
\end{equation}
\item  For all $\fp \mid p$, we have
\begin{equation} \label{e:ut}
\rho_f|_{G_{\fp}} \sim \mat{ \eta_{\fp}}{ 0}{ *}{ \xi_\fp}
\end{equation}
where $\eta_{\fp}, \xi_\fp : G_{\fp} \rightarrow E^*$ are characters with $\eta_\fp$ unramified and satisfying 
$\eta_{\fp}(\Frob_\fp) = c(\fp, f_p)$. Here $\Frob_\fp$ denotes a lifting to $G_{\fp}^{\ab}$ of the Frobenius element on the maximal unramified extension of $F_{\fp}$.
\end{itemize}

We now define
\[
\widetilde{\rho}_f = \rho_f \otimes_E K_f : G_F \longrightarrow \GL_2(K_f).
\] 

We need to show that $\widetilde\rho_f$ satisfies an analogue
of \eqref{e:ut} at $\fl$ in Case 1. We first justify the
weight-independent choice of $\fl$ made in
\S\ref{ss:hecke}.

\begin{lemma}\label{l:nonscalarprime}
Assume Case 1. There is a set of positive density of primes
$\fl\nmid\fn p$ whose Frobenius in $G$ is $c$, such that
$\rho_f(\Frob_\fl)$ is nonscalar for every $p$-ordinary
cuspidal Hilbert newform $f$ of parallel weight $k>1$,
level dividing $\fn\fP'$, and nebentypus an odd character
of $G$. In particular, the same prime $\fl$ may be used
for all such weights.
\end{lemma}

\begin{proof}
Hida's finiteness and control theorems at fixed tame level give
a bound $D$, independent of $k$, for the degrees
$[E_f:\Q_p]$ of the $p$-adic coefficient fields of these
forms; see \cite{hida-ordinary}*{\S3, Theorem~3.2}.
Indeed, the relevant ordinary Hecke algebras are finite over
their weight algebras, so their characteristic-zero arithmetic
fibers have uniformly bounded dimensions. The finite-order
characters occurring here take values in a fixed coefficient
field. These results apply also at $p=2$; see
\cite{hida-ordinary}*{Proposition~2.1}.
For this bound, one may impose a fixed auxiliary tame level
to meet the smoothness hypothesis in that reference; the
original ordinary spaces inject into the spaces at the
refined level.

For each $f$, choose a $G_F$-stable lattice in $E_f^2$,
and let $\cO_f$ be the ring of integers of $E_f$.
Reducing the resulting integral representation gives
\[
\rho_{f,p^2}:G_F\longrightarrow
\GL_2(\cO_f/p^2\cO_f).
\]
Let $L_f$ be the fixed field of its kernel. Then $L_f/F$
is Galois, unramified outside the primes dividing $\fn p$,
and
\[
[L_f:F]
\leq \bigl|\GL_2(\cO_f/p^2\cO_f)\bigr|
\leq p^{8D}.
\]
Hermite--Minkowski finiteness therefore implies that only
finitely many distinct fields $L_f$ occur, as $f$ and
$k$ vary. Let $L$ be their compositum with $H$.
This is a finite Galois extension of $F$, independent of
$k$.

Fix a complex conjugation $c_\infty\in G_F$. By the
\u{C}ebotarev density theorem, the primes $\fl\nmid\fn p$
whose Frobenius conjugacy class in $\Gal(L/F)$ is that of
$c_\infty|_L$ form a set of positive density. Such a prime
has Frobenius $c$ in $G$. Since every $\rho_f$ is
totally odd, it also satisfies
\[
\tr\bigl(\rho_f(\Frob_\fl)\bigr)
\equiv \tr\bigl(\rho_f(c_\infty)\bigr)
=0\pmod{p^2\cO_f}.
\]
If $\rho_f(\Frob_\fl)=aI_2$, preservation of the lattice
implies $a\in\cO_f^\times$. But
$2a\notin p^2\cO_f$, contradicting the trace congruence.
This argument includes $p=2$, where the reduction is
modulo $4$.
\end{proof}

Every $f\in\overline{\cM}$ whose level is not divisible by
$\fl$ has level dividing $\fn\fP'$ and is $p$-ordinary.
Thus the prime already chosen in \S\ref{ss:hecke} satisfies
the nonscalarity condition needed below, independently of the
subsequent choice of $k$.

\begin{lemma}  \label{l:localell}
In Case 1
we have
\begin{equation} \label{e:ut2}
\widetilde{\rho}_f|_{G_{\fl}} \sim \left(\begin{array}{cc} \eta_{\fl} & 0 \\ * &  \xi_\fl \end{array} \right)
\end{equation}
where $\eta_{\fl}, \xi_\fl : G_{\fl} \rightarrow K_f^*$ are characters with $\eta_\fl$ unramified and satisfying $\eta_{\fl}(\Frob_\fl) = U_\fl$. Here $K_f$ is as in (\ref{e:kfdef}).
 \end{lemma}
 
 \begin{proof}  If $\fl$ divides the level of $f$, then since $\fl^2$ does {\em not} divide the level of $f$ and the nebentypus $\psi$ is unramified at $\fl$, 
 the desired form (\ref{e:ut2}) follows from a result of Carayol \cite{carayol}.  Indeed, the
 local automorphic representation associated to $f$ is an unramified twist of Steinberg (see \cite{casselman} or \cite{lw}*{Proposition 2.8(2)}), and the Galois representation corresponding to this under local Langlands is described in  \cite{carayol}*{\S6.6}.
 
 Now suppose that $\fl$ does not divide the level of $f$. Then the representation $\rho_f$ is unramified at $\fl$, and equation (\ref{e:charpoly}) is satisfied for $\fq = \fl$.  We must therefore simply show that $\widetilde{\rho}_f(\Frob_\fl)$ is conjugate to a matrix of the form $\mat{U_\fl}{0}{*}{*}$ over $\GL_2(K_f)$.  This follows from the following elementary result in linear algebra.
 \end{proof}
 
 \begin{lemma} \label{l:linalg}
 Let $E$ be a field and $A \in M_2(E)$.  Suppose that $A$ has minimal polynomial $g(x)= x^2 - cx + d$
 and let $K = E[U]/g(U)$.  Then $A$ is conjugate via $GL_2(K)$ to a matrix of the form  $\mat{U}{0}{*}{*}$.
 \end{lemma}
 
 \begin{proof}  Let $U' = c - U$ denote the other root in the factorization $g(x) = (x - U)(x-U')$ over $K$.
 We must simply show that $A$ has an eigenvector  $v \in K^2$ with eigenvalue $U'$ such that $v$ can be completed to a basis of $K^2$.  The rational canonical form shows that $A$ is conjugate via $\GL_2(E)$ to $\mat{0}{1}{-d}{c}$.  Let $\{v_1, v_2\}$ be the associated basis.  Then $v = v_1 + v_2 U'$ is the desired vector.
 \end{proof}
 
 \begin{remark} Lemmas~\ref{l:localell} and \ref{l:linalg} make it clear why we must insist that $\rho_f(\Frob_\fl)$ is not scalar in Case 1 when $\fl$ does not divide the level of $f$---a scalar matrix is only conjugate to itself.  
 Interestingly, the situation where $\rho(\Frob_\fl)$ has one eigenvalue with multiplicity 2 but is not scalar is conjectured to  never occur. It is curious and somewhat fortuitous that we can still handle this conjecturally vacuous case via Lemma~\ref{l:linalg}, and that we {\em cannot} handle via similar means the case where $\rho(\Frob_\fl)$ is scalar (which {\em does} occur), but that we can avoid the scalar case by the \u{C}ebotarev argument of Lemma \ref{l:nonscalarprime}.
 \end{remark}
 
We define $\rho$ by putting the $\widetilde{\rho}_f$ together:
\begin{equation} \label{e:rho}
\rho = \prod_{f \in \overline{\cM}} \widetilde{\rho}_f \colon G_F \longrightarrow \GL_2(K).
\end{equation}
 By construction, the representation $\rho$ satisfies the following properties:
\begin{itemize}
\item $\rho$ is unramified outside $\widetilde{\fn} p $.
\item For all primes $\fq \nmid \widetilde{\fn} p $, the characteristic polynomial of $\rho(\Frob_{\fq})$ is given by 
\begin{equation} \label{e:charpolyq}
\text{char}(\rho(\Frob_{\fq}))(x) = x^2 - T_{\fq} x + \bpsi(\fq)\epsilon_{\cyc}^{k-1}(\fq).
\end{equation}
\item For all primes $\fp \mid p$, and for $\fp = \fl$ in Case 1, we have 
\[
\rho|_{G_{\fp}} \sim \mat{ \eta_{\fp}}{0}{*}{\xi_\fp},
\]
where $\eta_{\fp}, \xi_\fp \colon G_{\fp} \lra \widetilde{\T}^*$ are characters with $\eta_\fp$ unramified and  
$\eta_{\fp}(\Frob_\fp) = U_{\fp}$. 
\end{itemize}
Next we show that $\rho$ satisfies all the properties required by Theorem~\ref{t:rl}, as specified further in Theorem~\ref{t:galrep}.
 {Recall that $(\Sigma, \Sigma')$ are the minimal pair defined in Lemma~\ref{l:minimal}.  
 In particular
 \[ \Sigma =  \begin{cases}
\{ \fl\} & \text{ Case 1}, \\
\{ \fp \mid p: \fp \in S_{\ram} \} & \text{ Case 2}.
\end{cases} \]
 Let
\[
\cP = \{\fp \mid p, \fp \not \in \Sigma\} = \{\fp \mid p \colon \fp \not\in S_{\ram}\}.
\]
}

\begin{theorem}  Let $\chi = \epsilon_{\cyc}^{k-1}$. The representation $\rho$ satisfies the following.
\begin{enumerate}
\item For  $\sigma \in G_F,$  the characteristic polynomial $P_{\rho(\sigma)}(x)$ lies in $\T[x]$.  Furthermore we have
\[
P_{\rho(\sigma)}(x) \equiv (x - \chi(\sigma))(x - \bpsi(\sigma)) \pmod{I}.
\]
\item Let $K_0 = 
\red(K)$ denote the maximal reduced quotient of $K$. Write $K_0 = \prod_{i=1}^{m} k_i$ as a product of fields.  For every projection $K \rightarrow K_0 \rightarrow k_i$, the projection of $\rho$ to $\GL_2(k_i)$ is an irreducible representation of $G_F$ over $k_i$.
\item  For each $v \in \Sigma$, we have $\xi_v \equiv \bpsi|_{G_v} \pmod{\widetilde{I}}$.

\item   For all $v \in \cP$,  we have 
$(\xi_v)|_{I_v} = \chi|_{I_v}$.  
\end{enumerate}
\end{theorem}

\begin{proof}
We verify the statements in turn. 

\bigskip
(1) Since $\varphi(T_{\fq}) = \chi(\fq) + \bpsi(\fq)$ by (\ref{e:qcong}), we have 
\[
\text{char}(\rho(\Frob_{\fq}))(x) \equiv (x- \epsilon_{\cyc}^{k-1}(\Frob_{\fq}))(x - \bpsi(\Frob_{\fq})) \pmod{I}
\]
for all $\fq \nmid \widetilde{\fn} p$. Since the $\Frob_\fq$ are dense in $G_F$ by \u{C}ebotarev, we have 
$\text{char}(\rho(\sigma))(x) \in \T[x]$ for all $\sigma \in G_F$ and 
\[
\text{char}(\rho(\sigma))(x) \equiv (x - \epsilon_{\cyc}^{k-1}(\sigma)) (x - \bpsi(\sigma)) \pmod{I}.
\]

\bigskip
(2) The projections of $\rho$ to the field factors of the reduced quotient of $K$ are just the representations $\rho_f$, which are irreducible.

\bigskip
(3) For $v \in \Sigma$, we have  $\xi_v = \chi \bpsi \eta_v^{-1}$.
  We must show that $\chi|_{G_v} \equiv \eta_v \pmod{\widetilde{I}}$.  In Case 1, we have $\Sigma = \{\fl\}$.  
  Then $\eta_\fl$ and $\chi$ are both unramified at $\fl$.  Furthermore $\eta_\fl(\Frob_\fl) = U_\fl$.  
  If $\fP' \neq \fP''$, then $U_{\fl} \equiv \chi(\fl) \pmod{\widetilde{I}}$ by (\ref{e:lcong}), so we are done. 
  If $\fP' = \fP''$ then (\ref{e:lcong}) states $U_{\fl} \equiv 1 \pmod{\widetilde{I}}$.  
  But in this case we have $x=1$ so $\Theta^\#/2^n \in \widetilde{I}$, whence the congruence $\chi \equiv 1 \pmod{\Theta^\#/2^n}$  specified in (\ref{e:epscong}) yields the desired result. 
   In Case 2, we have $\Sigma = \{\fp \mid \fP\}$ and a similar argument holds. 
    We again have $x=1$ so $\chi \equiv 1 \pmod{\widetilde{I}}$, and (\ref{e:pcong}) shows that $\eta_\fp \equiv 1 \pmod{\widetilde{I}}$ as well.

\bigskip
(4) For $v \in \cP$, we again have $\xi_v = \chi \bpsi \eta_v^{-1}$, which equals $\chi$ on $I_v$ since $\eta_v$ and $\bpsi$ are unramified at $v$.
\end{proof}

We conclude this section by noting that in the application of 
Theorems~\ref{t:rl} and \ref{t:galrep} we choose for $\fp \in \cP$ the elements $\sigma_\fp \in G_\fp$ to be lifts of $\Frob_\fp \in G_\fp^{\ab}$.  Then 
\[ \xi_\fp(\sigma_\fp) - \chi(\sigma_\fp) = \epsilon_{\cyc}^{k-1}(\sigma_\fp) (\bpsi(\fp) U_\fp^{-1} - 1) \]
is a unit multiple of $(U_\fp - \bpsi(\fp))$.  Hence the implication
\[  y \in R \text{ and } \Big(\prod_{v \in \cP} (\xi_v(\sigma_v) - \chi(\sigma_v))\Big) y \in \widetilde{I} \Longrightarrow 
y \in (\Theta^\#/2^n)R \]
required by Theorem~\ref{t:galrep} follows from the last bullet point of Theorem~\ref{t:heckehom}. This concludes the proof that the Hecke algebras $\T_\fm \subset \widetilde{\T}_\fm$, the homomorphism $\varphi$, and the Galois representation $\rho$ satisfy the properties used in \S\ref{s:sr}.

\section{Construction of Generalized Ritter--Weiss Modules}\label{s:constr-rw}

In \S\S\ref{s:constr-rw}--\ref{s:duality}, we construct our generalized Ritter--Weiss modules $\nabla_{S}^{T}(H)$ and establish the properties stated in \S\ref{s:properties-rw}. 
{
Throughout \S\ref{s:constr-rw}, $H/F$ is a finite Galois extension of number fields with Galois group $G$. In Proposition~\ref{prop:quad-res}, which demonstrates the local quadratic presentation of $\nabla_S^T(H)$ under appropriate conditions, we assume that $G$ is abelian. 
The assumptions relevant to our applications to Brumer--Stark that $F$ is totally real,  $H$ is CM, and $G$ is
abelian  are imposed only in \S\S\ref{s:cmcase}--\ref{ss:changing-sets}.}

We will use the language of class formations; 
for the benefit of the reader, we recall this terminology in \S\ref{ss:class-formations}.

\subsection{Class modules and class formations}\label{ss:class-formations}
We recall the basic properties of class modules following \cite{nsw}. Let $G$ be a finite group. A discrete $G$-module $C$ is a {\em class module} if for every subgroup $U \subset G$, we have
\begin{enumerate}
\item $H^1(U, C) = 0$,
\item $H^2(U, C)$ is cyclic of order $|U|$.
\end{enumerate}

For a class module $C$, a choice of generator $\gamma_G \in H^2(G,C)$ is known as a {\em fundamental class}. For any $G$-module $C$, an element $\gamma \in H^2(G,C)$ defines a 2-extension of the form \begin{equation}\label{e:two-extn} 0 \lra C \lra C(\gamma) \lra \Z[G] \lra \Z \lra 0. \end{equation}
{Below $\hat{H}$ refers to Tate cohomology groups.}
\begin{theorem}[Tate]\label{thm:class-mod} Let $C$ be a discrete $G$-module and $\gamma \in H^2(G, C)$.  The following are equivalent:
\begin{enumerate}
\item $C$ is a class module with fundamental class $\gamma$.
\item $C(\gamma)$ is $G$-cohomologically trivial.
\item  For all $i \in \Z$ and all subgroups $H \subset G$, the cup product on Tate cohomology  
\[ \cup \res_{G,H}(\gamma) : \widehat{H}^i(H,\Z) \lra \widehat{H}^{i+2}(H,C) \] is an isomorphism.
\end{enumerate} 
\end{theorem}

In particular, if $C$ is a class module for $G$, cup-product with $\gamma$ defines the  {\em Nakayama map}:
\[ G^{\mathrm{ab}} = \widehat{H}^{-2}(G, \Z) \cong \widehat{H}^0(G, C) = C^G/N_G(C). \] The {\em reciprocity map} is then defined to be
\[ \mathrm{rec} \colon C^G \lra C^G/N_G(C) \cong G^{\mathrm{ab}}. \]

Now, consider a profinite group $G$ and a discrete $G$-module $C$. The pair $(G, C)$ is called a {\em class formation} if, for all open subgroups $H' \subset H$ of $G$ with $H'$ normal in $H$, the $H/H'$-module $C^{H'}$ is a class module for $H/H'$ equipped with a choice of fundamental class $\gamma_{H/H'}$. The fundamental classes are required to satisfy certain compatibilities under inflation and restriction (\cite{nsw}, Definition 3.1.8). When working with a fixed class formation, we will abbreviate $C(H/H') := C(\gamma_{H/H'})$.

For a profinite class formation $(G,C)$ and an open subgroup $H$ of $G$, there is a reciprocity map $\mathrm{rec} \colon C^H \lra H^{\mathrm{ab}}$, defined as the inverse limit of the reciprocity maps $C^{H} \lra (H/H')^{\mathrm{ab}}$ over open normal subgroups $H' \subset H$. The reciprocity map induces an isomorphism \[ \mathrm{rec} \colon \widehat{(C^H)}_{\mathrm{norm}} \cong H^{\mathrm{ab}}, \] where $\widehat{(C^H)}_{\mathrm{norm}} := \varprojlim_{H'} C^H/N_{H/H'}(C^{H'})$ is the completion in the norm topology on $C^H$.

We will also consider, for any open normal subgroup $H \subset G$, the 2-extension
\begin{equation}\label{e:norm-cplt-ext} 0 \lra H^{\mathrm{ab}} \lra E({G/H}) \lra \Z[G/H] \lra \Z \lra 0 \end{equation}
corresponding to pushing out the 2-extension \[ 0 \lra C^H \lra C(G/H) \lra \Z[G/H] \lra \Z \lra 0 \] along the map $C^H \lra H^{\mathrm{ab}}$.
The module $E({G/H})$ is topologized using the profinite topology on $H^{\mathrm{ab}}$ and the discrete topology on $\Z[G/H]$.

An \emph{inclusion of class formations}
$(G',C')\longrightarrow(G,C)$ consists of an inclusion of
profinite groups $G'\subset G$, together with homomorphisms
\[
j_U:(C')^{G'\cap U}\longrightarrow C^U
\]
for every open subgroup $U\subset G$.
These maps are compatible with conjugation: for $g\in G'$,
\[
j_{gUg^{-1}}(gx)=g\,j_U(x).
\]
They are also compatible with norms: for open subgroups
$U'\subset U$ with $U'$ normal in $U$,
\[
j_U\circ N_{(G'\cap U)/(G'\cap U')}
=
N_{U/U'}\circ j_{U'}.
\]
We require compatibility with reciprocity,
\[
\operatorname{rec}_U\circ j_U
=
\iota_{U,\mathrm{ab}}\circ
\operatorname{rec}'_{G'\cap U},
\]
where
\[
\iota_{U,\mathrm{ab}}:
(G'\cap U)^{\mathrm{ab}}\longrightarrow U^{\mathrm{ab}}
\]
is induced by inclusion.

Finally, for all open subgroups $H'\subset H$ of $G$
with $H'$ normal in $H$, we require
\[
\operatorname{res}(\gamma_{H/H'})
=
(j_{H'})_*
\bigl(\gamma_{(G'\cap H)/(G'\cap H')}\bigr)
\]
in
\[
H^2\bigl((G'\cap H)/(G'\cap H'),C^{H'}\bigr),
\]
where $\operatorname{res}$ is restriction along the natural
inclusion
\[
(G'\cap H)/(G'\cap H')\hookrightarrow H/H'.
\]

For an open normal subgroup $H \subset G$, let $H' = G' \cap H$. The methods of \cite{rw}*{\S2 and \S4} define a map $C'(G'/H') \lra C(G/H)$ such that the following diagram commutes:
\[\begin{tikzcd}
(C')^{H'} \arrow{r}\arrow{d}{j_H} & C'(G'/H') \arrow{r}\arrow{d} & \Z[G'/H'] \arrow{d} \\
C^H \arrow{r}& C(G/H) \arrow{r} & \Z[G/H].
\end{tikzcd}\]
We obtain a commutative diagram
\[\begin{tikzcd}
(H')^{\mathrm{ab}} \arrow{r}\arrow{d}& E(G'/H') \arrow{r}\arrow{d} & \Z[G'/H'] \arrow{d} \\
H^{\mathrm{ab}} \arrow{r}& E(G/H) \arrow{r} & \Z[G/H]. 
\end{tikzcd}\]

\subsection{Class field theory}\label{ss:cft}

% Now let $F$ be a number field and $H$ a finite Galois extension with $G = \Gal(H/F)$.  To maintain generality, we do not impose further conditions on $F$ and $H$.}
The main results of class field theory can be summarized in the results below. {In the following we fix a place $w$ of $H$ above a place $v$ of $F$.}
\begin{theorem}\label{thm:cft} We have
\begin{enumerate}
\item For $F_v$ a local field, $(G_{F_v}, \bigcup_{H_w/F_v} H_w^*)$ is a class formation.
\item For $F$ a global field, $(G_F, \bigcup_{H/F} \A_H^*/H^*)$ is a class formation.
\item There is an inclusion of class formations \[ (G_{F_v}, \bigcup_{H_w/F_v} H_w^*) \lra (G_{F}, \bigcup_{H/F} \A_H^*/H^*). \]
\end{enumerate}
For $U=G_L$, where $L/F$ is finite and $w$ is the
place selected by the fixed local embedding, the coefficient
map $j_U$ is the natural map
\[
L_w^\times\longrightarrow\A_L^\times/L^\times
\]
obtained by inserting a local element in the $w$-component and putting 1's in all other components.
\end{theorem}
\begin{proof}

 For (1), (2), see \cite{artin-tate}.  For (3), see \cite{cassels-frohlich}*{pp.\ 195--196}. 
\end{proof}

We will also need the following well-known facts about the local and global reciprocity maps.
\begin{theorem}\label{thm:rec} The following hold.
\begin{enumerate}
\item For $H_w$ an archimedean local field, we have an exact sequence
\[ \begin{tikzcd}
0 \ar[r] & (H_w^*)^{\circ} \ar[r] &  H_w^* \ar[r,"{\mathrm{rec}_v}"] &  G_{H_w}^{\mathrm{ab}} \ar[r] & 0, 
\end{tikzcd} \]
where $(H_w^*)^{\circ}$ denotes the connected component of the identity.
\item For $H_w$ a nonarchimedean local field, we have an exact sequence
\[ \begin{tikzcd}
 0 \ar[r] & H_w^* \ar[r, "{\mathrm{rec}_w}"] & G_{H_w}^{\mathrm{ab}} \ar[r] & \widehat{\Z}/\Z \ar[r] & 0. 
 \end{tikzcd} \]
\item For $H$ a number field,  the following diagram commutes:
\[ \begin{tikzcd}
\prod'_{w} H_w^* \arrow{r}\arrow{d}{\mathrm{rec}_w}& \A_H^*/H^* \arrow{d} \\
 \prod'_{w} G_{H_w}^{\mathrm{ab}} \arrow{r}& G_{H}^{\mathrm{ab}}.
\end{tikzcd} \]
Let $S$ be a finite set of places of $H$ containing the archimedean places. 
Restricting the local terms from $H_w^*$ to $\cO_w^*$ for $w \not\in S$,  we have a commutative diagram with exact rows: 
\[ \begin{tikzcd}
\prod_{w \in S} H_w^* \prod_{w \not\in S} \cO_w^* \arrow{r}\arrow{d}{\mathrm{rec}_w}& \A_H^*/H^* \arrow{r}\arrow{d}{\mathrm{rec}}& \Cl_S(H) \arrow{r}\arrow{d}{=}& 0 \\
 \prod_{w \in S} \rec(H_w^*)\prod_{w \not \in S} \mathrm{rec}(\cO_w^*) \arrow{r}& G_{H}^{\mathrm{ab}} \arrow{r}& \Cl_S(H) \arrow{r}& 0. 
\end{tikzcd} \]
\end{enumerate}
\end{theorem}

{For our finite Galois extension $H/F$ with $G = \Gal(H/F)$,}
 the fundamental class $\gamma_G \in H^2(G, \A_H^*/H^* )$ defines the {\em global Tate sequence}
\begin{equation} \label{e:ts}
0 \lra \A_H^*/H^*  \lra V \lra B \lra \Z \lra 0,
\end{equation}
an exact sequence of $\Z[G]$-modules, where $B = \Z[G]$. For $w$ a place of $H$ with decomposition group $G_w \subset G$, the fundamental class $\gamma_{G_w} \in H^2(G_w, H_w^*)$ defines the {\em local Tate sequence}
\begin{equation} \label{e:lts}
0 \lra H_w^* \lra V_w \lra B_w \lra \Z \lra 0,
\end{equation}
an exact sequence of $G_w$-modules, where $B_w = \Z[G_w]$. 
{
Define $W_w = V_w/\cO_w^*$.  
 Let $\Delta G_w$ denote the augmentation ideal of $\Z[G_w]$, and let $I_w \subset G_w$ denote the inertia subgroup at $w$. The norm $\N(I_w) = \sum_{\tau \in I_w} \tau$ defines a map \[ \N(I_w) \colon \Z[G_w/I_w] \longrightarrow \Z[G_w]. \] Let $\sigma_w \in G_w/I_w$ denote the arithmetic Frobenius. 
The module $W_w$ is shown in \cite{rw}*{\S3} to have the following description:
\begin{equation}W_w \cong \{ (x, y) \in \Delta G_{w} \times \Z[G_{w}/I_{w}] \mid \overline{x} = (1 - \sigma_w^{-1})y \}.
\label{e:wisom}
\end{equation}
The map $W_w \rightarrow B_w = \Z[G_w]$ is given by projection onto the first coordinate.
In particular, when $w$ is unramified, projection onto the second coordinate yields an isomorphism $W_w \cong \Z[G_w]$,
and the map \begin{equation} \label{e:wb}
 W_w \cong \Z[G_w] \lra B_w \cong \Z[G_w] \end{equation} is given by multiplication by $(1 - \sigma_w^{-1})$.}

There is a map from the local Tate sequence to the global Tate sequence:
\begin{equation*}\begin{tikzcd}
V_w \arrow{r} \arrow{d} & B_w \arrow{d} \\
V \arrow{r} & B
\end{tikzcd}\end{equation*}
 (see \S\ref{ss:class-formations}). This is a diagram of $\Z[G_w]$-modules, giving rise to a diagram of $\Z[G]$-modules
\begin{equation}\label{e:local-global}
\begin{tikzcd}
\Ind_{G_w}^G(V_w) \arrow{r} \arrow{d} & \Ind_{G_w}^G(B_w) \arrow{d} \\
V \arrow{r} & B,
\end{tikzcd}
\end{equation}
where the induced complex only depends on the place $v$ of $F$ under $w$.

As in Theorem \ref{thm:rec}, there are compatible {\em reciprocity maps}:
\[ \begin{tikzcd}
H_w^* \arrow{r}{\mathrm{rec}_w} \arrow{d}  & G_{H_w}^{\mathrm{ab}} \arrow{d} \\
\A_H^*/H^* \arrow{r}{\mathrm{rec}} & G_{H}^{\mathrm{ab}}. 
\end{tikzcd}\]
If we push out the local and global Tate sequences along these reciprocity maps, we obtain prodiscrete modules $V_{w,2}$ and $V_2$ sitting in a commutative diagram
\[ \begin{tikzcd}
0 \ar[r] & G_{H_w}^{\mathrm{ab}} \ar[r] \ar[d] & V_{w,2}  \ar[d] \ar[r] & B_w \ar[d] \ar[r] & \Z \ar[r]  \ar[d] & 0 \\
0 \ar[r] & G_{H}^{\mathrm{ab}}  \ar[r] &  V_2 \ar[r] & B \ar[r] & \Z \ar[r] & 0,
\end{tikzcd}
 \]

\subsection{The Ritter--Weiss complex} \label{ss:rwc}

For each place $v$ of $F$ fix a place $w \mid v$ of $H$. We will define, for each place $v$ of $F$, a complex of $\Z[G_w]$-modules. Unless mentioned otherwise, all the two term complexes below are in degree 0 and 1 {with the differential map going from degree 0 to degree 1}. 

\subsubsection{Global complex}

The global Tate sequence corresponds to the two-term complex \begin{equation} V \lra B, \end{equation} where $B = \Z[G]$. By Theorems \ref{thm:class-mod} and \ref{thm:cft}, $V$ is a cohomologically trivial $G$-module.  Of course, the same is true for $B$ as well. 

\subsubsection{Local complexes} \label{ss:lc}

{
For the rest of \S\ref{ss:rwc}--\ref{ss:changing-sets}, fix an RW pair $(S,T)$.  We will impose extra conditions as needed.}
\bigskip

\textbf{Case $v \in S$.}
The local Tate sequence corresponds to the two-term complex \begin{equation} \label{e:rwlis}
V_w \lra B_w \end{equation} 
{appearing in (\ref{e:lts}),} where $B_w = \Z[G_w]$. As in the global case, both $V_w$ and $B_w$ are cohomologically trivial $G_w$-modules.

\bigskip
\textbf{Case $v \in T$ (finite)}.
Define $U_w = \cO^*_{w,1}$ (the 1-units) for $w \in T_H$ finite. Consider the two-term complex \begin{equation}U_w \lra 0.\end{equation} Note that $U_{w}$ is $G_w$-cohomologically trivial if $T$ is tame, and $U_w \otimes_{\Z} \Z_{(p)}$ is $G_w$-cohomologically trivial if $T$ is $p$-tame (see \cite{dk}*{Lemma A.4}; {tame and $p$-tame are defined in \S \ref{s:properties-rw} above}).

\bigskip
\textbf{Case $v \in T$ (infinite)}.
Define $U_w = H_w$ for $w \in T_H$ infinite. Consider the two-term complex
\begin{equation} \label{e:rwlit}
U_w \lra 0.
\end{equation}
This maps to $V_w \lra B_w$ (and hence to $V \lra B$) via $H_w \xrightarrow{\exp} H_w^* \lra V_w$. It is easy to see that $U_w$ ($= \R$ or $\C$) is cohomologically trivial for $G_w$ ($=  \{ 1 \} $ or $\Z/2\Z$).

\bigskip
\textbf{Case $v \notin S \cup T$ (ramified).}
{
Consider the module $W_w = V_w/\cO_w^*$ and the Ritter-Weiss description stated in (\ref{e:wisom}).
There are two maps \begin{align*} \pi_1, \pi_2 \colon W_w \lra \Z[G_w],
\pi_1((x,y)) = x, \qquad \pi_2((x,y)) = \N(I_w)y.
\end{align*}
{Composing with $V_w \rightarrow W_w$,} we obtain two maps $\pi_1, \pi_2 \colon V_w \lra \Z[G_w]$.}
 Consider the two-term complex
\begin{equation}
V_w \xrightarrow{(\pi_1, \pi_2)} \Z[G_w]^2.
\end{equation}
The terms of this complex are $G_w$-cohomologically trivial, and it has \[ H^0 = \cO_w^*, \qquad  H^1 = W_w^* = \Hom_{\Z}(W_w, \Z) \](see \cite{dk}*{(138)--(139)}).

\bigskip

\textbf{Case $v \notin S \cup T$ (unramified).}
When $H_w/F_v$ is unramified, the complex \begin{equation}\cO_w^* \lra 0 \label{e:owo}
\end{equation} is quasi-isomorphic to the complex \begin{equation}V_w \lra W_w, \label{e:vw}
\end{equation} and the terms of either complex are $G_w$-cohomologically trivial (see \cite{dk}*{Lemma A.4}).

\subsection{Construction of the Ritter--Weiss complex} \label{ss:construction}

Let $S'$ be a finite set of  places of $F$ unramified in $H$ such that $S'$ is disjoint from $S \cup T$. Let \[ S'' = S \cup S' \cup \{v \colon v \notin (S \cup T), v \text{ ramified in } H\}.\]  We define \begin{align*}
V_{\mathrm{\loc}} &:= \prod_{v \in S''} \Ind_{G_w}^G(V_w) \times \prod_{v \in T} \Ind_{G_w}^G(U_w) \times \prod_{v \notin (S'' \cup T)} \Ind_{G_w}^G(\cO_w^*), \\
B_{\mathrm{\loc}} &:=  \prod_{v \in S} \Ind_{G_w}^G(B_w)  \times \prod_{v \in S'} \Ind_{G_w}^G(W_w) \times \prod_{v \in (S'' \setminus (S \cup S'))} \Ind_{G_w}^G(\Z[G_w]^2),
\end{align*}
and the complexes
\begin{align*}
C^*_{\mathrm{loc}} &: V_{\mathrm{loc}} \lra B_{\mathrm{loc}}, \\
C^*_{\mathrm{global}} &: V \lra B.
\end{align*}
Note that $C^*_{\mathrm{loc}}$ is simply the product, over all places $w$ of $H$, of the inductions of the local complexes defined above, each of which comes with a map to $C^*_{\mathrm{global}}$ (see (\ref{e:local-global})). 
The set $S'$ determines which of the two quasi-isomorphic complexes to use for unramified places $v \notin S \cup T$.  Thus we have a map \begin{equation}
\theta \colon C^*_{\mathrm{loc}} \lra C^*_{\mathrm{global}}.
\end{equation}

{
\begin{remark} \label{r:blocfree} By definition, the module $B$ in Tate's sequence (\ref{e:ts}) is $\Z[G]$.
Reviewing the definitions and applying the observation preceding (\ref{e:wb}) that $W_w \cong \Z[G_w]$ when $w$ is unramified, we see that $B_{\mathrm{loc}}$ is free of finite rank over $\Z[G]$. 
\end{remark}

\begin{remark} \label{r:ct}
The module $V$ in Tate's sequence (\ref{e:ts}) is $G$-cohomologically trivial.  By our construction in \S\ref{ss:lc}, as we described case by case, the module $V_{\mathrm{loc}}$ is $G$-cohomologically trivial if $T$ is tame.  The same is true after tensoring with $\Z_{(p)}$ if $T$ is $p$-tame.
\end{remark}
}

We define the {\em Ritter--Weiss complex}:
\begin{equation} C^*_{\mathrm{RW}, S'} = \mathrm{cone}(\theta)[-1]. \end{equation}

\begin{lemma}
Up to canonical quasi-isomorphism, $C^*_{\mathrm{RW}, S'}$ is independent of the auxiliary set $S'$.
\end{lemma}
\begin{proof}
This follows immediately from the fact that the two complexes (\ref{e:owo}) and (\ref{e:vw}), which the set $S'$ chooses between, are quasi-isomorphic.
\end{proof}

The {\em generalized Ritter--Weiss module} is defined to be \begin{equation} \label{e:rwm}
\nabla_S^T(H) := H^1(C^*_{\mathrm{RW}}).\end{equation}
Define the {\em logarithmic $(S,T)$-units}: \begin{equation} U_S^T(H) := \left\{ (x_w,y) \in \Big(\prod_{w \in (T_H)_{\infty}} H_w\Big) \times \cO_{H, S,T}^* \colon \exp(x_w) = \sigma_w(y) \right\}. \end{equation}

\begin{lemma} \label{l:nablaseq}
The cohomology of $C^*_{\mathrm{RW}}$ is
\begin{align*}
H^0(C^*_{\mathrm{RW}}) & = U_S^T(H), \\
H^1(C^*_{\mathrm{RW}}) &=\nabla_S^T(H),\\
  H^2(C^*_{\mathrm{RW}}) &= \coker\Big(\prod_{w \in S_H} \Z \times \prod_{\substack{w \notin (S \cup T) \\ w \text{ ramified}}} W_w^* \lra \Z\Big).
  \end{align*} 
The map $W_w^* \rightarrow \Z$ is the dual of the map $\Z \rightarrow W_w$ given by \[ a \mapsto (0,a \cdot \N(G_w/I_w)), \quad \N(G_w/I_w) = \sum_{\tau \in G_w/I_w} \tau. \]
If $(S,T)$ is complete (see \S\ref{s:properties-rw}), we have 
\[ H^2(C^*_{\mathrm{RW}}) = \begin{cases} \Z & \text{ if } S = \emptyset,  \\ 0 & \text{ otherwise.} \end{cases}\]
  
There is a short exact sequence
\begin{equation} \label{e:nablaseq}
  0 \lra \Cl_S^T(H) \lra \nabla_S^T(H) \lra X \lra 0 \end{equation}
 where \begin{equation} \label{e:xdefgen}
  X = \ker\Big(\prod_{w \in S_H} \Z \times \prod_{\substack{w \notin (S \cup T) \\ w \text{ ramified}}} W_w^* \lra \Z\Big). 
  \end{equation}
{If $(S,T)$ is complete we have $X = X_{H,S}$, the set of degree 0 divisors on the places of $H$ above $S$.}

\end{lemma}
\begin{proof}
We may assume $S' = \emptyset.$
There is an exact triangle 
\[ 
C^*_{\mathrm{loc}} \lra C^*_{\mathrm{global}} \lra C^*_{\mathrm{RW}}[1].
\] 
Computing cohomologies of the first two complexes, the associated long exact sequence in cohomology becomes
\[ \hspace{-1.3in} \begin{tikzcd}
\  & 0 \ar[r] & H^0(C^*_{\mathrm{RW}}) \ar[dll, out=-10, in=170] & \  \\ 
  \displaystyle\prod_{w \in S_H} \!\! H_w^* \times \!\!\! \prod_{w \in T_H} \!\! U_{w} \times \!\!\!\!\!\!\!\! \prod_{w \notin (S_H \cup T_H)} \!\!\!\!\!\!\!\!\! \cO_{w}^* \ar[r] & \A_H^*/H^*  \ar[r]
 & H^1(C^*_{\mathrm{RW}}) \ar[dll, out=-10, in=170] & \  \\ 
  \displaystyle\prod_{w \in S_H} \Z \times \!\!\!\!\!\prod_{\substack{w \notin (S \cup T) \\ w \text{ ramified}}} \!\!\!\!\!\! W_w^* \ar[r] &  \Z   \ar[r]
 & H^2(C^*_{\mathrm{RW}}) \ar[r] & 0. 
 \end{tikzcd} \]
The calculation of $H^*(C^*_{\mathrm{RW}})$ follows, as does the short exact sequence (\ref{e:nablaseq}).  The stated calculations under the assumption that {$(S,T)$ is complete} are immediate.
\end{proof}

\begin{lemma}\label{lem:coinv}
Let $(S,T)$ be a complete RW pair with $S\neq\emptyset$ and $T$ tame. Let $K$ be a subfield of $H$ containing $F$ such that $K/F$ is Galois. Then \[ \nabla_S^T(H)_{\Gal(H/K)} = \nabla_S^T(K), \] where the left side denotes $\Gal(H/K)$-coinvariants.

The same result holds after tensoring with $\Z_{(p)}$ if $T$ is only assumed to be $p$-tame rather than tame.
\end{lemma}
\begin{proof}
First, note that the Ritter--Weiss complex for $K$ is isomorphic to the coinvariants of the complex for $H$: \begin{equation} \label{e:coinv} C^*_{\mathrm{RW}}(H)_{\Gal(H/K)} = C^*_{\mathrm{RW}}(K).\end{equation}
This follows from Lemma \ref{lem:coinflation} applied to the local and global Tate sequences 
{
(note that the proof of the lemma is independent of the results of this section).}

By definition, {$\nabla_S^T(H)$ is equal to the cohomology of the sequence
\[
V_{\mathrm{loc}} \lra V \oplus \prod_{v \in S} \Ind_{G_w}^G(B_w)  \lra B.
\]
}
 
Now suppose $S \neq \emptyset$ and fix $v_0 \in S$. Since $\Ind_{G_{v_0}}^G(B_{w_0}) \lra B$ is an isomorphism, 
{
there is a canonical isomorphism
\[ \ker(V \oplus \prod_{v \in S} \Ind_{G_w}^G(B_w)  \lra B) \cong V \oplus \prod_{v \in S\setminus v_0} \Ind_{G_w}^G(B_w)
\]
obtained by forgetting the component at $v_0$.
We therefore obtain the following presentation of  $\nabla_S^T(H)$:
\[
 V_{\mathrm{loc}} \lra  V \oplus \prod_{v \in S \setminus v_0} \Ind_{G_w}^G(B_w)  \lra \nabla_S^T(H) \lra 0.
\]
In view of (\ref{e:coinv}), the result now follows from right exactness of taking co-invariants.}
\end{proof}

\subsection{Quadratic presentation} \label{s:qp}
{Let $R$ be a commutative ring.  An $R$-module $M$ is quadratically presented if there is an exact sequence \[ R^n \lra R^n \lra M \lra 0. \] 
The $R$-module $M$ is called locally quadratically presented if the localization $M_\fp$ is quadratically presented over $R_\fp$ for every prime ideal $\fp \subset R$.
For the remainder of \S\ref{s:constr-rw} we assume that $G$ is abelian so that $R=\Z[G]$ is a commutative ring.
In this subsection we prove that the Ritter--Weiss  module $\nabla_S^T(H)$ is locally quadratically presented under appropriate conditions.}

\begin{lemma} \label{l:sp}
We have
\begin{enumerate}
\item $V^{\mathrm{loc}} \lra V$ is surjective when $\cup_{v \in S'} G_v = G$ and $\Cl_{S'}^T(H) = 0$.
\item $B^{\mathrm{loc}} \lra B$ is surjective when $\cup_{v \in S'} G_v = G$ and $S \neq \emptyset$.
\end{enumerate}
\end{lemma}
\begin{proof}
See \cite{rw} Page 162 or \cite{dk} \S A.1.
\end{proof}

Now suppose that  $S \neq \emptyset$ and that $S'$ is sufficiently large to satisfy both conditions in Lemma~\ref{l:sp}(1). Consider {the complex} \begin{equation} \label{e:vtbt} V^{\theta} \lra B^{\theta}, \end{equation} where 
\[ V^{\theta} = \ker(V_{\mathrm{loc}} \lra V), \qquad B^{\theta} = \ker(B_{\mathrm{loc}} \lra B). \] 
{
Since the shifted cone of a degreewise surjective map of complexes is
quasi-isomorphic to its kernel, the complex (\ref{e:vtbt}) is
quasi-isomorphic to $C^*_{\mathrm{RW}}$. }

\begin{corollary} \label{c:uvbn}
{Suppose $S\neq\emptyset$, $S'$ is sufficiently large
to satisfy both conditions in Lemma~\ref{l:sp}(1), and that $T$ is torsion-free and tame (see \S\ref{s:rwm} and \S\ref{s:properties-rw} for the definitions of these properties).}
Then the $\Z[G]$-module $V^{\theta}$ is finitely generated and projective, 
the $\Z[G]$-module $B^{\theta}$ is finitely generated and free, 
and there is an exact sequence
\[ 0 \lra U_S^T(H) \lra V^{\theta} \lra B^{\theta} \lra \nabla_S^T(H) \lra 0. \]
If we retain the other assumptions but only assume that $T$ is $p$-tame rather than tame, then the same results hold after tensoring with $\Z_{(p)}$ and working over $\Z_{(p)}[G]$.
\end{corollary}

\begin{proof}
 {The exact sequence follows from Lemma~\ref{l:nablaseq} since the complex (\ref{e:vtbt}) is
quasi-isomorphic to $C^*_{\mathrm{RW}}$. 

Recall from Remark~\ref{r:blocfree} that $B_{\mathrm{loc}}$ is free of finite rank over $\Z[G]$.
Since $S\neq\emptyset$, choose $v\in S$. The summand
\[
\Ind_{G_{w}}^G(B_{w}) \cong \Z[G] \subset B_{\mathrm{loc}}
\]
maps isomorphically onto $B=\mathbf Z[G]$ via the local component of $\theta$ at $v$, which we denote $\theta_v$.
Write $B_{\mathrm{loc}} = \Ind_{G_{w}}^G(B_{w}) \oplus B_{\mathrm{loc}}'$ where $ B_{\mathrm{loc}}'$ is finite free over $\Z[G]$.
Then $B_{\mathrm{loc}}' \cong B^\theta$ via 
\[
x\longmapsto (-\theta_{v}^{-1}\theta(x),x).
\]
Therefore $B^\theta$ is free of finite rank over $\Z[G]$.

As noted in Remark~\ref{r:ct}, the module $V_{\mathrm{loc}}$ is $G$-cohomologically trivial when $T$ is tame.  Therefore, in this case, 
Lemma~\ref{l:sp}(1) implies that the module $V^{\theta}$ is $G$-cohomologically trivial, as it is the kernel of a surjective map of $G$-cohomologically trivial modules. 

 When $T$ is torsion-free, $V^{\theta}$ is $\Z$-torsion-free. 
Since the module $V^{\theta}$ is
 finitely generated, $\Z$-torsion-free, and $G$-cohomologically trivial, it is a projective $\Z[G]$-module by a result of Nakayama {\cite{nakayama}*{Theorem 1}}.

If $T$ is only $p$-tame, the same argument applies after tensoring with  $\mathbf Z_{(p)}$.  Note that Nakayama's criterion is known to hold over any Dedekind ring such as $\Z_{(p)}$; see for instance \cite{kock}*{Prop.~1.4}.}
\end{proof}

\begin{prop}\label{prop:quad-res}  Assume that $G = \Gal(H/F)$ is abelian.
Let $(S,T)$ be a complete RW pair with
$S\supset S_\infty$ and $T$ torsion-free. Assume  that $T$ is tame.
Then 
\[ V^{\theta} \lra B^{\theta} \lra \nabla_S^T(H) \lra 0 \]
is a locally quadratic presentation of $\nabla_S^T(H)$,
 i.e.\ the modules $V^{\theta}$ and  $B^{\theta}$  
 are finitely generated projective $\Z[G]$-modules of the same constant rank.
If $T$ is only $p$-tame rather than tame, while the other
assumptions are retained, the same conclusion holds over
$\mathbf Z_{(p)}[G]$ after tensoring the displayed sequence with
$\mathbf Z_{(p)}$.
\end{prop}

\begin{proof}
{
Let $A=\Z$ in the tame case and $A=\Z_{(p)}$ in the $p$-tame case, and
put $\Lambda=A[G]$. By Corollary~\ref{c:uvbn}, $B^\theta\otimes A$ is
free of some finite rank $r$ over $\Lambda$, while $V^\theta\otimes A$
is finitely generated and projective over $\Lambda$. It remains to show
that $V^\theta\otimes A$ has constant rank $r$.
The rank function of a finitely generated projective module is locally
constant on $\Spec \Lambda$, so it is enough to check the rank on each
connected component. Let $e\in \Lambda$ be the idempotent corresponding
to such a component. Since $\Lambda$ is finite free over $A$, the component $e\Lambda$ is $A$-torsion-free and hence
$(e\Lambda)\otimes_A \Q \neq 0.$
Thus every connected component of $\Spec \Lambda$ meets the generic fiber.
After extending scalars to $\overline{\Q}$, we have
\[
\overline{\Q}[G]\simeq \prod_{\chi\in \widehat G}\overline{\Q},
\]
since $G$ is abelian. The idempotent $e$ selects a subset
$J_e\subset \widehat G$. Therefore the rank of $V^\theta\otimes A$ on
the component corresponding to $e$ is determined by the dimensions of
$(V^\theta\otimes \overline{\Q})^\chi$ for $\chi\in J_e$.
It therefore suffices to show that
\begin{equation} \label{e:vqbq}
\dim_{\overline{\mathbf Q}} (V^\theta\otimes \overline{\mathbf Q})^\chi
=
\dim_{\overline{\mathbf Q}} (B^\theta\otimes \overline{\mathbf Q})^\chi
\end{equation}
for all characters $\chi$ of $G$.

Consider the exact sequence {from Corollary \ref{c:uvbn}, noting that \[ U_S^T(H) = \cO_{H,S,T}^* \] when $S \supset S_\infty$:}
\begin{equation} \label{e:fourterm}
0 \lra \cO_{H,S,T}^* \lra V^{\theta} \lra B^{\theta} \lra \nabla_S^T(H) \lra 0. 
\end{equation}
Tensoring
with $\overline{\mathbf Q}$ and taking $\chi$-components, the equality (\ref{e:vqbq}) is
equivalent to
\[
\dim_{\overline{\mathbf Q}} (\mathcal O_{H,S,T}^\ast\otimes
\overline{\mathbf Q})^\chi
=
\dim_{\overline{\mathbf Q}} (\nabla_S^T(H)\otimes
\overline{\mathbf Q})^\chi.
\]
Since $(S,T)$ is complete, the exact sequence (\ref{e:nablaseq}) with $X = X_{H,S}$
implies
\[
\nabla_S^T(H)\otimes\overline{\mathbf Q}
\simeq
X_{H,S}\otimes\overline{\mathbf Q}.
\]
On the other hand, the congruence conditions at $T$ only change the
$S$-unit group by finite index, so Dirichlet's unit theorem gives an
isomorphism of $\overline{\mathbf Q}[G]$-modules
\[
\mathcal O_{H,S,T}^\ast\otimes\overline{\mathbf Q}
\simeq
X_{H,S}\otimes\overline{\mathbf Q}.
\]
This proves that $V^\theta$ and $B^\theta$ have the same local rank as desired.}
\end{proof}

\subsection{$\nabla_S^T(H)$ and Galois cohomology} \label{s:gc}
{
In this subsection, we assume that $(S,T)$ is a complete RW pair.}  Our goal is to provide an interpretation of $\nabla_S^T(H)$ as representing a certain functor in Galois cohomology.

We consider the following variant of $C_{\mathrm{RW}}^*$, using the modified Tate sequences introduced in \S \ref{ss:cft}: \begin{align*}
V_{\mathrm{loc},2} &:= \prod_{v \in S} \Ind_{G_w}^G(V_{w,2}) \times \prod_{\substack{v \in T \\  v \text{ finite}}} \Ind_{G_w}^G(U_w) \times \prod_{v \notin (S \cup T)} \Ind_{G_w}^G(\cO_w^*), \\
B_{\mathrm{loc},2} &:=  \prod_{v \in S} \Ind_{G_w}^G(B_w) ,\\
C^*_{\mathrm{loc},2} &: V_{\mathrm{loc},2} \lra B_{\mathrm{loc},2}, \\
C^*_{\mathrm{global},2} &: V_2 \lra B, \\
C^*_{\mathrm{RW}, 2} &:= \mathrm{cone}(C^*_{\mathrm{loc},2} \lra C^*_{\mathrm{global},2})[-1].
\end{align*}
Combining the local and global reciprocity maps, we have a map \begin{equation}\mathrm{rec} \colon C^*_{\mathrm{RW}} \lra C^*_{\mathrm{RW},2}. \end{equation} Note that we can omit the infinite places in $T$ from $V_{\mathrm{loc},2}$ as $\mathrm{rec}_{w}(U_w) = 0$ for $w \in (T_H)_{\infty}$. 

\begin{prop}\label{prop:rw-rec}
The reciprocity map defines an isomorphism \[ \nabla_S^T(H) \cong H^1(C^*_{\mathrm{RW},2}). \]
\end{prop}
\begin{proof}
We have an exact sequence
\[  \prod_{w \in S_H} G_{H_w}^{\mathrm{ab}} \times \prod_{w \in T_H} U_{w} \times \!\!\!\!\!\! \prod_{w \notin (S_H \cup T_H)} \!\!\!\!\!\! \cO_{w}^* \lra G_H^{\mathrm{ab}} \lra H^1(C^*_{\mathrm{RW},2}) \lra \prod_{w \in S_H} \Z \lra \Z \]
from which we deduce the short exact sequence
\[  0 \lra \Cl_S^T(H) \lra H^1(C^*_{\mathrm{RW},2}) \lra {X_{H,S}} \lra 0 \]
as in (\ref{e:nablaseq}).
Therefore the canonical map $\nabla_S^T(H) = H^1(C^*_{\mathrm{RW}}) \lra H^1(C^*_{\mathrm{RW},2})$ is an isomorphism.
\end{proof}

Proposition \ref{prop:rw-rec} implies that {$\nabla_S^T(H)$ is naturally isomorphic to the cohomology of the sequence
\begin{equation}\label{e:present-nabla-bad}
V_{\mathrm{loc},2} \lra V_2 \oplus \prod_{v \in S} \Ind_{G_w}^G(B_w)  \lra B.
 \end{equation}
 
Now assume $S \neq \emptyset$ and fix $v_0 \in S$. Arguing as in the proof of Lemma~\ref{lem:coinv}, we obtain the following presentation of  
$\nabla_S^T(H)$:
\begin{equation}\label{e:present-nabla}
 V_{\mathrm{loc},2} \lra  V_2 \oplus \prod_{v \in S \setminus v_0} \Ind_{G_w}^G(B_w)  \lra \nabla_S^T(H) \lra 0.
 \end{equation}}

 \begin{theorem}\label{t:cocycle} {Let $(S,T)$ be a complete RW pair with $S \neq \emptyset$.}
Fix $v_0 \in S$, and let $A$ be a prodiscrete $G$-module. Then $\Hom_{\mathrm{cts},G}(\nabla_S^T(H), A)$ may be identified with the set of tuples $(\kappa, (x_v)_{v \in S \setminus \{v_0\}})$, where $\kappa \colon G_F \lra A$ is a 1-cocycle and $x_v \in A$, such that:
 \begin{enumerate}
 \item for $v \in S \setminus \{v_0\}$, $\kappa(\sigma_v) = (\sigma_v - 1)x_v$ for all $\sigma_v \in G_{F_v}$,
 \item $\kappa(G_{F_{v_0}}) = 0$,
 \item for finite $w \in T_H$, $\kappa|_{G_{H_w}}(U_w) = 0$,
 \item for $w \notin (S_H \cup T_H)$, $\kappa|_{G_{H_w}}(\cO_{H_w}^*) = 0$, i.e. $\kappa|_{G_{H_w}}$ is unramified.
 \end{enumerate}
 
The image of the map $\nabla_S^T(H) \lra A$ corresponding to $(\kappa, (x_v))$ is generated over $\Z[G]$ by $\kappa(G_{F}) \cup \{x_v \mid  v \in S \setminus \{ v_0 \} \}$.
\end{theorem} 

The proof of this theorem relies on a technical result on class formations, Theorem \ref{thm:desc-profinite}, whose proof will be the topic of \S\ref{s:duality}.

 \begin{proof}
The presentation (\ref{e:present-nabla}) implies that $\Hom_{\mathrm{cts},G}(\nabla_S^T(H), A)$ is the kernel of {the canonical map}:
 \begin{align*}
 & \Hom_{G, \mathrm{cts}}(V_2 \oplus \prod_{v \in S \setminus v_0} \Ind_{G_v}^G(B_w), A) \lra \\
 & \ \ \ \ \ \Hom_{G, \mathrm{cts}}\Big(\prod_{v \in S} \Ind_{G_w}^G(V_{w,2}) \oplus \prod_{v \in T \text{ finite}} \Ind_{G_w}^G(U_w) \oplus \prod_{v \notin (S \cup T)} \Ind_{G_w}^G(\cO_w^*), A\Big). 
 \end{align*} 
 
By Theorem \ref{thm:desc-profinite}, this is the same as the kernel of
 \begin{align*}
  & Z^1(G_{F}, A) \oplus \prod_{v \in S \setminus v_0} C^0(G_{F_v}, A) \lra \\ 
  & \ \ \ \ \  \prod_{v \in S} Z^1(G_{F_v}, A) \oplus \prod_{v \in T \text{ finite}} \Hom_{G_w, \mathrm{cts}}(U_w, A) \oplus \prod_{v \notin (S \cup T)} \Hom_{G_w, \mathrm{cts}}(\cO_w^*, A). \end{align*}
As the maps here are the ones arising from restriction and boundary maps in group cohomology, we obtain the first part of the theorem. 
 
It remains to show the claim regarding the image of $\nabla_S^T(H) \lra A$. This image is generated over $\Z$ by the images of 
\[ f \colon V_2  \lra A \] and \[ g \colon \bigoplus_{v \in S \setminus  v_0 } \Ind_{G_w}^G(B_w) \cong \bigoplus_{v \in S \setminus v_0} \Z[G] \lra A. \] The map $g$ sends $1$ in the $v$-component to $x_v \in A$. The map $f$ is determined by the factorization \[ \kappa \colon \Z[G_{F}] \xrightarrow{\kappa_{\mathrm{univ}}} V_2 \xrightarrow{f} A. \] Since $\kappa_{\mathrm{univ}}$ is surjective (see Theorem \ref{thm:desc-profinite}), we conclude that the image of $\nabla_S^T(H) \lra A$ equals the $\Z[G]$-submodule of $A$ generated by $\kappa(G_{F})$ and the elements $x_v$.
 \end{proof}
 
\begin{remark}
If we assume that $H/F$ is tamely ramified at $v \in T$, then condition (3) is equivalent to $\kappa|_{G_{F_v}}$ being tamely ramified. If we assume that $H/F$ is $p$-tamely ramified at $v$ and that $A$ is a pro-$p$ group, then condition (3) is equivalent to $\kappa|_{G_{F_v}}$ being tamely ramified. 
\end{remark} 
 
 \begin{remark}
Condition (4) is equivalent to $\kappa|_{G_{F_v}}$ being unramified (using the fact that $H/F$ is unramified at $v \notin (S \cup T)$). 
 \end{remark}

\subsection{The minus side in the  CM abelian case} \label{s:cmcase}

{
In this subsection, we suppose that $F$ is totally real and $H$ is a CM abelian extension of $F$.  Our eventual goal is to prove that after projecting to the minus side, a far weaker condition than $S \supset S_\infty$ is sufficient to imply local quadratic presentation of the Ritter--Weiss module (assuming the other conditions of Theorem~\ref{prop:quad-res} are maintained); indeed, this condition can be replaced by $S \neq \emptyset$.

We recall that for any $\Z[G]$-module $M$, we have defined \[ M_- = M/(1+c)M, \] where $c$ is the complex conjugation of $G = \Gal(H/F)$.
Our first task is to prove that under a mild hypothesis, the fundamental exact sequence (\ref{e:nablaseq}) involving $\nabla_S^T(H)$ remains exact after projecting to the minus side.}

\begin{prop}\label{p:exactness}
Suppose that $F$ is totally real and that $H/F$ is an abelian CM
extension. Let $(S,T)$ be a complete RW pair, and suppose that $S$
contains a place $v$ such that $c\in G_{v}$. Then
\[
\lrar_1^{\Z[G]}(X_{H,S},\Z[G]_-)=0.
\]
Consequently,
\[
0\longrightarrow
\Cl_S^T(H)_-
\longrightarrow
\nabla_S^T(H)_-
\longrightarrow
(X_{H,S})_-
\longrightarrow0
\]
is exact.
\end{prop}

Before we prove this proposition, we need a lemma on  $\lrar_{*}^{\ZG}(\cdot, N)$ where  $N$ = $\ZG_-$ or $\ZG_+$. We will suppress the ring $\ZG$ in our notation for $\lrar$. 
\begin{lemma} \label{l:tor}  {Suppose that $F$ is totally real, $H$ is  CM, and $G = \Gal(H/F)$ is abelian.} The following hold.
	\begin{enumerate}
		\item For $i\ge 1$ and any $\ZG$-module $M$, we have
	\begin{align} \label{e:twotors}
	\lrar_i(M, \ZG_+) &= \lrar_{i+1}(M, \ZG_-), \\
	  \lrar_i(M, \ZG_-) &= \lrar_{i+1}(M, \ZG_+). \label{e:twotors2}
	\end{align}
		\item For every $\ZG$-module $M$ such that $M[2]=0$, we have
		\begin{align*}
		\lrar_1(M, \ZG_+) &= M_-[2], \\ \lrar_1(M, \ZG_-) &= M_+[2].
		\end{align*}
	\end{enumerate}
\end{lemma}
\begin{proof}
	First, notice that there are two short exact sequences of $\ZG$-modules
	\begin{equation} \label{e:pm}
\begin{tikzcd}
	 0  \ar[r]  & \ZG_+ \ar[r,"c+1"]  & \ZG \ar[r] & \ZG_- \ar[r] & 0,
	 \end{tikzcd}
	 \end{equation}
	and
	\begin{equation} \label{e:mp}
	\begin{tikzcd}
	 0 \ar[r] & \ZG_-  \ar[r,"c-1"]  &  \ZG \ar[r] & \ZG_+ \ar[r] & 0,
\end{tikzcd}
	 \end{equation}
	where the inclusion maps are multiplication by $c+1$ and $c-1$ respectively. To verify that the first of these maps is injective, let $x = \sum \alpha_g g\in \ZG$ and note
	 \[ x(c+1) = 0 \iff \alpha_g+ \alpha_{gc} = 0 \text{ for all } g \in G \iff x = (1-c) \sum_{g\in G/\langle c \rangle} \alpha_g g. \]
	Similarly for the second map.
	
	If we apply $\otimes_{\ZG} M$ to the sequence (\ref{e:pm}), we obtain a long exact sequence 
	\begin{equation} \label{e:les}
	   \cdots \lra \lrar_1(M, \ZG_-) \lra M\otimes \ZG_+ \lra    M\otimes \ZG   \lra M\otimes \ZG_- \lra 0. \end{equation}
	Since $\ZG$ is projective and hence flat, we obtain  (\ref{e:twotors}).
	Similarly, we also obtain  (\ref{e:twotors2}) using (\ref{e:mp}). This finishes the proof of the first statement of the lemma.
	
	The sequence (\ref{e:les}) gives
\[
\operatorname{Tor}^{\ZG}_1(M,\ZG_-)
=\ker\bigl(M_+\xrightarrow{\,1+c\,}M\bigr).
\]
If $(1+c)m=0$, then $2\overline m=0$ in $M_+$.
Conversely, if $2\overline m=0$, write
$2m=(c-1)n$. Applying $1+c$ gives
$2(1+c)m=0$, and hence $(1+c)m=0$, since $M[2]=0$.
Thus the displayed kernel is $M_+[2]$.
The analogous argument using \eqref{e:mp} gives
$\operatorname{Tor}^{\ZG}_1(M,\ZG_+)=M_-[2]$.
\end{proof}

\begin{proof}[Proof of Proposition \ref{p:exactness}]
Since $(S,T)$ is complete, $X=X_{H,S}$, and we have the standard
exact sequence
\[
0\longrightarrow X_{H,S}
\longrightarrow Y_{H,S}
\longrightarrow \Z
\longrightarrow0,
\]
where
\[ Y=
Y_{H,S}=\bigoplus_{w\in S_H}\Z.
\]
Applying $\otimes \ZG_-$ to the sequence, we obtain a long exact sequence
	\begin{equation} \label{e:torseq}
\hspace{-0.95in}	\begin{tikzcd}
\cdots \ar[r] &   \lrar_2(Y, \ZG_-)   \ar[r] & \lrar_2(\Z, \ZG_-) \ar[r] 
 & \lrar_1(X, \ZG_-)  \ar[dll, out=-10, in=170] & \  \\ 
 &   \lrar_1(Y, \ZG_-) \ar[r]  & \cdots. 
 \end{tikzcd} 
	 \end{equation} 
	Writing $Y_v = \Ind_{G_v}^G \Z$ for the local components of $Y$, we  have \[ \lrar_1(Y, \ZG_-) = \oplus_{v\in S} \lrar_1(Y_v, \ZG_-) = \oplus_{v\in S} (Y_v)_+[2]. \] If $c\in G_v$, then $c$ acts trivially on $Y_v$, therefore \[ (Y_v)_+[2] = Y_v[2] = 0. \] If $c\notin G_v$, then \[ (Y_v)_+ = Y_v/ (c-1)Y_v \] is still torsion-free, so again $(Y_v)_+[2] = 0$. Therefore $\lrar_1(Y, \ZG_-) = 0$. In view of (\ref{e:torseq}), it remains to prove that
	\begin{equation} 
	\lrar_2(Y, \ZG_-)  \lra \lrar_2(\Z, \ZG_-) \label{e:tory}
	\end{equation} 
	is surjective.
	
		We first note that by Lemma~\ref{l:tor},  we have \[ \lrar_2(\Z, \ZG_-) = \lrar_1(\Z, \ZG_+) = \Z/2\Z.\]
Next we compute $\lrar_2(Y, \ZG_-) = \oplus_{v\in S} (Y_v)_-[2]$ as above. If $c\in G_v$, then $c$ acts trivially and $(Y_v)_-$ is just $Y_v/2Y_v$.  The component of (\ref{e:tory})
	is the surjective map $Y_v/2Y_v \lra \Z/2\Z$ induced by the surjection $Y_v \lra \Z$.  This concludes the proof. (Incidentally, we note that the condition $c \in G_v$ for some $v \in S$ is necessary, since if $c\notin G_v$, then 
	$(Y_v)_-$ is torsion-free, thus $(Y_v)_-[2]= 0$.) 
\end{proof}

\begin{remark}
The completeness hypothesis in Proposition~\ref{p:exactness} can be
weakened. Retaining the other hypotheses, the same conclusion holds
with $X_{H,S}$ replaced by the module $X$ of \ref{e:xdefgen}
provided that
$c\notin I_w$
for every finite ramified place $v\notin S\cup T$ and every $w\mid v$,
where $I_w\subset G_w$ denotes the inertia subgroup. In the
noncomplete case, the module $Y$ appearing in the proof contains, in
addition to $Y_{H,S}$, local summands associated with the finite
ramified places outside $S\cup T$. The displayed condition ensures
that the $\langle c\rangle$-coinvariants of each such summand are
$\Z$-torsion-free, which is the only additional point needed in the
proof.

The completeness hypothesis is the special case in which there are no
finite ramified places outside $S\cup T$.
\end{remark}

We now turn to the question of local quadratic presentation on the minus side in the CM setting without assuming $S \supset S_\infty$.  The following is an analogue of Corollary~\ref{c:uvbn}.

\begin{corollary} \label{c:vbminus}
Let $F$ be totally real and $H$  a CM abelian extension of $F$.  Suppose $S\neq\emptyset$, $S'$ is sufficiently large
to satisfy both conditions in Lemma~\ref{l:sp}(1), and that $T$ is torsion-free and tame.
The $\Z[G]_{-}$-modules $(V^{\theta})^{-}$ and $B^{\theta}_{-}$ are finitely generated and projective, and there is an exact sequence
\[ 0 \lra U_S^T(H)^{-} \lra (V^{\theta})^{-} \lra B^{\theta}_{-} \lra \nabla_S^T(H)_{-} \lra 0. \]
{If $T$ is only $p$-tame, the same statements hold after tensoring all
modules with $\mathbf Z_{(p)}$.}
\end{corollary}

\begin{proof}
 The result follows from Corollary~\ref{c:uvbn} using the fact that $(\cdot)^{-}$ is left-exact, $(\cdot)_{-}$ is right-exact, and $N^{-} \cong N_{-}$ for $G$-cohomologically trivial $\Z[G]$-modules $N$. {Note that this last isomorphism is not the canonical composition $N^- \rightarrow N \rightarrow N_-$, but rather the map \[ (1-c)x \mapsto x. \] $G$-cohomological triviality implies that \[ N^- = (1-c)N, \qquad \ker(1-c) = \im(1+c), \] so this map is well-defined. }
\end{proof}

 The following two lemmas are an adaptation of the results of \cite{nakayama} to the setting of $\Z[G]_{-}$-modules.  
Write $G = G_p \times G'$ where $G_p$ is the $p$-Sylow subgroup of $G$ and $G'$ is the subgroup of prime-to-$p$ order elements of $G$.
In the statement of part (1) below, note that if $p \neq 2$ then the composition $\Z_p[G_p] \lra \Z_p[G] \lra \Z_p[G]_-$ is injective.

\begin{lemma}\label{lem:rel-proj}  Suppose that $F$ is totally real and $H$ is a CM abelian extension of $F$.
We have:
\begin{enumerate}
\item If $p \neq 2$, a $\Z_p[G]_{-}$-module $B$ is projective if and only if it is projective with respect to the subring $\Z_p[G_p]$.
\item If $p = 2$, a $\Z_p[G]_{-}$-module $B$ is projective if and only if it is projective with respect to the subring $\Z_p[G_p]_{-}$.
\end{enumerate}
\end{lemma}
\begin{proof}
(1) $p \neq 2$.  As $\Z_p[G]_-$ is free as a $\Z_p[G_p]$-module, it is clear that $B$ is a projective $\Z_p[G_p]$ module if it is a projective $\Z_p[G]_-$-module.  For the other direction,
 note that we have an isomorphism of $\Z_p[G]_-$ modules
\begin{equation} \label{e:zpgb}
 \Z_p[G]_{-} \otimes_{\Z_p[G_p]} B \cong \Z_p[G']_- \otimes_{\Z_p} B, 
 \end{equation}
where
\begin{itemize}
\item $G$ acts on the left on the left side
\item $G$ acts via $G'$ on the left factor and $G_p$ on the right factor on the right side.
\end{itemize}
If $B$ is projective over $\Z_p[G_p]$, then the left side of (\ref{e:zpgb}) is projective over $\Z_p[G]_-$.  
As a $\Z_p[G]_-$-module, $B$ is a direct summand  of the right side of (\ref{e:zpgb}), since the size of $G'$ is coprime to $p$. 
Indeed, the map
\[ b \mapsto \frac{1}{|G'|} \sum_{g \in G'} g \otimes g^{-1}b
\]
from $B$ to the right side of (\ref{e:zpgb}) is $G$-equivariant, with a splitting given by $g \otimes b \mapsto gb$.
As a direct summand of a projective module is projective, the result follows.

(2) When $p = 2$ the proof is the same using the isomorphism 
\[ \Z_p[G]_{-}  \otimes_{\Z_p[G_p]_{-}} B \cong \Z_p[G/G_p] \otimes_{\Z_p} B. \]
\end{proof}

\begin{lemma}\label{l:nakayama} Suppose that $F$ is totally real and $H$ is a CM abelian extension of $F$.
Suppose  we have a short exact sequence of $\Z[G]_-$-modules
\[ 0 \lra P \lra Q \lra A \lra 0, \]
with $P$ and $Q$ projective and 
 $A$ finitely-generated and $\Z$-torsion-free.  Then $A$ is a projective $\Z[G]_{-}$-module.
\end{lemma}
\begin{proof}
It will suffice to show that, for all primes $p$, $A \otimes \Z_p$ is a projective $\Z_p[G]_{-}$ module. 
By Lemma \ref{lem:rel-proj}, we may instead show that $A \otimes \Z_p$ is a projective $\Z_p[G_p]_{-}$-module (when $p = 2$) or a projective $\Z_p[G_p]$-module (when $p \neq 2$).  

When $p = 2$, the ring $\Z[G_2]_{-} \otimes \Z/2\Z$ is isomorphic to $(\Z/2\Z)[G_2/\langle c \rangle]$. Since $A$ is $2$-torsion-free, the sequence $0 \lra P/2P \lra Q/2Q \lra A/2A \lra 0$ is exact. The module $P$ and $Q$ are projective as $\Z_2[G_2]_{-}$-modules by Lemma \ref{lem:rel-proj}, and so $P/2P$ and $Q/2Q$ are projective $(\Z/2\Z)[G_2/\langle c \rangle]$-modules. Thus we find that $A/2A$ is a cohomologically trivial $G_2/\langle c \rangle$-module, and hence is a free  $(\Z/2\Z)[G_2/\langle c \rangle]$-module by 
{\cite{nakayamaI}*{pg.~42}}.  By Nakayama's Lemma, a $(\Z/2\Z)[G_2/\langle c \rangle]$-basis of $A/2A$ lifts to $\Z_2[G_2]_{-}$-generators of $A \otimes \Z_2$. This gives an exact sequence \[ 0 \lra K \lra \oplus_{i=1}^n \Z_2[G_2]_{-} \lra A \otimes \Z_2 \lra 0. \] Since $A \otimes \Z_2$ is 2-torsion-free, we may tensor with $\Z/2\Z$ to conclude that $K/2K = 0$, and hence by Nakayama's Lemma $K = 0$. Therefore $A \otimes \Z_2$ is a free $\Z_2[G_2]_{-}$-module.

When $p \neq 2$,  we can simply apply the results of \cite{nakayama} directly to see that the $G_p$-cohomologically trivial,  finitely-generated, $p$-torsion-free module $A \otimes \Z_p$ is a free $\Z_p[G_p]$-module.
\end{proof}

 Consider the exact sequence of Corollary~\ref{c:vbminus}:
\begin{equation} \label{e:uvbn}
0 \lra U_S^T(H)^{-} \lra (V^{\theta})^{-} \lra (B^{\theta})_{-} \lra (\nabla_S^T)_{-} \lra 0. \end{equation}
We now quotient $(V^\theta)^-$ by the logarithmic lattice in its kernel in order to obtain a quadratic presentation.
{Write $T_{H,\infty}$ for the set of all complex places of $H$ above places in $T$.} We have an exact sequence 
\[ 
0 \lra \bigoplus_{\substack{v \in T_{H,\infty}}} 2 \pi i \Z \lra U_S^T(H) \lra \cO_{H, S,T}^* \lra 0. 
\] 
Taking minus-parts, we get
\[ 
0 \lra \bigoplus_{\substack{v \in T_{H, \infty}}}  2 \pi i \Z \lra U_S^T(H)^{-} \lra (\cO_{H, S,T}^*)^{-} \lra 0, \]
where right-exactness follows from $H^1(\langle c \rangle, \Z) = 0$.

Define 
\begin{equation} \label{e:vpdef}
V' := (V^{\theta})^{-}/ \bigoplus_{\substack{v \in T_{H, \infty}}}  2 \pi i \Z. 
\end{equation} We obtain an exact sequence
\begin{equation} 0 \lra (\cO_{H, S,T}^*)^{-} \lra V' \lra (B^{\theta})_{-} \lra  \nabla_S^T(H)_{-} \lra 0. \label{e:vp}
\end{equation}

{
\begin{prop} \label{prop:quad-resCM} Suppose $F$ is totally real and $H$ is a CM abelian extension of $F$.  
Let $(S,T)$ be a complete RW pair with
$S\neq\emptyset$ and $T$ torsion-free. Assume either that $T$ is tame, or work after tensoring
with $\mathbf Z_{(p)}$ and assume that $T$ is $p$-tame.
Then 
\[ V'\lra (B^{\theta})_{-} \lra (\nabla_S^T)_{-} \lra 0\]
is a locally quadratic presentation of $\nabla_S^T(H)_-$.
\end{prop}
}
\begin{proof} 
First we observe that since $B^\theta$ is a free $\Z[G]$-module, the quotient $(B^\theta)_-$ is a free $\Z[G]_-$-module.
Next, note that
\[ 
\bigoplus_{\substack{v \in T_{H, \infty}}} 2 \pi i \Z \cong \bigoplus_{v \in T_{\infty}} \Z[G]_{-},
\]
{where $T_{\infty}$ is the set of  archimedean places in $T$.} Therefore $V'$ is a quotient of projective $\Z[G]_{-}$-modules.  Furthermore $V'$ is torsion-free by (\ref{e:vp}) since  {$T$ is torsion-free}. By Lemma \ref{l:nakayama}, this implies that $V'$ is itself projective {as a $\Z[G]_-$-module}. 
In the $p$-tame case, the same argument is applied after tensoring with
$\mathbf Z_{(p)}$; we use the localized analogue of Lemma~\ref{l:nakayama},
whose proof is identical with $\mathbf Z_{(p)}$ in place of $\mathbf Z$.

{
As in the proof of Proposition~\ref{prop:quad-res}, to complete the proof of local quadratic presentation it remains to show that}
$\dim_{\overline{\Q}} (\cO_{H,S,T}^* \otimes \overline{\Q})^{\chi}$ equals \[ \dim_{\overline{\Q}} (\nabla_S^T(H) \otimes \overline{\Q})^{\chi} = \dim_{\overline{\Q}} (X_{H,S} \otimes \overline{\Q})^{\chi} \] for all odd $\chi$. As \[ \cO_{H,S,T}^* \otimes \overline{\Q} \cong X_{H,S \cup T_{\infty}} \otimes \overline{\Q} \] by Dirichlet's Unit Theorem, the claim follows from the fact that \[ K := X_{H,S \cup T_{\infty}}/X_{H,S} \cong \bigoplus_{w \in T_{H, \infty}} \!\! \Z \] satisfies $K \otimes \overline{\Q}[G]_{-} = 0$.
\end{proof}

{
\begin{remark}
 Retaining the assumption that $H/F$ is an abelian CM extension of a totally real field: if  $T$ is torsion-free and tame, but $(S,T)$ is not complete, we still may construct the module $V'$ and obtain the exact sequence (\ref{e:vp}).
\end{remark}

Propositions~\ref{prop:quad-res} and~\ref{prop:quad-resCM} complete the proof of Theorem~\ref{thm:quad}.
}

\subsection{Fitting Ideals and Transposes}\label{ss:changing-sets}

In this subsection we continue to assume that $F$ is totally
real and $H/F$ is a CM abelian extension. Let $(S,T)$ be an
RW pair with $S\neq\emptyset$ and $T$ torsion-free. We assume
that $T$ is tame, or work after tensoring with $\Z_{(p)}$
and assume that $T$ is $p$-tame. Thus
Theorem~\ref{thm:quad}(2) provides the presentation
\eqref{e:vp} by finitely generated projective modules.
When $(S,T)$ is complete, this presentation is locally
quadratic by Proposition~\ref{prop:quad-resCM}.

We establish some elementary properties of the Fitting ideal of $\nabla_{S}^{T}(H)_-$ as certain places are moved in or out of $S$ and $T$. 
\begin{lemma} \label{l:fittinfty}
{Assume the hypotheses of Proposition \ref{prop:quad-resCM}}, with $T$ tame. Then for any infinite place $v \in T$, we have \[ \Fitt_{\Z[G]_{-}}(\nabla_{S \cup \{ v \}}^{T\setminus \{v\}}(H)_{-}) = 2 \Fitt_{\Z[G]_{-}}(\nabla_S^T(H)_{-}). \]
We have the corresponding identity after tensoring with
$\mathbf Z_{(p)}$ if $T$ is only $p$-tame rather than tame.
\end{lemma}
\begin{proof}
{
Let
\[
\xymatrix{ C_{S,T}^\bullet:\quad V'_{S,T} \ar[r]^{ \ \ \ d_1} & (B^\theta_{S,T})_-
}
\]
and
\[
\xymatrix{ C_{S \cup \{v\},T \setminus \{v\}}^\bullet:\quad V'_{S \cup \{v\},T \setminus \{v\}} \ar[r]^{\ \ \ \ \ \ \ d_2} & (B^\theta_{S \cup \{v\},T \setminus \{v\}})_- }
\]
be the complexes yielding locally quadratic presentations of $\nabla_S^T(H)_-$ and
$\nabla_{S \cup \{v\}}^{T \setminus \{v\}}(H)_-$, respectively, obtained in Proposition \ref{prop:quad-resCM}. We obtain an
injective map of complexes
\[
C_{S,T}^\bullet \longrightarrow C_{S \cup \{v\},T \setminus \{v\}}^\bullet.
\]

 First we compute the quotient of the first terms in the complexes.
Let $w\mid v$, so $G_w=\mathrm{Gal}(\mathbf C/\mathbf R)=\{1,c\}$.
The local term when $v\in T$ is \[ \mathbf C^-/2\pi i\mathbf Z = \R i /2\pi i\Z, \]
while after moving $v$ to $S$, the local term is
$V_w^-$  (see (\ref{e:rwlis}), (\ref{e:rwlit}), and (\ref{e:vpdef}) to recall the definitions of our local complexes at the infinite places).
The exponential map identifies
\[
\mathbf \R i/2\pi i\mathbf Z \cong (\mathbf C^\times)^- = \text{unit circle}.
\]
Let $d_w:V_w\to\Z[G_w]$ be the Tate differential, whose
image is $\Delta G_w=\Z[G_w]^-$. Since $V_w$ is
cohomologically trivial,
\[
d_w(V_w^-)=(1-c)\Delta G_w=2\Z[G_w]^-.
\]
Its kernel on $V_w^-$ is $(\C^\times)^-$, so half the
Tate differential gives an isomorphism
\[
V_w^-/(\C^\times)^-
\xrightarrow{\ \frac12d_w\ }\Z[G_w]^-.
\]
The map on the second terms is split, with cokernel
$\Ind_{G_w}^G\Z[G_w]_-$.
Recall that the differential to minus coinvariants uses the
inverse-norm identification
$(1-c)b\mapsto\overline b$ from Corollary~\ref{c:vbminus}.
Under the displayed normalization, the local map on cokernels is
therefore the inclusion-induced map
\[
\Z[G_w]^-\longrightarrow\Z[G_w]_-.
\]
It sends the generator $1-c$ to $2$.
Thus, after choosing generators, the quotient differential is
multiplication by $2$.

To summarize, after passing from these local calculations to the global $G$-modules, we have a commutative diagram:
\[
\xymatrix{
V'_{S,T} \ar[r]^{d_1} \ar[d] &
(B^\theta_{S,T})_- \ar[d] \\
V'_{S \cup \{v\},T \setminus \{v\}} \ar[r]^{d_2} \ar[d] &
(B^\theta_{S \cup \{v\},T \setminus \{v\}})_- \ar[d] \\
\Z[G]_- \ar[r]^-{\times 2} &
\Z[G]_-.
}
\]
Since $\Z[G]_-$ is $\Z$-torsion-free, the quotient complex
\[
\Z[G]_- \xrightarrow{\times 2} \Z[G]_-
\]
has $H^0=0$. Therefore the above short exact sequence of complexes gives a short
exact sequence of $\Z[G]_-$-modules
\[
0\lra \nabla_S^T(H)_-
\lra \nabla_{S\cup\{v\}}^{T\setminus\{v\}}(H)_-
\lra \Z[G]_-/2\Z[G]_-
\to 0.
\]
Since $\Z[G]_-/2\Z[G]_-$ is quadratically presented with
Fitting ideal $(2)$, \cite{dk}*{Lemma 2.6} gives \[
\operatorname{Fitt}_{\Z[G]_-}
(\nabla_{S\cup\{v\}}^{T\setminus\{v\}}(H)_-)
=
2\operatorname{Fitt}_{\Z[G]_-}
(\nabla_S^T(H)_-).
\]
The $p$-tame case is identical after tensoring with $\Z_{(p)}$.
}
\end{proof}

The following lemma is well-known, but since our construction of $\nabla$ is different than what has previously appeared in the literature, we include it for completeness.

\begin{lemma}
{Assume the hypotheses of Proposition \ref{prop:quad-resCM}}, with $T$ tame. Let $v \notin S \cup T$ be a finite place of $F$ that is unramified in $H/F$. Then
\begin{enumerate}
\item $\Fitt_{\Z[G]_-}(\nabla_{S \cup \{ v \}}^{T}(H)_-) = (1- \sigma_v^{-1} )\Fitt_{\Z[G]_-}(\nabla_S^T(H)_-)$
\item $\Fitt_{\Z[G]_-}(\nabla_{S}^{T \cup \{ v \}}(H)_-) = (1 - \sigma_v^{-1} \N(v)) \Fitt_{\Z[G]_-}(\nabla_S^T(H)_-)$
\end{enumerate}
The corresponding identities hold over $\Z_{(p)}[G]_-$ after
tensoring all modules with $\Z_{(p)}$ if $T$ is only $p$-tame
rather than tame.
\end{lemma}
\begin{proof}  {
We argue as in Lemma~\ref{l:fittinfty} using the complex 
\[
\xymatrix{
C_{S,T}^\bullet:\quad V'_{S,T} \ar[r]^{\ \ \ d_1} & (B^\theta_{S,T})_-
}
\]
yielding the locally quadratic presentation of $\nabla_S^T(H)_-$ produced in Proposition \ref{prop:quad-resCM}.  We assume that $v \not\in S'$.

 (1) Consider the map between the relevant local complexes: \[ (\cO_w^* \lra 0) \lra (V_w \lra B_w)\] which has cokernel
  \[ V_w/\cO_w^* = W_w \cong \Z[G_w] \lra B_w = \Z[G_w] . \]  
 As described in (\ref{e:wb}), this map is given by multiplication by $1 - \sigma_w^{-1}$.

We therefore obtain a diagram as below:
\[
\xymatrix{
V'_{S,T} \ar[r]^{d_1} \ar[d] &
(B^\theta_{S,T})_- \ar[d] \\
V'_{S \cup \{v\}} \ar[r]^{d_2} \ar[d] &
(B^\theta_{S \cup \{v\}})_- \ar[d] \\
\Z[G]_- \ar[r]^-{1- \sigma_v^{-1}} &
\Z[G]_-
}
\]
Although the complex in the bottom row need not be acyclic in degree $1$, it gives a quadratic
presentation whose Fitting ideal is generated by $1-\sigma_v^{-1}$. Localizing
at each prime of $\Z[G]_-$ and choosing splittings of the presentation modules,
the matrix for the presentation $d_2$ of
$\nabla_{S\cup\{v\}}^T(H)_-$ is block upper-triangular with diagonal blocks
$d_1$ and $1-\sigma_v^{-1}$. Hence
\[ \Fitt_{\Z[G]_-}(\nabla_{S \cup \{ v \}}^{T}(H)_-) = (1- \sigma_v^{-1} )\Fitt_{\Z[G]_-}(\nabla_S^T(H)_-) \] follows.

 (2) Consider the map between the relevant local complexes: \[ (U_w \lra 0) \lra (\cO_{H_w}^* \lra 0). \] We find that $C^*_{S,T \cup \{v\}} \lra C^*_{S,T}$ has cokernel $\Ind_{G_w}^G(\F_w^*) \lra 0$.
We obtain a commutative diagram:
\[
\xymatrix{
V'_{S,T\cup\{v\}} \ar[r]^{d_2} \ar[d]_{f} &
(B^\theta_{S,T\cup\{v\}})_- \ar[d] \\
V'_{S,T} \ar[r]^{d_1} \ar[d] &
(B^\theta_{S,T})_- \ar[d] \\
(\Ind_{G_w}^G \mathbf F_w^\ast)_- \ar[r] &
0
}
\]
The $B$ modules on the right are equal, and the map between them is the identity.  After passing to a localization where the four modules in the upper square are free, we therefore find 
\[ \det(d_2) = \det(f) \det(d_1).\]
After localizing at a prime of $\Z[G]_-$, the determinant of $f$ generates
the localized Fitting ideal of $(\Ind_{G_w}^G\mathbf F_w^\ast)_-$.
Since
 \[ \Fitt_{\Z[G]_-}((\Ind_{G_w}^G\F_w^*)_-) = (1 - \sigma_v^{-1}\N(v)), \] we obtain the desired result
\[ \Fitt_{\Z[G]_-}(\nabla_{S}^{T \cup \{ v \}}(H)_-) = (1 - \sigma_v^{-1} \N(v)) \Fitt_{\Z[G]_-}(\nabla_S^T(H)_-). \]
}
\end{proof}

We record a general description of the transpose that will be used in
the proof of Lemma~\ref{l:causel}.
We equip the $\Z[G]_-$-linear dual
$\Hom_{\Z[G]_-}(M,\Z[G]_-)$ with the contragredient action
\[
(g\cdot f)(m)=g^{-1}f(m).
\]
Thus it is the $\#$-twist of the usual scalar action on this
Hom module. We use the same contragredient convention on
$\Z$-linear and Pontryagin duals. With these actions, the
coefficient pairing on $\Z[G]_-$ gives a natural isomorphism
\begin{equation} \label{e:functoreq}
\Hom_{\Z[G]_-}(M,\Z[G]_-)\simeq\Hom_\Z(M,\Z).
\end{equation}

\begin{lemma} \label{l:transprop}
 {Suppose that $F$ is totally real and $H$ is a CM abelian extension of $F$.  Suppose that $S\neq\emptyset$ and $T$ is torsion-free and tame. Consider the presentation (\ref{e:vp})} of $ \nabla_S^T(H)_{-}$.  Let $\Sel_S^T(H)_-$ denote the transpose of $ \nabla_S^T(H)_{-}$ over $\Z[G]_-$ associated to {(\ref{e:vp})} in the sense of Jannsen, i.e., the module defined by the exact sequence
\[  \Hom_{\Z[G]_-}((B^{\theta})_{-} , \Z[G]_-) \lra  \Hom_{\Z[G]_-}(V' , \Z[G]_-) \lra \Sel_S^T(H)_- \lra 0. \]
Then we have a short exact sequence
\begin{equation} \label{e:selseq}
  0 \lra  \bigl((\nabla_S^T(H)_-)_{\mathrm{tors}}\bigr)^\vee \lra \Sel_S^T(H)_- \lra \Hom_{\Z}((\cO_{H,S,T}^*)^{-}, \Z) \lra 0. 
  \end{equation}
\end{lemma}

\begin{proof} Split up (\ref{e:vp}) into two short exact sequences
\[ 0 \lra (\cO_{H,S,T}^*)^{-}   \lra  V'  \lra V'/(\cO_{H,S,T}^*)^{-} \lra 0, \]
\[ 0 \lra V'/(\cO_{H,S,T}^*)^{-} \lra  (B^\theta)_- \lra  \nabla_S^T(H)_{-} \lra 0. \]
We use the identification of functors on the category of $\Z[G]_-$-modules given in (\ref{e:functoreq}).
Applying $ \cF(-) = \Hom_{\Z}(-, \Z)$ to the two sequences above, we obtain
\begin{equation} \label{e:f1}
 0 \lra \cF(V'/(\cO_{H,S,T}^*)^{-}) \lra \cF(V') \lra \cF( (\cO_{H,S,T}^*)^{-}) \lra 0 \end{equation}
since $V'/(\cO_{H,S,T}^*)^{-} \subset (B^\theta)_-$ is a free $\Z$-module, and
\begin{equation} \label{e:f2}
 \cF((B^\theta)_-) \lra \cF(V'/(\cO_{H,S,T}^*)^{-}) \lra \Ext^1_\Z(\nabla_S^T(H)_-, \Z) \lra 0.
\end{equation}
Modding out the first two nontrivial terms of (\ref{e:f1}) by the image of  $ \cF((B^\theta)_-)$ in each, applying (\ref{e:f2}), and noting that 
\[ \Ext^1_\Z(\nabla_S^T(H)_-, \Z) \cong (\nabla_S^T(H)_-)_{\mathrm{tors}}^\vee, \]
we obtain the desired exact sequence (\ref{e:selseq}).
\end{proof}

We now prove Lemma~\ref{l:causel}, whose statement we recall.
Fix a real place $w$ of $F$, and put
\[
T^\dagger=T\cup(S_\infty\setminus\{w\}).
\]  Lemma~\ref{l:causel} states that 
there is a natural surjective homomorphism of $\Z[G]_-$-modules
\[
\Sel_{\{w\}}^{T^\dagger}(H)_-
\longtwoheadrightarrow
\bigl(\Cl_{\mathrm{au}}^T(H)\bigr)^\vee.
\]

\begin{proof}[Proof of Lemma~\ref{l:causel}]
Choose anti-invariant ideals
\[
\fa_1,\ldots,\fa_r\in I^T(H)^-
\]
whose classes generate $\Cl_{\mathrm{au}}^T(H)$. Choose the auxiliary
unramified set $S'$ in the Ritter--Weiss construction sufficiently large
that it satisfies the usual conditions, that the corresponding depleted
$T$-ray class group is trivial, and that it contains every unramified
place below a prime in the support of some $\fa_i$. Put
\[
\mathcal S=S_{\ram}\cup S'.
\]
Let $I_{\mathcal S}(H)\subset I^T(H)$ denote the subgroup
of fractional ideals supported above $\mathcal S$.
Then
\[
\mathcal C_{\mathcal S}^T(H)
:=
I_{\mathcal S}(H)^-
\big/
\{(u):u\in(\cO_{H,\mathcal S,T}^*)^-\}
\]
maps isomorphically onto $\Cl_{\mathrm{au}}^T(H)$. The map is
injective because an anti-unit generating an ideal supported above
$\mathcal S$ belongs to $(\cO_{H,\mathcal S,T}^*)^-$.  It is
surjective by the choice of the ideals $\fa_i$. Equivalently,
\[
\mathcal C_{\mathcal S}^T(H)
\cong
\coker\left(
(\cO_{H,\mathcal S,T}^*)^-
\xrightarrow{\operatorname{ord}}
Y_{H,\mathcal S}^-
\right).
\]

We use the chosen set $S'$ in the Ritter--Weiss presentation for
$(\{w\},T^\dagger)$. Put
\[
W_{\loc}=\im(V_{\loc}\longrightarrow B_{\loc}),
\qquad
W=\Delta G,
\qquad
W^\theta=\ker(W_{\loc}\longrightarrow W).
\]
Let $W'_{\loc}$ be obtained from $W_{\loc}$ by replacing the component
$W_v$ at each finite $v\in\mathcal S$ by $\Delta G_v$, leaving the
component $\Delta G_w$ at $w$ unchanged, and put
\[
W'^\theta=\ker(W'_{\loc}\longrightarrow W).
\]

The local exact sequences
\[
0\longrightarrow\Z
\longrightarrow W_v
\longrightarrow\Delta G_v
\longrightarrow0
\qquad (v\in\mathcal S)
\]
and the usual snake-lemma argument give exact sequences
\begin{equation}\label{e:wwprime}
0\longrightarrow Y_{H,\mathcal S}
\longrightarrow W^\theta
\longrightarrow W'^\theta
\longrightarrow0
\end{equation}
and
\begin{equation}\label{e:vwprime}
0\longrightarrow\widetilde{\cO}_{H,\mathcal S,T}^*
\longrightarrow V^\theta
\longrightarrow W'^\theta
\longrightarrow0,
\end{equation}
where $\widetilde{\cO}_{H,\mathcal S,T}^*$ is the group of logarithmic
lifts of the elements of $\cO_{H,\mathcal S,T}^*$ at the places above
$S_\infty\setminus\{w\}$. In particular, writing
\[
L_\infty=
\bigoplus_{\widetilde v\mid S_\infty\setminus\{w\}}2\pi i\Z,
\]
we have
\[
0\longrightarrow L_\infty
\longrightarrow
(\widetilde{\cO}_{H,\mathcal S,T}^*)^-
\longrightarrow
(\cO_{H,\mathcal S,T}^*)^-
\longrightarrow0.
\]
The last map is surjective because an anti-invariant element has
absolute value $1$ at every complex place and hence admits purely
imaginary logarithms.

Since
\[
V'=(V^\theta)^-/L_\infty,
\]
taking anti-invariants in \eqref{e:wwprime} and
\eqref{e:vwprime} gives injections
\[
\frac{(W^\theta)^-}{Y_{H,\mathcal S}^-}
\hookrightarrow (W'^\theta)^-,
\qquad
\frac{V'}{(\cO_{H,\mathcal S,T}^*)^-}
\hookrightarrow (W'^\theta)^-.
\tag{$\dagger$}
\]
The maps in \eqref{e:wwprime} and \eqref{e:vwprime} are compatible, and
the map
\[
(\cO_{H,\mathcal S,T}^*)^-
\longrightarrow Y_{H,\mathcal S}^-
\]
that they induce is the equivariant valuation map.

Using
\[
(B^\theta)^-\xrightarrow{\sim}(B^\theta)_-,
\qquad
(1-c)b\longmapsto\overline b,
\]
we regard $(W^\theta)^-$ as a submodule of $(B^\theta)_-$. The map
$V'\to(B^\theta)_-$ factors through $(W^\theta)^-$, and we have a
commutative square
\[
\begin{CD}
(\cO_{H,\mathcal S,T}^*)^- @>>> V'\\
@V{\operatorname{ord}}VV                @VVV\\
Y_{H,\mathcal S}^-              @>>> (W^\theta)^-.
\end{CD}
\tag{$\ddagger$}
\]
The injections in $(\dagger)$ imply that this square is cartesian.
Indeed, if an element of $V'$ maps into $Y_{H,\mathcal S}^-$, then its
image in $(W'^\theta)^-$ is zero, so its class modulo
$(\cO_{H,\mathcal S,T}^*)^-$ is zero.

It follows from $(\ddagger)$ that the lower horizontal map induces an
injection
\[
\mathcal C_{\mathcal S}^T(H)
\hookrightarrow
\coker\left(
V'\longrightarrow(W^\theta)^-
\right).
\]
Since $(W^\theta)^-\subset(B^\theta)_-$, the latter cokernel injects
into $\nabla_{\{w\}}^{T^\dagger}(H)_-$. Thus
\[
\Cl_{\mathrm{au}}^T(H)
\cong\mathcal C_{\mathcal S}^T(H)
\hookrightarrow
\bigl(\nabla_{\{w\}}^{T^\dagger}(H)_-\bigr)_{\mathrm{tors}},
\]
where the image lies in the torsion subgroup because
$\Cl_{\mathrm{au}}^T(H)$ is finite.

Finally,
\[
\cO_{H,\{w\},T^\dagger}^*
=
\cO_{H,S_\infty,T}^*.
\]
Complex conjugation acts trivially on this group, and it is torsion-free
by the hypothesis on $T$. Hence
$(\cO_{H,\{w\},T^\dagger}^*)^-$ is trivial.
Lemma~\ref{l:transprop} therefore gives
\[
\Sel_{\{w\}}^{T^\dagger}(H)_-
\cong
\bigl(\nabla_{\{w\}}^{T^\dagger}(H)_-\bigr)_{\mathrm{tors}}^\vee.
\]
Dualizing the preceding injection proves the asserted surjection.

The construction is compatible with enlarging $S'$, so the resulting
map is independent of the sufficiently large auxiliary set.
\end{proof}

\section{A Duality Theorem for Class Formations}\label{s:duality}

 {In this section $G$ is any profinite group.}  Consider a class formation $(G, C)$ and an open normal subgroup $H$. We will prove that, 
 given a prodiscrete $\Z[G/H]$-module $A$, the complex \[ \Hom_{\Z[G/H]}(\Z[G/H], A) \lra \Hom_{\Z[G/H],\text{cts}}(E({G/H}), A), \] dualizing the modified Tate extension (\ref{e:norm-cplt-ext}), is quasi-isomorphic to $\tau_{\leq 1} C^*(G, A)$, the truncation of Galois cohomology. This result was used in the proof of Theorem \ref{t:cocycle} above, which gave an interpretation of the module $\nabla_S^T(H)$ in terms of Galois cohomology.
Our results in this section may be of independent interest.

\subsection{Duality ($G$ finite)}

Consider a class formation $(G,C)$ with $G$ finite. Recall the extension (\ref{e:norm-cplt-ext})  obtained by pushing out the extension (\ref{e:two-extn}) along $C^H \lra C^H/N_H(C)\ \cong H^{\mathrm{ab}}$: 
\begin{equation} 0 \lra C^H/N_H(C) \lra E({G/H}) \lra \Z[G/H] \lra \Z \lra 0,  \label{e:egh}
\end{equation}
where $ E({G/H}) := C({G/H})/N_H(C)$.

Let \[ C^{\mathrm{bar}}_{*}: \Z[G] \longleftarrow \Z[G^2] \longleftarrow \cdots \] be the bar resolution of $\Z$ by $G$-modules. For a discrete $G/H$-module $A$,
\begin{itemize}
\item $C^*(H, A) := \Hom_{\Z[H]}(C^{\mathrm{bar}}_{*}, A)$ is a complex of $G/H$-modules whose cohomology equals $H^*(H, A)$ as a $G/H$-module.
\item $C^*(G, A) := \Hom_{\Z[G]}(C^{\mathrm{bar}}_{*}, A)$ is a complex whose cohomology equals $H^*(G,A)$.
\end{itemize}

We will restrict our attention to the truncation $\tau_{\leq 1}C^*(G,A)$, i.e.
\begin{equation}
\tau_{\leq 1}C^*(G,A) \colon C^0(G, A) \lra Z^1(G,A).
\end{equation}

\begin{theorem}\label{thm:finite-descent}

Consider a class formation $(G, C)$ with $G$ finite, a normal subgroup $H \subset G$, and a discrete $G/H$-module $A$. \begin{enumerate}
\item The first cohomology group of the complex 
\[ 
\Hom(\Z[G/H], A) \lra \Hom(E({G/H}), A)
\]
is $\Hom(C^H/N_H(C), A)$. Furthermore, there is a quasi-isomorphism from the truncation $\tau_{\leq 1} C^*(H, A)$ to 
\[
\Hom(\Z[G/H], A) \lra \Hom(E({G/H}), A)
\]
such that the map induced on the first cohomology groups 
\[
H^1(H,A) = \Hom(H^{\mathrm{ab}},A) \lra \Hom(C^H/N_H(C), A) 
\] 
is dual to the reciprocity map $C^H/N_H(C) \lra H^{\mathrm{ab}}$. 
\item The truncation $\tau_{\leq 1} C^*(G,A)$ is quasi-isomorphic to
\[ \Hom_{\Z[G/H]}(\Z[G/H], A) \lra \Hom_{\Z[G/H]}(E({G/H}), A). \]
\item The inhomogeneous 1-cocycle $\Z[G] \lra E({G/H})$ corresponding to the universal element $\kappa_{\mathrm{univ}} \in Z^1(G, E({G/H}))$ is surjective.
\end{enumerate}
\end{theorem}

We begin with some preparatory lemmas.

\begin{lemma}\label{lem:coinflation}
Consider a class formation $(G,C)$. For $H \subset G$ a normal subgroup, we have $C(G)_H \cong C({G/H}).$ In particular, $C(G)_G \cong C^G.$ There is a commutative diagram
with exact rows \[ \begin{tikzcd} 0 \arrow{r} & C \arrow{r} \arrow{d}& C(G) \arrow{r} \arrow{d}& \Z[G] \arrow{r} \arrow{d}& \Z \arrow{r} \arrow{d}& 0 \\
0 \arrow{r}& C^H \arrow{r} & C({G/H}) \arrow{r} & \Z[G/H] \arrow{r} & \Z \arrow{r} & 0, 
\end{tikzcd} \]
where the vertical map on the left term is the norm $N_H$ and the map $\Z[G] \lra \Z[G/H]$ is the canonical projection.
\end{lemma}
\begin{proof}
Write $I_G \subset \Z[G]$ for the augmentation ideal. Multiplying by the norm $N_H = \sum_{h \in H} h$ yields a commutative diagram  
\begin{equation} \label{e:snake}
\begin{tikzcd}
\ & C_H \arrow{r} \ar[d,"N_H"] & C(G)_H  \arrow{r} \ar[d,"N_H"] &(I_G)_H  \arrow{r}\ar[d,"N_H"] & 0 \\
0  \arrow{r} & C^H  \arrow{r} & C(G)^H  \arrow{r} & (I_G)^H. & \\
\end{tikzcd}
\end{equation}
Let $K(G, H)$ denote the kernel of $N_H\colon (I_G)_H \lra (I_G)^H$.
The snake lemma applied to (\ref{e:snake})  yields a map 
\begin{equation} \label{e:kdef}
 K(G, H) \lra C^H/N_H(C). \end{equation} Since $C(G)$ is $H$-cohomologically trivial, the middle vertical arrow in (\ref{e:snake}) is an isomorphism, whence (\ref{e:kdef}) is an isomorphism as well.

Note that $(I_H)_H \cong H^{\mathrm{ab}}$ is contained in 
$K(G,H).$ 
In fact, these are equal: since both $\Z[G]$ and $\Z[H]$ are $H$-cohomologically trivial, we have \[ \widehat{H}^{-1}(H, I_G) \cong \widehat{H}^{-2}(H, \Z) \cong \widehat{H}^{-1}(H, I_H). \]

 Therefore the image of the norm map $(I_G)_H \lra (I_G)^H$ is isomorphic to \[ I_G/(\Z[G]I_H) \cong I_{G/H}. \] This proves that $C(G)_H$ is an extension \[ 0 \lra C^H \lra C(G)_H \lra I_{G/H} \lra 0. \]

It remains to verify that $C({G/H}) \cong C(G)_H$. There is a commutative diagram with exact rows
\[ \begin{tikzcd} 0 \arrow{r} & C \arrow{r} \arrow[d,"N_H"] & C(G) \arrow{r} \arrow{d}& \Z[G] \arrow{r} \arrow{d}& \Z \arrow{r} \arrow{d}& 0 \\
0 \arrow{r}& C^H \arrow{r} \arrow{d}& C(G)_H \arrow{r} \arrow{d}& \Z[G]_H \arrow{r} \arrow{d}& \Z \arrow{r} \arrow{d}& 0 \\
0 \arrow{r} & C \arrow{r} & A \arrow{r} & \Z[G]_H \arrow{r} & \Z \arrow{r} & 0.
\end{tikzcd} \]

The map from the first row to the second row on the middle two terms is the canonical projection. 
The map from the second row to third row is the pushout along the  inclusion $C^H \subset C$.
The bottom two rows imply that $C(G)_H \lra \Z[G]_H$, viewed as an element of $\Ext^2_{\Z[G/H]}(\Z, C^H)$, inflates to the extension $A \lra \Z[G]_H$ under the map \[ \Ext^2_{\Z[G/H]}(\Z, C^H) \lra \Ext^2_{\Z[G]}(\Z, C). \] On the other hand, this diagram implies that $A \lra \Z[G]_H$ is \emph{equivalent}, as a 2-extension, to the pushout of the first row along the map \[ C \xrightarrow{N_H} C^H \subset C. \]

The map \[ N_H \colon \Ext^2_{\Z[G]}(\Z, C) \lra \Ext^2_{\Z[G]}(\Z, C) \] equals multiplication by $|H|$, as can be verified using a dimension-shifting argument, since $N_H$ gives this multiplication on degree $0$.

This implies that the 2-extension $C(G)_H \lra \Z[G]_H$, when inflated to an element of $H^2(G, C)$, equals $|H| \gamma_G$. This is the same as the inflation of $\gamma_{G/H}$. Since inflation \[ H^2(G/H, C^H) \lra H^2(G, C) \] is an injection, this implies that the 2-extension $C(G)_H \lra \Z[G]_H$ is equivalent, as a 2-extension of $(G/H)$-modules, to $C({G/H}) \lra \Z[G/H]$. Moreover, this implies that the 1-extension \[ 0 \lra C^H \lra C(G)_H \lra I_{G/H} \lra 0\] is isomorphic to 
\[0 \lra C^H \lra C({G/H}) \lra I_{G/H} \lra 0,\] and hence that $ C(G)_H \cong C({G/H})$.
\end{proof}

\begin{lemma}\label{lem:fund-class-recip-2}
The isomorphisms $C(G)_G \cong C^G$ and $(I_G)_G \cong G^{\mathrm{ab}}$ identify the map \[ C(G)_G \lra (I_G)_G\] with the reciprocity map $C^G \lra G^{\mathrm{ab}}$.
\end{lemma}
\begin{proof}
 The exact sequence in Tate cohomology associated to the short exact sequence
\[ 0 \lra C \lra C(G) \lra I_G \lra 0 \]
 yields a boundary map \[ \delta \colon \widehat{H}^{-1}(G, I_G) \cong \widehat{H}^0(G, C). \] Note that $\widehat{H}^{-1}(G, I_G) = (I_G)_G$ and $\widehat{H}^0(G, C) = C^G/N_G(C)$. After identifying $(I_G)_G \cong G^{\mathrm{ab}}$, the map $\delta$ is by definition the Nakayama map.

By the construction of boundary maps in Tate cohomology, $\delta$ has the following direct description: take an element of $(I_G)_G$, lift it to $C(G)$, apply the norm map to get an element of $C(G)^G$, observe that it is contained in $C^G$, and project to $C^G/N_G(C)$.

The reciprocity map is therefore given by \[ \begin{tikzcd}
C^G \cong C(G)^G \ar[r, "\sim" inner sep = 0.3mm, "N_G^{-1}"'] & C(G)_G \twoheadrightarrow (I_G)_G. \end{tikzcd} \] The composition of the first two maps is the identification $C^G \cong  C(G)_G$ from Lemma \ref{lem:coinflation}.
\end{proof}

We can now prove the duality theorem.
\begin{proof}[Proof of Theorem \ref{thm:finite-descent}]
(1)  Since $\Z$ and $I_{G/H}$ are $\Z$-free, the sequence (\ref{e:egh}) remains exact after applying $\Hom( -, A)$.
It follows that \begin{equation} \label{e:homa}
 \Hom(\Z[G/H], A) \lra \Hom(E({G/H}), A) \end{equation}
 has cohomology 
\begin{align*}
 H^0 &= \Hom(\Z, A) = H^0(H, A), \\
  H^1 &= \Hom(C^H/N_H(C), A) = \Hom(H^{\mathrm{ab}}, A) = H^1(H, A),
 \end{align*}
 where the last equalities on each line hold since $A$ has trivial $H$-action.
    Therefore the  complex (\ref{e:homa}) has the same cohomology as
    $\tau_{\leq 1} C^*(H, A)$, and our goal is to show that these two complexes are in fact quasi-isomorphic as complexes of $\Z[G/H]$-modules.

Using the lifting property for projective $\Z[G]$-modules, we may choose a commutative diagram 
\[ \begin{tikzcd}
\Z[G^3] \arrow{r} \arrow{d} & \Z[G^2] \arrow{r} \arrow{d} & \Z[G]  \arrow{d}{=}\\
C \arrow{r} & C(G) \arrow{r} & \Z[G] \\ \end{tikzcd} \]

Applying $\Hom_{\Z[H]}(\cdot, A)$ to the first complex gives the first few terms of the complex $C^*(H,A)$. We obtain a map of chain complexes of $\Z[G/H]$-modules:
\[ \tau_{\leq 1}\Hom_{\Z[H]}(C \lra C(G) \lra \Z[G], A) \lra \tau_{\leq 1} C^*(H, A). \]

By Lemma \ref{lem:coinflation}, this becomes:
\[ \tau_{\leq 1} \Hom(C_H \lra C({G/H}) \lra \Z[G/H], A) \lra \tau_{\leq 1} C^*(H, A). \]
The map $C_H \lra C({G/H})$ factors as $C_H \xrightarrow{N_H} C^H \subset C({G/H})$.
Therefore we obtain a map of complexes:
\[ \Hom(E({G/H}) \lra \Z[G/H], A) \lra \tau_{\leq 1} C^*(H, A). \]

It remains to check that the induced map on $H^1$,
\[ \Hom(C^H/N_H(C), A) \lra H^1(H,A), \]
is the inverse of the dual of the reciprocity map, and hence
is an isomorphism.
We can reduce to the case $G = H$, using the map
\[ (C \lra C(H) \lra \Z[H]) \lra (C \lra C(G) \lra \Z[G]) \]
considered as complexes of $\Z[H]$-modules.

The partial resolution $0 \lra I_H \lra \Z[H] \lra \Z \lra 0$ of $\Z$ induces a map \[ H^1(H, A) \lra \Hom_{\Z[H]}(I_H, A) = \Hom_{\Z[H]}(I_H, A) = \Hom(H^{\mathrm{ab}}, A) \] which is the usual isomorphism. This partial resolution can be extended to $C(H) \lra \Z[H] \lra \Z \lra 0$, inducing a map \[ H^1(H, A) \cong \Hom_{\Z[H]}(I_H, A) \lra \Hom_{\Z[H]}(C(H), A).\] By Lemma \ref{lem:fund-class-recip-2}, the second map is dual to the reciprocity map.

(2) The argument of (1) implies that $\tau_{\leq 1} C^*(G,A)$ is quasi-isomorphic to \[ \tau_{\leq 1} \Hom_{\Z[G]}(C \lra C(G) \lra \Z[G], A). \]
Using Lemma \ref{lem:coinflation}, we find that the following sequence is exact:
\[ 0 \lra \Hom(E({G/H}), A) \lra \Hom_{\Z[H]}(C({G}), A) \lra \Hom_{\Z[H]}(C, A). \]
Hence, by left-exactness of $H^0(G/H, \cdot)$, the following sequence is exact as well: \[ 0 \lra \Hom_{\Z[G/H]}(E({G/H}), A) \lra \Hom_{\Z[G]}(C({G}), A) \lra \Hom_{\Z[G]}(C, A). \]
Therefore
\[
\tau_{\leq 1}\Hom_{\Z[G]}
(C\longrightarrow C(G)\longrightarrow\Z[G],A)
\]
is isomorphic to
\[
\Hom_{\Z[G/H]}
(E(G/H)\longrightarrow\Z[G/H],A).
\]

(3) The inhomogeneous 1-cocycle $\kappa_{\mathrm{univ}} \colon G \lra E(G/H)$ in $Z^1(G, E(G/H))$ corresponds to the homogeneous 1-cocycle $\Z[G^2] \lra  E(G/H)$ occuring in: \
\[ \begin{tikzcd}
\Z[G^3] \arrow{r}\arrow{d}& \Z[G^2] \arrow{r}\arrow{d}& \Z[G] \arrow{d} \\
H^{\mathrm{ab}} \arrow{r} & E(G/H) \arrow{r}& \Z[G/H]   \\
\end{tikzcd}\]

Since we have an exact sequence \[ 0 \lra H^{\mathrm{ab}} \lra E(G/H) \lra I_{G/H} \lra 0, \] we must show that $f \colon \Z[G^2] \lra I_{G/H}$ is surjective, and that $\ker(f)$ surjects onto $H^{\mathrm{ab}}$.

The map \[ \Z[G^2] \lra I_G \lra \Z[G/H] \] has image equal to $I_{G/H}$. The $H$-submodule $\Z[H^2] \subset \Z[G^2]$ is contained in $\ker(f)$. The restriction of $\Z[G^2] \lra E(G/H)$ to $\Z[H^2]$ is a homogeneous 1-cocycle $\Z[H^2] \lra E(G/H)$ whose image is $H^{\mathrm{ab}}$, corresponding to \[ \id \in \Hom(H^{\mathrm{ab}}, H^{\mathrm{ab}}) = H^1(H, H^{\mathrm{ab}}). \] 
The result follows.
\end{proof}

\subsection{Functoriality}

\begin{lemma}[Functoriality]\label{lem:compatibility-finite} Suppose we have a map of class formations $(G', C') \lra (G,C)$. 
Let $H$ be a normal subgroup of $G$, and let $H' = G' \cap H$.
Let $A$ be a $\Z[G/H]$-module. The restriction map \[ Z^1(G, A) \lra Z^1(G', A) \] is dual to $E({G'/H'}) \lra E({G/H}).$
\end{lemma}
\begin{proof}
The argument of Theorem \ref{thm:finite-descent} proves that the complex
\begin{equation} \label{e:bar1}
 0 \lra H^{\mathrm{ab}} \lra \Z[G^2]_H/\Z[G^3]_H \lra \Z[G/H] \lra \Z \lra 0 \end{equation} is isomorphic to 
\begin{equation} \label{e:bar2}
 0 \lra C^H/N_H(C) \lra E(G/H) \lra \Z[G/H]\lra \Z \lra 0. 
 \end{equation}
 {This isomorphism is obtained by taking the lift of the bar
resolution to the Tate 2-extension used in the proof of
Theorem~\ref{thm:finite-descent}, then passing to $H$-coinvariants. The quotient
of the bar complex gives (\ref{e:bar1})
while Lemma~\ref{lem:coinflation} identifies the corresponding quotient of
the Tate extension with (\ref{e:bar2}).
Under this identification the left terms are matched by the inverse of the
reciprocity map, by Lemma~\ref{lem:fund-class-recip-2}. This gives the claimed
isomorphism of complexes.}

Since we have a commutative diagram \[ \begin{tikzcd}
(H')^{\mathrm{ab}} \arrow{r}  \arrow{d} & (C')^{H'}/N_{H'}(C') \arrow{d} \\
 H^{\mathrm{ab}} \arrow{r} & C^H/N_H(C), \end{tikzcd} \]
we obtain a commutative diagram
 \[ \begin{tikzcd}
\Z[(G')^2]_{H'}/\Z[(G')^3]_{H'} \arrow{r}  \arrow{d} & E(G'/H') \arrow{d} \\
\Z[G^2]_H/\Z[G^3]_H \arrow{r} & E(G/H). \end{tikzcd} \]
 
Applying $\Hom_{\Z[G'/H']}( \cdot, A)$ to the first row  and $\Hom_{\Z[G/H]}( \cdot, A)$ to the second, we obtain the required compatibility.
\end{proof}

\subsection{Duality theorem ($G$ profinite)} 
Recall that a prodiscrete group is an inverse limit of discrete groups with surjective transition maps. Note that, for a class formation $(G,C)$, the $\Z[G/H]$-module $E({G/H})$ is prodiscrete.

\begin{theorem}\label{thm:desc-profinite}
Let $(G,C)$ be a class formation and let $H$ be an open normal subgroup of $G$. For any prodiscrete $\Z[G/H]$-module $A$:

\begin{enumerate}
\item The complex $\tau_{\leq 1} C^*(H, A)$ is quasi-isomorphic to \[ \Hom(\Z[G/H], A) \lra \Hom_{\mathrm{cts}}(E({G/H}), A). \]
\item The complex $\tau_{\leq 1} C^*(G,A)$ is quasi-isomorphic to
\[ \Hom_{\Z[G/H]}(\Z[G/H], A) \lra \Hom_{\Z[G/H], \mathrm{cts}}(E({G/H}), A). \]
\item The inhomogeneous 1-cocycle $\Z[G] \lra E({G/H})$ corresponding to the universal element $\kappa_{\mathrm{univ}} \in Z^1(G, E({G/H}))$ is surjective.
\end{enumerate}
\end{theorem}
\begin{proof}

(1) We may assume that $A$ is discrete, since \[ \varprojlim_n \Hom_{\mathrm{cts}}(E({G/H}), A/A_n) \cong \Hom_{\mathrm{cts}}(E({G/H}), \varprojlim_n A/A_n) = \Hom_{\cts}(E({G/H}), A)\] and \[ \varprojlim_n Z^1(H, A/A_n) \cong Z^1(H, \varprojlim_n A/A_n) \cong Z^1(H,A). \] 

Since $A$ is discrete, $Z^1(H, A) = \varinjlim_{H'} Z^1(H/H', A)$ as $H'$ ranges over open normal subgroups of $G$ contained in $H$, 
 and \[ \Hom_{\mathrm{cts}}(E({G/H}), A) = \varinjlim_{H'} \Hom(E({(G/H')/(H/H')}), A). \]
As $E({(G/H')/(H/H')})$ comes from the \emph{finite} class formation $(G/H', C^{H'})$, we may apply Theorem \ref{thm:finite-descent} to deduce that \[ \Hom(E({(G/H')/(H/H')}), A) \cong Z^1(H/H', A). \]

(2) and (3) follow from Theorem \ref{thm:finite-descent} via a similar argument.
\end{proof}

Suppose we have a map of class formations $(G', C') \lra (G,C)$. Let $H$ be an open normal subgroup of $G$, and let $H' = G' \cap H$.
\begin{lemma}[Functoriality]\label{lem:compatibility}
For any prodiscrete $\Z[G/H]$-module $A$, the map \[ Z^1(H, A) \lra Z^1(H', A) \] is the continuous dual $\Hom_{\mathrm{cts}}(\cdot, A)$ of the map $E({G'/H'}) \lra E({G/H})$.
\end{lemma}
\begin{proof}
This follows from the analogous result in the finite case, Lemma \ref{lem:compatibility-finite}.
\end{proof}

This concludes the construction of the generalized Ritter--Weiss module $\nabla_S^T(H)$ and the proof of its properties that were used in \S\ref{s:sr} to prove the Brumer--Stark conjecture.  What remains is the proof of Ribet's Lemma (Theorem~\ref{t:rl}).

\section{Ribet's Lemma}\label{s:cocycle}

In the remainder of the paper (\S\ref{s:cocycle}--\S\ref{s:cc}), we prove Theorem~\ref{t:rl}, Ribet's Lemma in the residually indistinguishable case.  This result strengthens Theorem~\ref{t:rli} stated in the introduction by incorporating local conditions.
We also relax the assumption that $\T$ is reduced, though we still need an assumption on the ring $K$ in which our representation $\rho$ lands.

\begin{remark} We have included the assumption of residual indistinguishability \[ \chi \equiv \psi \pmod{\fm} \] in the statement of Theorem~\ref{t:rl} since that is the salient case for this paper.  The theorem  remains true without that assumption and can be proven when $\chi \not\equiv \psi \pmod{\fm}$ using the ``$b$-coefficient'' of the representation $\rho$ in the appropriate basis, as indicated in the introduction.  The theorem was essentially proven this way in the residually distinguishable case in \cite{dk}*{\S 9}, though it was not stated in  this precise form.
\end{remark}

\begin{remark} \label{remark:top} In practical applications, the finite generation of $\rho(\T[\cG])$ {as a $\T$-module}
 will follow from topological considerations.  
As mentioned in the proof of Theorem~\ref{t:galrep}, in the setting of our proof of Brumer--Stark, the result follows since $K$ is a finite dimensional vector space over $\Q_p$ and hence has a canonical $p$-adic topology, with respect to which $\rho$ is continuous. 

Here is another sufficient condition. Suppose that $\cG$ is a
topological group, that the projections of
$\rho$ to the field factors of $K_0=\red(K)$ are absolutely
irreducible, and that $\tr\circ\rho:\cG\to\T$ is continuous
for the $\fm$-adic topology. Put
$\mathscr R=\rho(\T[\cG])$.
Absolute irreducibility and Nakayama's lemma imply
$K\mathscr R=M_2(K)$. Choose $g_1,\ldots,g_N\in\cG$
such that the matrices $\rho(g_i)$ span $M_2(K)$ over $K$.
The map
\[
\mathscr R\longrightarrow\T^N,\qquad
b\longmapsto\bigl(\tr(b\rho(g_i))\bigr)_{i=1}^N
\]
is injective, by the nondegeneracy of the matrix trace pairing.
Since $\T$ is Noetherian, $\mathscr R$ is therefore finitely
generated over $\T$. Moreover, the composite of $\rho$
with this map is continuous, since its coordinates are
$g\mapsto\tr(\rho(gg_i))$. By the Artin--Rees lemma, the
subspace topology on $\mathscr R\subset\T^N$ is its
$\fm$-adic topology, so $\rho$ is continuous for that topology.
For related finiteness results for completed Cayley--Hamilton
algebras, see \cite{chenevier}*{Lemma~4.5}.

The remarks above are useful for applications to the $\Lambda$-adic setting, i.e.\ when considering Galois representations associated to $p$-adic families of modular forms.  \end{remark}

\begin{remark}\label{remark:cont}
In the setting of Theorem~\ref{t:rl}, suppose that $\cG$ is a
topological group and give the finitely generated $\T$-module
$\rho(\T[\cG])$ its $\fm$-adic topology. If $\rho$ is
continuous for this topology and
$\chi,\psi:\cG\to\T^*$ are continuous for the $\fm$-adic
topology, then the constructed cocycle $\kappa$ is continuous.
In particular, this applies to our Brumer--Stark construction,
where $\chi$ and $\psi$ are continuous Galois characters.
\end{remark}

\subsection{Construction of the module $M$}

Our setting is as in the statement of Theorem~\ref{t:rl}.
Extend $\rho$ to a $\T$-algebra homomorphism
\[ 
\T[\cG] \longrightarrow M_2(K). 
\] 
Similarly extend $\chi, \psi$ to $\T$-algebra homomorphisms $\T[\cG] \longrightarrow \T$.
It is well-known that the congruence (\ref{e:redpoly}) extends to all $t \in \T[\cG]$:  
\begin{equation} \label{e:redpoly2}
P_{\rho(t)}(x) \equiv (x - \chi(t))(x - \psi(t)) \pmod{I}.
\end{equation}
See for instance \cite{ow}*{Lemma 3.1}.

Define two $\T$-submodules of $M_2(K)$:
 \[ \Delta_\chi = \{ \rho(t) - \chi(t) \colon t \in \T[\cG]\}, \qquad \Delta_{\psi} =   \{ \rho( t ) - \psi(t)\colon t \in \T[\cG]\}. \]
  Here $\chi(t)$ and $\psi(t)$ denote scalar matrices. {Then $\Delta_\chi$ and $\Delta_\psi$ are $\rho(\T[\cG])$-bimodules, }
and in particular $\Delta_{\chi} \Delta_{\psi} \subset \Delta_{\psi}$, where the module on the left represents the $\T$-module generated by all products $tt'$ with $t \in \Delta_\chi, t' \in \Delta_\psi$.    We then define the $\T$-module
\begin{equation} M_0 = \Delta_{\psi}/ \Delta_{\chi} \Delta_{\psi}.
\end{equation}
Note that $M_0$ is finitely generated over $\T$ from the assumption that $\rho(\T[\cG])$ is finitely generated, together with the fact that $\T$ is Noetherian.
 
\begin{lemma} The map $\kappa \colon \cG \longrightarrow \Delta_{\psi}/\Delta_{\chi}\Delta_{\psi}$ given by
\[
 g \mapsto \psi^{-1}(g)(\rho(g) - \psi(g))
\]
defines a cocycle $\kappa \in Z^1(\cG, M_0(\chi\psi^{-1}))$.
\end{lemma}
\begin{proof}
We check the cocycle condition:
\begin{align*}
\kappa(g_1g_2) \ - & \ \kappa(g_1) - \psi^{-1}\chi(g_1)\kappa(g_2) \\
= & \ \kappa(g_1g_2) - \psi^{-1}(g_1)(\rho(g_1) - \psi(g_1)) - \psi^{-1}\chi(g_1) \psi^{-1}(g_2)(\rho(g_2) - \psi(g_2)) \\
= & \  \psi^{-1}(g_1g_2)(\rho(g_1g_2) - \psi(g_1g_2)) \\
& \ \ \ \ \  - \psi^{-1}(g_1) \rho(g_1) + 1 - \psi^{-1}(g_1g_2)\chi(g_1) \rho(g_2) + \psi^{-1}(g_1)\chi(g_1) \\ 
= & \  \psi^{-1}(g_1g_2)(\rho(g_1) - \chi(g_1))(\rho(g_2) - \psi(g_2)).
\end{align*}
The last item vanishes in the quotient defining $M_0$.
\end{proof}  

Define \[ M_1 = M_0 \oplus \bigoplus_{v \in \Sigma \setminus \{v_0\}} \T y_v, \]
i.e.\ the direct sum of $M_0$ with the free $\T$-module on the set $\Sigma \setminus \{v_0\}$. 
Since $M_0$ is finitely generated over $\T$ and $\Sigma$ is finite, $M_1$ is finitely generated as well.
Let $Q \subset M_1$ denote the $\T$-submodule generated by the relations we must impose for $\kappa$ to satisfy the desired properties, namely:
\begin{itemize}
\item $\kappa(\sigma)$ for $\sigma \in \cG_{v_0}$ (in the case that $\Sigma$ is nonempty).
\item $\kappa(\sigma) - (\chi \psi^{-1}(\sigma)-1) y_v $ for $\sigma \in \cG_v$, $v \in \Sigma \setminus \{v_0\}$. 
\item $\kappa(\sigma)$ for $\sigma \in \cI_v$, $v \in \cP$.
\end{itemize}
We then define
\[ M = M_1/Q. \]
Note that $M$ inherits the property of being finitely generated over $\T$ from $M_1$.
The first {four} bullet points required of the module $M$ in Theorem~\ref{t:rl} are clearly satisfied.  

The remainder of the paper is taken up in proving the last bullet point of Theorem~\ref{t:rl}, which states
\begin{equation} \label{e:lbp}
 \left(\prod_{v \in \cP} (\xi_v(\sigma_v) - \chi(\sigma_v))\right)\Fitt_{\T}(M) \subset \widetilde{I}.  
 \end{equation}
  
\subsection{Relation Matrices for Fitting Ideal: the original matrix $D_0$} \label{s:fi}

{
We begin by choosing  auxiliary data  to write down a presentation of $M$.
These choices will remain fixed throughout \S\S\ref{s:cocycle}--\ref{s:formal}.

If $\Sigma$ is nonempty, we fix the basis for the representation $\rho$ corresponding to (\ref{e:rholocal}) for the place $v_0$, i.e.\ we choose a basis such that $\rho|_{\cG_{v_0}}$ is lower triangular. 
If $\Sigma$ is empty, choose any basis. For $\sigma\in \cG$, write
\[
\rho(\sigma)-\psi(\sigma)
=
\begin{pmatrix}
a(\sigma) & b(\sigma)\\
c(\sigma) & d(\sigma)
\end{pmatrix}.
\]

Next choose elements $g_1,\dots,g_r\in \cG$ such that the elements
\begin{equation} \label{e:rhoidef}
\rho_i:=\rho(g_i)-\psi(g_i)
\qquad (1\le i\le r)
\end{equation}
generate $\Delta_\psi$ as a $\T$-module, and such that
\[
\rho(g_i)-\chi(g_i)
\qquad (1\le i\le r)
\]
generate $\Delta_\chi$ as a $\T$-module. 
We allow repetitions in the list $g_1,\ldots,g_r$.
Enlarging this list if necessary, choose distinct indices
$j_v$, one for each $v\in\cP$, such that
\begin{equation} \label{e:jvdef}
g_{j_v}=\sigma_v.
\end{equation}
The index $j_v$ is reserved for the distinguished occurrence at
$v$, even if the same group element occurs elsewhere in the
family.

Write
\[
\rho_i=
\begin{pmatrix}
a_i & b_i\\
c_i & d_i
\end{pmatrix}.
\]

The module $M$ is generated over $\T$ by the images of
\[
\rho_1,\dots,\rho_r
\]
together with the images of the elements $y_v$ for $v\in\Sigma\setminus
\{v_0\}$. We use this ordered generating set throughout. Write 
\[ s = |\Sigma \setminus \{v_0\}|, \quad  \text{and }s=0 \text{ if } \Sigma \text{ is empty.} \]
For each $\sigma \in \cG$ choose coefficients $\alpha_{\sigma, i} \in \T$ such that
\[ \rho(\sigma) - \psi(\sigma) = \sum_{i=1}^r \alpha_{\sigma,i} \rho_i. \]

With respect to the ordered generating set
\[
\rho_1,\dots,\rho_r,\{y_v\}_{v\in\Sigma\setminus\{v_0\}}
\]
of size $r+s$, the relations needed to give a presentation of $M$ are of the following five types.}
\begin{itemize}
\item[I.] Relations among the $\rho_i$ that already hold in $\Delta_\psi$, before any quotient is taken.  Suppose
\begin{equation} \label{e:epsilondef}
\sum_{k=1}^r \epsilon_k \rho_k = 0  \qquad \text{ with } \epsilon_k \in \T.
\end{equation}
We have the relation
\begin{equation} \label{e:epsrel}
(\epsilon_1, \dots, \epsilon_r, 0, \dots, 0)  \qquad (s \text{ zeros}).
\end{equation}
\item[II.]  Relations arising from the quotient by $\Delta_\chi \Delta_\psi$.  Write 
\begin{equation} \label{e:nudef}
 \nu_i = \psi(g_i) - \chi(g_i) {\ \in \T } \end{equation} 
with the $g_i$ as in (\ref{e:rhoidef}). 
  For each $1 \le i, j \le r$, write
\begin{equation} \label{e:deltadef}
(\rho_i + \nu_i) \rho_j = \sum_{k = 1}^r \delta_{ijk} \rho_k
\end{equation}
 with $\delta_{ijk} \in \T$. 
We have the relation
\begin{equation} \label{e:delrel}
(\delta_{ij1}, \dots, \delta_{ijr}, 0, \dots, 0). 
\end{equation}
\item[III.] Relations arising from the quotient by $\kappa(\cG_{v_0})$.  For each $\sigma \in \cG_{v_0}$, we have the relation
\begin{equation} \label{e:alpharel}
(\alpha_{\sigma,1}, \dots, \alpha_{\sigma,r}, 0, \dots, 0). 
\end{equation}
\item[IV.] Relations arising from the quotient by $\kappa(\cI_v)$ for $v \in \cP$.  For each $\sigma \in \cI_v$, we have the relation
\begin{equation} \label{e:betarel}
(\alpha_{\sigma,1}, \dots, \alpha_{\sigma,r}, 0, \dots, 0). 
\end{equation}
\item[V.] Relations arising from the quotient by the conditions at $v \in \Sigma \setminus \{v_0\}$.  For each $\sigma \in \cG_v$, we have the relation
\begin{equation} \label{e:gammarel}
(\alpha_{\sigma,1}, \dots, \alpha_{\sigma,r}, 0, \dots, 0, \psi(\sigma) - \chi(\sigma), 0). 
\end{equation}
Here the entry $\psi(\sigma) - \chi(\sigma)$ occurs in the coordinate corresponding to $y_v$.
\end{itemize}

{
Thus $\Fitt_\T(M)$ is generated by the determinants of all
square matrices $D_0$ of size $r+s$ whose rows are chosen from the
relations of types (\ref{e:epsrel})--(\ref{e:gammarel}) above.
}

{
\subsection{Roadmap and notation for \S\ref{s:cocycle}--\S\ref{s:cc}}

The proof of Theorem \ref{t:rl} proceeds by transforming the
relation matrix $D_0$ through several auxiliary matrices. We summarize the notation
here.

Starting from a square relation matrix $D_0$ whose determinant contributes
to $\operatorname{Fitt}_{\T}(M)$, we perform the following steps:
\[
D_0
\;\leadsto\;
D
\;\leadsto\;
E
\;\leadsto\;
E'.
\]
The passage $D_0\leadsto D$ enlarges the chosen generating set of
$\Delta_\psi$; this does not change the determinant:
\begin{equation} \label{e:dzd}
\det(D) = \det(D_0).
\end{equation}
 The passage
$D\leadsto E$ adds rows and columns corresponding to the local conditions
at the primes in $\cP$, and gives
\begin{equation} \label{e:epd}
\det(E)=\left(\prod_{v\in \cP}(\xi_v(\sigma_v)-\chi(\sigma_v))\right)\det(D).
\end{equation}
The matrix $E'$ is obtained from $E$ by altering certain entries. The purpose of this alteration
is that
\begin{equation} \label{e:dep}
\det(E')=0.
\end{equation}
Thus, to prove the desired Fitting ideal containment, it remains to show
that
\[
\det(E')-\det(E)\in \widetilde I.
\]

In \S\ref{s:trr} we introduce a polynomial ring $R$, formal matrices
$\mathbf E,\mathbf E'\in M_{r+s+t}(R)$, and a specialization map
$\pi\colon R \lra K$, so that
\begin{equation} \label{e:peep}
\pi(\mathbf E)\equiv E \pmod{\widetilde I},
\qquad
\pi(\mathbf E')=E'.
\end{equation}
The remainder of the proof shows that
\[
\pi(\det(\mathbf E')-\det(\mathbf E))\in \widetilde I,
\]
which in view of (\ref{e:dzd})--(\ref{e:peep}) yields
\[ \left(\prod_{v\in \cP}(\xi_v(\sigma_v)-\chi(\sigma_v))\right)\det(D_0) \in \widetilde{I}. \]
As $D_0$ is an arbitrary relation matrix, we obtain the desired result (\ref{e:lbp}).

\begin{table}[h]
\centering
\caption{Notation used in \S\ref{s:cocycle}--\S\ref{s:cc}.}
\begin{tabularx}{\textwidth}{c|X}
Notation & \text Meaning \\ \hline
$D_0$ & An arbitrary square relation matrix contributing to  $\Fitt_{\T}(M)$. \\
$D$ & The enlarged relation matrix obtained from $D_0$
 by adding generators $\rho_\sigma.$   We have $\det(D)=\det(D_0)$.\\
 $E$ & The matrix obtained from $D$ by adding the local rows and columns indexed by $\cP$.\\
$E'$ & The alteration of $E$ whose determinant is 0.\\
$R$ & The polynomial ring of formal variables introduced in \S7.\\
$\pi:R\lra K$ & The specialization map sending formal variables to the actual coefficients of $\rho$.\\
$\mathbf E,\mathbf E'$ & Formal lifts of $E,E'$ with entries in $R$.\\
$J$ & The relation ideal in $R$ encoding the formal matrix identities.\\
$A$ & The trace--determinant subring of $R$; later identified with invariants.\\
$J'$ & The $\B$-stable subideal of $J$ used in the cohomological vanishing argument arising from the $b$-coefficients of the relations defining $J$.\\
$\mathcal \cC^\bullet$ & A complex of $\B$-modules resolving $R/J'$; corestricting its degree-one differential gives the exact complex ending in $J'$.
\\
$\mathcal \cD^\bullet$ & A complex of $\B$-acyclic modules receiving a map from $\mathcal \cC^\bullet$
 and satisfying $H^0(\mathcal \cD^\bullet)=R/J$.
\end{tabularx}
\end{table}
}

\subsection{Enlarging the relation matrix: $D_0 \leadsto D$} \label{s:ns}

Suppose we are given a matrix $D_0$ as above.
If, for some $v\in\Sigma\setminus\{v_0\}$, the matrix $D_0$
has no type V row associated to $v$, then its $y_v$-column
is zero and $\det(D_0)=0$. Such matrices already satisfy the
desired containment, so we may assume that every such $v$
has at least one associated type V row.

In what follows, it is notationally convenient if we replace our given matrix $D_0$ with another matrix $D$ that has the same determinant.
Namely, we may choose a larger set of generators of $\Delta_{\psi}$ so that, for all $\sigma \in  \cG$ contributing rows of types III, IV, or V to $D_0$, we include as one of our generators $\rho_i$ the element 
\[ {\rho_\sigma := \rho(\sigma) - \psi(\sigma).} \]
We introduce a fresh generator index for each occurrence of a
row of type III, IV, or V in $D_0$, even if its group element
already occurs among the initial generators or in another row.
In particular, all these new indices are distinct from the
reserved indices $j_v$.

This change replaces the matrix $D_0$ with a matrix $D$ defined as follows:
\begin{itemize} 
\item The new matrix $D$ has one additional column associated to each row of type III, IV, or V in $D_0$, i.e. to each new generator $\rho_\sigma$.
\item Each row of type I or II in $D_0$ has a corresponding row in $D$, with 0's in the new columns.
\item  Each row of type III, IV, or V in $D_0$ gets replaced by 2 rows in $D$:
\begin{itemize}
\item One row $(\alpha_{\sigma,1}, \dots, \alpha_{\sigma,r}, 0, \dots, 0, -1, 0 , 0 , \dots, 0)$, with the $-1$ in the new column corresponding to $\rho_\sigma$. Note that this is a  row of type I.
\item One row \begin{equation} \label{e:newtype}
(0, 0 , \dots, 1, 0 , \dots , 0), \end{equation}
 with the $1$ in the $\rho_\sigma$ column if we have a row of type III or IV, or 
one row \begin{equation} \label{e:newtypeV} (0,  0 , \dots, 1, 0 , \dots, 0, \psi(\sigma) - \chi(\sigma), 0, \dots, 0), \end{equation} with the $1$ in the $\rho_\sigma$ column and $\psi(\sigma) - \chi(\sigma)$ in the $y_v$ column if we have a row of type V.
\end{itemize}
\end{itemize}

 It is easy to see that the new matrix $D$ has the same determinant as the old one, but the form of the rows of type III, IV, or V has been greatly simplified.
Going forward, we will ignore $D_0$ and only work with the matrix $D$. In particular, $r$ denotes the number of generators corresponding to $D$ (i.e., the $r$ for $D_0$ plus the number of new generators included).  For the matrix $D$, the rows of type III or IV all have the form (\ref{e:newtype}), and the rows of type V all have the form (\ref{e:newtypeV}).

\subsection{Adding local factors: $D \leadsto E$}  {Let $t = |\cP|$ and write} $\cP = \{ v_1, \dots, v_t\}$.
{Recall that in the statement of Theorem~\ref{t:rl}, one chooses an element $\sigma_{v_i} \in \cG_{v_i}$ for each $v_i \in \cP$.}
 For notational simplicity, we write $\sigma_i = \sigma_{v_i}$ and define
\begin{equation} \label{e:zidef}
{z_{\sigma_{i}} = \xi_{v_i}(\sigma_{i}) - \chi(\sigma_{i})  \in \widetilde{\T}, \qquad  1 \leq i \leq t}. 
\end{equation}

In the row associated to $v_i$, we use the expression
\[
\rho(\sigma_i)-\psi(\sigma_i)=\rho_{j_{v_i}},
\]
where the index $j_{v_i}$ was defined in (\ref{e:jvdef}).
Thus the coefficient in the $\rho_j$-column of this row is
$\delta_{j_{v_i},j}$, where $\delta$ denotes the Kronecker delta.

  Associated to our matrix $D$ from \S\ref{s:ns}, we define an auxiliary $(r + s + t) \times (r+s + t)$ matrix $E$ with coefficients in {$\widetilde{\T}$} as follows. 

\begin{equation} E =  \left( \ 
\begin{tikzcd}[cramped, sep=small, row sep=-3pt]
z_{\sigma_1} &0 &[-30pt] \cdots  &[-30pt]  0 & \delta_{j_{v_{1}},1} & \cdots & \delta_{j_{v_{1}},r} & 0 & \cdots & 0 \\
0  & z_{\sigma_2} & \cdots & 0 & \delta_{j_{v_{2}},1} & \cdots & \delta_{j_{v_{2}},r} & 0 & \cdots & 0 \\
    &      & \ddots &  &   & \vdots &   &  & \vdots &  \\
    0  & 0 & \cdots &  z_{\sigma_t} & \delta_{j_{v_t},1} & \cdots & \delta_{j_{v_t},r} & 0 & \cdots & 0 \\
 \ \\
 \ \\
    &     &  \!\!\!\!\! \scaleto{0_{(r+s)\times t}}{24pt} & & & &  \scaleto{D}{24pt} & & & 
 \end{tikzcd} \  \right)
 \end{equation}

\bigskip

Clearly we have \[ \det(E) = (\prod_{i =1}^t z_{\sigma_i}) \det(D). \] {So our main goal, which is to prove Theorem \ref{t:rl} by establishing the inclusion (\ref{e:lbp}), is to show that} \[ \det(E) \in \widetilde{I}. \]

To prove this, we will define an alteration of the matrix $E$ yielding a new square matrix $E'$ of dimension $(r+s+t)$ with coefficients in $K = \Frac(\widetilde{\T})$, rather than $\widetilde{\T}$.  The motivation behind these alterations will become clear in Lemma~\ref{l:detzero} below.
For the first $t$ rows of $E$, we define
\begin{equation} \label{e:xidef}
x_{\sigma_{i}} = \xi_{v_i}(\sigma_{i}) - \psi(\sigma_{i})
\end{equation}
and we  replace the element $z_{\sigma_i}$ by 
\[ x_{\sigma_i}-a(\sigma_i) = x_{\sigma_i}-a_{j_{v_i}}. \]

\subsection{Altering to create a singular matrix: $E \leadsto E'$}

{We alter the last $r+s$ rows of $E$ to form a matrix $E'$, based on the type of relation appearing in the corresponding row of $D$:}
\begin{itemize}
\item[I.] No change.
\item[II.]  With notation as in (\ref{e:delrel}), if $i\neq j$, replace $\delta_{ijj}$ by $\delta_{ijj} - a_i - \nu_i$, replace $\delta_{iji}$ by $\delta_{iji} - d_j$, and leave the other $\delta_{ijk}$ unchanged.  If $i=j$, replace $\delta_{iii}$ by $\delta_{iii}-a_i-\nu_i-d_i$.
\item[III.] No change.
\item[IV.] The place {$v \in \cP$ corresponds to one of the first $t$ columns of $E$.  Change the entry in this row and column from 0 to 
 \[ \xi_v(\sigma) - \psi(\sigma) - a(\sigma), \]
 where $\sigma \in \cI_v$ indexes the relation.}
\item[V.] Replace $\psi(\sigma) - \chi(\sigma)$ by 
\[ \xi_v(\sigma) - \psi(\sigma) - a(\sigma), \]
 {where $\sigma \in \cG_v$ indexes the relation.}
\end{itemize}

The purpose of these replacements is that with these small changes, the matrix $E'$ obtained has vanishing determinant.

\begin{lemma} We have $\det(E') = 0.$ \label{l:detzero}
\end{lemma}

\begin{proof}  By assumption, the ring $K$ is a product of local rings, and the maximal ideal of each of these local rings is principal.  It suffices to prove the result on each such local factor, so we assume that $K$ is a local ring with maximal ideal generated by an element $\pi$.  We will prove the result by demonstrating an appropriate vector in the kernel of $E'$.

{In (\ref{e:rholocal}),} we are given that for each place $v \in \Sigma \sqcup \cP$, there exists a change of basis matrix \[ M_v = \mat{A_v}{B_v}{C_v}{D_v} \in \GL_2(K) \] such that
\begin{equation} \label{e:changebasis}
  \mat{a(\sigma)+ \psi(\sigma)}{b(\sigma)}{c(\sigma)}{d(\sigma) + \psi(\sigma)}M_v = M_v\mat{\eta_v(\sigma)}{0}{*}{\xi_v(\sigma)}  \end{equation}
for all $\sigma \in \cG_v$.  Equating the upper right entries of this matrix equation, we obtain 
\begin{equation} \label{e:db}
 D_v b(\sigma) = B_v(\xi_v(\sigma) - \psi(\sigma) - a(\sigma)). 
 \end{equation} 
Let us for the moment assume that each $D_v$ is invertible in $K$.  
Writing $\cP = \{v_1, \dots, v_t\} $ and $\Sigma \setminus \{v_0\} = \{w_1, \dots, w_s \}$, 
we define a (column) vector of length $r+s+t$ as follows:
\[ w = \left(- \frac{B_{v_1}}{D_{v_1}}, \dotsc, -\frac{B_{v_t}}{D_{v_t}}, b_1, b_2, \dots, b_r, -\frac{B_{w_1}}{D_{w_1}}, \dotsc, -\frac{B_{w_s}}{D_{w_s}}\right). \]
The fact that $E' w = 0$ follows from the definition of the alterations made {to define $E'$ from $E$}.  For example, for the first $t$ rows, the dot product of the $i$th row of $E'$ with the vector $w$ is
\begin{align*}
& \ (\xi_{v_i}(\sigma_{i}) - \psi(\sigma_{i}) - a(\sigma_{i}))\left(-\frac{B_{v_i}}{D_{v_i}}\right) + b_{j_{v_i}} \\
   = & \  (\xi_{v_i}(\sigma_{i}) - \psi(\sigma_{i}) - a(\sigma_{i}))\left(-\frac{B_{v_i}}{D_{v_i}}\right) + b(\sigma_{i}) \\
  = & \ 0 
  \end{align*}
    by (\ref{e:db}). 
 The orthogonality with the other rows is similar, and we content ourselves with brief explanations: 
  \begin{itemize}
  \item For rows of type IV or V, we again use (\ref{e:db}).
  \item For rows of type I, we use
  \[ \sum_{j=1}^r \epsilon_{ij} b_j = 0. \]
The same  holds for rows of type III---this is where our chosen basis for $\rho$ is used (recall we chose the basis such that $\rho|_{\cG_{v_0}} $ is lower triangular).
\item For rows of type II, we use the equation
\[ (a_i+\nu_i)b_j + b_i d_j = \sum_{k=1}^r \delta_{ijk} b_k \]
arising from the ``$b$''-component of the equation (\ref{e:deltadef}) defining the $\delta_{ijk}$.
  \end{itemize}
  The second bulleted assumption in the statement of Theorem~\ref{t:rl} regarding the irreducibility of the projections of $\rho$ implies that the elements $b_1, \dots, b_r$ generate $K$ as a $K$-module.  Therefore the equation $E'w = 0$ implies that $\det(E') = 0$ as desired.
  
  This argument is easily adapted in the case that some $D_{v_i}$ is not a unit in $K$.  Recall we have reduced to the setting that $K$ is a local ring with maximal ideal generated by an element $\pi$ (which is necessarily nilpotent).  For each $v \in (\Sigma \setminus \{v_0\}) \sqcup \cP$, we can write $D_{v} = u_v \pi^{e_v}$ for some unit $u_v \in K^*$ and integer $e_v \ge 0$.   Let $m$ denote the maximal $e_v$ (if some $D_{v} =0$, then we let 
$m = e_v$ be the minimal exponent $k$ such that $\pi^k = 0$).  Equation (\ref{e:db}) yields
\[  \pi^m b(\sigma) = \frac{B_{v_i}}{u_{v_i}} \pi^{m - e_{v_i}}(\xi_{v_i}(\sigma) - \psi(\sigma) - a(\sigma)). 
\]
We therefore define
\begin{align*}
 w = \Big( & - \frac{B_{v_1}}{u_{v_1}} \pi^{m-e_{v_1}}, \dotsc, -\frac{B_{v_t}}{u_{v_t}}\pi^{m-e_{v_t}}, b_1 \pi^m, b_2 \pi^m, \dots, b_r \pi^m, \\
 & -\frac{B_{w_1}}{u_{w_1}}\pi^{m-e_{w_1}}, \dotsc, -\frac{B_{w_s}}{u_{w_s}}\pi^{m-e_{w_s}}\Big). 
 \end{align*}
The equation $E'w = 0$ follows as before.  Now for the $v$ yielding the maximal $e_v = m$, the component of $w$ is 
$-B_v/u_v \in K^*$.  Here we use the fact that $D_v \not\in K^* \Longrightarrow B_v \in K^*$ since $M_v \in \GL_2(K)$.  Therefore $E' w = 0$ again implies $\det(E') = 0$.
\end{proof}

Our goal is to show that $\det(E) \in \widetilde{I}$, and we have shown that $\det(E') = 0$.  It therefore suffices to show that 
$\det(E') - \det(E) \in \widetilde{I}$.  In other words, {we must show that }the alterations used to pass from $E$ to $E'$ are small enough to leave the determinant unchanged modulo $\widetilde{I}$. Let us motivate our strategy to prove this with an example.

\subsection{Example: the trace identity in a small case} \label{s:example}

Suppose $r=2$ and that $\Sigma \sqcup \cP$ is empty.   We consider a matrix with only rows of type II, namely

\[ E = D = \mat{\delta_{121}}{\delta_{122}}{\delta_{211}}{\delta_{212}},  \text{ whence }  E' = \mat{\delta_{121} - d_2}{\delta_{122} - a_1 - \nu_1}{\delta_{211} - a_2 - \nu_2}{\delta_{212} - d_1}. \]
As shorthand, write $t_i = a_i + d_i$ for the trace of $\rho_i$ and $t_{12}$ for the trace of $\rho_1 \rho_2$.
By multilinearity of the determinant, we have \begin{align}
\det(E') - \det(E) &= - \det\mat{d_2}{a_1+\nu_1}{\delta_{211}}{\delta_{212}} - \det \mat{\delta_{121}}{\delta_{122}}{a_2+\nu_2}{d_1} \\
& \ \ \ \ + \det \mat{d_2}{a_1 + \nu_1}{a_2+\nu_2}{d_1} \nonumber \\
& = ( t_1 + \nu_1) \delta_{211} - (d_1 \delta_{211} + d_2 \delta_{212})  +   \nonumber \\
& \ \  \ \  + (t_2 + \nu_2) \delta_{122} - (d_1 \delta_{121} + d_2 \delta_{122} )  +   \nonumber\\
&  \ \ \ \   +  (d_1d_2 - (a_1 + \nu_1)(a_2 + \nu_2))   \nonumber \\
\begin{split}
&= ( t_1 + \nu_1) \delta_{211} - (c_2b_1 + (d_2 + \nu_2)d_1)  +   \label{e:exdel}  \\
& \ \  \ \  + (t_2 + \nu_2) \delta_{122} - (c_1b_2 + (d_1 + \nu_1)d_2 )  +   \\
&  \ \ \ \   +  (d_1d_2 - (a_1 + \nu_1)(a_2 + \nu_2))    
\end{split} \\
&= ( t_1 + \nu_1) \delta_{211} + (t_2 + \nu_2) \delta_{122} - (t_{12} + t_1 \nu_2 + t_2 \nu_1 + \nu_1\nu_2 ).
\label{e:example}
 \end{align}
 Equation (\ref{e:exdel}) uses the definition (\ref{e:deltadef}) of the $\delta$'s.
By the fundamental assumption concerning the representation $\rho$ (namely, the first bullet point in the statement of Theorem~\ref{t:rl}), we know that $\tr(\rho(\sigma)) - \psi(\sigma) - \chi(\sigma) \in I$. Recalling the definition of $\nu_i$ from (\ref{e:nudef}), it directly follows that $t_i + \nu_i \in I$.  A similar elementary computation that is generalized below also shows that $t_{12} + t_1 \nu_2 + t_2 \nu_1 + \nu_1\nu_2 \in I$. Hence (\ref{e:example}) lies in $I$ as desired.

\subsection{Trace identities modulo $I$} \label{s:ti}

We now return to the general setting.  We prove that certain expressions involving traces and determinants, generalizing the elements in parentheses in (\ref{e:example}), lie in $I$.

\begin{lemma}\label{l:dets} We have $\det(\Delta_\psi) \subset I$.
\end{lemma}
\begin{proof}
{For all $t \in \T[\cG]$,} we have \[ \det(\rho(t) - \psi(t)) = P_{\rho(t)}(\psi(t)) \equiv 0 \pmod{I} \] by (\ref{e:redpoly2}).
\end{proof}

Let $\T\langle X_1, \dots, X_r \rangle$ denote the polynomial algebra over $\T$ in $r$ noncommuting variables $X_i$.  Define a $\T$-algebra homomorphism \[ V \colon \T\langle X_1, \dots, X_r \rangle \longrightarrow \T, \qquad X_i \mapsto -\nu_i, \]
where $ \nu_i = \psi(g_i) - \chi(g_i)$ as in (\ref{e:nudef}).
\begin{lemma}\label{l:tr-char}
For any $f \in  \T\langle X_1, \dots, X_r \rangle$ with constant term 0, we have {the congruence
\[ \tr(f(\rho_1, \dotsc, \rho_r)) \equiv V(f) \pmod{I} \]
in $\T$.}
\end{lemma}
\begin{proof}
It suffices to prove the result for a monomial
\[
X_{a_1}X_{a_2}\cdots X_{a_n},
\qquad
n\geq1,\quad 1\leq a_k\leq r.
\]
To simplify notation, put
\[
\begin{gathered}
g_k=g_{a_k},\qquad
\rho_k=\rho(g_k)-\psi(g_k),\\[2pt]
\chi_k=\chi(g_k),\qquad
\psi_k=\psi(g_k),\qquad
\nu_k=\psi_k-\chi_k.
\end{gathered}
\]
We calculate the following, where the congruence holds modulo $I$.  All products below are taken in the original order of the factors, and
the omitted higher-order sums are indexed by strictly increasing
tuples of indices.
\begin{align}
\tr\Big(\prod_{i=1}^n \rho_{i}\Big) &= \tr\Big(\prod_{i=1}^n (\rho(g_{i}) - \psi(g_{i}))\Big) \nonumber \\
&= \tr\Big(\prod_{i=1}^n \rho(g_{i})\Big) - \sum_{i =1}^n \psi_{i} \tr\Big(\prod_{k \neq i} \rho(g_{k})\Big)  \nonumber \\
& \ \ \ \ \  +  \sum_{1 \leq i < j \leq n} \psi_{i}\psi_{j} \tr\Big(\prod_{k \neq i, j} \rho(g_{k})\Big) + \cdots + (-1)^n 2 \prod_{i=1}^n \psi_{i}  \nonumber \\
\begin{split}
&\equiv \Big(\prod_{i=1}^n \chi_{i} + \prod_{i=1}^n \psi_{i}\Big) - \sum_{i=1}^n \psi_{i} \Big(\prod_{k \neq i} \chi_{k} + \prod_{k \neq i} \psi_{k}\Big)   \\
 & \ \ \ \ \ + \sum_{1 \leq i < j \leq n} \psi_{i} \psi_{j} \Big(\prod_{k \neq i, j} \chi_{k} + \prod_{k \neq i, j} \psi_{k}\Big) + \cdots + (-1)^n 2\prod_{i=1}^n \psi_{i}. \label{e:ti}
\end{split}
\end{align}

For $0\le k<n$, the terms selecting exactly $k$ scalar factors contribute $(-1)^k\binom nk\prod_i\psi_i$; the final term contributes $2(-1)^n\prod_i\psi_i$.
 Taking into account the signs, we see that the total contribution of $\prod_i \psi_i$ to the sum is $(1-1)^n + (-1)^n = (-1)^n$. 
The expression (\ref{e:ti}) becomes
{\begin{align*}
 \prod_{i=1}^n \chi_{i} & - \sum_{i=1}^n \psi_{i} \prod_{k \neq i} \chi_{k}  + \sum_{1 \leq i < j \leq n} \psi_{i} \psi_{j} \prod_{k \neq i, j} \chi_{k} + \cdots + (-1)^n \prod_{i=1}^n \psi_{i} \\
 & = \prod_{i=1}^n (\chi_{i} - \psi_{i}) = \prod_{i=1}^n (-\nu_{i}). 
\end{align*}
}
This is the desired result.
\end{proof}

\section{Formal Variables} \label{s:formal}

{
The output of \S\ref{s:cocycle} was that, to prove Theorem \ref{t:rl}, it
remains to show
\[
\det(E')-\det(E)\in \widetilde I.
\]
The example of \S\ref{s:example} illustrates how such a containment can be
proved by rewriting $\det(E')-\det(E)$ in terms of traces and determinants,
to which the identities of \S\ref{s:ti} apply. In the general case, however,
writing down such an identity explicitly seems intractable.

The purpose of this section is to replace this concrete calculation by a
formal one. We introduce a polynomial algebra $R$ and a specialization
homomorphism
\[
\pi\colon R \longrightarrow K=\Frac(\widetilde{\T})
\]
so that the matrices $E$ and $E'$ are obtained by specializing formal
matrices with entries in $R$. We will then reduce the desired containment
to a statement in the quotient $R/J$, where invariant theory can be applied.}

\subsection{The ring $R$} \label{s:trr}

We first clarify the indexing used in the formal construction.
Local subscripts $\sigma$ label occurrences: one for each row of
type III, IV, or V of $D$, and one distinguished occurrence at each
$v\in\cP$. Each occurrence is attached to a group element and to
its own matrix-generator index. Evaluations such as $\rho(\sigma)$
and $\xi_v(\sigma)$, and membership statements such as
$\sigma\in\cI_v$, refer to the attached group element. Distinct
occurrences have distinct assigned indices, even when their
attached group elements coincide. For the distinguished occurrence
at $v\in\cP$, the assigned index is $j_v$.

We define the polynomial ring
\begin{equation} \label{e:r0def}
R_0 =
\Z[
\boldsymbol{\nu}_i,
\boldsymbol{\epsilon}_i,
\boldsymbol{\delta}_{ijk},
\boldsymbol{x}_\sigma
],
\end{equation}
whose independent variables and specialization homomorphism
\[ \pi:R_0\lra K \] are specified as follows:
\begin{itemize}
\item
For $1\leq i\leq r$, we introduce a variable
$\boldsymbol{\nu}_i$ with
\[
\pi(\boldsymbol{\nu}_i)=\nu_i,
\]
where $\nu_i$ is defined in \eqref{e:nudef}.

\item
For each row of type I in $D$, we introduce a separate set of
variables
\[
\boldsymbol{\epsilon}_1,\ldots,\boldsymbol{\epsilon}_r
\]
with $\pi(\boldsymbol{\epsilon}_i)=\epsilon_i$, where the
$\epsilon_i$ are the coefficients of that row.
For each row of type II associated to a pair $(i,j)$, we introduce
a separate set of variables
\[
\boldsymbol{\delta}_{ij1},\ldots,\boldsymbol{\delta}_{ijr}
\]
with $\pi(\boldsymbol{\delta}_{ijk})=\delta_{ijk}$.
The row labels are suppressed in this notation: different rows
have different sets of variables, even if two type II rows are
associated to the same pair $(i,j)$.

\item
For each occurrence $\sigma$ corresponding to a row of type IV
or V in $D$, or to one of the first $t$ rows of $E$, we
introduce a variable $\boldsymbol{x}_\sigma$. If the occurrence
is associated to $v$, we set
\[
\pi(\boldsymbol{x}_\sigma)=\xi_v(\sigma)-\psi(\sigma).
\]
In particular, there is a separate variable
$\boldsymbol{x}_{\sigma_v}$ for the distinguished occurrence
at each $v\in\cP$.
\end{itemize}

We next define
\begin{equation} \label{e:r1def}
R_1 =
R_0[\boldsymbol{a}_i,\boldsymbol{b}_i,
    \boldsymbol{c}_i,\boldsymbol{d}_i]_{i=1}^r
\end{equation}
and extend $\pi$ to $R_1$ by
\[
\pi(\boldsymbol{a}_i)=a_i,\qquad
\pi(\boldsymbol{b}_i)=b_i,\qquad
\pi(\boldsymbol{c}_i)=c_i,\qquad
\pi(\boldsymbol{d}_i)=d_i,
\]
where $a_i,b_i,c_i,d_i$ are the entries of $\rho_i$.

For a local occurrence $\sigma$ with assigned matrix-generator
index $i$, we use the abbreviations
\[
\boldsymbol{a}_\sigma=\boldsymbol{a}_i,\qquad
\boldsymbol{b}_\sigma=\boldsymbol{b}_i,\qquad
\boldsymbol{c}_\sigma=\boldsymbol{c}_i,\qquad
\boldsymbol{d}_\sigma=\boldsymbol{d}_i,\qquad
\boldsymbol{\nu}_\sigma=\boldsymbol{\nu}_i.
\]
Thus, for the distinguished occurrence at $v\in\cP$,
\[
\boldsymbol{a}_{\sigma_v}=\boldsymbol{a}_{j_v},
\qquad
\boldsymbol{\nu}_{\sigma_v}=\boldsymbol{\nu}_{j_v},
\]
and similarly for the other matrix entries. These are abbreviations
for already defined variables, not additional relations. Variables
attached to distinct indices remain independent even when their
images under $\pi$ coincide.

Finally, let $\cB_{v_0}$ denote the set of labels of the type III
rows of $D$, with their assigned supplemental generator indices.
If $\Sigma=\emptyset$, we take $\cB_{v_0}=\emptyset$.
Define
\begin{equation} \label{e:rdef}
R =
R_1/(\boldsymbol{b}_\sigma:\sigma\in\cB_{v_0}).
\end{equation}
Our chosen basis makes $\rho|_{\cG_{v_0}}$ lower triangular, so
$\pi(\boldsymbol{b}_\sigma)=0$ for every
$\sigma\in\cB_{v_0}$. Consequently, $\pi$ descends to a
homomorphism
\[
\pi:R\lra K.
\]

{
\subsection{Formal lifts: $E, E' \leadsto \boldsymbol{E}, \boldsymbol{E}'$} \label{s:edef}
We next define matrices \[ \boldsymbol{E}, \boldsymbol{E}' \in M_{r+s+t}(R) \] whose images under $\pi$ are related to $E, E'$, respectively, as follows:}
\[
\pi(\mathbf E)-E \in M_{r+s+t}(\widetilde I),
\qquad
\pi(\mathbf E')=E' \]
  The matrix $\boldsymbol{E}$ is defined from $E$ simply by making every entry bold, with three exceptions. 
\begin{itemize}
\item In the first $t$ rows of $E$ we replace each 
$z_{\sigma_i}$ {appearing as an entry of $E$ by $\boldsymbol{x}_{\sigma_i} + \boldsymbol{\nu_{\sigma_i}}$ as the corresponding entry of $\boldsymbol{E}$}. 
\item For the rows of type IV, we insert  $\boldsymbol{x}_\sigma + \boldsymbol{\nu}_\sigma$ into the column of $\boldsymbol{E}$ corresponding to $v$. 
\item For the rows of type V, we replace each  $\psi(\sigma) - \chi(\sigma)$ {appearing as an entry of $E$ by $\boldsymbol{x}_\sigma + \boldsymbol{\nu}_\sigma$ as the corresponding entry of $\boldsymbol{E}$}.
\end{itemize}
To motivate these replacements, note that
\[ \psi(\sigma) - \chi(\sigma) = \nu_\sigma \equiv -\tr(\rho_{\sigma}) \pmod{I}, \]
which leads to the following result.
\begin{lemma} \label{l:bolde}
$\pi(\boldsymbol{E}) - E$ has coefficients in $\widetilde{I}$. 
\end{lemma}
\begin{proof}
For the  entries of $\boldsymbol{E}$ defined from $E$  by making every entry bold, the corresponding entry of $\pi(\boldsymbol{E}) - E$ is 0.
We must check that for the exceptions to this rule outlined in the three bullet points above, the corresponding entry of $\pi(\boldsymbol{E}) - E$ still lands in $\widetilde{I}$.   
For the first bullet point, this follows directly from the definitions of $\pi$, and $\nu_i$ and $z_i$ given in equations (\ref{e:nudef}) and (\ref{e:zidef}), respectively. 
For the rows of type IV, we use the assumption that $(\xi_v)|_{\cI_v} \equiv \chi|_{\cI_v} \pmod{\widetilde{I}}$. For the rows of type V, we use the assumption $\xi_v \equiv \psi|_{\cG_v} \pmod{\widetilde{I}}$.
 \end{proof}

Since $E$ has coefficients in $\widetilde{\T}$, the fact that  $\pi(\boldsymbol{E}) - E$ has coefficients in $\widetilde{I}$ implies that 
\begin{equation} \label{e:econg}
 \pi(\det(\boldsymbol{E})) \equiv \det(E) \pmod{\widetilde{I}}. \end{equation}

We next define a matrix $\boldsymbol{E'}$ from $\boldsymbol{E}$ by enacting the same alterations used to define $E'$ from $E$. More precisely:
\begin{itemize} \item For the first $t$ rows of $\boldsymbol{E}$, the rows of type IV, and the rows of type V, replace
 \[\boldsymbol{x}_\sigma + \boldsymbol{\nu}_\sigma  \quad \text{ by } \quad 
 \boldsymbol{x}_\sigma - \boldsymbol{a}_{\sigma}. \]
  \item For rows of type I or III, no change.
 \item For rows of type II, if $i\neq j$ replace $\boldsymbol{\delta}_{ijj}$ by $\boldsymbol{\delta}_{ijj} - \boldsymbol{a}_i - \boldsymbol{\nu}_i$, replace $\boldsymbol{\delta}_{iji}$ by $\boldsymbol{\delta}_{iji} - \boldsymbol{d}_j$, and leave the other $\boldsymbol{\delta}_{ijk}$ unchanged.  If $i=j$, 
 replace $\boldsymbol{\delta}_{iii}$ by $\boldsymbol{\delta}_{iii} - \boldsymbol{a}_i - \boldsymbol{\nu}_i - \boldsymbol{d}_i.$
\end{itemize}
We then have $\pi(\boldsymbol{E}') = E'$.  Combining with (\ref{e:econg}), in order to prove our desired result $\det(E') - \det(E) \in \widetilde{I}$, it  suffices to prove that \[ \pi(\det(\boldsymbol{E}') - \det(\boldsymbol{E})) \in \widetilde{I}. \]

\subsection{Relation Ideal} \label{s:ri}

We now define the polynomial relations that allow us to reduce $\det(\boldsymbol{E}') - \det(\boldsymbol{E})$ to an expression involving traces and determinants, as in the example of \S\ref{s:example}.  
Define \[ \boldsymbol{\rho_i} = \mat{\boldsymbol{a}_i}{\boldsymbol{b}_i}{\boldsymbol{c}_i}{\boldsymbol{d}_i} \in M_2(R). \]
We let $J \subset R$ be the ideal generated by the following:
\begin{itemize}
\item The 4 matrix coefficients of 
\begin{equation} 
\sum_{i=1}^r \boldsymbol{\epsilon}_i \boldsymbol{\rho}_i \label{e:b1}
\end{equation}
 for each row of type I in $D$.
\item The 4 matrix coefficients of 
\begin{equation} \label{e:b2}
(\boldsymbol{\rho}_i + \boldsymbol{\nu}_i) \boldsymbol{\rho}_j - \sum_{k = 1}^r \boldsymbol{\delta}_{ijk} \boldsymbol{\rho}_k
\end{equation} for each row of type II in $D$.
 \item For $v\in\cP$, let $\cB_v$ be the disjoint union of
the distinguished label $\sigma_v$, assigned to the initial
index $j_v$, and the labels of the type IV rows of $D$
associated to $v$. For $v\in\Sigma\setminus\{v_0\}$, let
$\cB_v$ be the set of labels of the type V rows of $D$
associated to $v$. The row labels carry the fresh generator
indices introduced in \S\ref{s:ns}. These indexing sets are
pairwise disjoint, and distinct labels need not represent
distinct group elements.

For $v \in (\Sigma\setminus\{v_0\}) \sqcup \cP$ and each $\sigma \in \mathcal{B}_v$, recall the variables  $\boldsymbol{a}_\sigma, \boldsymbol{b}_\sigma, \boldsymbol{c}_\sigma, \boldsymbol{d}_\sigma$ corresponding to the entries of $\boldsymbol{\rho}_{\sigma} (=\boldsymbol{\rho}_i$ for some $i$).
For each pair of distinct $\sigma, \tau \in \mathcal{B}_v$ we include as generators for $J$ the 4 matrix coefficients of
\begin{equation} \label{e:b4}
 \mat{A(\sigma, \tau)}{ B(\sigma, \tau) }{C(\sigma,\tau)}{D(\sigma, \tau)}
\end{equation}
where 
\begin{equation} \label{e:ABCD}
\begin{split}
A(\sigma, \tau) & = \boldsymbol{b}_\sigma\boldsymbol{c}_\tau - (\boldsymbol{x}_{\tau} - \boldsymbol{d}_\tau)(\boldsymbol{x}_{\sigma} - \boldsymbol{a}_\sigma),
 \\
B(\sigma, \tau) & = \boldsymbol{b}_\sigma(\boldsymbol{x}_{\tau} - \boldsymbol{a}_\tau) - 
\boldsymbol{b}_\tau(\boldsymbol{x}_{\sigma} - \boldsymbol{a}_\sigma), \\
C(\sigma, \tau) & = \boldsymbol{c}_\sigma(\boldsymbol{x}_{\tau} - \boldsymbol{d}_\tau) - 
\boldsymbol{c}_\tau(\boldsymbol{x}_{\sigma} - \boldsymbol{d}_\sigma),  \\
D(\sigma, \tau) & = \boldsymbol{b}_\tau\boldsymbol{c}_\sigma - (\boldsymbol{x}_{\sigma} - \boldsymbol{d}_\sigma)(\boldsymbol{x}_{\tau} - \boldsymbol{a}_\tau).
\end{split}
\end{equation}
Of course, the elements $D(\sigma, \tau)$ are redundant since $D(\sigma, \tau) = A(\tau, \sigma)$, but our reasons for describing these generators in terms of the matrix (\ref{e:b4}) will become apparent when we prove Lemma~\ref{l:stable} below.
\end{itemize}

\begin{lemma} We have $\pi(J) =0$. \end{lemma}

\begin{proof}
The fact that the elements in (\ref{e:b1}) and (\ref{e:b2}) map to 0 under $\pi$ is clear since these are precisely the relations defining $\epsilon_i, \delta_{ijk}$ in $\T$.  For the generators in (\ref{e:b4}) we compute as follows.  
\begin{equation} \label{e:pibst}
 \pi(B(\sigma, \tau)) = b(\sigma)(\xi_v(\tau) - \psi(\tau) - a(\tau)) - b(\tau)(\xi_v(\sigma) - \psi(\sigma) - a(\sigma)).
\end{equation}
Using (\ref{e:db}), it is easy to see that the expression (\ref{e:pibst}) vanishes if it is multiplied by $B_v$ or by $D_v$.
Since $B_v$ and $D_v$ generate $K$ as a $K$-module, it follows that (\ref{e:pibst}) vanishes.  The proof that $A(\sigma, \tau), C(\sigma, \tau),$ and $D(\sigma, \tau)$ also vanish under $\pi$ is similar and uses the equation
\[
 B_v c(\sigma) = D_v(\xi_v(\sigma) - \psi(\sigma) - d(\sigma))
\]
obtained from equating the lower right entries of (\ref{e:changebasis}).
\end{proof}

\subsection{Subring of traces and determinants} \label{s:subring}

In order to apply the results on traces and determinants from \S\ref{s:ti}, we define a certain subring $A \subset R$ as follows.
Let $A_1 \subset R_1$ denote the sub $R_0$-algebra generated by the traces and determinants of all matrices in the noncommutative $\Z$-algebra generated by the matrices $\boldsymbol{\rho}_i$:
\[ A_1 = R_0[\boldsymbol{a}_i + \boldsymbol{d}_i, \boldsymbol{a}_i\boldsymbol{d}_i - \boldsymbol{b}_i\boldsymbol{c}_i, 
\boldsymbol{a}_i\boldsymbol{a}_j + \boldsymbol{b}_i\boldsymbol{c}_j + \boldsymbol{c}_i\boldsymbol{b}_j + \boldsymbol{d}_i\boldsymbol{d}_j,  \dotsc ]. \]
Let $A_0$ denote the image of $A_1$ in $R = R_1/(\boldsymbol{b}_{\tau} \colon \tau \in \cB_{v_0})$. 
Let $A \subset R$ denote the ring generated over $A_0$ by the $\boldsymbol{d}_\tau$ for $\tau \in \cB_{v_0}$.
Denote by $\overline{A}$ the image of $A$ in $R/J$.  
{
The goal of the remainder of this subsection is to show:
\begin{prop} \label{p:goalepe}
 In order to prove the desired result 
 \[ \det(E') - \det(E) \in \widetilde{I}, \]
 it suffices to prove that the image of 
\[ \boldsymbol{e} := \det(\boldsymbol{E}') - \det(\boldsymbol{E}) \] in $R/J$ lies in the subring $\overline{A}$, the image of $A$ in $R/J$.
\end{prop}
}

For this, it is important that $\boldsymbol{e}$ lies in the following ideal of $R$:
\[ I_R = \langle \boldsymbol{a}_i + \boldsymbol{\nu}_i,   \boldsymbol{b}_i,  \boldsymbol{c}_i,  \boldsymbol{d}_i   \rangle. \]
The fact that  $\boldsymbol{e}$ lies in $I_R$ follows from multilinearity of the determinant since every entry of 
$\boldsymbol{E}' - \boldsymbol{E}$ lies in $I_R$.

\begin{lemma} We have $\pi(A \cap (I_R, J)) \subset \widetilde{I}$. \label{l:pia}
\end{lemma}

\begin{proof}
Let  $\cB_{v_0} = \{\tau_1, \dots, \tau_h\}$ be the set defined in the paragraph preceding (\ref{e:rdef}).
Define a $R_0$-algebra homomorphism
\[ \boldsymbol{V} \colon R_0 \langle X_1, \dots, X_r \rangle \longrightarrow R_0, \qquad  X_i \mapsto -\boldsymbol{\nu_i}. \]
Any element of $A$ can be written as a polynomial with coefficients in $R_0$ in the expressions
\begin{align}
&  \tr(f(\boldsymbol{\rho}_1, \dots, \boldsymbol{\rho}_r)) - \boldsymbol{V}(f), \label{e:trv} \\
&  \det(f(\boldsymbol{\rho}_1, \dots, \boldsymbol{\rho}_r)), \label{e:detf} \\
& \text{ and } \boldsymbol{d}_{\tau_1}, \dots, \boldsymbol{d}_{\tau_h}, \label{e:ds} \end{align}
  as $f$ ranges over polynomials with no constant term.  As we now explain, the elements in (\ref{e:trv})--(\ref{e:ds}) have images under $\pi$ lying in $\widetilde{I}$.  This holds by  Lemma \ref{l:tr-char} for (\ref{e:trv}) and 
  by Lemma~\ref{l:dets} for (\ref{e:detf}).  For (\ref{e:ds})
 we note that {
 the choice of basis made at the outset of \S\ref{s:fi} (such that (\ref{e:rholocal}) holds for the place $v_0$) ensures that
  \[ \pi(\boldsymbol{d}_{\tau_i}) = \xi_{v_0}(\tau_i) - \psi(\tau_i). \]
By assumption (\ref{e:Sigmaassumption}) this element lies in $\widetilde{I}$.}
  It is also easy to check that each  expression (\ref{e:trv})--(\ref{e:ds}) lies in $I_R$.  For example, we note that modulo $I_R$ we have $\boldsymbol{\rho}_i \equiv \mat{-\boldsymbol{\nu}_i}{0}{0}{0}$,
  from which it follows that  (\ref{e:trv}) lies in $I_R$. The computation that (\ref{e:detf})  lies in $I_R$ is similar.  The elements in (\ref{e:ds}) lie in $I_R$ by definition.
  
  To conclude, it therefore suffices to prove that $\pi(R_0 \cap (I_R, J)) \subset \widetilde{I}$.  To calculate the intersection $R_0 \cap (I_R, J)$, we note that
  \[ R / I_R \cong R_0, \] 
  with $\boldsymbol{a}_i \mapsto - \boldsymbol{\nu_i}$ and $\boldsymbol{b}_i, \boldsymbol{c}_i, \boldsymbol{d}_i \mapsto
  0$.  
  To compute the quotient by $J$ we then mod out by the elements of $J$ with these substitutions made, and obtain
  \[ R/(I_R, J) \cong R_0/I_{R_0}
  \]
  where $I_{R_0} \subset R_0$ is the ideal generated by
  \begin{itemize}
 \item  $\sum_{i=1}^r \boldsymbol{\epsilon}_i \boldsymbol{\nu}_i$ for each row of type I,
 \item $\sum_{k=1}^r \boldsymbol{\delta}_{ijk} \boldsymbol{\nu}_k$ for each row of type II,
\item $(\boldsymbol{x}_\sigma + \boldsymbol{\nu}_{\sigma}) \boldsymbol{x}_\tau$ for each distinct $\sigma, \tau \in \cB_v$, for each $v \in (\Sigma \setminus \{v_0\}) \sqcup \cP.$
\end{itemize}
Thus $R_0 \cap (I_R, J) = I_{R_0}$ and it remains to check that each of the bulleted expressions maps to $\widetilde{I}$ under $\pi$.  For instance, for a row of type I we have
\[ 0 = \tr\Big(\sum_{i=1}^r  {\epsilon}_i {\rho}_i \Big)
\equiv - \sum_{i=1}^r {\epsilon}_i {\nu}_i \pmod{I} \]
by Lemma~\ref{l:tr-char}.  The proof for the second bulleted item is similar.  For the third, we note that 
\begin{equation} \label{e:pisc}
 \pi((\boldsymbol{x}_\sigma + \boldsymbol{\nu}_{\sigma}) \boldsymbol{x}_\tau) = (\xi_v(\sigma) - \chi(\sigma))(\xi_v(\tau) - \psi(\tau)). 
 \end{equation}
All terms in \eqref{e:pisc} lie in $\widetilde{\T}$.  If
$v\in\Sigma$, then the second factor in \eqref{e:pisc} lies in
$\widetilde I$ by \eqref{e:Sigmaassumption}.

Now suppose that $v\in\cP$.  If $\sigma\in\cI_v$, then the first
factor in \eqref{e:pisc} lies in $\widetilde I$ by
\eqref{e:Passumption}.  The only remaining case is
\[
\sigma=\sigma_v,\qquad \tau\in\cI_v.
\]
By \eqref{e:Passumption} and the trace congruence following from
\eqref{e:redpoly}, we have
\begin{equation} \label{e:taucong}
\xi_v(\tau)\equiv\chi(\tau),
\qquad
\eta_v(\tau)\equiv\psi(\tau)
\pmod{\widetilde I}.
\end{equation}
Applying the same trace congruence to $\sigma_v$ and
$\sigma_v\tau$ respectively gives
\begin{align}
\eta_v(\sigma_v)+\xi_v(\sigma_v)
& \equiv
\chi(\sigma_v)+\psi(\sigma_v)
\pmod{\widetilde I}, \label{e:ex1} \\
\eta_v(\sigma_v)\eta_v(\tau)
+\xi_v(\sigma_v)\xi_v(\tau)
& \equiv
\chi(\sigma_v)\chi(\tau)
+\psi(\sigma_v)\psi(\tau)
\pmod{\widetilde I}. \label{e:ex2}
\end{align}
Subtracting $\psi(\tau)$ times (\ref{e:ex1}) from (\ref{e:ex2}),
and using (\ref{e:taucong}) yields
\[
\bigl(\xi_v(\sigma_v)-\chi(\sigma_v)\bigr)
\bigl(\chi(\tau)-\psi(\tau)\bigr)
\in\widetilde I.
\]
Since $\xi_v(\tau)-\chi(\tau)\in\widetilde I$, it follows that
\[
\bigl(\xi_v(\sigma_v)-\chi(\sigma_v)\bigr)
\bigl(\xi_v(\tau)-\psi(\tau)\bigr)
\in\widetilde I.
\]
Thus every expression in the third bullet defining $I_{R_0}$ maps
into $\widetilde I$, which concludes the proof.
\end{proof}

We can now establish:

\begin{proof}[Proof of Proposition~\ref{p:goalepe}]
 If the image of $\boldsymbol{e}$ in $R/J$ lies in $\overline{A}$, then we may write $\boldsymbol{e} = a+ j$ with $a \in A$ and $j \in J$.  Then $a = \boldsymbol{e} - j \in A \cap (I_R, J)$ since $\boldsymbol{e} \in I_R$, so Lemma~\ref{l:pia} implies that $\pi(a) \in \widetilde{I}$.  Since $\pi(j) = 0$, we obtain $\pi(\boldsymbol{e}) \in \widetilde{I}$ as well.  We showed in \S\ref{s:trr} that $\pi(\boldsymbol{e}) \equiv \det(E') - \det(E) \pmod{\widetilde{I}}$.  The result follows.
\end{proof}

{
Thus our remaining task is purely invariant-theoretic: prove that the image of
\[
\boldsymbol{e}=\det(\mathbf E')-\det(\mathbf E)
\]
in $R/J$ lies in the trace--determinant subring $A$.
}

\section{Matrix Invariant Theory and Rational Cohomology} \label{s:mit}

{
As mentioned above, we do not know how to use direct computation to show that the image of $\det(\boldsymbol{E}') - \det(\boldsymbol{E})$ in $R/J$ lies in $\overline{A}$.  Instead, we identify the subring $A \subset R$ as the ring of invariants of a Borel subgroup $\B$ of the algebraic group 
$\G = \GL_2\!/\Z$ acting by simultaneous conjugation on the matrices $\boldsymbol{\rho}_i$.  
The element $\boldsymbol{e} \in R$ is not itself $\B$-invariant, but we show
that its image in $R/J$ is $\B$-invariant.
We must then show that the resulting invariant lifts from 
$R/J$ to an invariant element of $R$, and hence an element of $A$; this will prove the desired result that $\boldsymbol{e}$ lies in $\overline{A}$. }

The obstruction to lifting invariants is a class in {\em rational cohomology}, or alternately {\em Hochschild cohomology}. For the convenience of the reader we review some basic definitions.  As mentioned earlier, the word ``rational" is used in the sense of rational maps, not the rational numbers; we work integrally over $\Z$ throughout.

\subsection{Rational cohomology}
Let $\G$ denote the algebraic group $\GL_2\!/\Z$, i.e.
\[ \G = \Spec \mathfrak{G}, \qquad \mathfrak{G} = \Z[a, b, c, d, (ad-bc)^{-1}] \]
endowed with its usual structure as a group scheme given by a comultiplication
\[ c\colon \mathfrak{G} \longrightarrow \mathfrak{G} \otimes_{\Z} \mathfrak{G} \]
and a counit $e \colon \mathfrak{G} \lra \Z$.
 The lower triangular Borel $\B \subset \G$ will play an important role in this study:
\[ \B = \Spec \mathfrak{B}, \qquad \mathfrak{B} = \Z[x, y, z, (xz)^{-1}], \]
with $\mathfrak{G} \rightarrow \mathfrak{B}$ given by $\mat{a}{b}{c}{d} \mapsto \mat{x}{0}{y}{z}$.

\begin{definition}[\cite{jantzen}*{I.2.7-I.2.8}]  A {\em $\G$-module} is a $\Z$-module $V$ endowed with a \emph{comodule map} $w \colon V \longrightarrow \mathfrak{G} \otimes_{\Z} V$, a map of $\Z$-modules satisfying:
\begin{itemize}
\item $(1 \otimes w) \circ w = (c \otimes 1) \circ w$,
\item $(e \otimes 1) \circ w = id$.
\end{itemize}

Equivalently, a $\G$-module is a $\Z$-module $V$ endowed with a functorial action, for each commutative ring $A$, of $\G(A) = \GL_2(A)$ on $A \otimes_{\Z} V$. Consider the universal element \[ g =\mat{a}{b}{c}{d} \in \G(\mathfrak{G}),\] corresponding to $id$ under the isomorphism $\G(\mathfrak{G}) \cong \Hom(\mathfrak{G}, \mathfrak{G})$. The action of $g$ on $\mathfrak{G} \otimes_{\Z} V$ defines a map $\mathfrak{G} \otimes_{\Z} V  \lra \mathfrak{G} \otimes_{\Z} V$ whose restriction to $1 \otimes V$ is $w$. 
\

  We say that $V$ is a $\G$-module over a commutative ring $R$ if $V$ has a structure of $R$-module that commutes with the $\G$-action.
\end{definition}

The center $\G_m \subset \G$ decomposes any $\G$-module $V$ into a sum of \emph{homogeneous} $\G$-modules: $V = \oplus_{n \in \Z} V^{(n)}$, where \[ V^{(n)} = \{v \in V\colon \mat{x}{0}{0}{x} v = x^n \otimes v\}. \]

\bigskip

\begin{example} \label{e:adjoint}
 We denote by $\cA = \Z A \oplus \Z B \oplus \Z C \oplus \Z D$ the $\G$-module given by the adjoint representation, i.e. for $g = \mat{a}{b}{c}{d}$ the universal element,
\[ \mat{g \cdot A}{g \cdot B}{g \cdot C}{g \cdot D} = \mat{a}{b}{c}{d}^{-1} \mat{A}{B}{C}{D}  \mat{a}{b}{c}{d}. \]

Note that for $g = \mat{x}{0}{y}{z} \in \B$, we have \[ g \cdot A = A + \frac{y}{x} B, \quad g \cdot B = \frac{z}{x} B. \] In particular, 
{
\[ \mathcal{B} := \Z B \quad \text{ and }  \quad \mathcal{V} :=  \Z A \oplus \Z B \]}
 are $\B$-submodules of $\cA$.
For future reference, we also record \[ g \cdot C = -\frac{y}{z}A - \frac{y^2}{xz}B + \frac{x}{z}C + \frac{y}{z} D, \qquad g \cdot D = -\frac{y}{x} B + D. \]

For a $\B$-module $M$, let 
\begin{equation} \label{e:twistdef}
M(n) := M \otimes \mathcal{B}^{\otimes n}. \end{equation} Note that $\bigwedge^2 \mathcal{V} \cong \Z(1)$ and 
$\mathcal{V}^*{= \Hom_{\Z}(\mathcal{V}, \Z)} \cong \mathcal{V}(-1)$. 

\end{example}

\begin{definition}[\cite{jantzen}*{I. 4}]  If $V$ is a $\G$-module, define the invariants
\begin{align*} H^0(\G, V) = V^{\G} &= \{v \in V \colon v \mapsto 1 \otimes v \text{ under the coaction map } V \rightarrow \mathfrak{G} \otimes V \} \\
&= \{ v \in V \colon \text{for all commutative rings $A$ and $g \in \G(A)$, }  \\
& \ \ \ \ \ \ \ \ \ \ \ \ \ \ \ \text{ we have } g \cdot v = v\}.  \end{align*}
The rational cohomology groups $H^i(\G, -)$ are the right derived functors of $H^0(\G, -)$.  \end{definition}

All of the definitions above carry through with $\G$ replaced by $\B$.
The following is proven in \cite{cpsv}*{Theorem 2.1} for semisimple groups over algebraically closed fields in characteristic $p$. However, these hypotheses can be removed using the techniques of \cite{jantzen}.
\begin{theorem} \label{t:bg} Let $V$ be a $\G$-module that is $\Z$-flat.  The restriction map \[ H^i(\G, V) \lra H^i(\B, V) \] is an isomorphism for all $i \ge 0$. 
\end{theorem}

\begin{proof} 
For a field $k$, the isomorphism $H^i(\G_k, V_k) \cong H^i(\B_k, V_k)$ is \cite{jantzen}*{II.4.7}. {
The result over $\Z$ follows from this.
Let  $C_G^\bullet(V)$ and $C_B^\bullet(V)$ be complexes computing the rational cohomology of $G$ and $B$ with coefficients in $V$, 
and let
 $\cC^\bullet = \mathrm{Cone}\big(C_G^\bullet(V)\lra C_B^\bullet(V)\big)$. 
By the universal coefficient theorem \cite{jantzen}*{I.4.18} (which we can apply since $V$ is $\Z$-flat), for any field $k$ we have
\[
\cC^\bullet \otimes_{\mathbf Z}^{\mathbf L} k \simeq 
\mathrm{Cone}\big(C_{G_k}^\bullet(V_k)\lra C_{B_k}^\bullet(V_k)\big),
\]
which is acyclic since the result holds over fields. Since $\Z$ has 
global dimension 1, vanishing after derived tensoring with every field 
forces $\cC^\bullet$ itself to be acyclic.
 (See also the proof of \cite{jantzen}*{II.4.5} for a similar base-change argument).}
\end{proof}

\begin{remark} Theorem~\ref{t:bg} is our motivation for using rational cohomology.  The analogous statement in ordinary group cohomology for the groups $\G(\Z) = \GL_2(\Z)$ and $\B(\Z)$ is in general false.
\end{remark}

\subsection{Acyclic $\B$-modules}

Given a $\B$-module $W$, the algebraic induction $\Ind_{\B}^{\G}(W)$ is defined to be the $\G$-module $(\mathfrak{G} \otimes_{\Z} W)^{\B}$, where the action of $\B$ on $\mathfrak{G}$ comes from right-multiplication on $\G$ (\cite{jantzen} I.3.3).

We identify the weight $(a,b)$ with the character
$\operatorname{diag}(x,z)\mapsto x^a z^b$, and use the usual
notion of dominance $a\geq b$. With our lower triangular
Borel, the induced module of dominant weight $\lambda$ is
$\Ind_\B^\G(\lambda)$.

\begin{definition}
A $\G$-module $V$ admits a \emph{finite good filtration} if
there is a finite filtration by $\G$-submodules
\[
0=V_0\subset V_1\subset\cdots\subset V_m=V
\]
such that each quotient $V_i/V_{i-1}$ is isomorphic to
$\Ind_\B^\G(\lambda_i)$ for a dominant weight $\lambda_i$,
regarded as a rank-one $\B$-module.

More generally, we call a $\G$-module $V$ \emph{good} if it
is a filtered union of $\G$-submodules admitting finite good
filtrations.
\end{definition}
 
\begin{remark}
The modules $\Ind_\B^\G(\lambda)$ are $\Z$-flat.
Extensions and filtered colimits of flat modules are flat, so
every good $\G$-module is $\Z$-flat.

Rational cohomology of $\G$ and $\B$ commutes with filtered
colimits \cite{jantzen}*{I.4.17}. We can use this to pass from
finite good filtrations to general good modules.
In particular, if
\[
V=\bigoplus_{d\geq0}V_d
\]
is a grading by $\G$-submodules and every $V_d$ admits a finite
good filtration, then $V$ is good, by taking the union of the
finite partial sums. This grading need not be the grading by
central characters.
\end{remark}
 
 The following result is proven in \cite{cpsv} in characteristic $p$. The proof easily generalizes to $\G = \GL_2\!/\Z$ using results of \cite{jantzen}.
\begin{theorem}\label{t:good-acyclic}
Every good $\G$-module $V$ is acyclic:
\[
H^i(\G,V)=0\qquad(i>0).
\]
\end{theorem} 
\begin{proof}
Since rational cohomology commutes with filtered colimits, we
may assume that $V$ admits a finite good filtration.
Induction on its length, using the long exact cohomology
sequence, reduces the assertion to
\[
H^i(\G,\Ind_\B^\G(\lambda))=0
\qquad(i>0)
\]
for every dominant weight $\lambda$.
The line bundle on $\G/\B\simeq\mathbf P^1$ associated to
$\lambda$ is $\mathcal O(d)$ for some $d\geq0$.
Consequently,
\[
R^b\Ind_\B^\G(\lambda)=0\qquad(b>0).
\]
The Grothendieck spectral sequence
\cite{jantzen}*{I.4.5}
\[
H^a\bigl(\G,R^b\Ind_\B^\G(\lambda)\bigr)
\;\Rightarrow\;
H^{a+b}(\B,\lambda)
\]
therefore gives
\[
H^i(\G,\Ind_\B^\G(\lambda))
\cong H^i(\B,\lambda).
\]
 Since $\lambda$ is $\Z$-flat, the universal coefficient theorem applied as in the proof of Theorem~\ref{t:bg} 
reduces the vanishing of $H^i(\B,\lambda)$ to the case of fields. Over an
algebraically closed field this is \cite{cpsv}*{2.2}, and the general field case follows
by extension of scalars.
\end{proof}

The following result is well-known (see \cite{cp}*{Theorem 11.1}, {which works over $\Z$ and calls good filtrations ``$\nabla$-filtrations'').}
\begin{theorem}\label{t:tensor-good}
If $V$ and $W$ are good $\G$-modules, then so is $V \otimes_{\Z} W$.
\end{theorem} 
 
  The cited result gives the assertion for finite good filtrations.
For general good modules, write $V$ and $W$ as filtered unions
of submodules with finite good filtrations. Their tensor products
embed in $V\otimes_\Z W$, by $\Z$-flatness, and form a filtered
union equal to $V\otimes_\Z W$. The general assertion follows.

\begin{corollary}\label{c:adj-good}
The adjoint representation $\cA$ is a good $\G$-module, as are its tensor powers $\cA^{\otimes k}$.
\end{corollary}
\begin{proof}
We have
\[
\mathrm{std}\cong\Ind_\B^\G(1,0),\qquad
\mathrm{std}^{\vee}\cong\Ind_\B^\G(0,-1).
\]
Hence $\cA\cong\mathrm{std}\otimes\mathrm{std}^{\vee}$ is good
by Theorem~\ref{t:tensor-good}.
\end{proof}

We also need the following result (\cite{cp} Corollary 11.7):
\begin{theorem}\label{t:r-good}
Every component of fixed total polynomial degree in
\[
\Z[\boldsymbol a_i,\boldsymbol b_i,
   \boldsymbol c_i,\boldsymbol d_i]_{i=1,\ldots,n}
\]
admits a finite good filtration for the action of $\G$ by
simultaneous conjugation. In particular, the entire polynomial
ring is a good $\G$-module.
\end{theorem}

We also use finite good filtrations over a field \(k\),
with \(\Z,\G,\B\) replaced by \(k,\G_k,\B_k\).
For \(a\geq b\), the \emph{Weyl module} associated to
\((a,b)\) is
\[
\Delta_k(a,b)
:=\bigl(\Ind_{\B_k}^{\G_k}(-b,-a)\bigr)^\vee
\cong
\bigl(\Sym^{a-b}(\mathrm{std}_k)\bigr)^\vee\otimes\det\nolimits^a.
\]
A \emph{Weyl filtration} is a finite filtration by
\(\G_k\)-submodules with Weyl modules as successive quotients.

For a finite-dimensional rational \(\G_k\)-module \(M\),
its \emph{good filtration dimension} \(\operatorname{gfd}(M)\)
is the least \(d\geq0\) for which there is an exact sequence
\[
0\longrightarrow M\longrightarrow Q^0\longrightarrow\cdots
\longrightarrow Q^d\longrightarrow0
\]
with every \(Q^s\) admitting a finite good filtration.
Its \emph{Weyl filtration dimension} \(\operatorname{wfd}(M)\)
is defined in the same way using exact sequences
\[
0\longrightarrow P_d\longrightarrow\cdots\longrightarrow P_0
\longrightarrow M\longrightarrow0
\]
with every \(P_s\) admitting a Weyl filtration.
These are the conventions of \cite{p}*{Definition~2.2}.
Dualizing gives
\(\operatorname{gfd}(M^\vee)=\operatorname{wfd}(M)\).
We will use the vanishing
\[
\Ext_{\G_k}^q(M,N)=0
\qquad
\bigl(q>\operatorname{wfd}(M)+\operatorname{gfd}(N)\bigr)
\]
from \cite{p}*{Lemma~2.3}.
In particular, positive-degree Ext groups vanish when the
first argument has a Weyl filtration and the second has a
good filtration.

\begin{lemma} \label{l:vanishing}
 Let $V$ be a good $\G$-module. 
  For any $j \geq 0$, we have \[ H^i(\B, V(j)) = 0 \] for $i > j$.
\end{lemma}
\begin{proof}
 Since twisting and rational cohomology commute with filtered
colimits, we may assume that $V$ admits a finite good filtration.
If $j = 0$, this follows from Theorem \ref{t:good-acyclic} and Theorem \ref{t:bg}. 
{If $j \geq 1$, 
we first show that for $i \ge 1$, \begin{equation} \label{e:groth}
 H^i(\B, V(j)) \cong H^{i-1}(\G, V\otimes H^1(\mathbf{P}^1, \cO(-2j))) \end{equation} (see \cite{jantzen}*{I.4.5 and I.4.8}). 
This isomorphism follows from the Grothendieck spectral sequence
\[
E_2^{a,b}
=
H^a\bigl(\G, R^b\Ind_\B^\G(V(j))\bigr)
\;\Rightarrow\;
H^{a+b}(\B,V(j)).
\]
The  line bundle on $\G/\B \simeq \mathbf{P}^1$ associated to the $\B$-module $\cB$ that defines our twisting is $\cO(-2)$. Thus
\[
R^0\Ind_\B^\G(V(j)) = 0,
\qquad
R^1\Ind_\B^\G(V(j))
\cong
V \otimes H^1(\mathbf{P}^1,\mathcal{O}(-2j)),
\]
and $R^b\Ind_\B^\G(V(j))=0$ for $b>1$. Here we use the tensor identity for derived induction:
\[
R^b\Ind_\B^\G(\operatorname{Res}_\B^\G V\otimes M)
\cong
V\otimes R^b\Ind_\B^\G(M)
\] since $V$ is $\Z$-flat.
Thus the spectral sequence
degenerates to give
\[
H^i(\B,V(j))
\cong
H^{i-1}\bigl(\G, V \otimes H^1(\mathbf{P}^1,\mathcal{O}(-2j))\bigr).
\]
Serre duality on $\mathbf P^1$ yields an isomorphism of $\G$-modules \begin{equation} \label{e:sd} H^1(\mathbf P^1,\mathcal O(-2j)) \cong D_j := \Sym^{2j-2}(\mathrm{std})^{\vee} \otimes {\det}^{j-1}. \end{equation} 

Combining (\ref{e:groth}) and (\ref{e:sd}), to prove the lemma we need to show that \[ H^{i-1}(\G, V \otimes D_j) = 0 \] when $i > j$.} By \cite{jantzen}*{I.4.18}, it suffices to prove that for all primes $p > 0$, we have \[ H^{i-1}(\GL_{2}\!/{\F_p}, V \otimes D_j  \otimes \F_p) = 0\]  when $i > j$. 
By \cite{p}*{Proposition~5.6}, applied with \(n=2\),
and duality, we have the following bound on the good filtration dimension (gfd):
\[
\operatorname{gfd}(D_j\otimes\F_p)
\leq
\left\lfloor\frac{2j-2}{p}\right\rfloor
\leq j-1.
\]
Here we use that tensoring with a power of the determinant
preserves good filtration dimension.
The Friedlander--Parshall vanishing criterion
\cite{fp}*{\S3} therefore gives the required vanishing.
\end{proof}

\begin{lemma}\label{l:higher-cohom}
Let $V$ be a good $\G$-module. We have:
\begin{enumerate}
\item $
H^1(\B,V(1))\cong H^0(\G,V).$
\item For $i>1$ and every prime $p>2$, 
$H^i(\B,V(i))\otimes\F_p=0.
$
\item For $i>1$, put
\[
V_2=V\otimes\F_2,
\qquad
\cA_2=\cA\otimes\F_2,
\]
where $\cA$ denotes the adjoint representation as in Example~\ref{e:adjoint}.
There is a natural isomorphism
\[
H^i(\B,V(i))\otimes\F_2
\cong
\coker\left(
\Hom_{\G_{\F_2}}(\cA_2,V_2)
\xrightarrow{\operatorname{res}}
H^0(\G_{\F_2},V_2)
\right),
\]
where $\operatorname{res}$ is restriction along the inclusion $\F_2\longrightarrow\cA_2$
of the scalar matrices. Moreover,
\[
H^0(\G_{\F_2},V_2)
\cong
H^0(\G,V)\otimes\F_2.
\]
\end{enumerate}

The isomorphism in {\rm (3)} is functorial in $V$ and is compatible
with multiplication by invariant scalars.
\end{lemma}

\begin{proof}
All the functors in the statement commute with filtered colimits
in $V$. For the $\Hom$-term in {\rm (3)}, this follows because
$\cA_2$ is finite-dimensional and taking invariants commutes
with filtered colimits. We may therefore assume that $V$
admits a finite good filtration.

As in Lemma~\ref{l:vanishing}, for every $m\geq1$ and every
$q\geq1$, we have
\begin{equation}\label{e:higher-cohom-D}
H^q(\B,V(m))
\cong
H^{q-1}(\G,V\otimes D_m),
\end{equation}
where $D_m$ is as in \eqref{e:sd}.
Since $D_1=\Z$, taking $m=q=1$ proves {\rm (1)}.

Now let $i>1$ and $p>2$. By the universal coefficient theorem,
\[
H^i(\B,V(i))\otimes\F_p
\hookrightarrow
H^{i-1}\bigl(
\G_{\F_p},
(V\otimes D_i)\otimes\F_p
\bigr).
\]
The module $V\otimes\F_p$ has a good filtration, while
\[
\operatorname{gfd}(D_i\otimes\F_p)
\leq
\left\lfloor\frac{2i-2}{p}\right\rfloor
<i-1.
\]
The Friedlander--Parshall vanishing criterion used in the proof of
Lemma~\ref{l:vanishing} therefore implies that
the group on the right is zero. This proves {\rm (2)}.

It remains to consider $p=2$. Taking $m=i$ and $q=i+1$ in
\eqref{e:higher-cohom-D}, and applying
Lemma~\ref{l:vanishing}, gives
\[
H^i(\G,V\otimes D_i)
\cong
H^{i+1}(\B,V(i))
=0.
\]
The universal coefficient theorem and tensor--Hom adjunction therefore
give
\begin{align}
H^i(\B,V(i))\otimes\F_2
&\cong
H^{i-1}\bigl(
\G_{\F_2},
V_2\otimes(D_i\otimes\F_2)
\bigr) \nonumber\\
&\cong
\Ext_{\G_{\F_2}}^{i-1}
\bigl(D_i^\vee\otimes\F_2,V_2\bigr).
\label{e:higher-cohom-ext}
\end{align}

For \(r\geq1\), the module \(D_{r+1}^\vee\otimes\F_2\)
is the induced module of highest weight \((r,-r)\).
By \cite{p}*{Lemma~5.2}, all the weights \((r,-r)\),
\(r\geq0\), are regular in characteristic \(2\), since
\(2r+1\not\equiv0\pmod2\).
For every \(r\geq2\), the dimension-shifting isomorphism
at the end of the proof of \cite{p}*{Theorem~4.2},
specialized to \(\GL_2\) in characteristic \(2\) using
\cite{p}*{Lemmas~3.5 and~5.1}, gives a natural isomorphism
\[
\Ext_{\G_{\F_2}}^r
\bigl(D_{r+1}^\vee\otimes\F_2,V_2\bigr)
\cong
\Ext_{\G_{\F_2}}^{r-1}
\bigl(D_r^\vee\otimes\F_2,V_2\bigr).
\]
Indeed, the two outer terms in Parker's exact sequence still
vanish with \(V_2\) as the second argument, by
\cite{p}*{Lemma~2.3}, since \(V_2\) has a good filtration.
The construction is compatible with extension of scalars,
so the isomorphism holds over \(\F_2\).
Iterating yields
\begin{equation}\label{e:extf2}
\Ext_{\G_{\F_2}}^{i-1}
\bigl(D_i^\vee\otimes\F_2,V_2\bigr)
\cong
\Ext_{\G_{\F_2}}^1
\bigl(D_2^\vee\otimes\F_2,V_2\bigr).
\end{equation}
For \(i=2\), this is simply the identity map.
The application of \cite{p}*{Lemma~4.1} in this
dimension-shifting argument requires \(r-1\geq1\).
It therefore does not cover \(r=1\), so we compute the
remaining \(\Ext^1\)-group directly.

Using
\[
\mathrm{std}^\vee
\cong
\mathrm{std}\otimes(\det)^{-1},
\]
the standard exact sequence
\[
0\longrightarrow
\bigwedge\nolimits^2\mathrm{std}
\longrightarrow
\mathrm{std}\otimes\mathrm{std}
\longrightarrow
\Sym^2(\mathrm{std})
\longrightarrow0
\]
gives, after reduction modulo $2$ and twisting by $\det^{-1}$,
an exact sequence
\begin{equation}\label{e:adjoint-D2}
0\longrightarrow
\F_2
\longrightarrow
\cA_2
\longrightarrow
D_2^\vee\otimes\F_2
\longrightarrow0.
\end{equation}
The first map identifies $\F_2$ with the scalar matrices in
$\cA_2$.

The trace pairing identifies $\cA_2$ with its linear dual.
Dualizing \eqref{e:adjoint-D2} therefore gives an exact sequence
\[
0\longrightarrow
D_2\otimes\F_2
\longrightarrow
\cA_2
\longrightarrow
\F_2
\longrightarrow0.
\]
Both outer terms are Weyl modules, so $\cA_2$ has a Weyl
filtration. Since $V_2$ has a good filtration,
\[
\Ext_{\G_{\F_2}}^q(\cA_2,V_2)=0
\qquad(q>0).
\]
Applying $\Hom_{\G_{\F_2}}(-,V_2)$ to
\eqref{e:adjoint-D2} therefore gives
\[
\Ext_{\G_{\F_2}}^1
\bigl(D_2^\vee\otimes\F_2,V_2\bigr)
\cong
\coker\left(
\Hom_{\G_{\F_2}}(\cA_2,V_2)
\xrightarrow{\operatorname{res}}
H^0(\G_{\F_2},V_2)
\right).
\]
Together with \eqref{e:higher-cohom-ext} and \eqref{e:extf2}, this proves
the first assertion of {\rm (3)}.

Finally, $H^1(\G,V)=0$ by Theorem~\ref{t:good-acyclic}. The universal
coefficient sequence therefore gives
\[
H^0(\G_{\F_2},V_2)
\cong
H^0(\G,V)\otimes\F_2.
\]
All the isomorphisms above arise from natural adjunctions, connecting
homomorphisms, and long exact sequences. They are therefore functorial
in $V$ and compatible with multiplication by invariant scalars.
\end{proof}

\begin{corollary}\label{c:base-change}
Let $S_1\subset S_2$ be an inclusion of commutative
$\G$-algebras equipped with nonnegative gradings preserved
by $\G$, and suppose that the inclusion respects these
gradings. Assume that every graded piece of $S_1$ and $S_2$
is finitely generated over $\Z$ and admits a finite good
filtration. Then the natural map
\[
H^i(\B,S_1(i))
\otimes_{H^0(\B,S_1)}
H^0(\B,S_2)
\longrightarrow
H^i(\B,S_2(i))
\]
is surjective for every $i\geq0$.
\end{corollary}
\begin{proof}
The case $i=0$ is automatic, and the case $i=1$ follows from
Lemma~\ref{l:higher-cohom}(1).

Suppose that $i>1$. The map in the statement preserves the
grading, and its target is finitely generated over $\Z$
in each degree. The same is therefore true of its cokernel.
It suffices to prove surjectivity after tensoring with
$\F_p$ for every prime $p$. By
Lemma~\ref{l:higher-cohom}(2), both sides vanish for $p>2$.

It remains to consider $p=2$. By Theorem~\ref{t:bg}, we identify
\[
H^0(\B,S_j)=H^0(\G,S_j).
\]
Put
\[
\cA_2=\cA\otimes\F_2,\qquad
S_{j,2}=S_j\otimes\F_2,\qquad
R_j=H^0(\G,S_j),
\]
\[
P_j=\Hom_{\G_{\F_2}}(\cA_2,S_{j,2}),\qquad
Q_j=H^0(\G_{\F_2},S_{j,2}),
\]
and
\[
C_j=H^i(\B,S_j(i))\otimes\F_2
\qquad(j=1,2).
\]
By Lemma~\ref{l:higher-cohom}(3), there are functorial right-exact
sequences
\[
P_j\longrightarrow Q_j\longrightarrow C_j\longrightarrow0.
\]
Tensoring the sequence for $j=1$ over $R_1$ with $R_2$, and
using the compatibility with multiplication by invariant scalars,
gives a commutative diagram whose columns are right exact:
\[
\xymatrix@C=4em@R=1.5em{
P_1\otimes_{R_1}R_2 \ar[r] \ar[d] &
P_2 \ar[d] \\
Q_1\otimes_{R_1}R_2 \ar[r]^-{\sim} \ar[d] &
Q_2 \ar[d] \\
C_1\otimes_{R_1}R_2 \ar[r] \ar[d] &
C_2 \ar[d] \\
0 \ar[r] & 0.
}
\]
There is a canonical identification
\[
C_1\otimes_{R_1}R_2
\cong
\left(
H^i(\B,S_1(i))\otimes_{R_1}R_2
\right)\otimes\F_2,
\]
under which the bottom horizontal arrow is the reduction modulo $2$
of the base-change map in the statement.
Since $S_1$ and $S_2$ are good, Theorem~\ref{t:good-acyclic}
and the universal coefficient sequence give
\[
Q_j
=
H^0(\G_{\F_2},S_{j,2})
\cong
R_j\otimes\F_2.
\]
Thus the middle horizontal arrow in the diagram is indeed an
isomorphism. It follows from the right exactness of the columns that
\[
C_1\otimes_{R_1}R_2\longrightarrow C_2
\]
is surjective, which is the desired assertion after tensoring with
$\F_2$.\end{proof}

\begin{lemma}\label{l:nice-res}
If a $\B$-module $V$ admits a finite resolution \[ \begin{tikzcd}
 V_n(n) \ar[r,"f_n"] & V_{n-1}(n-1)  \ar[r,"f_{n-1}"] & \cdots \ar[r, "f_2"] & V_1(1) \ar[r,"f_1"] & V_0 \ar[r,"f_0"] & V
 \end{tikzcd} \] where each $V_j$ is a good $\G$-module {that is $\Z$-flat,} then $V$ is an acyclic $\B$-module.
\end{lemma}
\begin{proof}
Put $W_0=V$, $W_{n+1}=0$, and, for $1\leq j\leq n$, put
\[
W_j=\ker(f_{j-1})=\operatorname{im}(f_j).
\]
For $j=0,\ldots,n$, we then have a short exact sequence
\[
0\longrightarrow W_{j+1}
\longrightarrow V_j(j)
\longrightarrow W_j
\longrightarrow0.
\]
Taking the associated long exact sequences in $\B$-cohomology, a downward induction using Lemma~\ref{l:vanishing} shows that $H^i(\B, W_j)$ vanishes for $i > j$. 
For $j=0$, this is the desired result.
\end{proof}

\begin{theorem}\label{t:r-acyclic}
Consider the $\B$-module \[ S = \Z[\boldsymbol{a}_i, \boldsymbol{b}_i, \boldsymbol{c}_i, \boldsymbol{d}_i]_{i=1,\ldots,n}/(\boldsymbol{b}_i \colon 1 \leq i \leq k), \] where $\B$ acts by simultaneous conjugation. Let $V$ be a good $\G$-module such that $V$ is $\Z$-flat.  Then $V \otimes_{\Z} S$ is an acyclic $\B$-module.
\end{theorem}
\begin{proof}
Let $S_0 = \Z[\boldsymbol{a}_i, \boldsymbol{b}_i, \boldsymbol{c}_i, \boldsymbol{d}_i]_{i=1,\ldots,n}$. The elements $\boldsymbol{b}_1,\ldots, \boldsymbol{b}_k$ form a regular sequence in $S_0$, hence the Koszul complex for this sequence yields an exact sequence of $\B$-modules:
\[ \left(\bigwedge^k_{S_0} \bigoplus_{i=1}^k S_0\right)(k) \lra \cdots \lra \left(\bigoplus_{i=1}^k S_0\right)(1) \lra S_0. \]
This gives a resolution of $S = S_0/(\boldsymbol{b}_1, \dotsc, \boldsymbol{b}_k)$. If we tensor this complex with $V$, we obtain an exact complex of $\B$-modules resolving $V \otimes_\Z S$. The $k$th term of this complex is of the form $\left(\bigoplus V \otimes S_0 \right)(k)$. As $V \otimes S_0$ is a good $\G$-module by Theorem \ref{t:r-good} and Theorem \ref{t:tensor-good}, we may apply Lemma \ref{l:nice-res} to conclude that $V \otimes S$ is an acyclic $\B$-module.
\end{proof}

In particular, this theorem may be applied to the ring $R$ defined in \S\ref{s:formal}, with the action of $\B$ defined in \S\ref{s:invt-theory}.

\subsection{Matrix Invariant Theory}\label{s:invt-theory}

The following important classical result is the subject of the beautiful book \cite{cp} by de Concini and Procesi.  {Theorem~1.10 of {\em loc.\ cit.}\ is the statement above with $R_0 = \Z$ and invariants taken for the group $\G(\Z) = \GL_2(\Z)$.}

{
\begin{theorem}\label{t:invariants}Let $\Lambda_0$ be a $\Z$-flat commutative ring and let \[ \Lambda = \Lambda_0[a_i, b_i, c_i, d_i]_{i=1}^{\lambda}. \] Endow $\Lambda$ with a $\G$-action defined by simultaneous conjugation on the matrices $\rho_i = \mat{a_i}{b_i}{c_i}{d_i}$. The ring $H^0(\G, \Lambda) \subset \Lambda$ is generated as an algebra over $\Lambda_0$ by the traces and determinants of all matrices in the algebra generated by the $\rho_i$.
\end{theorem}
}

 Since $\G$-invariance is stronger than $\G(\Z)$-invariance, the statement above follows for $\Lambda_0 = \Z$ from 
 \cite{cp}*{Theorem~1.10}.  Extension to arbitrary $\Z$-flat algebras $\Lambda_0$ on which $\G$ acts trivially is then immediate.

We now apply Theorem~\ref{t:invariants}  to our setting.
Recall the rings  $R_0, R_1$, and $R$ defined in (\ref{e:r0def})--(\ref{e:rdef}) above and the subring $A \subset R$ defined at the start of \S\ref{s:subring}. {
In particular we have \[ R_1 = R_0[\boldsymbol{a}_i, \boldsymbol{b}_i, \boldsymbol{c}_i, \boldsymbol{d}_i]_{i=1}^r \quad \text{and } \quad
 R = R_1/(\boldsymbol{b}_{\tau} \colon \tau \in \mathcal{B}_{v_0}), \]
where the $\boldsymbol{b}_{\tau}$ are certain generators $\boldsymbol{b}_i$ of $R_1$, indexed by a set $\mathcal{B}_{v_0}$ corresponding to the rows of $D$ of type III.}

Let $\B$ act on $R_1 = R_0[\boldsymbol{a}_i, \boldsymbol{b}_i, \boldsymbol{c}_i, \boldsymbol{d}_i]_{i=1}^r $ by simultaneous conjugation on the matrices $\boldsymbol{\rho}_i$. In particular, $\B$ acts trivially on $R_0.$ This action preserves the ideal $(\boldsymbol{b}_{\tau} \colon \tau \in \mathcal{B}_{v_0})$ (see Example \ref{e:adjoint}), and hence descends to an action on $R = R_1/(\boldsymbol{b}_{\tau} \colon \tau \in \mathcal{B}_{v_0}).$

\begin{corollary}\label{c:a} 
We have $H^0(\B, R) = A$.
\end{corollary}

\begin{proof}
It is convenient to write $R$  in the form
\begin{equation} \label{e:taunotation}
 R = R_0[\boldsymbol{a}_i, \boldsymbol{b}_i, \boldsymbol{c}_i, \boldsymbol{d}_i]_{i = 1}^t[\boldsymbol{a}_\tau,\boldsymbol{b}_\tau, \boldsymbol{c}_\tau, \boldsymbol{d}_\tau]/(\boldsymbol{b}_{\tau}) \end{equation}
{
where we exclude the variables corresponding to the entries of $\boldsymbol{\rho}_{\tau}$ with $\tau \in \mathcal{B}_{v_0}$ from the list of variables indexed by $i$, and $\tau$ ranges over $\mathcal{B}_{v_0}$.}
We first handle the case $t=0$, i.e.\ 
\[ R = R_0[\boldsymbol{a}_\tau,\boldsymbol{b}_\tau, \boldsymbol{c}_\tau, \boldsymbol{d}_\tau]/(\boldsymbol{b}_{\tau}). \]
We claim that in this case $H^0(\B, R)$ is the image of $R_0[\boldsymbol{a}_\tau,\boldsymbol{d}_\tau]$ in $R$.  Let $\D \subset \B$ denote the torus of diagonal matrices.  There is a $\D$-module isomorphism
\[ R \cong R_0[\boldsymbol{a}_\tau, \boldsymbol{c}_\tau, \boldsymbol{d}_\tau], \qquad  \boldsymbol{b}_\tau \mapsto 0. \]
Since $\D$ acts trivially on $\boldsymbol{a}_\tau$ and $\boldsymbol{d}_\tau$ but scales $\boldsymbol{c}_\tau$ non-trivially, it is clear that 
$H^0(\D, R)$ is the image of  $R_0[\boldsymbol{a}_\tau,\boldsymbol{d}_\tau]$.  Since the elements $\boldsymbol{a}_\tau$ and $\boldsymbol{d}_\tau$ are $\B$-invariant in $R$, the claim follows.

Now we return to the case of general $t$.
Let \[ S = R_0[\boldsymbol{a}_i, \boldsymbol{b}_i, \boldsymbol{c}_i, \boldsymbol{d}_i],  \text{ so that } R_1 = S[\boldsymbol{a}_\tau,\boldsymbol{b}_\tau, \boldsymbol{c}_\tau, \boldsymbol{d}_\tau]. \] Let $\fb = (\boldsymbol{b}_{\tau} \colon \tau \in \mathcal{B}_{v_0})$.  By Theorem \ref{t:invariants},  the image of \[ H^0(\G, R_1) = H^0(\B,R_1) \]  in $H^0(\B, R_1/\fb)$ is the ring $A_0$ defined in \S\ref{s:subring}.  

Let \[ R'_1 =R_0[\boldsymbol{a}_\tau,\boldsymbol{b}_\tau, \boldsymbol{c}_\tau, \boldsymbol{d}_\tau] \subset R_1, \]
 and let \[ \fb' = R'_1 \cap \fb.\] 
 By the claim above for $t=0$, we have \[ H^0(\B, R'_1/\fb') = R_0[\boldsymbol{a}_\tau,\boldsymbol{d}_\tau]. \]
 The ring $A$ is by definition generated over $A_0$ by the $\boldsymbol{d}_\tau$'s and hence also by the larger generating set of $\boldsymbol{a}_\tau$'s and $\boldsymbol{d}_\tau$'s (note $\boldsymbol{a}_\tau + \boldsymbol{d}_\tau \in A_0$). Thus, to prove that $H^0(\B, R_1/\fb) = A$,  it will suffice to show that the map \begin{equation}\label{e:base-change} H^0(\B, R'_1/\fb') \otimes_{H^0(\B, R'_1)} H^0(\B, R_1) \lra H^0(\B, R_1/\fb)\end{equation} is surjective.

 Equip $R'_1$ and $R_1$ with total polynomial degree, assigning
degree $1$ to every polynomial generator, including those
of $R_0$. The $\G$-actions preserve these gradings, and
$R'_1\subset R_1$ is a graded inclusion. Each graded piece is
finitely generated over $\Z$ and admits a finite good
filtration, by Theorem~\ref{t:r-good} and the trivial action on
$R_0$. Thus this inclusion satisfies the hypotheses of
Corollary~\ref{c:base-change}.
 
 Since  $\{\boldsymbol{b}_\tau \}$ is a regular sequence in $R_1$, we may provide a $\B$-module resolution of $R_1/\fb$ via a Koszul complex:
\[ \left(\bigwedge\nolimits^k_{R_1}  \left(\bigoplus_{i=1}^k R_1\right)\right)(k) \lra \cdots \lra  \bigoplus_{i=1}^k R_1(1) \lra R_1 \lra R_1/\fb \lra 0.  \]

Here $k$ denotes the number of $\tau$'s, and the map $\bigoplus_{i=1}^k R_1(1) \lra R_1$ sends each standard basis vector on the left to a $\boldsymbol{b}_\tau$.
Let \[ C_j =  \left(\bigwedge\nolimits^j_{R_1}  \left(\bigoplus_{i=1}^k R_1\right)\right)(j) \]
denote the terms of this complex.
 There is an analogous Koszul complex for $R'_1/\fb'$, which we denote $C'_j$. We have $C_{\bullet} \cong C'_{\bullet} \otimes_{R'_1} R_1$.

{Viewing $C_\bullet$ as a cochain complex via $C^{-i} = C_i$, the
hypercohomology spectral sequence (cf. \cite{weibel}*{\S5.7}) gives} a
  2nd-quadrant spectral sequence \[ E_1^{-i,j} = H^j(\B,  C_i) \implies H^{j-i}(\B, R_1/\fb). \]  Since $R_1$ is a good $\G$-module (Theorem \ref{t:r-good}),  Lemma \ref{l:vanishing} yields \[ H^j(\B, C_i) = 0 \quad \text{ for }j > i. \]
This vanishing gives a finite filtration
\[
0=F_{-1}\subset F_0\subset\cdots\subset F_k
=H^0(\B,R_1/\fb)
\]
with $F_i/F_{i-1}=E_\infty^{-i,i}$.
Every differential leaving the diagonal has target above the
diagonal, where the spectral sequence vanishes. Consequently
$E_\infty^{-i,i}$ is a natural quotient of
$E_1^{-i,i}=H^i(\B,C_i)$.
The corresponding spectral sequence for $C'_\bullet$ gives
such a filtration $F'_i$ on $H^0(\B,R'_1/\fb')$.

The map of Koszul complexes and multiplication by
$H^0(\B,R_1)$ respect these filtrations and their spectral
sequences. By Corollary~\ref{c:base-change},
\[
H^i(\B,C'_i)\otimes_{H^0(\B,R'_1)}H^0(\B,R_1)
\longrightarrow H^i(\B,C_i)
\]
is surjective for each $i$. Passing to the natural diagonal
quotients therefore gives surjections
\[
(F'_i/F'_{i-1})\otimes_{H^0(\B,R'_1)}H^0(\B,R_1)
\longrightarrow F_i/F_{i-1}.
\]
Right exactness of tensor product and induction on this finite
filtration now imply that \eqref{e:base-change} is surjective.
\end{proof}

\subsection{Roadmap}

We now return to the proof of our main result, Theorem~\ref{t:rl}, picking up the argument of \S\ref{s:formal}.
Recall that we had defined an element
\[ \boldsymbol{e} = \det(\boldsymbol{E}') - \det(\boldsymbol{E}). \] {By Proposition~\ref{p:goalepe},} our goal is to show that 
 its image in $R/J$, denoted $\overline{\boldsymbol{e}}$, lies in  $\overline{A} \subset R/J$. {Here $J \subset R$ is the ideal defined in \S\ref{s:ri} and } $\overline{A}$ is the image of the subring $A \subset R$ in $R/J$.
 Our strategy is as follows.  We will first prove that the ideal $J$ is stable under the action of $\B$.
We will then show that 
\[ \overline{\boldsymbol{e}} \in H^0(\B, R/J). \]
The long exact sequence in rational cohomology associated to 
\[ 0 \longrightarrow J \longrightarrow R \longrightarrow R/J \longrightarrow 0 \]
 yields an exact sequence
\begin{equation} \label{e:les2}
 \begin{tikzcd}
H^0(\B, R) = A \ar[r] & H^0(\B, R/J) \ar[r, "c_J"] & H^1(\B, J).  \end{tikzcd}
\end{equation} 
The equality on the left is Corollary~\ref{c:a}.  Let $\beta \in H^1(\B, J)$ denote the image of $\overline{\boldsymbol{e}}$ under the connecting homomorphism $c_J$.
In view of (\ref{e:les2}), the desired result $\overline{\boldsymbol{e}} \in \overline{A}$ will follow if we can show that $\beta =0$.

For this, we will define a certain $\B$-stable subideal $J' \subset J$ and show that in fact 
\[ \overline{\boldsymbol{e}}\in H^0(\B, R/J'). \]
(This is a slight abuse of notation, as here $\overline{\boldsymbol{e}}$ denotes the reduction of $\boldsymbol{e}$ modulo $J'$.)

We let $\alpha = c_{J'}(\overline{\boldsymbol{e}}) \in H^1(\B, J')$ defined as above.  Then $\beta = \iota_*(\alpha)$ where $\iota \colon J' \rightarrow J$ is the inclusion and $\iota_*$ is the induced map on rational cohomology.  To conclude, we will prove that the map
\[ \iota_*\colon H^1(\B, J') \longrightarrow H^1(\B, J) \]
vanishes.  Therefore $\beta = 0$, and our result follows.  

\subsection{Invariance} \label{s:invariance}
{
The ideal $J$ was defined as the ideal generated by all 4 coefficients of the matrices in (\ref{e:b1})--(\ref{e:b4}).  We now define
 $J' \subset J$ to be the subideal generated by just the ``$b$'' coefficients of the matrices in (\ref{e:b1})--(\ref{e:b4}).}  To be precise, $J'$ is generated by:
\begin{itemize}
\item The elements $
\sum_{i=1}^r \boldsymbol{\epsilon}_i \boldsymbol{b}_i 
$
 for each row of type I in $D$.
\item The elements 
$
(\boldsymbol{a}_i + \boldsymbol{\nu_i})\boldsymbol{b}_j 
+ \boldsymbol{b}_i \boldsymbol{d}_j - \sum_{k = 1}^r \boldsymbol{\delta}_{ijk} \boldsymbol{b}_k
$ for each row of type II in $D$.
\item For each $v \in (\Sigma \setminus \{v_0\}) \sqcup \cP$ and distinct $\sigma, \tau \in \cB_v$, the elements
\[ B(\sigma, \tau) = \boldsymbol{b}_\sigma(\boldsymbol{x}_{\tau} - \boldsymbol{a}_\tau) - 
\boldsymbol{b}_\tau(\boldsymbol{x}_{\sigma} - \boldsymbol{a}_\sigma). \]
\end{itemize}

\begin{lemma} \label{l:stable}
The ideals $J'$ and $J$ are stable under the action of $\B$.  Moreover, {for the ideal $J$,}  the $\Z$-module spanned by each {set of 4 matrix coefficients} in (\ref{e:b1})--(\ref{e:b4}) is isomorphic as a $\B$-module to a copy of the adjoint representation $\cA$ (see Example~\ref{e:adjoint}).
\end{lemma}

\begin{proof}
This is clear for the relations in (\ref{e:b1}) and (\ref{e:b2}) since these relations are linear combinations of products of the $\boldsymbol{\rho}_i$, with coefficients in $R_0$ (on which $\B$ acts trivially), and the definition of our action is by simultaneous conjugation on  $\boldsymbol{\rho}_i$. We must verify this for (\ref{e:b4}) by direct computation.  The Borel $\B$ is generated by its torus and unipotent subgroups:
\begin{equation} \label{e:tu}
 \left\{ \sigma_{x,y} = \mat{x}{0}{0}{y} \right\}, \qquad \left\{ \tau_x = \mat{1}{0}{x}{1} \right\}. \end{equation}

Acting on the adjoint $\cA$, the element $\sigma_{x,y}$ fixes $A$ and $D$, scales $B$ by $y/x$ and scales $C$ by $x/y$.  From the definitions in (\ref{e:ABCD}), it follows that $\sigma_{x,y}$ acts the same way on $A(\sigma, \tau), \dots, D(\sigma, \tau)$.  Similarly, $\tau_x$ has the following action:
\[ \mat{A}{B}{C}{D} \mapsto \mat{A + Bx}{B}{C + (D-A)x - Bx^2 }{D-Bx},
\]
and one checks that the action on the matrix (\ref{e:b4}) is the same.  For instance,
\[
\tau_x B(\sigma, \tau)  = \boldsymbol{b}_\sigma(\boldsymbol{x}_{\tau} - \boldsymbol{a}_\tau- \boldsymbol{b}_\tau x) - 
\boldsymbol{b}_\tau(\boldsymbol{x}_{\sigma} - \boldsymbol{a}_\sigma - \boldsymbol{b}_\sigma x)
 = B(\sigma, \tau).
\]
The verification of the other 3 coefficients is similar and left to the reader.
The stability of $J'$ under $\B$ follows since the $b$-coefficient of the adjoint is stable under $\B$.
\end{proof}

The goal of the rest of this subsection is to prove the following.

\begin{lemma} We have $\overline{\boldsymbol{e}} \in H^0(\B, R/J')$. \label{l:ebar}
\end{lemma}

{
The matrices $\boldsymbol{E}$ and $\boldsymbol{E'}$ were defined in \S\ref{s:edef}, and $\boldsymbol{e} = \det(\boldsymbol{E}') - \det(\boldsymbol{E}).$}
 The matrix $\boldsymbol{E}$ has coefficients in $R_0$, on which $\B$ acts trivially, so we must show that 
\[ \overline{\det(\boldsymbol{E}')} \in H^0(\B, R/J'). \]
Again we use the fact that $\B$ is generated by the $\sigma_{x,y}$ and $\tau_x$ defined in (\ref{e:tu}).  The matrix $\boldsymbol{E}'$ has coefficients in $R_0[\boldsymbol{a}_i, \boldsymbol{d}_i]$ on which  $\sigma_{x,y}$ acts trivially.  Therefore we must only consider the action of $\tau_x$.

{
From the definition of  $\boldsymbol{E'}$ given in \S\ref{s:edef},} 
one sees that the first $t$ rows of $\boldsymbol{E}'$, the rows of type IV, and the rows of type V all have the same shape.  To each such row is attached a $v \in (\Sigma \setminus\{v_0\}) \sqcup \cP$ and a $\sigma \in \cG_v$.  
In the column corresponding to $v$ we have
$\boldsymbol x_\sigma-\boldsymbol a_\sigma$, in the
matrix-generator column corresponding to the label $\sigma$
we have $1$, and all other entries are zero. Therefore to streamline the exposition we will treat these rows the same and call them ``of type L.''  So $\boldsymbol{E}'$ has rows of type I, II, III, or L.

Conjugation by $\tau_x$ sends $\boldsymbol{a}_i \mapsto \boldsymbol{a}_i + \boldsymbol{b}_i x$ and $\boldsymbol{d}_i \mapsto \boldsymbol{d}_i - \boldsymbol{b}_i x$.  This action fixes rows of type I or III, but is non-trivial on rows of type II or L.  We will handle the rows of type II first, and then the rows of type L one $v$ at a time.  
Write $(\Sigma \setminus\{v_0\}) \sqcup \cP= \{v_1, \dots, v_{s+t} \}$.
For $i = 0, \dots, s+t$, let 
$\boldsymbol{E}'_i$  be the matrix obtained from $\boldsymbol{E}'$, 
but with the rows of type II and rows of type L associated to $v_j$ for $j \le i$ replaced by their image under $\tau_x$.
In particular $\boldsymbol{E}'_{s+t} = \tau_x(\boldsymbol{E}').$

\begin{lemma} We have $\det(\boldsymbol{E}'_0) \equiv \det(\boldsymbol{E}') \pmod{J'}$. \label{l:e0}
\end{lemma}

\begin{proof}
Suppose we have a row $w$ of type II associated to a pair $(i,j)$.  Then 
\[ \tau_x(w) - w = (0, 0, \dots, \boldsymbol{b}_j x, 0, \dots, -\boldsymbol{b}_i x, 0, \dots, 0), \]
where there is $\boldsymbol{b}_j x$ in the $\boldsymbol{\rho}_i$ slot and $ -\boldsymbol{b}_i x$ in the 
$\boldsymbol{\rho}_j$ slot.
The difference $\det(\boldsymbol{E}'_0) - \det(\boldsymbol{E}')$ is a linear combination (with coefficients $\pm 1$) 
of determinants of all matrices $X$ obtained by starting with $\boldsymbol{E}'$ and replacing some nonempty subset of the rows $w$ of type II with $\tau_x(w) - w$.

Suppose we are given such a matrix $X$ and row $w$.  We want to show that $\det(X) \equiv 0 \pmod{J'}.$  
Note that if $i = j$ then $\tau_x(w) - w = 0$ and hence $\det(X)=0$.  So we assume $i \neq j$.
We make the following alterations to $X$ which do not change the determinant:
\begin{itemize}
\item We replace $\tau_x(w) - w$ by 
\[ (0, 0, \dots,  x, 0, \dots, -x, 0, \dots, 0). \]
At the same time, in every row other than $w$ we multiply the $\boldsymbol{\rho}_j$-coordinate by $\boldsymbol{b}_j$ and the $\boldsymbol{\rho}_i$-coordinate by $\boldsymbol{b}_i$. 
\item We then add the new  $\boldsymbol{\rho}_j$ column to the new $\boldsymbol{\rho}_i$-column.
\item For each $1 \le k \le r, k \neq i, j$, we add $\boldsymbol{b}_k$ times the $\boldsymbol{\rho}_k$-column to the 
new $\boldsymbol{\rho}_i$-column.
\end{itemize}
Let us make some observations about the matrix $X$ that results from these changes:
\begin{itemize}
\item In rows of type I, II, or III, the
$\boldsymbol{\rho}_i$-coordinate lies in $J'$.
More precisely, it equals the associated generator of $J'$
for a row of type I, and the negative of the associated
generator for a row of type II that was not replaced by a
difference row. It is zero for rows of type III and for all
the type II difference rows, including the distinguished
row $w$.
\item In rows of type L, we have $\boldsymbol{b}_\sigma$ in the $\boldsymbol{\rho}_i$-coordinate and 
$\boldsymbol{x}_\sigma - \boldsymbol{a}_\sigma$ in the $v$-coordinate.  The $2 \times 2$ determinant of any two such rows (associated to the same $v$) is one of the generators of $J'$.
\end{itemize}
As we now explain, the fact that $\det(X) \in J'$ follows from these two facts.  We compute $\det(X)$  as the sum
\begin{equation} \label{e:detm}
\sum_{\pi} \sgn(\pi) \prod_{w} X_{w, \pi(w)} \end{equation}
as $\pi$ ranges over all bijections between the set of rows and columns of $X$, and $w$ runs over the set of rows of $X$.
For a given $\pi$, consider the $w$ such that $\pi(w)$ is the  $\boldsymbol{\rho}_i$-column.  If $w$ is of type I, II, or III, then the contribution to $\det(X)$ lies in $J'$ by the first bulleted observation.  Now suppose  $w$ is of type L, associated to some 
$v \in (\Sigma \setminus\{v_0\}) \sqcup \cP$ and $\sigma \in \cG_v$.
The only rows of $X$ with nonzero component in the $v$-column are those of type L associated to $v$.  So we only obtain a nonzero contribution to $\det(X)$ if there is such a row $w'$, say attached to $\tau \in \cG_v$, and 
 $\pi(w') = $ the $v$-column.  In this case we consider the permutation $\pi'$ obtained from $\pi$ by swapping the roles of $w$ and $w'$.  So $\pi'(w) = v$-column, $\pi'(w') = \boldsymbol{\rho}_i$-column, and $\pi'(u) = \pi(u)$ for $u \neq w, w'$.  Then {$\sgn(\pi') = - \sgn(\pi)$} and by the second bulleted observation we find
 \[ \prod_{u} X_{u, \pi(u)} - \prod_{u} X_{u, \pi'(u)} = (\boldsymbol{b}_\sigma(\boldsymbol{x}_\tau - \boldsymbol{a}_\tau) - 
 \boldsymbol{b}_\tau(\boldsymbol{x}_\sigma - \boldsymbol{a}_\sigma))
  \prod_{u \neq w, w'} X_{u, \pi'(u)} \in J'. \]
 
Pairing off the permutations $(\pi, \pi')$ in this way shows that $\det(X) \in J'$. This concludes the proof.
\end{proof}

By a similar argument, we can show:

\begin{lemma}  For $i = 1, \dots, s+t$, we have $\det(\boldsymbol{E}'_i) \equiv \det(\boldsymbol{E}'_{i-1}) \pmod{J'}$.  \label{l:ei}
\end{lemma}

\begin{proof} Note that $\boldsymbol{E}'_i$ and $\boldsymbol{E}'_{i-1}$ only differ in the rows of type L associated to $v_i$, and only in the column associated to $v_i$.  In these entries, $\boldsymbol{E}'_{i-1}$ has the value $\boldsymbol{x}_\tau - \boldsymbol{a}_\tau$ and $\boldsymbol{E}'_i$
has the value $\boldsymbol{x}_\tau - \boldsymbol{a}_\tau - \boldsymbol{b}_\tau x$.  

We therefore modify $\boldsymbol{E}'_i$ by adding, for each $j=1, \dotsc, r$, the $\boldsymbol{\rho}_j$-column scaled by $\boldsymbol{b}_j x$ to the $v_i$-column.  The resulting matrix $X$ satisfies $\det(X) = \det(\boldsymbol{E}'_i)$.  Furthermore
the entries of $X$ and $\boldsymbol{E}'_{i-1}$ are equal except in the column associated to $v_i$ for the rows {\em not} associated to $v_i$.  In the $v_i$-column we have:
\begin{itemize}
\item For rows of type I, $x$ times the associated generator
of $J'$; for rows of type II, $-x$ times the associated
generator of $J'$; and for rows of type III, $0$.
\item For rows of type L associated to $v_j \neq v_i$, the element $\boldsymbol{b}_\sigma x.$
\end{itemize}
Again we compute $\det(X)$ and $\det(\boldsymbol{E}'_{i-1})$ via (\ref{e:detm}), this time focusing on the row $w$ such that $\pi(w) = v_i$-column.  

If $w$ is of type I, II, or III, then the first bulleted point shows that the contribution to $\det(X)$ is 0 modulo $J'$. The same is true for 
$\det(\boldsymbol{E}'_{i-1})$ since the $(w, \pi(w))$-entry of $\boldsymbol{E}'_{i-1}$ is 0.

  If $w$ is of type L associated to $v_i$, then the $(u, \pi(u))$-entries of $X$ and $\boldsymbol{E}'_{i-1}$ are equal for all rows $u$, so the contributions to the determinants are the same. 
  
   Finally, if $w$ is of type L associated to some other place $v_j \neq v_i$, then the $(w, \pi(w))$-entry of $\boldsymbol{E}'_{i-1}$ is 0, so we must show that the contribution to $\det(X)$ is also $0$ modulo $J'$.  By the second bullet point above, the $(w, v_i)$-entry of $X$ is 
   $\boldsymbol{b}_\sigma x$ and the $(w, v_j)$-entry is either $\boldsymbol{x}_\sigma - \boldsymbol{a}_\sigma$ or 
   $\boldsymbol{x}_\sigma - \boldsymbol{a}_\sigma - \boldsymbol{b}_\sigma x$, 
   depending on whether $j > i$ or $j \le i-1$.   Now, the only other rows of $X$ that have nonzero $v_j$-coordinate are those of type L associated to $v_j$. So the contribution of $\pi$ to $\det(X)$ is 0 unless $\pi(w') = v_j$-column for some other $w'$ of type L associated to $v_j$, say attached to the element $\tau \in \cG_{v_j}$.  The entries of $X$ in the $w'$ row are the same as those described for the $w$ row, with $\sigma$ replaced by $\tau$.  In both the cases $j > i$ or $j \le i-1$, the $2 \times 2$ determinant given by the $(w, w')$ rows and $(v_i, v_j)$ columns is equal to 
\begin{equation}\label{e:baj}
x\Bigl(
\boldsymbol{b}_\sigma
  (\boldsymbol{x}_\tau-\boldsymbol{a}_\tau)
-
\boldsymbol{b}_\tau
  (\boldsymbol{x}_\sigma-\boldsymbol{a}_\sigma)
\Bigr)\in J'.
\end{equation}
 Once again we define $\pi'$ to be the permutation  $\pi$ with the roles of $w, w'$ swapped.
As before, (\ref{e:baj}) implies that the sum of the contributions of $\pi$ and $\pi'$ to $\det(X)$ is 0 modulo $J'$.  Of course, the contribution of $\pi'$ to $\det(\boldsymbol{E}'_{i-1})$ is also 0.
Pairing off $(\pi, \pi')$ in this way concludes the proof.
\end{proof}

Lemmas~\ref{l:e0} and~\ref{l:ei} imply the desired result
\[ \det(\tau_x(\boldsymbol{E}')) = \det(\boldsymbol{E}'_{s+t}) \equiv \det(\boldsymbol{E}') \pmod{J'}. \]
This concludes the proof of Lemma~\ref{l:ebar}.

\subsection{A cascade of cohomology classes}
 
 In this section, we will assume the existence of a certain morphism of complexes (see Theorem~\ref{t:comm} below) and use it to prove the following.
 
 \begin{theorem} \label{t:push}
 Let $\iota \colon J' \rightarrow J$ be the inclusion, let $j \ge 1$, and let  
\[ \iota_*\colon H^j(\B, J') \longrightarrow H^j(\B, J) \]
be the induced map on rational cohomology groups.  Then $\iota_* = 0$.
 \end{theorem}

 \begin{theorem}\label{t:comm}
  There exists a morphism of complexes of $\B$-modules:
\begin{equation} \label{e:Koszul2} \begin{tikzcd}
0 \ar[r]  & \cC^r  \ar[r,"f_r"] \ar[d,"\iota_r"]
 & \cC^{r-1}   \ar[r, "f_{r-1}"] \ar[d,"\iota_{r-1}"]
 &  \cdots \ar[r,"f_2"] & \cC^1 \ar[r,"f_1"]  \ar[d,"\iota_1"] &  \cC^0 = J' \ar[d, "\iota_0",shift right=5]  \ar[r]
 & 0
 \\
 0 \ar[r] & \cD^r \ar[r,"g_r"] 
 & \cD^{r-1}   \ar[r, "g_{r-1}"]
 &  \cdots \ar[r,"g_2"] &  \cD^1 \ar[r,"g_1"] &  \cD^0 = J  \ar[r]
 & 0
 \end{tikzcd}
 \end{equation}
such that $\iota_0:J'\longrightarrow J$ is the inclusion,
the complex $\cC^\bullet$ is exact, and every module
$\cD^i$, for $i>0$, is an acyclic $\B$-module.
 \end{theorem}
 
 We assume Theorem \ref{t:comm} for now and use it to prove Theorem~\ref{t:push}.
 
 \begin{proof}[Proof of Theorem~\ref{t:push}]  
 We extend both complexes in \eqref{e:Koszul2} by zero outside
degrees $0,\ldots,r$. In particular, we use the zero maps
\[
f_0:\cC^0\longrightarrow0,\qquad
f_{r+1}:0\longrightarrow\cC^r,\qquad
g_{r+1}:0\longrightarrow\cD^r.
\]
The image of the map $f_1$ is $J'$.   Let \[ \alpha_j \in H^j(\B, J') = H^j(\B, \im(f_1)). \]
Our goal is to show $\iota_{0, *} \alpha_j = 0$.  Applying the coboundary in 
the long exact sequence in cohomology associated to 
\[ 0 \lra \ker(f_1) \lra \cC^1 \lra \im(f_1) \lra 0 \] to $\alpha_j$ 
 yields a class $\alpha_{j+1} \in H^{j+1}(\B, \ker(f_1))$. The class $\alpha_{j+1}$ represents the obstruction to lifting $\alpha_j$ to a class in $H^j(\B, \cC^1)$.  Writing $\ker(f_1) = \im(f_2)$, we can view \[ \alpha_{j+1} \in H^{j+1}(\B, \im(f_2)). \] Repeat the process above, using the coboundary in the long exact sequence associated to 
\[ 0 \lra \ker(f_2) \lra \cC^2 \lra \im(f_2) \lra 0 \]
to obtain $\alpha_{j+2} \in H^{j+2}(\B, \ker(f_2))$.  Continuing in this way we obtain 
\[ \alpha_{j+i} \in H^{j+i}(\B, \ker(f_{i})) = H^{j+i}(\B, \im(f_{i+1})) \] for $i = 0, \dots, r$.  Note $\alpha_{j+r} = 0$ since $\ker(f_r)=0$.

 For $0\leq i\leq r$, commutativity of \eqref{e:Koszul2} gives
a restriction
\[
\iota_i:
\operatorname{im}(f_{i+1})
\longrightarrow
\operatorname{im}(g_{i+1}).
\]
Using these restricted maps, define
\[
\beta_{j+i}
=
\iota_{i,*}\alpha_{j+i}
\in
H^{j+i}(\B,\operatorname{im}(g_{i+1})).
\]
The image of $\beta_j$ under the map
\[
H^j(\B,\operatorname{im}(g_1))
\longrightarrow H^j(\B,J)
\]
is $\iota_*(\alpha_j)$, so it suffices to prove that
$\beta_j=0$.

For $0\leq i<r$, naturality of the connecting homomorphisms
shows that the obstruction to lifting $\beta_{j+i}$ to
$H^{j+i}(\B,\cD^{i+1})$ is the image of $\beta_{j+i+1}$
under the map induced by
\[
\operatorname{im}(g_{i+2})
\hookrightarrow
\ker(g_{i+1}).
\] Now, $\beta_{j+r} = 0$ since $\alpha_{j+r}=0$, and hence we conclude that $\beta_{j+r-1}$ lifts to a class in $H^{j+r-1}(\B, \cD^r)$.  However, $\cD^r$ is $\B$-acyclic, so $H^{j+r-1}(\B, \cD^r) =0$ and hence $\beta_{j+r-1}=0$.  Therefore, $\beta_{j+r-2}$ lifts to a class in  $H^{j+r-2}(\B, \cD^{r-1})$; again this cohomology group vanishes so $\beta_{j+r-2} = 0$.  This downward cascading continues and we obtain $\beta_i = 0 $ for all $i= j, \dotsc, j+r$.  In particular $\beta_j = 0$ as desired.
\end{proof}

In the remainder of the paper we will prove Theorem~\ref{t:comm}.  This will conclude the proof of Theorem~\ref{t:rl}.

\section{Construction of Complexes} \label{s:cc}

The only result from \S\ref{s:mit} that remains to be proven is
Theorem~\ref{t:comm}, which asserts the existence of a morphism of
complexes
\[
\mathcal C^\bullet\longrightarrow\mathcal D^\bullet
\]
with $\mathcal C^0=J'$, $\mathcal D^0=J$, and degree-zero map
the inclusion $J'\hookrightarrow J$, such that
$\mathcal C^\bullet$ is exact and $\mathcal D^i$ is an acyclic
$\B$-module for every $i>0$. In this section, we construct
such a morphism.

\subsection{Generators of $J'$} \label{s:generators}

Recall the ring $R$ defined in (\ref{e:r0def})--(\ref{e:rdef}) and the ideal $J' \subset R$ defined in \S\ref{s:invariance}.
We distinguish three types of generators of $J'$:
\begin{enumerate}
\item[i.] The elements $
\sum_{i=1}^r \boldsymbol{\epsilon}_i \boldsymbol{b}_i 
$
 for each row of type I in $D$.
\item[ii.] The elements 
$
(\boldsymbol{a}_i + \boldsymbol{\nu}_i)\boldsymbol{b}_j 
+ \boldsymbol{b}_i \boldsymbol{d}_j - \sum_{k = 1}^r \boldsymbol{\delta}_{ijk} \boldsymbol{b}_k
$ for each row of type II in $D$.
\item[iii.] For each $v \in (\Sigma \setminus \{v_0\}) \sqcup \cP$ and distinct $\sigma, \tau \in \cB_v$, the elements
\[ B(\sigma, \tau) = \boldsymbol{b}_\sigma\boldsymbol{b}'_{\tau} - 
\boldsymbol{b}_\tau\boldsymbol{b}'_{\sigma}, \]
where \begin{equation} \label{e:bpdef}
\boldsymbol{b}'_{\sigma} := \boldsymbol{x}_{\sigma} - \boldsymbol{a}_{\sigma}. \end{equation}
\end{enumerate}

Let $k$ be the number of initial generators chosen in
\eqref{e:rhoidef}, including the $|\cP|$ reserved generators
with indices $j_v$. The supplemental generators introduced
in \S\ref{s:ns} have indices $k+1,\ldots,r$.
Reorder only these supplemental indices so that the type III
indices are $\ell+1,\ldots,r$, whose
$\boldsymbol b_i$-variables vanish in $R$.
The surviving variables are therefore
$\boldsymbol b_1,\ldots,\boldsymbol b_\ell$.

The local blocks use the reserved initial variables
$\boldsymbol b_{j_v}$, $v\in\cP$, together with the
supplemental variables associated to type IV and V rows.
The remaining initial variables lie outside all local blocks.

Returning now to our generators of the ideal $J'$ described above, we may view the generators of types {i and ii} as linear forms in the variables $\boldsymbol{b}_1,\ldots, \boldsymbol{b}_{\ell}$: \begin{equation} L_i = \sum_{j=1}^{\ell} \boldsymbol{V}_{ij} \boldsymbol{b}_j \in R. \end{equation} The coefficients $\boldsymbol{V}_{ij} \in R$ are \emph{generic}, in that we can present $R$ as 
\[ R = S_0[\boldsymbol{a}_i, \boldsymbol{c}_i, \boldsymbol{d}_i]_{i=1,\ldots,r}[\boldsymbol{V}_{ij}][\boldsymbol{b}_1,\ldots, \boldsymbol{b}_{\ell}], \] 
Here the coefficients $\boldsymbol V_{ij}$ and the subring
$S_0\subset R_0$ are described as follows. For a row of type I,
the coefficient of $\boldsymbol b_j$ is
$\boldsymbol\epsilon_j$. For a row of type II associated to
the pair $(u,v)$, the coefficient of $\boldsymbol b_j$ is
\[
-\boldsymbol\delta_{uvj}
+\delta_{jv}(\boldsymbol a_u+\boldsymbol\nu_u)
+\delta_{ju}\boldsymbol d_v,
\qquad 1\leq j\leq\ell,
\]
where the unbolded symbols $\delta_{jv}$ and $\delta_{ju}$
are Kronecker deltas. In particular, when $u=v$, the
coefficient of $\boldsymbol b_u$ is
\[
\boldsymbol a_u+\boldsymbol\nu_u+\boldsymbol d_u
-\boldsymbol\delta_{uuu}.
\]
We obtain $S_0$ by deleting from the polynomial generators of
$R_0$, separately for each row, the variables
$\boldsymbol\epsilon_j$ for a type I row and
$\boldsymbol\delta_{uvj}$ for a type II row, with
$1\leq j\leq\ell$. Replacing these deleted variables by the
corresponding coefficients $\boldsymbol V_{ij}$ is an
invertible change of variables over
\[
S_0[\boldsymbol a_i,\boldsymbol c_i,\boldsymbol d_i]_{i=1}^r,
\]
which gives the displayed presentation of $R$.

The coefficients of type I and the unshifted coefficients
$-\boldsymbol\delta_{uvj}$ are fixed by $\B$. The diagonal
coefficient
$\boldsymbol a_u+\boldsymbol\nu_u+\boldsymbol d_u
-\boldsymbol\delta_{uuu}$ is also fixed, since
$\boldsymbol a_u+\boldsymbol d_u$ is a trace.
For $u\neq v$, the coefficients involving
$\boldsymbol a_u$ or $\boldsymbol d_v$ need not be fixed.
Thus, for constructing the resolutions of
Theorem~\ref{t:comm}, we may treat the generators of types
i and ii as generic linear forms
$\sum_j\boldsymbol V_{ij}\boldsymbol b_j$;
for proving $\B$-acyclicity, we retain their descriptions
in terms of the original type I and II relations.

We may do a further change of variables and present $R$ as 
\begin{equation}\label{e:present-R}
R = T_0[\boldsymbol{b}'_{\sigma} , \boldsymbol{V}_{ij}][\boldsymbol{b}_1,\ldots,\boldsymbol{b}_{\ell}],
\end{equation}
for $T_0$ a subring of $R_0[\boldsymbol{a}_i, \boldsymbol{c}_i, \boldsymbol{d}_i]_{i=1,\ldots,r}$.  Here the $\boldsymbol{b}'_{\sigma}$ are defined in (\ref{e:bpdef}).

\subsection{Motivation via regular sequences}

For each $v \in (\Sigma \setminus \{ v_0 \}) \sqcup \cP$, fix $\sigma_v \in \mathcal{B}_v$. Over the ring $R' = R[(\boldsymbol{b}'_{\sigma_v})^{-1}]_{v \in (\Sigma \setminus \{ v_0 \}) \sqcup \cP}$, we may choose a smaller set of elements to generate the extended ideal $R' J'$, by removing certain generators of {type iii}. Indeed, the following elements generate $R'J'$:
\begin{itemize}
\item the linear forms $\sum \boldsymbol{V}_{ij} \boldsymbol{b}_j$ {of types i or ii},
\item for each $v \in (\Sigma \setminus \{ v_0 \}) \sqcup \cP$, $B(\sigma_v,\tau) =  \boldsymbol{b}'_{\tau}\boldsymbol{b}_{\sigma_v} - 
\boldsymbol{b}'_{\sigma_v} \boldsymbol{b}_\tau$ for all $\tau \in \cB_v \setminus \{ \sigma_v \}$.
\end{itemize}

Indeed, for $\sigma \in \mathcal{B}_v$, we have 
\[ B(\sigma, \tau) = (\boldsymbol{b}'_{\sigma_v})^{-1}( \boldsymbol{b}'_{\sigma} B(\sigma_v, \tau) - \boldsymbol{b}'_{\tau} B(\sigma_v, \sigma)). \]
These bulleted elements form a regular sequence in $R'$, a fact whose proof we omit as it will not end up being used. To obtain an exact complex of $R'$-modules resolving $R'/R' J'$ as required in Theorem~\ref{t:comm}, we would simply take the Koszul complex on this regular sequence.

However, to prove Theorem \ref{t:comm} we need to construct a resolution of $R/J'$, not of $R'/R'J'$. Unfortunately, the ring $R/J'$ cannot typically be resolved using a Koszul complex, as the ideal $J'$ is not typically generated by a regular sequence. Nevertheless, it is still possible to resolve $R/J'$ using a generalization of the Koszul complex due to Buchsbaum--Rim \cite{b}, which we now describe.

\subsection{Buchsbaum--Rim complexes}
For this section, let $R$ denote any commutative ring. Associated to an $R$-linear map \[ f \colon V \lra W,  \qquad V = R^n,\ W = R^m, m \leq n, \] Buchsbaum--Rim \cite{b} define two complexes $R(f)$ and $R(\det(f))$. {These complexes have a degree-decreasing differential.}

\subsubsection{Special cases} \label{s:special}

\noindent \textbf{Case $m = 1$.}
For a map $f \colon V \lra R$, \[ R(f) = R(\det(f)) = \text{Koszul complex on the elements } f(e_1),\ldots, f(e_n) \in R. \] The degree $k$ term of the complex $R(f)$ is $\bigwedge^k V$, and in degrees $1$ and $0$ the complex is given by $\begin{tikzcd} V \ar[r,"f"] & R \end{tikzcd}$.  Here and throughout, all alternating powers are understood to be over $R$.

\bigskip

\noindent\textbf{Case $m = 2$.}
Consider a map $f \colon V \lra W$, with \[ V = R^n = \bigoplus_i R e_i, \qquad W = R^2 = R w_1 \oplus R w_2. \] Write $f(e_i) = b_i w_1 + b_i' w_2$, and {define 
\[ r_{ij} := b_i b'_j - b_j b'_i, \qquad d_{ijk} = r_{ij} e_k + r_{jk} e_i + r_{ki} e_j.\]
}  The complex $R(f)$ begins as follows:
\begin{equation}
\begin{tikzcd}
 \cdots \ar[r] & \bigwedge^2(W^*) \otimes \bigwedge^3 V \ar[r] &  V \ar[r, "f"] & W.
\end{tikzcd}
\end{equation}
{
The image of the map 
\[ \bigwedge^2(W^*) \otimes \bigwedge^3 V \lra V, \qquad 
(w_1^* \wedge w_2^*)\otimes(e_i \wedge e_j \wedge e_k) \mapsto d_{ijk}\]is generated as an $R$-module by  the $d_{ijk}$.}
The complex $R(\det(f))$ begins as follows:
\begin{equation}
\begin{tikzcd}
 \cdots \ar[r] & \Big(W^* \otimes \bigwedge^3 V\Big) \oplus \Big(\bigwedge^2(W^*) \otimes \bigwedge^4 V\Big) \ar[r] & \bigwedge^2 V \ar[r, "\det(f)"] & \bigwedge^2 W.
\end{tikzcd}
\end{equation}

For $k \geq 1$, the degree $k$ term of the complex $R(\det(f))$ is:
\begin{equation} \bigoplus_{s_i \in \{1, 2\}} \Big(\bigwedge\nolimits^{s_1} W^*\Big) \otimes \cdots \otimes \Big(\bigwedge\nolimits^{s_{k-1}} W^*\Big) \otimes \bigwedge\nolimits^{2 + \sum s_i} V. \end{equation}

The image of the map $\det(f) \colon \bigwedge^2(W^*) \otimes \bigwedge^2 V \lra R$ is generated over $R$ by the $2 \times 2$ minors of the $2 \times n$ matrix corresponding to $f$.

The complex $R(\det(f))$ is functorial in the following sense. Given another $R$-linear map \[ f' \colon V' \lra W, \quad V' = R^{n'},\] as well as an 
$R$-linear map $g \colon V \lra V'$ such that $f' \circ g = f$, there is a chain homomorphism \[ \det(g) \colon R(\det(f)) \lra R(\det(f')) \] defined using the map $\bigwedge^i g \colon \bigwedge^i V \lra \bigwedge^i V'$.  

\subsubsection{Exactness of BR-complexes}
For $m \leq k \leq n$, define $f(k) \colon R^k \lra R^m$ by restricting $f$ to the first $k$ summands of $R^n$. 
\begin{definition}[\cite{b}] \label{d:regular}
The map $f$ is {\em regular} if for all $m \leq k \leq n$, we have $H_1(R(f(k))) = 0$. 
\end{definition}

{\begin{remark} \label{r:regular} As described in \S\ref{s:special} above, when $m=1$ the complex $R(f(k))$ is the precisely the Koszul complex associated to the elements $f(e_1), \dots, f(e_k)$. The vanishing of $H_1(R(f(k)))$ for all $1 \le k \le n$ as stated in Definition~\ref{d:regular}
is equivalent to the usual definition of $f(e_1), \dots, f(e_n)$ forming a regular sequence in $R$.
 Indeed, arguing
inductively on $k$, one has
\[ H_1(K(x_1,\dots,x_k))
\cong
\frac{\{a\in R: ax_k\in (x_1,\dots,x_{k-1})\}}
{(x_1,\dots,x_{k-1})},
\]
provided $x_1,\dots,x_{k-1}$ is already regular.
See Lemma~\ref{l:reg} for an analogous statement when $m=2$. \end{remark}}

The following result of Buchsbaum--Rim generalizes the well-known exactness of the Koszul complex associated to a regular sequence.
\begin{theorem}[\cite{b}]  The complexes $R(f)$ and $R(\det(f))$ satisfy the following.
\begin{enumerate}
\item If $f$ is regular then $H_i(R(f)) = 0$ for all $i \geq 1$.
\item If $H_i(R(f)) = 0$ for $i \geq j$, then $H_i(R(\det(f))) = 0$ for $i \geq j$.
\end{enumerate}
\end{theorem}

\begin{corollary}
If $f \colon R^n \lra R^m$ is regular then $R(\det(f))$ is a resolution of $R/\im(\det(f))$.
\end{corollary}

We will apply this result to multiple complexes $R(\det(f))$ and then take the tensor product, necessitating the following lemma.

\begin{lemma}\label{l:tensor}
Consider a ring $R$ with maps $f_i \colon R^{n_i} \lra R^{m_i}$, $1\leq i \leq k$, such that $m_i \leq n_i$. Define ideals $J_i = \im(\det(f_i))$, and let $J = J_1 + \cdots + J_{k}$.
 Suppose that $f_i \otimes R/(J_1 + \cdots + J_{i-1})$ is regular for all $i \leq k$. Then the tensor product over $R$ of these complexes, $\bigotimes_{j=1}^k R(\det(f_j))$, is a resolution of $R/J$.
\end{lemma}
\begin{proof}
By induction on $k$, we may suppose that the natural map \[ \otimes_{j=1}^{k-1} R(\det(f_j)) \lra R/\Big(\sum_{j=1}^{k-1} J_j\Big)\] is a quasi-isomorphism. Then, since $R(\det(f_k))$ is a complex of free $R$-modules, the Tor spectral sequence {
(see \cite{stacks-tor}*{\S15.62})} implies that 
\[ \otimes_{j=1}^{k} R(\det(f_j)) \lra R/\Big(\sum_{j=1}^{k-1} J_j\Big) \otimes R(\det(f_k)) \] is also a quasi-isomorphism. By the assumption on regularity, the codomain has zero homology in positive degrees and the natural map \[ R/\Big(\sum_{j=1}^{k-1} J_j\Big) \otimes R(\det(f_k)) \lra R/\Big(\sum_{j=1}^{k} J_j\Big) \] is a quasi-isomorphism. 
\end{proof}

\subsection{Outline of proof of Theorem \ref{t:comm}} \label{s:outline}

{
Recall that the ideal $J$ is the
full relation ideal from \S\ref{s:ri}. The smaller ideal $J'\subset J$ is only
$\B$-stable, not $\G$-stable; it is generated by the $b$-coefficients of the
relations. 
}

In \S\ref{s:generators}, we described generators for $J'$ of the form:
\begin{itemize}
\item $L_i=\sum_{j=1}^{\ell}\boldsymbol V_{ij}\boldsymbol b_j,
 1\leq i\leq n,$
\item $B(\sigma, \tau) = \boldsymbol{b}_\sigma\boldsymbol{b}'_{\tau} - 
\boldsymbol{b}_\tau\boldsymbol{b}'_{\sigma}$ for $v \in (\Sigma \setminus \{v_0\}) \sqcup \cP,$ distinct {$\sigma, \tau \in \cB_{v}$}.
\end{itemize}
Recall also the $\mathbf{B}$-modules $\cB$ and $\cV$ defined in Example \ref{e:adjoint}, and the notation
   $R(1) = \cB \otimes_{\Z} R.$ We also write $\cV_R = \cV \otimes_{\Z} R$; this is $\cong R^2$ as an $R$-module.
   
To handle the first collection of generators bulleted above, we define the $\B$-equivariant map
\begin{equation} \label{e:fdef}
 f \colon R(1)^n = \bigoplus_{i=1}^n R(1)e_i \lra R, \qquad f(e_i) = L_i.  
 \end{equation}
For the second collection of generators bulleted above, we define for each $v \in (\Sigma \setminus \{v_0\}) \sqcup \cP$ the 
$\B$-equivariant map 
\begin{equation} f_v \colon \bigoplus_{\sigma \in \cB_{v}} R(1)e_{\sigma} \lra \mathcal{V}_R, \qquad
f_v(e_{\sigma}) = A \otimes \boldsymbol{b}_{\sigma} + B \otimes\boldsymbol{b}'_{\sigma}.
\end{equation}

We consider the Buchsbaum--Rim complexes associated to the maps $f, f_v$:
\begin{itemize}
\item The Koszul complex $R(f)$ is given in degrees $1$ and $0$ by the map $f$, whose image is the ideal $J'_0 := (\{\sum \boldsymbol{V}_{ij} \boldsymbol{b}_j\}_i)$.
\item The Buchsbaum--Rim complex $R(\det(f_v))(-1)$ (note the twist) is given in degrees $1$ and $0$ by the map \[ \det(f_v) \colon \Big(\bigwedge^2 \oplus_{\sigma} R(1)\Big) \otimes R(-1) \lra \bigwedge^2\mathcal{V}_R \otimes R(-1) \cong R. \] The image of $\det(f_v)$ is the ideal $J'_v := (B(\sigma,\tau) \colon \sigma, \tau \in \mathcal{B}_v) \subset R$. 
\end{itemize}

If $\cB_v$ consists of a single label, then $J'_v=0$, and
we interpret $R(\det(f_v))(-1)$ as $R$ concentrated in
degree zero.

The tensor product of these complexes, \begin{equation} \cC^\bullet := R(f) \otimes_R \bigotimes_{v \in (\Sigma \setminus \{ v_0 \}) \sqcup \cP} R(\det(f_v))(-1), \label{e:cdef}
 \end{equation} is a complex of $\B$-modules such that $H^0(\cC^\bullet) = R/(J'_0 + \sum_v J'_v) = R/J'$.

To prove Theorem \ref{t:comm} it will suffice to prove:
\begin{prop}\label{p:resolve-j} The complex $\cC^\bullet$ satisfies the following.
\begin{enumerate}
\item $H^i(\cC^\bullet) = 0$ for $i > 0$, i.e. $\cC^\bullet$ is a resolution of $R/J'$ and the image of $\cC^1 \lra \cC^0 = R$ is the ideal $J'$.
 \item $\cC^\bullet$ maps to a complex $\cD^\bullet$, with
degree-zero map $\operatorname{id}_R$, such that the image of
$\cD^1\lra\cD^0=R$ is $J$, and $\cD^i$ is an acyclic
$\B$-module for every $i\geq0$.
\end{enumerate}
\end{prop}

To obtain the complexes of Theorem~\ref{t:comm}, replace the
degree-zero terms $\cC^0=\cD^0=R$ by the images $J'$ and
$J$, respectively, and corestrict the degree-one differentials.
Leave all positive-degree terms and maps unchanged. The resulting
top complex is exact, the positive-degree terms of the bottom
complex remain $\B$-acyclic, and the degree-zero map is the
inclusion $J'\hookrightarrow J$.

We now explain how this proposition follows from the results of \S\ref{ss:extend-f}-\S\ref{ss:tensor-exact}. In the next two sections (\S \ref{ss:extend-f} and  \S \ref{ss:extend-fv}), we will extend the maps $f$ and $f_v$ to maps \[ \widetilde{f} \colon \bigoplus_{i=1}^n \mathcal{A}_R \lra R,\qquad \widetilde{f}_v \colon \bigoplus_{\sigma \in \mathcal{B}_v} \mathcal{V}_R \lra \mathcal{V}_R. \]
Here $\cA_R = \cA \otimes_\Z R$ where $\cA$ denotes the adjoint representation
defined in example \ref{e:adjoint}. 
There are natural maps $R(f) \lra R(\widetilde{f})$ and $R(\det(f_v))(-1) \lra R(\det(\widetilde{f}_v))(-1)$. We will show that these maps factor through subcomplexes \[ R(\widetilde{f})' \subset R(\widetilde{f}), \qquad R(\det(\widetilde{f}_v))(-1)' \subset R(\det(\widetilde{f}_v))(-1) \] consisting of good $\G$-modules tensored with the $\B$-module $R$. 

Define \[ J_0 := \text{ideal of } R  \text{ generated by the coefficients of (\ref{e:b1}) and (\ref{e:b2}),} \] and for $v \in (\Sigma \setminus \{v_0\}) \sqcup \cP$ define \[ J_v :=  (A(\sigma,\tau),\ B(\sigma,\tau),\ C(\sigma,\tau),\ D(\sigma,\tau) \colon \sigma, \tau \in \mathcal{B}_v,\  \sigma \neq \tau) \subset R. \] Then by definition {(see \S\ref{s:ri}),}  \[ J = J_0 + \sum_{v \in (\Sigma \setminus \{v_0\}) \sqcup \cP} J_v. \]

We will show that $H^0(R(\widetilde{f})') = R/J_0$ and $H^0(R(\det(\widetilde{f}_v))'(-1)) = R/J_v$.
Define \begin{equation}\cD^\bullet := R(\widetilde{f})' \otimes_R \bigotimes_{v \in (\Sigma \setminus \{ v_0 \}) \sqcup \cP} R(\det(\widetilde{f}_v))'(-1).  \label{e:ddef}
\end{equation} It follows that \[ H^0(\cD^\bullet) = R/(J_0 + \sum_{v \in (\Sigma \setminus \{v_0\}) \sqcup \cP} J_v) = R/J, \] i.e. that the image of the map $\cD^1 \lra \cD^0 = R$ is $J$. As all terms of $\cD^\bullet$ consist of good $\G$-modules tensored with $R$, we conclude by Theorem \ref{t:r-acyclic} that all terms of $\cD^\bullet$ are acyclic $\B$-modules.

It will then remain to show that $H^i(\cC^\bullet) = 0$ for $i > 0$. This will be proved using Proposition \ref{p:br-exact} below, building on the work of sections \S\ref{ss:br-regular}-\S\ref{ss:tensor-exact}.

\subsection{Extension of $f$}\label{ss:extend-f}

We define a map \[ \widetilde{f} \colon \cA_R^n =  \bigoplus_{i=1}^n \cA_R e_i \lra R \] extending the map $f$ defined in (\ref{e:fdef}) as follows (here {\em extension} refers to the natural inclusion $R(1) = \cB_R \subset \cA_R$, see Example~\ref{e:adjoint}).
\begin{itemize}
\item For a type I row indexed by $i$, write its original
matrix relation as
\[
\sum_{j=1}^r\boldsymbol\epsilon_{i,j}\boldsymbol\rho_j.
\]
We set
\begin{align*}
\widetilde f(Ae_i)&=\sum_{j=1}^r
  \boldsymbol\epsilon_{i,j}\boldsymbol a_j,&
\widetilde f(Be_i)&=\sum_{j=1}^r
  \boldsymbol\epsilon_{i,j}\boldsymbol b_j=L_i,\\
\widetilde f(Ce_i)&=\sum_{j=1}^r
  \boldsymbol\epsilon_{i,j}\boldsymbol c_j,&
\widetilde f(De_i)&=\sum_{j=1}^r
  \boldsymbol\epsilon_{i,j}\boldsymbol d_j.
\end{align*}
The extra row subscript on $\boldsymbol\epsilon$ is only to
distinguish the separate sets of variables attached to different
rows. The sums retain the terms with $j>\ell$: only their
$\boldsymbol b_j$ coefficients vanish in $R$.

\item For linear forms $L_i$ corresponding to relations of type II, i.e.\ to the $b$-coefficient of $(\boldsymbol{\rho}_i + \boldsymbol{\nu}_i) \boldsymbol{\rho}_j - \sum_{k = 1}^r \boldsymbol{\delta}_{ijk} \boldsymbol{\rho}_k$, we similarly extend $f$ using the other 3 coefficients.
 
\end{itemize}
 
The Koszul complex $R(\widetilde{f})$ is a complex of $\B$-modules. By the functoriality of the Koszul complex we obtain a $\B$-equivariant map $R(f) \lra R(\widetilde{f})$.

\begin{lemma} Expanding the terms of the complex $R(\widetilde{f})$ via the Kunneth formula, 
the image of $R(f) \lra R(\widetilde{f})$ is contained in the subcomplex $R(\widetilde{f})'$ whose $k$th term is 
\[\bigoplus_{\# I = k} \bigotimes_{i \in I} (\cA e_i) \otimes R. \]
\end{lemma}
\begin{proof}
The $k$th degree term of the Koszul complex $R(f)$ is given by \[\bigwedge^k \Big(\bigoplus_i R(1) e_i \Big) \cong \bigoplus_{\# I = k} \bigotimes_{i \in I} R(1), \] 
since $R(1)$ is free of rank 1 over $R$. The image of this map is contained in the subcomplex $R(\widetilde{f})'$ whose $k$th degree term is identified by the Kunneth formula with \[ \bigoplus_{\# I = k} \bigotimes_{i \in I} \cA_R e_i, \]
i.e.\ the submodule of $\bigwedge^k\bigoplus \cA_R e_i$ generated by wedge products \[ v_{i_1} e_{i_1} \wedge \cdots \wedge v_{i_k} e_{i_k}, \quad v_{i_j} \in \mathcal{A}_R,\] such that $i_1 < \cdots < i_k$.
\end{proof}

\begin{lemma}\label{l:adj-r-acyclic}
The terms of $R(\widetilde{f})'$ are the tensor product of good $\G$-modules with $R$.
\end{lemma}
\begin{proof}
This follows directly from Corollary \ref{c:adj-good}.
\end{proof}

\subsection{Extension of $f_v$}\label{ss:extend-fv}

For $\sigma \in \cB_v$, recall the notation $\boldsymbol{b}'_\sigma = \boldsymbol{x}_\sigma - \boldsymbol{a}_\sigma$. Define
now also $\boldsymbol{c}'_\sigma = \boldsymbol{x}_\sigma - \boldsymbol{d}_\sigma$.
We define \[ \widetilde{f}_v \colon \bigoplus_{\sigma \in \cB_v} \mathcal{V}_R e_{\sigma} \lra \mathcal{V}_R\] by \begin{align*}
\widetilde{f}_v(A e_{\sigma}) &= A \otimes \boldsymbol{c}'_\sigma + B \otimes \boldsymbol{c}_\sigma, \\
\widetilde{f}_v(B e_{\sigma}) &= A \otimes \boldsymbol{b}_\sigma + B \otimes \boldsymbol{b}'_\sigma. 
\end{align*}

\begin{lemma}
The map $\widetilde{f}_v$ is $\B$-equivariant.
\end{lemma}
\begin{proof}

Let $g = \mat{x}{0}{y}{z}$. We have
\begin{align*}
 g \cdot \boldsymbol{b}_\sigma &= \frac{z}{x} \boldsymbol{b}_\sigma, \\
 g \cdot \boldsymbol{b}'_\sigma &= g \cdot (\boldsymbol{x}_\sigma - \boldsymbol{a}_\sigma) = \boldsymbol{x}_\sigma - \boldsymbol{a}_\sigma - \frac{y}{x} \boldsymbol{b}_\sigma = \boldsymbol{b}'_\sigma - \frac{y}{x} \boldsymbol{b}_\sigma \\
 g \cdot \boldsymbol{c}_\sigma &= -\frac{y}{z}\boldsymbol{a}_\sigma - \frac{y^2}{xz}\boldsymbol{b}_\sigma + \frac{x}{z}\boldsymbol{c}_\sigma + \frac{y}{z} \boldsymbol{d}_\sigma\\
 g \cdot \boldsymbol{c}'_\sigma &= g \cdot (\boldsymbol{x}_\sigma - \boldsymbol{d}_\sigma) = \boldsymbol{x}_\sigma - \boldsymbol{d}_\sigma + \frac{y}{x} \boldsymbol{b}_\sigma =  \boldsymbol{c}'_\sigma + \frac{y}{x} \boldsymbol{b}_\sigma. 
 \end{align*}
We compute
\begin{align*}
 g \cdot \left(A \otimes \boldsymbol{b}_\sigma + B \otimes \boldsymbol{b}'_\sigma\right)  &= \left(A + \frac{y}{x} B\right)  \otimes \left(\frac{z}{x} \boldsymbol{b}_\sigma\right)  + \left(\frac{z}{x} B\right)  \otimes  \left(\boldsymbol{b}'_\sigma - \frac{y}{x} \boldsymbol{b}_\sigma\right)  \\
 &= \frac{z}{x}\left( A \otimes \boldsymbol{b}_\sigma + B \otimes \boldsymbol{b}'_\sigma\right).
 \end{align*}
We also have
\begin{align*}
g \cdot\left(A \otimes \boldsymbol{c}'_\sigma + B \otimes \boldsymbol{c}_\sigma \right) &=  \left(A + \frac{y}{x} B\right)  \otimes \left(\boldsymbol{c}'_\sigma + \frac{y}{x} \boldsymbol{b}_\sigma\right)  +\\
& \ \ \ \  \left(\frac{z}{x} B\right)  \otimes \left(-\frac{y}{z}\boldsymbol{a}_\sigma - \frac{y^2}{xz}\boldsymbol{b}_\sigma + \frac{x}{z}\boldsymbol{c}_\sigma + \frac{y}{z} \boldsymbol{d}_\sigma\right)  \\
 &= A \otimes \left(  \boldsymbol{c}'_\sigma + \frac{y}{x} \boldsymbol{b}_\sigma\right)  \\ & \ \ \ +  B\otimes \left( \frac{y}{x}\left(\boldsymbol{x}_\sigma - \boldsymbol{d}_\sigma\right)  + \frac{y^2}{x^2} \boldsymbol{b}_\sigma - \frac{y}{x} \boldsymbol{a}_\sigma - \frac{y^2}{x^2}\boldsymbol{b}_\sigma + \boldsymbol{c}_\sigma  + \frac{y}{x} \boldsymbol{d}_\sigma\right) \\
 &= \left(A \otimes \boldsymbol{c}'_\sigma + B \otimes \boldsymbol{c}_\sigma \right)  + \frac{y}{x}  \left( A \otimes \boldsymbol{b}_\sigma + B \otimes \boldsymbol{b}'_\sigma\right).
\end{align*}
These two identities show that
\[
g\cdot\widetilde f_v(Be_\sigma)
=\widetilde f_v(g\cdot Be_\sigma)
\quad\text{and}\quad
g\cdot\widetilde f_v(Ae_\sigma)
=\widetilde f_v(g\cdot Ae_\sigma).
\]
The result follows.
\end{proof}

The BR-complex $R(\det(f_v))$ is a complex of $\B$-modules whose degree $0$ term is \[ \bigwedge\nolimits^2 \mathcal{V}_R \cong R(1),\] and whose degree $k \geq 1$ term is:
\begin{equation} \bigoplus_{s_i \in \{1, 2\}} \Big(\bigwedge\nolimits^{s_1} \mathcal{V}_R^*\Big) \otimes \cdots \otimes \Big(\bigwedge\nolimits^{s_{k-1}} \mathcal{V}_R^*\Big) \otimes \bigwedge\nolimits^{2 + \sum s_i} \Big(\bigoplus_{\sigma \in \cB_v} R(1)e_{\sigma}\Big). \end{equation}
The BR-complex $R(\det(\widetilde{f}_v))$ is the same but with $\bigwedge^{k} (\bigoplus_{\sigma} R(1)e_{\sigma})$ replaced by $\bigwedge^{k} (\bigoplus_{\sigma} \mathcal{V}_Re_{\sigma})$. The map $R(\det(f_v)) \lra R(\det(\widetilde{f}_v))$ has image contained in the subcomplex $R(\det(\widetilde{f}_v))'$ where each $\bigwedge^k (\bigoplus_{\sigma} \mathcal{V}_Re_{\sigma})$ is replaced by the $R$-submodule generated by wedge products $v_{\sigma_1} e_{\sigma_1} \wedge \cdots \wedge v_{\sigma_{k}} e_{\sigma_k}$, $v_i \in \mathcal{V}_R$, such that $\sigma_{i} \neq \sigma_{j}$ for $i \neq j$. This submodule is isomorphic to a direct sum of modules $\bigotimes^k_R \mathcal{V}_R$.

\begin{lemma}
$H^0(R(\det(\widetilde{f}_v))'(-1)) = R/J_v$.
\end{lemma}
\begin{proof}

We are restricting the map \[ \det(\widetilde{f}_v) \colon \bigwedge^2 (\oplus_{\sigma} \mathcal{V}_R e_{\sigma}) \otimes R(-1) \lra R \] to the wedge-products $v_{\sigma} e_{\sigma} \wedge v_{\tau} e_{\tau}$, $v_{\sigma}, v_{\tau} \in \mathcal{V}_R$, with $\sigma \neq \tau$. We need to show that the image of this restricted map is $J_v$.
A direct computation using the definitions in (\ref{e:ABCD}) shows that the elements $B e_{\sigma} \wedge A e_{\tau}$, $B e_{\sigma} \wedge B e_{\tau}$, $A e_{\sigma} \wedge A e_{\tau}$, $A e_{\sigma} \wedge B e_{\tau}$, map to $\pm A(\sigma, \tau)$, $\pm B(\sigma, \tau)$, $\pm C(\sigma, \tau)$, $\pm D(\sigma,\tau)$. By definition, these 
latter elements generate $J_v$.  The result follows.
\end{proof}

\begin{lemma} \label{l:localgood}
The terms of $R(\det(\widetilde{f}_v))'(-1)$ are the tensor products of good $\G$-modules with $R$.
\end{lemma}
\begin{proof}

The degree $0$ term of $R(\det(\widetilde{f}_v))'(-1)$ is $R$.  
The degree $k \geq 1$ term of $R(\det(\widetilde{f}_v))'(-1)$
 is a direct sum of modules of the form \begin{equation}\label{e:good1}\left(\Big(\bigwedge\nolimits^{s_1} \mathcal{V}_R^*\Big) \otimes \cdots \otimes \Big(\bigwedge\nolimits^{s_{k-1}} \mathcal{V}_R^*\Big) \otimes \mathcal{V}_R^{\otimes 2 + \sum s_i}\right) \otimes R(-1).\end{equation}

 Setting $m = \sum (s_i - 1) = \Big(\sum s_i\Big) - (k-1)$, (\ref{e:good1}) is isomorphic to 
 \begin{equation} \label{e:good2}
 \left(\Big(\bigwedge\nolimits^2 \mathcal{V}_R^* \Big)^{\otimes m} \otimes (\mathcal{V}_R^*)^{\otimes (k-1- m)} \otimes \mathcal{V}_R^{\otimes(k+m+1)}\right) \otimes R(-1)
 \end{equation}
 
 Using the isomorphisms $\bigwedge^2 \mathcal{V} \cong \Z(1)$, $\mathcal{V} \otimes \mathcal{V}^* \cong \mathcal{A}$, and $\mathcal{V} \otimes \mathcal{V} = \mathcal{A}(1)$, 
  one readily checks that (\ref{e:good2}) is isomorphic to $\mathcal{A}^{\otimes k} \otimes R$.
The result now follows from Corollary~\ref{c:adj-good}.
\end{proof}

As described in \S\ref{s:outline} above, Lemmas \ref{l:adj-r-acyclic} and \ref{l:localgood} imply that the complex $\cD^\bullet$
defined in (\ref{e:ddef}) satisfies the property that $\cD^i$ is $\B$-acyclic for $i > 0$.  In other words, we have finished the proof of the second part of Proposition~\ref{p:resolve-j}.  It remains to prove the first part of the proposition.

\subsection{Regularity criterion for maps $f \colon R^n \lra R^2$}\label{ss:br-regular}

Let $R$ be a commutative ring. Given a map $f \colon R^n \lra R^2$, write $f(e_i) = (b_i, b_i')$, and define $r_{ij} = b_i b'_j - b_j b'_i$.

\begin{lemma}\label{l:reg}
Suppose that for $k = 2, \dots, n$, the congruence
\[ r_{1k} x \equiv 0 \pmod{(r_{12}, \dots, r_{1(k-1)})} \]
{for $x \in R$} implies that
\[ x \in (r_{ij} \colon i,j < k ). \]
Then the map $f \colon R^n \lra R^2$ is regular.
\end{lemma}
\begin{proof}
The goal is to show that the kernel of $f$ is generated as an $R$-module by \[ d_{ijk} = r_{ij}e_k + r_{jk}e_i + r_{ki}e_j. \] 

We proceed via induction. In the base case $n=2$, our assumption is that  $r_{12}$ is a non-zerodivisor, and this implies that $\ker(f) = 0$.  We proceed to the general case.

Any element $x = (x_1, \ldots, x_n) \in \ker(f)$ satisfies $\sum r_{1i} x_i = 0$. 
Hence \[ r_{1n}x_n \equiv 0 \mod (r_{12}, \ldots, r_{1(n-1)}). \] The assumption implies that  we may write $x_n = \sum_{i<j < n} \alpha_{ij} r_{ij}$ for some $\alpha_{ij} \in R$. Therefore the final coordinate of $(x_1, \ldots,x_n) - \sum \alpha_{ij} d_{ijn}$ is zero. Writing this vector as $y = (y_1, \dots, y_{n-1}, 0)$, we can apply the inductive hypothesis to $y$ and the restriction of $f$ to its first $n-1$ coordinates.  Hence $y$ is a linear combination of the $d_{ijk}$, so $x$ is as well.
\end{proof}

\begin{corollary} \label{c:genericb}
Let $R_0$ be a  commutative ring, and let \[ R = R_0[b'_1, \ldots, b'_n, b_1, \ldots, b_n]. \] The map $f \colon R^n \lra R^2$ defined by $f(e_i) = (b_i, b_i')$ is regular.
\end{corollary}
\begin{proof}
We will show that the criterion of Lemma \ref{l:reg} holds. 
Let \[ J_k = (\{ r_{ij} \}_{1 \leq i,j \leq k}) \subset R. \] 

First we write $V = b'_1$ and work in $R[V^{-1}]$.
Note that \[ R[V^{-1}]\cdot J_k = (r_{12},\ldots,r_{1k}), \] and that \begin{equation} \label{e:rv}
 R[V^{-1}]/R[V^{-1}]\cdot J_k \cong R_0[b_1', \dots, b_n', b_1, b_{k+1}, \dots, b_n][(b_1')^{-1}]. \end{equation} 
 It is clear that the image of $r_{1(k+1)}$ is a non-zerodivisor in (\ref{e:rv}), and hence that the criterion of Lemma \ref{l:reg} holds.
 
It remains to prove that $R \cap (R[V^{-1}] \cdot J_k) = J_k.$  
Therefore suppose that $x \in R$ and that the image of $x$ vanishes in $(\ref{e:rv})$.  The element $x$ is a sum terms that are coefficients in $R_0$ times monomials
\[ b_1^{a_1} \cdots b_n^{a_n}(b_1')^{a_1'} \cdots (b_n')^{a_n'} \]
for nonnegative integers $a_i, a_i'$.  {As we now explain, we can partition the collection of these monomials such that the sum of the terms in each part vanishes modulo $J_k$.  To this end, fix nonnegative integers $c_1, \dots, c_k$, nonnegative integers $a_i, a_i'$ for $k < i \le n$, and nonnegative integers $A$ and $C$.
Consider all the monomials with $a_i + a_i' = c_i$ for $i \le k$, $\sum a_i =A$, and $\sum a_i' = C.$ These are the parts of our partition.
Applying the isomorphism (\ref{e:rv}), it follows that the sum of the coefficients of the monomials in $x$ in each part vanishes.

For example, $x$ could be a linear combination of expressions of the form $b_2 b'_3 - b'_2 b_3$ or 
\[ d_1 b_2 b_3 b'_4 + d_2 b_2 b'_3 b_4 + d_3 b'_2 b_3 b_4 \quad \text{ with } d_1 + d_2 + d_3 = 0. \]

 It is easy to see that (without their coefficients) all such monomials are equivalent modulo $J_k$, and therefore that the sum of these terms (with their coefficients) vanishes modulo $J_k$.  Therefore $x \in J_k$ as desired.}
\end{proof}

\subsection{Application of regularity criteria to the resolution of $J'$}\label{ss:tensor-exact}

\begin{prop}\label{p:br-exact}
Let $R_0$ be a commutative ring, let
\[
\{1,\ldots,N\}=S_1\sqcup\cdots\sqcup S_h
\]
be a partition into nonempty sets, and let $q,m\geq0$ satisfy
$m\leq q+h$. Put
\[
R=R_0[b'_1,\ldots,b'_N,b_1,\ldots,b_{N+q},
       V_{ij}:1\leq i\leq m,\ 1\leq j\leq N+q].
\]
For $1\leq a\leq h$, define
\[
f_a:R^{S_a}\longrightarrow R^2,
\qquad f_a(e_j)=(b_j,b'_j),
\]
and define
\[
f_{h+1}:R^m\longrightarrow R,
\qquad f_{h+1}(e_i)=L_i:=\sum_{j=1}^{N+q}V_{ij}b_j.
\]
Let $I_a=\operatorname{im}(\det(f_a))$ and
$I=I_1+\cdots+I_{h+1}$. Then
\[
\bigotimes_{a=1}^{h+1}R(\det(f_a))
\]
is a resolution of $R/I$. For a singleton $S_a$, we interpret
$R(\det(f_a))$ as $R$ concentrated in degree zero.
\end{prop}

\begin{proof}
The determinantal blocks use disjoint variables. Hence
Corollary~\ref{c:genericb} and Lemma~\ref{l:tensor}
show that the tensor product of the first $h$ complexes resolves
$R/(I_1+\cdots+I_h)$.
We prove that the remaining Koszul complex is exact in positive
degrees, first over a field and then universally over $\Z$.

Suppose $R_0=K$ is a field, and put
\[
A=K[b'_1,\ldots,b'_N,b_1,\ldots,b_{N+q}]/(I_1+\cdots+I_h).
\]
Each block of size $d$ is a generic rank-one determinantal ring
of dimension $d+1$; it is Cohen--Macaulay
\cite{brw}*{Theorem~3.5(c)}.
The blocks have disjoint variables, so $A$ is Cohen--Macaulay
of dimension $N+h+q$. If $B=(b_1,\ldots,b_{N+q})\subset A$,
then $A/B=K[b'_1,\ldots,b'_N]$, and
\[
\operatorname{ht}(B)=h+q.
\]
For $1\leq a\leq m$, consider a prime of $A[V_{ij}]$
containing $(L_1,\ldots,L_a)$. If it contains $B$, its height
is at least $h+q\geq a$. Otherwise some $b_j$ is invertible
locally, and the equations $L_i=0$, $1\leq i\leq a$, eliminate
the distinct variables $V_{ij}$. They therefore form a regular
sequence in that localization, and the prime has height at least
$a$. Krull's height theorem gives the opposite bound for the
height of $(L_1,\ldots,L_a)$. Thus every such prefix has height
$a$, and Cohen--Macaulayness implies that
$L_1,\ldots,L_m$ is a regular sequence in $A[V_{ij}]$.
This proves the assertion over every field.

Now form the same complex over $\Z$, with all the displayed
variables indeterminate. Give the variables positive degree and
use the homogeneous shifts in the determinantal and Koszul
complexes. Each graded piece is a bounded complex of finitely
generated free abelian groups. Its base change to $\Q$ and to
every $\F_p$ has zero homology in positive degrees by the field
case. The universal coefficient sequence therefore implies that
the integral positive homology vanishes and that its degree-zero
homology is $\Z$-torsion-free. Hence the quotient is $\Z$-flat,
and tensoring this universal resolution with any $R_0$ preserves
exactness. This proves the result in general.
\end{proof}

We may now complete the proof of  Proposition \ref{p:resolve-j}, which as noted in \S\ref{s:outline} completes the proof of Theorem~\ref{t:comm} (and hence of Theorem~\ref{t:rl}).  The second part of the proposition was proven in \S\ref{ss:extend-fv}, so we must only demonstrate the first part.

To apply Proposition~\ref{p:br-exact}, we verify its numerical
hypothesis using the notation $k,r,\ell$ from
\S\ref{s:generators}. Recall also that
$s=|\Sigma\setminus\{v_0\}|$, with $s=0$ when
$\Sigma=\emptyset$, and $t=|\cP|$.

The matrix $D$ has $r+s$ rows. Exactly $r-k$ of these are
of types III, IV, or V, one for each supplemental generator.
Thus the number of linear generators $L_i$ is
\[
(r+s)-(r-k)=k+s.
\]
There are $\ell$ surviving $\boldsymbol b$-variables.
The local blocks contain the $\ell-k$ supplemental variables
associated to type IV and V rows, together with the $t$
reserved initial variables with indices $j_v$, $v\in\cP$.
Consequently, the number of variables in the local blocks is
\[
N=\ell-k+t,
\]
and the number outside those blocks is
\[
q=\ell-N=k-t.
\]
The $s+t$ local blocks are nonempty: the $\cP$-blocks
contain their distinguished labels, and we have already
discarded matrices $D_0$ having no type V row for some
$v\in\Sigma\setminus\{v_0\}$.

In the notation of Proposition~\ref{p:br-exact}, we therefore have
\[
m=k+s,\qquad h=s+t,\qquad
m=(k-t)+(s+t)=q+h.
\]
Thus the required inequality $m\leq q+h$ holds.
After relabeling the variables for this application, the
presentation \eqref{e:present-R} identifies $R$ and its
generators of $J'$ with the setup of that proposition.
It follows that the complex $\cC^\bullet$ defined in
\eqref{e:cdef} is a resolution of $R/J'$, proving
part {\rm (1)} of Proposition~\ref{p:resolve-j}.

\end{document}